\documentclass[a4paper,11pt,titlepage,twoside,openright]{book} 

\usepackage[T1]{fontenc}
\usepackage[latin1]{inputenc}
\usepackage[english, italian]{babel} 
\usepackage{cite} 
\usepackage{frontespizio} 

\usepackage{amsmath}
\usepackage{amssymb}
\usepackage{amsthm}
\usepackage{mathtools}
\usepackage{dsfont}
\usepackage{textcomp}
\usepackage{mathrsfs}

\usepackage{tikz-cd}
\usepackage{graphicx}
\usetikzlibrary{cd,arrows,chains,matrix,positioning,scopes}
\usepackage{adjustbox}
\usepackage{url}

\usepackage{pgfplots}
\pgfplotsset{/pgf/number format/use comma,compat=newest}
\usetikzlibrary{decorations.text}

\usepackage{hyperref}
\hypersetup{colorlinks,%
citecolor=black,%
filecolor=black,%
linkcolor=black,%
urlcolor=black}

\usepackage{fancyhdr} 

\usepackage{makeidx} 
\usepackage{nomencl} 

\newcommand\restr[2]{{
  \left.\kern-\nulldelimiterspace 
  #1 
  \vphantom{\big|} 
  \right|_{#2} 
  }}

{\left\lbrace\begin{array}{@{}l@{}}}%
{\end{array}\right.}

\newcommand{\N}{\mathbb{N}}
\newcommand{\Z}{\mathbb{Z}}

\renewcommand{\P}{\mathds{P}}
\newcommand{\V}{\mathds{V}}
\newcommand{\A}{\mathds{A}}
\renewcommand{\O}{\mathcal{O}}
\newcommand{\F}{\mathcal{F}}
\newcommand{\G}{\mathcal{G}}
\newcommand{\T}{\mathcal{T}}
\renewcommand{\L}{\mathcal{L}}

\DeclarePairedDelimiter{\intinf}{\lfloor}{\rfloor}

\DeclareMathOperator{\jac}{Jac}
\DeclareMathOperator{\Alb}{Alb}
\DeclareMathOperator{\alb}{alb}
\DeclareMathOperator{\Hom}{Hom}
\DeclareMathOperator{\Pic}{Pic}
\DeclareMathOperator{\ns}{NS}

\DeclareMathOperator{\im}{Im}
\DeclareMathOperator{\aut}{Aut}
\DeclareMathOperator{\bir}{Bir}
\DeclareMathOperator{\spec}{Spec}
\DeclareMathOperator{\specu}{\underline{Spec}}
\DeclareMathOperator{\proj}{Proj}
\DeclareMathOperator{\proju}{\underline{Proj}}

\DeclareMathOperator{\rk}{rk}
\DeclareMathOperator{\ext}{Ext}
\DeclareMathOperator{\cha}{char}
\DeclareMathOperator{\bcha}{\mathbf{char}}

\DeclareMathOperator{\ann}{Ann}
\DeclareMathOperator{\coker}{coker}
\DeclareMathOperator{\diff}{{D}iff}
\DeclareMathOperator{\iso}{Iso}
\DeclareMathOperator{\isou}{\underline{Iso}}
\DeclareMathOperator{\Fet}{Fet}
\DeclareMathOperator{\Fib}{Fib}

\theoremstyle{remark}
\newtheorem{remark}{Remark}[section]
\newtheorem{example}[remark]{Example}
\theoremstyle{definition}
\newtheorem{definition}[remark]{Definition}
\theoremstyle{plain}
\newtheorem{theorem}[remark]{Theorem}
\newtheorem{lemma}[remark]{Lemma}
\newtheorem{corollary}[remark]{Corollary}
\newtheorem{proposition}[remark]{Proposition}

\newcommand{\introthmname}{}
\newtheorem{introthminn}{\introthmname}
\newenvironment{introthm}[1]
  {\renewcommand{\introthmname}{#1}\begin{introthminn}}
  {\end{introthminn}}
	
	\swapnumbers

\makeindex 
\makenomenclature 

\begin{document}

\frontmatter 

\pagestyle{plain}

\begin{frontespizio} 
\Logo[6cm]{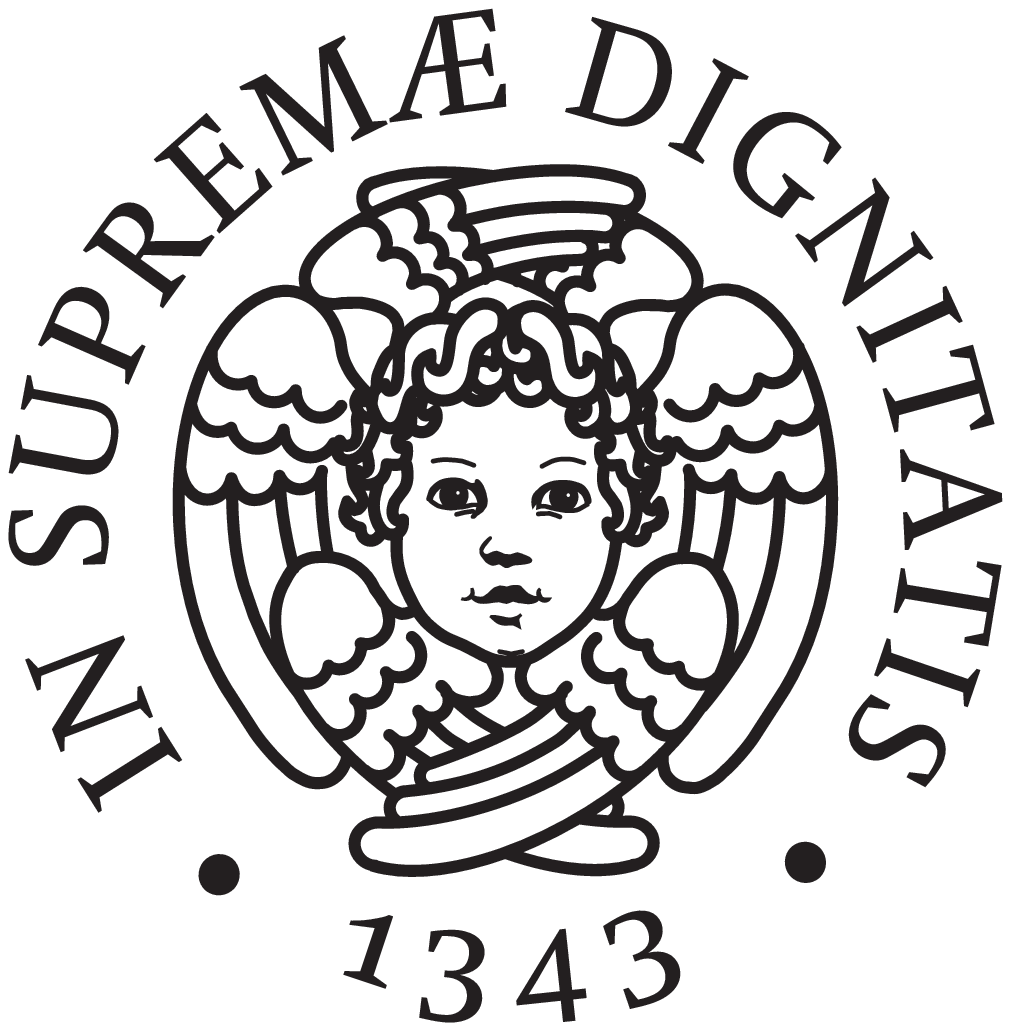}
\Istituzione{Universit\`a di Pisa}
\Dipartimento{Matematica}
\Corso[Dottorato di Ricerca]{Matematica}
\Titoletto{Tesi di Dottorato di Ricerca}
\Titolo{Surfaces close to the Severi lines}
\Candidato{Federico Cesare Giorgio Conti}
\Relatore{Prof.ssa Rita Pardini}
\NRelatore{Relatrice}{Relatori}
\Rientro{1cm}
\Piede{XXXIII Ciclo, Pisa 2021}
\end{frontespizio}

\begin{otherlanguage}{english} 

\chapter*{Abstract}

In this Thesis we study surfaces of general type with maximal Albanese dimension for which the quantity $K_X^2-4\chi(\O_X)-4(q-2)$ vanishes or is "small", that is surfaces close to the Severi lines. 
Over the complex numbers, it is known that a surface $X$, provided that $K_X^2<\frac{9}{2}\chi(\O_X)$, has to satisfy the inequality $K_X^2-4\chi(\O_X)-4(q-2)\geq 0$. We give a constructive and complete classification of surfaces for which equality holds: these are surfaces whose canonical model is a double cover of an Abelian surface ($q=2$) or of a product elliptic surface ($q\geq 3$) branched over an ample divisor with at most negligible singularities which intersects the elliptic fibre twice in the latter case.

We also prove, in the same hypothesis, that a surface $X$ with $K_X^2\neq 4\chi(\O_X)+4(q-2)$ satisfies $K_X^2\geq 4\chi(\O_X)+8(q-2)$ and we give a characterization of surfaces for which the equality holds. These are surfaces whose canonical model  is a double cover of an isotrivial smooth elliptic surface fibration branched over an ample divisor with at most negligible singularities whose intersection with the elliptic fibre is $4$.

Because these results are intimately related to theory of double covers, we see that their proof extend almost step by step to the case of any algebraically closed field of characteristic different from $2$.
We also give some partial results over algebraically closed fields of characteristic $2$ after a study of double covers in that case. 

\tableofcontents 

\pagestyle{fancy} 
\renewcommand{\chaptermark}[1]{\markboth{#1}{}}
\renewcommand{\sectionmark}[1]{\markright{\thesection\ #1}}
\fancyhf{}
\fancyhead[LE,RO]{\bfseries\thepage}
\fancyhead[RE]{\bfseries\footnotesize\nouppercase{\leftmark}}
\fancyhead[LO]{\bfseries\footnotesize\nouppercase{\rightmark}}

\phantomsection 
\addcontentsline{toc}{chapter}{List of Symbols} 
\printnomenclature 

\chapter{Introduction}

The aim of this thesis is the study and the classification of surfaces of general type with maximal Albanese dimension which lie on or close to the Severi lines, that is for which the quantity $K_X^2-4\chi(\O_X)-4(q-2)$ vanishes or is "small". 
\nomenclature{$K_X$}{the canonical divisor of a smooth, or at least Gorenstein, variety $X$}
\nomenclature{$\chi(\O_X)$}{the Euler characteristic of the structure sheaf of $X$}

More in detail, let $X$ be a minimal surface of general type with maximal Albanese dimension (recall that a surface is called of maximal Albanese dimension if its Albanese morphism is generically finite). We denote by $K_X$ the canonical divisor, by $\chi(\O_X)$ the Euler characteristic of the structure sheaf and by $q$ the dimension of the Albanese variety (which, in positive characteristic, may be smaller than the irregularity $q'=h^1(X,O_X)$). 
\nomenclature{$q(X)$}{the dimension of the Albanese variety of $X$}
\nomenclature{$q'(X)$}{the irregularity of a surface $X$, i.e. $q'(X)=h^1(X,\O_X)$}

The value we are interested in is strictly related to the so called Severi inequality (cf. \cite{par_sev} over the complex numbers and  \cite{yuan} for the general case), which states that a minimal surface of general type with maximal Albanese dimension satisfies
\begin{equation*}
K_X^2\geq 4\chi (\O_X).
\label{sev}
\end{equation*}
In \cite{barparsto} there is a characterization of complex surfaces for which the inequality above is indeed an equality, namely these are surfaces whose canonical model is a double cover of its Albanese variety branched over an ample divisor with at most negligible singularities (in particular $q=2$). There are many generalizations of the Severi inequality; in particular Lu and Zuo have proved in \cite{lu} a similar inequality involving also the irregularity $q$:  a complex surface of general type and maximal Albanese dimension satisfies
\begin{equation*}
K_X^2\geq\min\Bigl\{\frac{9}{2}\chi(\O_X),4\chi(\O_X)+4(q-2)\Bigr\}
\label{lusev}
\end{equation*}
or, equivalently, if $K_X^2<\frac{9}{2}\chi(\O_X)$ then $K_X^2\geq 4\chi(\O_X)+4(q-2)$.
They also give conditions for a surface to satisfy the equality
\begin{equation*}
\frac{9}{2}\chi(\O_X)>K_X^2=4\chi(\O_X)+4(q-2).
\label{uglu}
\end{equation*}
The condition $K_X^2<\frac{9}{2}\chi(\O_X)$ is necessary to prove that there exists an involution $i$ for which the Albanese morphism of $X$ is composed with $i$ (cf. \cite{lu} Theorem 3.1) which is central in their argument. 

The first original results of this thesis are obtained over the complex numbers: we give  a complete characterization of surfaces satisfying $K_X^2=4\chi(\O_X)+4(q-2)$ in case $q\geq 3$ and $K_X^2<\frac{9}{2}\chi(\O_X)$: in \cite{lu} it is proved that the canonical model of such a surface is a double cover of a smooth isotrivial elliptic surface branched over a divisor $R$ with at most negligible singularities. We prove here that this elliptic surface has to be a product $C\times E$ and we also determine the linear class of the branch divisor. 

\begin{introthm}{Theorem}[\ref{mioa}]
\label{introteo1}
Let $X$ be a minimal complex surface of general type with maximal Albanese dimension satisfying $q=q(X)\geq 3$ such that $K_X^2<\frac{9}{2}\chi(\O_X)$.
Then 
\begin{equation*}
K_X^2= 4\chi(\O_X) +4(q-2)
\end{equation*}
if and only if the canonical model of $X$ is isomorphic to a double cover of a product elliptic surface $Y=C\times E$ where $E$ is an elliptic curve and  $C$ is a curve of  genus $q-1$, whose branch divisor $R$ has at most negligible singularities and 
\[R\sim_{lin}C_1+C_2+\sum_{i=1}^{2d}E_i,\]
where $E_i$ (respectively $C_i$) is a fibre of the first projection (respectively the second projection) of $C\times E$ and $d>7(q-2)$. Moreover, we have that $\Alb(X)\simeq\Alb(Y)$ and, in particular, if $q\geq 3$, the Albanese variety of $X$ is not simple.
\nomenclature{$C$}{a curve of genus greater than $1$}
\nomenclature{$E$}{an elliptic curve}
\nomenclature{$R$}{the branch divisor of a finite morphism}
\nomenclature{$\Alb(X)$}{the Albanese variety of $X$}
\end{introthm} 


The second result is about surfaces that are not on the Severi lines but are close to them. We see that in this case $K_X^2\geq 4\chi(\O_X)+ 8(q-2)$ and we also give a characterization for surfaces that satisfy this equality.

\begin{introthm}{Theorem}[\ref{miob}]
\label{introteo2}
Let $X$ be a minimal complex surface of general type with maximal Albanese dimension with $K_X^2<\frac{9}{2}\chi(\O_X)$.
\begin{enumerate}
\item If $K^2_X>4\chi(\O_X)+4(q-2)$, then $K^2_X\geq 4\chi(\O_X)+8(q-2).$
\item If $q=2$ and $K^2_X>4\chi(\O_X)$, then $K^2_X\geq 4\chi(\O_X)+2.$
\item If $q\geq 3$, equality holds, i.e.
\[K^2_X=4\chi(\O_X)+8(q-2),\]
if and only if the canonical model of $X$ is isomorphic to a double cover of a smooth isotrivial elliptic surface fibration $Y$ over a curve $C$ of genus $q-1$, branched over a divisor $R$ with at most negligible singularities for which $K_Y.R=8(q-2)$. In particular, we have that $\Alb(X)\simeq\Alb(Y)$ and, if $q(X)\geq 3$, we have that the Albanese variety of $X$ is not simple.
\nomenclature{$D_1.D_2$}{the intersection product between two divisors $D_1$ and $D_2$}
\end{enumerate}
\end{introthm}

Next we try to recover similar results over algebraically closed fields of positive characteristic.
The Severi inequality was first proved over algebraically closed fields of any characteristic  in \cite{yuan} using different methods from those used over the complex numbers.
In \cite{gusunzhou} it is pointed out that the same proof used in \cite{par_sev} can be adjusted to be used over every algebraically closed field: in the same article it is proven a refined Severi inequality, i.e. for every surface $X$ of general type with maximal Albanese dimension
\[K_X^2\geq\Bigl(4+\min \{c(X),\frac{1}{3}\}\Bigr)\chi(\O_X)\]
holds, where $c(X)$ is a constant that only depends on the morphisms $X\to Y_i$ of degree two which are factors of the Albanese morphism $\alb_X\colon X\to \Alb(X)$ of $X$.
\nomenclature{$\alb_X$}{the Albanese morphism from $X$ to its Albanese variety}
With a slight improvement of Theorem  3.1 of \cite{gusunzhou}, we improve this inequality and, in parallel, we also show that, for a surface of general type $X$ with maximal Albanese dimension, if $K_X^2<\frac{9}{2}\chi(\O_X)$, there exists a morphism $X\to Y$ of degree $2$ which is a factor of the Albanese morphism. 
That is, we prove the two following Theorems.

\begin{introthm}{Theorem}[\ref{teo_albfact2}]
Let $X$ be a surface of general type with maximal Albanese dimension and suppose that $K_X^2<\frac{9}{2}\chi(\O_X)$. 
Then there exists a morphism of degree two $f\colon X\to Y$ to a normal surface $Y$ such that the Albanese morphism of $X$ factors through $f$, i.e. the following diagram commutes:
\begin{equation*}
\begin{tikzcd}
X\arrow{rr}{f}\arrow{dr}{\alb_X} & & Y\arrow{ld}{}\\
& \Alb(X). & 
\end{tikzcd}
\end{equation*}
\label{introteo3}
\end{introthm}

This Theorem is particularly interesting over a field of characteristic $2$, indeed for a field of characteristic different from $2$ is very similar to the case of complex numbers and is tacitly proved in \cite{gusunzhou}.

\begin{introthm}{Theorem}[\ref{teo_cxl}]
Let $X$ be a minimal surface of general type of maximal Albanese dimension.
Then we have
\begin{equation*}
K_X^2\geq\Bigl(4+\min \{c(X),\frac{1}{2}\}\Bigr)\chi(\O_X)
\label{eq_cxl}
\end{equation*}
for a constant $c(X)$.
\label{introteo4}
\end{introthm}

Once this has been done, we try to extend the results obtained in \cite{lu} and by us over the complex numbers in the study of surfaces of general type with maximal Albanese dimension which lie on or close to the Severi lines for fields of positive characteristic.
Unfortunately there some issues over fields of characteristic $2$, based on the fact that in this case flat double covers behave wildly,  that we were not able to solve: in particular we generalize Theorems \ref{introteo1} and \ref{introteo2} only over algebraically closed fields of characteristic different from $2$.

\begin{introthm}{Theorem}[\ref{teo_ultimmio+}]
\label{introteo5}
Let $X$ be a minimal surface of general type of maximal Albanese dimension over an algebraically closed field of characteristic different from $2$ and assume that $K_X^2<\frac{9}{2}\chi(\O_X)$.
Then we have that 
\begin{equation*}
K_X^2\geq 4\chi(\O_X)+4(q-2)
\label{eq:luzuoneq3}
\end{equation*}
where $q$ is the dimension of the Albanese variety of $X$.

In particular equality holds, i.e.
\begin{equation*}
K_X^2=4\chi(\O_X)+4(q-2)
\label{eq_teomioa+}
\end{equation*}
if and only if the canonical model of $X$ is isomorphic to a double cover of a product elliptic surface ($q\geq 3$) $Y=C\times E$ where $E$ is an elliptic curve and $C$ is a curve of genus $q-1$, whose branch divisor $R$ has a most negligible singularities and
\begin{equation*}
R\sim_{lin} C_1+C_2+\sum_{i=1}^{2d}E_i
\label{eq_teomioabranch}
\end{equation*}
where  $E_i$ (respectively $C_i$) is a fibre of the first projection (respectively the second projection) of $C\times E$ and $d>7(q-2)$ or the canonical model of $X$ is a double cover of an Abelian surfaces branched over an ample divisor ($q=2$). 
Moreover we have that $\Alb(X)=\Alb(Y)$ and $q(X)=q'(X)$, i.e. the Picard scheme of $X$ is reduced.
In particular, if $q\geq 3$, the Albanese variety of $X$ is not simple.
\end{introthm}

\begin{introthm}{Theorem}[\ref{teo_closesev+}]
\label{introteo6}
Let $X$ be a minimal surface of general type with maximal Albanese dimension with $K_X^2<\frac{9}{2}\chi(\O_X)$ over an algebraically closed field of characteristic different from $2$ and $q$ be the dimension of the Albanese variety $\Alb(X)$ of $X$.
\begin{enumerate}
\item If $K^2_X>4\chi(\O_X)+4(q-2)$, then $K^2_X\geq 4\chi(\O_X)+8(q-2).$
\item If $q=2$ and $K^2_X>4\chi(\O_X)$, then $K^2_X\geq 4\chi(\O_X)+2.$
\item If $q\geq 3$, equality holds, i.e.
\begin{equation*}
\label{eq_8q-2}
K^2_X=4\chi(\O_X)+8(q-2),
\end{equation*}
if and only if the canonical model of $X$ is isomorphic to a double cover of a smooth isotrivial elliptic surface fibration $Y$ over a curve $C$ of genus $q-1$, branched over a divisor $R$ with at most negligible singularities for which $K_Y.R=8(q-2)$. In particular, we have that $\Alb(X)\simeq\Alb(Y)$ and $q(X)=q'(X)$, i.e. the Picard variety of $X$ is reduced, and the Albanese variety of $X$ is not simple.
\end{enumerate}
\end{introthm}

As we have already observed, there are some problems which prevent us from proving these Theorems in characteristic $2$: anyway we would like to stress that in \cite{gusunzhou} there is a characterization of minimal surfaces of general type with maximal Albanese dimension satisfying $K_X^2=4\chi(\O_X)$ as surfaces whose Albanese morphism is a generically finite morphism of degree $2$, their proof works also over fields of characteristic $2$ and the Albanese morphism may well be inseparable.

The outline of the Thesis is the following.
In Chapter \ref{chap_preliminaries} we collect some preliminary results that we will need throughout all the dissertation.
A particular attention is given to the study of elliptic surfaces of maximal Albanese dimension: as pointed out in Theorems \ref{introteo1}, \ref{introteo2}, \ref{introteo5} and \ref{introteo6}, the surfaces of general type we are interested in turn out to be double covers of this kind of elliptic surfaces.
Specifically we show that every smooth elliptic surface fibration admits a base change such that it becomes a product surface (cf. Proposition \ref{propo_isoell} and Corollary \ref{coro_isoell}) and base change is not needed when there exists a section of the fibration.
Then a study of the invariant Picard group of a product elliptic surface under the action of an Abelian group is done in order to be able to construct examples of non-trivial smooth elliptic surface fibrations.
In this Chapter we also recall the classification of algebraic surfaces, the extension of the classical theory for smooth curves to the case of Gorenstein curves (with a particular attention to the case of integral curves on surfaces), the theory of semi-stable vector bundles on curves, some results on projective bundles and on the algebraic fundamental group.

In Chapter \ref{chap_double}, we collect all that is known about double covers of surfaces. 
When the characteristic is different from two, the argument is standard and well known (cf. \cite{barth} Section I.17).
Over fields of characteristic $2$ different phenomena arise (such as inseparable or non-splittable double covers) and the theory is more involved.
In particular the fact that the line bundle $L$ associated with a double cover is not necessarily an effective divisor (up to a multiple) is an obstruction to a possible extension of Theorems \ref{introteo1} and \ref{introteo2} over fields of characteristic $2$.

In Chapter \ref{chap_severi0} we focus on complex surfaces, we recall the work of \cite{lu} and we prove Theorems \ref{introteo1} and \ref{introteo2}.
Moreover we give constructive examples of surfaces satisfying equalities $K_X^2-4\chi(\O_X)=4(q-2)$ and $K_X^2-4\chi(\O_X)=8(q-2)$.

In Chapter \ref{chap_from0to+} we prove Theorems \ref{introteo3} and \ref{introteo4}, the former of which has to be intended as a motivation for a generalization of the results of Chapter \ref{chap_severi0} and of \cite{lu} over fields of any characteristic while the latter is simply a generalization of Theorem 5.14 of \cite{gusunzhou}.
In order to do so, we have to improve a Theorem on non-hyperelliptic fibred surfaces (\cite{gusunzhou} Theorem 3.1) and we do this in Theorem \ref{teo_slopegusunzhou}.

\begin{introthm}{Theorem}[\ref{teo_slopegusunzhou}]
\label{introteo7}
Let $f\colon X\to C$ be a non-hyperelliptic  surface fibration such that $K_{X/C}$ is nef and all the quotients appearing in the Harder-Narasimhan filtration of $f_*\omega_{X/C}$ are strongly semi-stable.
\nomenclature{$f_*F$}{the pushforward of a coherent sheaf $F$ via the morphism $f$}
 Assume that there exists a constant $0\leq c\leq\frac{1}{2}$ such that for every diagram
\[
\begin{tikzcd}
X\arrow{rr}{}\arrow[swap]{dr}{f} & & Y\arrow{dl}{\overline{f}}\\
& C & 
\end{tikzcd}
\]
 where $X\to Y$ is a morphism of degree $2$ with $Y$ a smooth surface, $\overline{g}\geq cg$ holds where $g$ and  $\overline{g}$ are the arithmetic genus of a general fibre of $f$ and $\overline{f}$ respectively: then
	\begin{equation}
	K_{X/C}^2\geq (4+c)\frac{g-1}{g+2}\deg(f_*\omega_{X/C}).
	\label{eq_gusunzhou2}
	\end{equation}
	\nomenclature{$g(C)$}{the arithmetic genus of a Gorenstein curve $C$}
\end{introthm}

Notice that this is a  generalization of Theorem 3.1 of \cite{gusunzhou} because we allow $c\in[0,\frac{1}{2}]$ instead of $c\in[0,\frac{1}{3}]$.
The proof of this result is rather technical and require the usage of some geometrical construction that have been used for the proof of other results, such as the slope and the Severi inequality.
That the same proofs used over the complex numbers for these inequalities, up to minor considerations which one has to do over fields of positive characteristic (mainly coming from the fact that the general fibre of a surfaces fibration need not be smooth), still work over fields of any characteristic have been pointed out (cf. \cite{gu2019slope} Section 2.1 for the slope inequality and \cite{gusunzhou} Section 4 for the Severi inequality): however, in order to get the reader more acquainted with the geometric objects used and make a more linear discussion so that they understand better the proof of Theorem \ref{teo_slopegusunzhou} and the results we derive from it, we thought it to be useful to recall most of these proofs here.

Finally, in Chapter \ref{chap_severi+} we see how to adapt the proofs of Chapter \ref{chap_severi0} for surfaces over fields of characteristic different from $2$, i.e. we prove Theorems \ref{introteo5} and \ref{introteo6}.
We also collect partial results we have obtained over fields of characteristic $2$ for Severi type inequalities and we give some examples also in this case.

\section*{Notations and conventions}
Throughout all the thesis by $k$ we mean an algebraically closed field.
In chapter \ref{chap_severi0} we work over the complex numbers.
By a variety we mean a projective integral scheme over $k$.
For a smooth surface $X$ we will indicate by $q(X)$ the dimension of the Albanese dimension and by $q'(X)$ its irregularity, i.e. $q'(X)=h^1(X,\O_X)$: clearly over the complex numbers these two quantities coincide and we will indicate it as $q(X)$.
For an integral Gorenstein curve $C$ we denote by $g(C)$ it's arithmetic genus, which coincides with it's geometric genus if and only if $C$ is smooth.
We use interchangeably the notion of line bundles and Cartier divisors and we use both additive and multiplicative notations.
When we consider a product $A\times B$ of two schemes we will always denote by $\pi_A\colon A\times B\to A$ and $\pi_B\colon A\times B\to B$ the first and the second projection. We will denote by $A$ and $B$ a general fibre of $\pi_B$ and $\pi_A$ respectively: furthermore, we will use $A_b$ and $B_a$ to indicate the fibre over $b\in B$ and $a\in A$ respectively.

\mainmatter

\chapter{Preliminaries}
\label{chap_preliminaries}
In this Chapter we present the preliminary materials needed in this thesis.
In Section \ref{sec_frob} we recall what is the $k$-linear Frobenius morphism and some of its properties.
In Section \ref{sec_sss} we present the theory of semi-stable vector bundles on curves and the Harder-Narasimhan filtration; moreover we define and give some properties of saturated subsheaves in every dimension.
In Section \ref{sec_projbundle} we collect some results on projective bundles, in particular we are interested in the relation between the nefness of some line bundles on $\P(E)$ and the numerical invariants of the Harder-Narasimhan filtration of $E$.
In Section \ref{sec_fundgroup} we define and look at some properties of the algebraic fundamental group: in particular we are interested in its relation with finite \'etale covers, in the Lefschetz hyperplane Theorem and the algebraic fundamental group of a projective bundle.
In Section \ref{sec_surfaces} we briefly recall the classification of algebraic surfaces.
In Section \ref{sec_curves} we collect some known results about non-smooth curves (especially integral ones) on surfaces: we need this in the study of surface fibrations, because in positive characteristic the general fibre may not be smooth but has to be integral.
In Section \ref{sec_ell} we focus on elliptic surfaces, i.e. surfaces admitting a surjective morphism to a smooth curve such that the general fibre is a smooth curve of genus $1$: as we will see, the surfaces lying close to the Severi lines (which is the main topic of this thesis) turn out to be double covers of elliptic surface, that's why we are interested in them.
In particular we show that every smooth elliptic surface fibration is trivial after a base change on the base: if there exists a section of the fibration, we even show that actually no base change is needed.
In Section \ref{sec_piccxe} we describe the formula of the Picard group of a very particular elliptic surface, i.e. $C\times E$ where $E$ is an elliptic curve and $C$ a curve of higher genus: we also focus on the case where a finite Abelian group $G$ acts freely on both $C$ and $E$, and see what can be said for the Picard group of the quotient $(C\times E)/G$.
In Section \ref{sec_gentype} we recall something about the geography of surfaces of general type, stressing the differences between the case of characteristic zero and the case of positive characteristic.

\section{Frobenius morphisms\index{Frobenius morphism}}
\label{sec_frob}
In this section we briefly recall what is the $k$-linear Frobenius morphism and some of its properties.

\begin{definition}
Let $X$ be a smooth scheme over an algebraically closed field $k$ of characteristic $p$.
\nomenclature{$k$}{an algebraically closed field} 
The absolute Frobenius morphism\index{Frobenius morphism! absolute} is the morphism which is topologically equal to the identity and at the level of functions it behaves as the p-th power morphism.
\end{definition}

\begin{remark}
Observe that the absolute Frobenius morphism commutes with respect to any morphism of schemes $f\colon X\to Y$, i.e. the following diagram is commutative:
\begin{equation}
\begin{tikzcd}
X\arrow{d}{f}\arrow{r}{F}&X\arrow{d}{f}\\
Y\arrow{r}{F}&Y.
\end{tikzcd}
\label{eq_commfrobabs}
\end{equation}
\label{rem_absfrobcomm}
\end{remark}

\begin{definition}
\label{def_kfrob}
Let $X$ be a smooth scheme and consider the following Cartesian diagram:
\begin{equation}
\label{eq_kfrob}
\begin{tikzcd}
X^{(-1)}\arrow{r}{} \arrow[dr, phantom, "\square"] \arrow{d}{\pi^{(-1)}} & X\arrow{d}{\pi}\\
\spec(k)\arrow{r}{F} & \spec(k),
\end{tikzcd}
\end{equation}
where the lower arrow is the absolute Frobenius and the right one is the structure morphism of $X$.

Then the morphism $F_k^{(0)}\colon X^{(0)}=X\to X^{(-1)}$, given by the following diagram, is called $k$-linear Frobenius morphism\index{Frobenius morphism!k-linear @ $k$-linear}:
\begin{equation}
\label{eq_kfrobdef}
\begin{tikzcd}
X\arrow[bend right=30]{ddr}{}\arrow[bend left=30]{drr}{F}\arrow{dr}{F_k^{(0)}} & & \\
& X^{(-1)}\arrow{r}{} \arrow[dr, phantom, "\square"] \arrow{d}{\pi^{(-1)}} & X\arrow{d}{\pi}\\
& \spec(k)\arrow{r}{F} & \spec(k).
\end{tikzcd}
\nomenclature{$F_k$}{the $k$-linear Frobenius morphism}
\end{equation}
Inductively, we define $F_k^{(-n)}\colon X^{(-n)}\to X^{(-n-1)}$ for every positive integer $n$, where, by definition, $X^{(-n-1)}=(X^{(-n)})^{(-1)}$.

Consider the following commutative diagram
\begin{equation}
\begin{tikzcd}
X\arrow{r}{F}\arrow{d}{\pi}&X\arrow{d}{\pi}\\
\spec(k)\arrow{r}{F}&\spec(k)
\end{tikzcd}
\label{eq_kfrobdef1}
\end{equation}
where the horizontal arrows are the absolute Frobenius morphisms of $X$ and $\spec(k)$ respectively.
Denote by $X^{(1)}$ the $k$-scheme defined by $X$ with structure morphism given by $\pi\circ F$: in this case $F$ becomes a $k$-linear morphism $F_k^{(1)}\colon X^{(1)}\to X$ which we still call the $k$-linear Frobenius morphism.
As above, we define inductively $F_k^{(n)}\colon X^{(n)}\to X^{(n-1)}$ for every positive integer $n$, where, by definition, $X^{(n)}=(X^{(n-1)})^{(1)}$.

We will usually write $F_k$ instead of $F_k^{(n)}$ because it will be clear by the context which $k$-linear Frobenius morphism we consider.
\end{definition}

\begin{remark}
Let $f\colon X\to Y$ be a morphism of $k$-schemes and consider the following commutative diagram
\begin{equation}
\begin{tikzcd}
X\arrow{dd}{f}\arrow{dr}{F_k}\arrow[bend left=25]{drr}{F} & &\\
& X^{(-1)}\arrow{r}{}\arrow[ddr, phantom, "\square"]& X\arrow{dd}{f}\\
Y\arrow[bend right=30]{dddr}{}\arrow[bend left=25]{drr}[near start]{F}\arrow{dr}{F_k} & & \\
& Y^{(-1)}\arrow[from=uu, crossing over, "f^{(-1)}" near end]\arrow{r}{} \arrow[ddr, phantom, "\square"] \arrow{dd}{\pi^{(-1)}} & Y\arrow{dd}{\pi}\\
& & \\
& \spec(k)\arrow{r}{F} & \spec(k):
\end{tikzcd}
\label{eq_functkfrob}
\end{equation}
this is clearly commutative, in particular $f^{(-1)}\circ F_k=F_k\circ f$.

Similarly, consider the following diagram
\begin{equation}
\begin{tikzcd}
X^{(1)}\arrow{r}{F_k}\arrow{d}{f^{(1)}} &X\arrow{d}{f}\\
Y^{(1)}\arrow{d}{}\arrow{r}{F_k}& Y\arrow{d}{}\\
\spec(k)\arrow{r}{F}&\spec(k):
\end{tikzcd}
\label{eq_functkfrob1}
\end{equation}
this is also commutative, in particular $F_k\circ f^{(1)}=f\circ F_k$.

It is immediate by definition that $Id_X^{(1)}=Id_{X^{(1)}}$, $Id_X^{(-1)}=Id_{X^{(-1)}}$, $f^{(1)}\circ g^{(1)}=(f\circ g)^{(1)}$ and $f^{(-1)}\circ g^{(-1)}=(f\circ g)^{(-1)}$.
We have thus shown how to make $\cdot^{(1)}$ and $\cdot^{(-1)}$ (and proceeding by induction $\cdot^{(n)}$ for any integer $n$) into a functor from $k$-schemes to $k$-schemes, such that its action on morphisms commutes with the $k$-linear Frobenius morphism.
\end{remark}

\begin{remark}
\label{rem_frobassrel}
Restricting to an affine open subset of $X$ we may assume that $X=\spec k[x_1,\ldots,x_n]/I$ (clearly here we are assuming $X$ being of finite type over $k$).
Both the absolute and the $k$-linear Frobenius are topologically equivalent to the identity: hence they are uniquely characterized by their action on functions.
The absolute morphism acts on functions as the p-th power morphism, while the $k$-linear Frobenius acts as a "$k$-linear p-th power".
Indeed the diagram \ref{eq_kfrobdef} in the affine chart becomes
\begin{equation}
\label{eq_kfrobdefaff}
\begin{tikzcd}
k[x_1,\ldots,x_n]/I\arrow[bend right=20,leftarrow]{ddr}{}\arrow[bend left=20,leftarrow]{drr}{F}\arrow[leftarrow]{dr}{F_k} & & \\
& k[x_1,\ldots,x_n]/I^{(p)}\arrow[leftarrow]{r}{} \arrow[leftarrow]{d}{} & k[x_1,\ldots,x_n]/I\arrow[leftarrow]{d}{}\\
& k\arrow[leftarrow]{r}{F} & k,
\end{tikzcd}
\end{equation}
where, if $I$ is generated by polynomials $f_i=\sum_Jf_{iJ}x^J$ (where for $J=(j_1,\ldots,j_n)$ we denote by $x^J=x_1^{j_1}\cdot\ldots\cdot x_n^{j_n}$), then $I^{(p)}$ is generated by $f_i^{(p)}=\sum_Jf_{iJ}^px^J$.
It is then clear that $F_k\colon k[x_1,\ldots,x_n]/I^{(p)}\to k[x_1,\ldots,x_n]/I$ is defined by $\sum_Ja_Jx^J\mapsto\sum_Ja_J(x^J)^p$.

Similarly, the diagram \ref{eq_kfrobdef1} in the affine chart becomes 
\begin{equation}
\begin{tikzcd}
k[x_1,\ldots,x_n]/I^{(1/p)}\arrow[leftarrow]{r}{F_k}\arrow[leftarrow]{d}{} &k[x_1,\ldots,x_n]/I\arrow[leftarrow]{d}{}\\
k\arrow[leftarrow]{r}{F} & k,
\end{tikzcd}
\label{eq_kfrobaffdef1}
\end{equation}
where the upper arrow is a morphism of $k$-algebra and $I^{(1/p)}$ is generated by $f_i^{(1/p)}=\sum_Jf_{iJ}^{1/p}x^J$ with the same notation as above (recall that over an algebraically closed field of characteristic $p$, the $p$-th root is uniquely determined).
It is then clear that $F_k\colon k[x_1,\ldots,x_n]/I\to k[x_1,\ldots,x_n]/I^{(1/p)}$ is defined by $\sum_Ja_Jx^J\mapsto\sum_Ja_J(x^J)^p$.

By this local computation, it is also clear that $(X^{(1)})^{(-1)}=X$, $(X^{(-1)})^{(1)}=X$ and, more generally, $(X^{(m)})^{(n)}=X^{(m+n)}$ for any integers $m$ and $n$.
\end{remark}

\begin{proposition}
\label{propo_factfrob}
A morphism $f\colon X\to Y$ between smooth varieties factors through the $k$-linear Frobenius $F_k\colon X\to X^{(-1)}$ if and only if the morphism $f^*\Omega^1_Y\to\Omega^1_X$ is trivial.
\nomenclature{$\Omega^1_X$}{the sheaf of K\"ahler differentials on $X$}
\nomenclature{$f^*F$}{the pull-back of a coherent sheaf $F$ via the morphism $f$}
\end{proposition}

\begin{proof}[Sketch of proof]
This follows almost immediately from the fact that for a regular ring $A$ the kernel of the universal derivation $d\colon A\to\Omega^1_A$ is $A^p$ (cf. \cite{brion} Theorem 1.3.4).

We restrict to affine open subsets such that $f$ is given by $f^*\colon B\to A$ where $B=k[y_1,\ldots,y_m]/J$ and $A=k[x_1,\ldots,x_n]/I$.
If $f^*$ factors as $k[y_1,\ldots,y_m]/J \to k[x_1,\ldots,x_n]/I^{(p)}\xrightarrow{F_k} k[x_1,\ldots,x_n]/I$ then clearly the cotangent map is trivial.

Conversely suppose that $f^*\colon\Omega^1_B\to\Omega^1_A$ is trivial.
In particular, for every $b\in B$, we have $0=f^*db=df^*b$, which we have already seen to be equivalent to $f^*b=a^p$ for a unique $a\in A$.
If we write $a=\sum_Ja_Jx^J$ we clearly see that we can factor $f^*\colon B\to A$ as $b\mapsto \sum_Ja_J^px^J\mapsto \sum_Ja_J^p(x^J)^p=f^*b$, where the last part is clearly the $k$-linear Frobenius.
One has to verify that these local morphisms glue together and this follows from the  uniqueness of the $k$-linear Frobenius morphism and of the p-th root in characteristic $p$.
\end{proof}

\begin{remark}
\label{rem_albfrob}
The morphism $\alb_X\colon X\to\Alb(X)$ does not factors through the $k$-linear Frobenius morphism.
\index{Albanese!variety} \index{Albanese!morphism}
Indeed, by a Theorem of Igusa \cite{igusa}, we know that the pull-back via the Albanese morphism of global differential 1-forms is injective: in particular $\alb_X^*\Omega^1_{\Alb(X)}\to \Omega^1_X$ can not be trivial.

Moreover, if the Albanese morphism $\alb_X\colon X\to \Alb(X)$ of a surface $X$ is inseparable, the image of $\alb_X^*\Omega^1_{\Alb(X)}\to \Omega^1_X$ is  a coherent sheaf of rank $1$: indeed, because the morphism $\alb_X$ is inseparable, it can not be of maximal rank (cf. \cite{liu} Lemma 6.1.13(b)).
\end{remark}

\section{Semi-stable sheaves}
\label{sec_sss}
%
We start this section defining saturated subsheaves and giving some properties that we will need in the study of foliations in Section \ref{sec_foliations}.
We will then  study semi-stable sheaves and the Harder-Narasimhan filtration of a vector bundle over a smooth curve.
In this Section, by $X$ we mean, unless otherwise stated, a smooth variety and by $C$ a smooth curve.

\begin{definition}
\label{def_torsion}
Let $E$ be a coherent sheaf on $X$, we denote by $T(E)$ the torsion subsheaf\index{subsheaf!torsion} \nomenclature{$T(E)$}{the torsion subsheaf of a coherent sheaf $E$} of $E$, that is the sheaf whose sections are the torsion sections of $E$.

A sheaf $E$ is said to be torsion-free\index{sheaf!torsion-free} if $T(E)=0$ or, equivalently, if for every $x\in X$, the stalk $E_x$ \nomenclature{$E_x$}{the stalk of a coherent sheaf $E$ at $x$} is a torsion-free $\O_{X,x}$-module\nomenclature{$\O_X$}{the structure sheaf of $X$}.
\end{definition}

\begin{definition}
\label{def_saturation}
Let $E$ be a torsion-free sheaf.
The saturation\index{subsheaf!saturation of} of a subsheaf $F\subset E$ is the minimal subsheaf $F'$ containing $F$ such that the quotient $E/F'$ is torsion-free or zero. 
If $F'=F$ we say that $F$ is saturated\index{subsheaf!saturated}.
\end{definition}

Clearly we have that $F'$ is the kernel of the morphism $E\to E/F\to (E/F)/T(E/F)$.

\begin{proposition}[\cite{hartref} Proposition 1.1]
\label{propo_satref}
A coherent sheaf $F$ is reflexive\index{sheaf!reflexive} (that is $\Hom(\Hom(F,\O_X),\O_X)\simeq F$) if and only if (at least locally) it can be included in an exact sequence
\[0\to F\to E\to E/F\to 0\]
where $E$ is locally free and $T(E/F)=0$.
\end{proposition}

Notice that Proposition \ref{propo_satref} implies that every saturated subsheaf of a vector bundle is reflexive. This, together with the next Proposition, will be important in the study of purely inseparable double covers of surfaces.

\begin{proposition}[\cite{hartref} Proposition 1.9]
\label{propo_ref}
Let $X$ be a variety which is locally factorial (i.e. all its local rings are unique factorization domains, e.g. $X$ smooth). Then every reflexive rank $1$ sheaf is actually a line bundle.
\end{proposition}

\begin{remark}
\label{rem_folfree}
Combining Propositions \ref{propo_satref} and \ref{propo_ref} we get that every rank $1$ saturated subsheaf of a vector bundle on a smooth variety is a line bundle. 
We will use this in the study of 1-foliations on surfaces.
\end{remark}

From now on we restrict to study coherent sheaves on a smooth curve $C$.
In this case, it is known (by the well-known fact that over a principal ideal domain a torsion-free module is free) that every torsion-free sheaf is locally free\index{sheaf!locally free|see{vector bundle}}, i.e. it is a vector bundle\index{vector bundle}. 
In particular, every subsheaf of a vector bundle is a vector bundle.
Moreover the saturation $F'$ of a subsheaf $F$ of a vector bundle $E$ is a subvector bundle, i.e. the quotient $E/F'$ is a vector bundle too. 

\begin{definition}
\label{def_deg}
Let $E$ be a coherent sheaf on $C$.
We define the degree\index{sheaf!degree of} of $E$ to be $\deg(E)=c_1(E)$\index{degree! of a coherent sheaf on a smooth curve} the first Chern class of $E$\index{sheaf!first Chern class of}.
\nomenclature{$\deg(E)$}{the degree of a coherent sheaf $E$ on a smooth curve}
\nomenclature{$c_1(E)$}{the first Chern class of a coherent sheaf $E$}
\end{definition}

Recall that the first Chern class of a vector bundle $V$ can be defined as the first Chern class of the determinant bundle $\det(V)$ of $V$, while for a general coherent sheaf $E$ by taking a finite locally free resolution of $E$ and using additivity of first Chern classes on exact sequences.

\begin{definition}
\label{def_slope}
Let $E$ be a coherent sheaf.
We define the slope\index{sheaf!slope of} of $E$ as
\[\mu(E)=\frac{\deg(E)}{\rk(E)},\]
\nomenclature{$\mu(E)$}{the slope of a coherent sheaf $E$ on a smooth curve}
where, denoting by $\eta$ the generic point of $X$, $\rk(E)$\index{sheaf!rank of} is the dimension of the vector space $E\otimes_{\O_X}k(\eta)$ over the field $k(\eta)$ of the generic point of $C$.
\nomenclature{$\rk(E)$}{the rank of a coherent sheaf $E$}
\end{definition}

We are now ready to define what is a semi-stable vector bundle.

\begin{definition}
\label{def_semi}
A vector bundle $E$ is semi-stable\index{vector bundle!semi-stable} if  $\mu(F)\leq\mu(E)$ for every non-trivial subsheaf $F$ of $E$. 
\end{definition}

\begin{lemma}
\label{lemma_equivalentss}
Let $E$ be a vector bundle. Then the following are equivalent:
\begin{enumerate}
	\item $E$ is semi-stable;
	\item for every saturated subsheaf $F$ of $E$ one has $\mu(F)\leq\mu(E)$;
	\item for every quotient vector bundle $G$ of $E$ one has $\mu(E)\leq\mu(G)$;
	\item for every quotient sheaf $G$ of $E$ one has $\mu(E)\leq\mu(G)$.
\end{enumerate}
\end{lemma}

\begin{proof}
Consider an exact sequence of sheaves
\begin{equation}
0\to F\to E\to G\to 0.
\label{eq_equivss}
\end{equation}
By additivity of the degree and the rank function, we have $\deg(E)=\deg(F)+\deg(G)$ and $\rk(E)=\rk(F)+\rk(G)$.
In particular $\rk(F)(\mu(E)-\mu(F))=\rk(G)(\mu(G)-\mu(E))$ from which one can easily derive the equivalence of conditions 2 and 3 and of conditions 1 and 4.

Clearly 1 implies 2 and the reverse follows from the fact that the saturation $F'$ of $F\subseteq E$ satisfies $\rk(F')=\rk(F)$ and $\deg(F)\leq\deg(F')$.
\end{proof}

\begin{definition}
\label{def_HN}
Let $E$ be a vector bundle. 
A Harder-Narasimhan\index{vector bundle!Harder-Narasimhan filtration of} filtration for $E$ is an increasing filtration
\[0=E_0\subset E_1\subset \ldots \subset E_n=E,\]
such that the subsequent quotients $E_i/E_{i-1}$ are semi-stable vector bundles with slopes $\mu_i=\mu(E_i/E_{i-1})$ satisfying
\[\mu_{max}(E)=\mu_1>\ldots>\mu_n=\mu_{min}(E).\]
\nomenclature{$\mu_{min}(E)$}{the minimal slope of the subsequent quotients of the Harder-Narasimhan filtration of a coherent sheaf $E$}
\nomenclature{$\mu_{max}(E)$}{the maximal slope of the subsequent quotients of the Harder-Narasimhan filtration of a coherent sheaf $E$}
\end{definition}

\begin{lemma}[\cite{huy} Lemma 1.3.5] 
\label{lem_destabilizing}
Let $E$ be a vector bundle. 
Then there is a subsheaf $F\subset E$ such that for all $G\subset E$ on has $\mu(F)\geq \mu(G)$ and, if equality occurs, then $G\subset F$.
$F$ is uniquely determined by these properties, is semi-stable, saturated (hence the quotient is a vector bundle) and is called the maximal destabilizing subsheaf\index{subsheaf!maximal destabilizing} of $E$.
\end{lemma}

By the previous Lemma one can easily derive the existence of a Harder-Narasimhan filtration.
One takes as $E_1$ the maximal destabilizing subsheaf of $E$ and then proceeds by induction on the quotient $E/E_1$.

\begin{theorem}[\cite{huy} Theorem 1.3.4]
\label{propo_HN}
Every vector bundle has a unique Harder-Narasimhan filtration.
\end{theorem}

\begin{corollary}
\label{coro_hom0}
If $F$ and $G$ are vector bundles with $\mu_{min}(F)>\mu_{max}(G)$, then $\Hom(F,G)=0$.
\end{corollary}

\begin{proof}
Suppose that $\phi\colon F \to G$ is a non-trivial morphism and let $0=F_0\subset F_1\subset\ldots\subset F_n=F$ and $0=G_0\subset G_1\subset\ldots\subset G_m=G$ be their Harder-Narasimhan filtrations.
Let $i$ and $j$ be the smallest indexes such that $\restr{\phi}{F_i}$ is not identically trivial and  $\restr{\phi}{F_i}(F_i)\subset G_j$.
This in particular gives rise to a non-trivial morphism $\widetilde{\phi}\colon F_i/F_{i-1}\to G_j/G_{j-1}$ where both $F_i/F_{i-1}$ and $G_j/G_{j-1}$ are  semi-stable vector bundles and $\mu(F_i/F_{i-1})\geq\mu_{min}(F)>\mu_{max}(G)\geq\mu(G_j/G_{j-1})$.
Let $E$ be the image of $\widetilde{\phi}$: by Lemma \ref{lemma_equivalentss} we have that $\mu(G_j/G_{j-1})\geq\mu(E)\geq\mu(F_i/F_{i-1})$ which is a contradiction.
\end{proof}

\begin{remark}
\label{rem_minmaxss}
It is possible to give another construction of the Harder-Narasimhan filtration, starting from the top instead of from the bottom.
Indeed it is possible to define the minimal destabilizing quotient\index{minimal destabilizing quotient sheaf} of $E$ as a quotient $G$ of $E$ such that $\mu(G)\leq\mu(F)$ for every other quotient $F$ of $E$ and equality occurs if and only if $F$ is a quotient of $G$ too.
Notice that $G$ turns out to be a semi-stable vector bundle.
In this case, one constructs the Harder-Narasimhan filtration taking $E_{n-1}=\ker(E\to G)$ and proceeding by induction on $E_{n-1}$.
The fact that one obtains the same filtration is ensured by the uniqueness part of Theorem \ref{propo_HN}.

It is then clear that $\mu_{min}(E)$ and $\mu_{max}(E)$ can be defined as the minimal slope of a quotient sheaf of $E$ and the maximal slope of a subsheaf of $E$.
In particular, for every subsheaf $F\subseteq E$ we have $\mu_{max}(F)\leq\mu_{max}(E)$ and, for every quotient sheaf $E\twoheadrightarrow G$ we have $\mu_{min}(E)\leq\mu_{min}(G)$.
\end{remark}

\begin{lemma}
\label{lemma_deghn}
Let $E$ be a vector bundle and $0=E_0\subset E_1 \subset\ldots\subset E_n=E$ be its Harder-Narasimhan filtration with slopes $\mu_1>\mu_2>\ldots>\mu_n$ and let $r_i=\rk E_i$.
Then
\begin{equation}
\deg(E)=\sum_{i=1}^nr_i(\mu_i-\mu_{i+1})
\label{eq_deghn}
\end{equation}
where, by definition, $\mu_{n+1}=0$.
\end{lemma}

\begin{proof}
For $i=2,\ldots,n$ we have the short exact sequence
\begin{equation}
0\to E_{i-1}\to E_i\to E_i/E_{i-1}\to 0,
\label{eq_shortdeg}
\end{equation}
 which gives $\deg(E_i)-\deg(E_{i-1})=\deg(E_i/E_{i-1})=\mu_i(r_i-r_{i-1})$.
By induction we obtain
\begin{equation}
\begin{split}
\deg&(E_n)=\\
=\deg(E_1)+(\deg(E_2)-\deg(E_1))+&\ldots+(\deg(E_n)-\deg(E_{n-1}))=\\
\mu_1r_1+\mu_2(r_2-r_1)+\ldots+\mu_n&(r_n-r_{n-1})=\sum_{i=1}^n r_i(\mu_i-\mu_{i+1})
\end{split}
\label{eq_degcan}
\end{equation}
and the Lemma is proven.
\end{proof}

\begin{corollary}
\label{coro_belloneq}
Let $E$ be a vector bundle and $0=E_0\subset E_1 \subset\ldots\subset E_n=E$ be its Harder-Narasimhan filtration with slopes $\mu_1>\mu_2>\ldots>\mu_n$ and let $r_i=\rk E_i$.
Let $0=F_0\subset F_1\subset\ldots\subset F_m=E$ be another filtration of $E$ and let $\widetilde{r}_j$ and $\widetilde{\mu}_j$ numbers such that $\mu_{min}(F_j)\geq\widetilde{\mu}_j$ and $\rk(F_j)\geq\widetilde{r}_j$ for all $j\in\{1,\ldots,m\}$ and $\widetilde{r}_m=\rk(F_m)=r_n=\rk(E)$.
Then we have
\begin{equation}
\deg(E)\geq\sum_{i=1}^m\widetilde{r}_j(\widetilde{\mu}_j-\widetilde{\mu}_{j+1})
\label{eq:}
\end{equation}  
\end{corollary}

\begin{proof}
For every $j\in\{1,\ldots,m\}$ there exists a unique $i\in\{1,\ldots,n\}$, call it $i_j$, such that $F_j\subseteq E_{i_j}$ and $F_j\nsubseteq E_{i_j-1}$: in particular we have a non-trivial morphism $F_j\to E_{i_j}/E_{i_j-1}$, which, thanks to Corollary \ref{coro_hom0}, implies that $\widetilde{\mu}_j\leq\mu_{min}(F_j)\leq \mu_{max}(E_{i_j}/E_{i_j-1})=\mu_{i_j}$.
We also clearly have that $\widetilde{r}_j\leq\rk(F_j)\leq r_{i_j}$.

\begin{figure}%
\centering
\begin{tikzpicture}[xscale=1, yscale=1, every node/.style={scale=1}]
	\begin{axis}[ticks=none, xmin=-5, xmax=5, ymin=0, ymax=7, axis lines=middle, xlabel=$\mu$, ylabel=$r$, enlargelimits]
	
	\draw[white!0] (0,0) -- (0,7);

	\fill [gray!20] (4,0) -- (4,1) -- (2,1) -- (2,0) -- (4,0);
	\fill [gray!20] (3,1) -- (3,3) -- (1,3) -- (1,1) -- (3,1);
	\fill [gray!20] (1,2) -- (1,4) -- (-2,4) -- (-2,2) -- (1,2);
	\fill [gray!20] (-2,3) -- (-2,4) -- (-3,4) -- (-3,3) -- (-2,3);
	\fill [gray!20] (-1,4) -- (-1,6) -- (-4,6) -- (-4,4) -- (-1,4);

	\draw (4,0) -- (4,1);
	
	\draw (4,1) -- (3,1);
	
	\draw (3,1) -- (3,3);
	
	\draw[dashed] (3,3) -- (1,3);
	
	\draw[dashed] (1,3) -- (1,4);
	
	\draw (1,4) -- (-1,4);
	
	\draw (-1,4) -- (-1,6);
	
	\draw (-1,6) -- (-4,6);
	
	\filldraw  (4,1) circle (1pt);
	\filldraw  (3,3) circle (1pt);
	\filldraw  (1,4) circle (1pt);
	\filldraw  (-1,6) circle (1pt);
	
	\node at (-0.9,6.4) {$(\mu_n,r_n)$};
	\node at (4.1,1.4) {$(\mu_1,r_1)$};
	\node at (1.1,4.4) {$(\mu_{n-1},r_{n-1})$};
	\node at (3.1,3.4) {$(\mu_2,r_2)$};

	\draw (2,0) -- (2,1);
	
	\draw (2,1) -- (1,1);
	
	\draw (1,1) -- (1,2);
	
	\draw[dashed] (1,2) -- (-2,2);
	
	\draw[dashed] (-2,2) -- (-2,3);
	
	\draw (-2,3) -- (-3,3);
	
	\draw (-3,3) -- (-3,4);
	
	\draw (-3,4) -- (-4,4);
	
	\draw (-4,4) -- (-4,6);
	
	\filldraw  (2,1) circle (1pt);
	\filldraw  (1,2) circle (1pt);
	\filldraw  (-2,3) circle (1pt);
	\filldraw  (-3,4) circle (1pt);
	\filldraw  (-4,6) circle (1pt);
	
	\node at (1,0.6) {$(\widetilde{\mu}_1,\widetilde{r}_1)$};
	\node at (0,1.6) {$(\widetilde{\mu}_2,\widetilde{r}_2)$};
	\node at (-4,2.6) {$(\widetilde{\mu}_{m-2},\widetilde{r}_{m-2})$};
	\node at (-4,3.6) {$(\widetilde{\mu}_{m-1},\widetilde{r}_{m-1})$};
	\node at (-4,6.4) {$(\widetilde{\mu}_m,\widetilde{r}_m)$};

	\end{axis}
\end{tikzpicture}
\caption{}%
\label{fig_hn}%
\end{figure}

Viewing the points $(\mu_i,r_i)$ and $(\widetilde{\mu}_j,\widetilde{r}_j)$ in a plane as in Figure \ref{fig_hn} we see that 
\begin{equation}
\deg(E)-\sum_{i=1}^m\widetilde{r}_j(\widetilde{\mu}_j-\widetilde{\mu}_{j+1})=\sum_{i=1}^nr_i(\mu_i-\mu_{i+1})-\sum_{i=1}^m\widetilde{r}_j(\widetilde{\mu}_j-\widetilde{\mu}_{j+1})
\end{equation}
is the area of the grey part, hence the Corollary follows.
\end{proof}

\begin{lemma}[cf. \cite{miyaoka} Proposition 3.2] 
\label{lemma_finitestable}
Let $f\colon C'\to C$ be a separable finite morphism of smooth projective curves. 
Then a vector bundle $E$ on $C$ is semi-stable if and only if $f^*E$ is.
In particular the pull-back of the Harder-Narasimhan filtration of $E$ is the Harder-Narasimhan filtration of $f^*E$.
\end{lemma}

\begin{proof}
Because the slope is multiplicative by the degree of $f$ under pull-back via a finite morphism $f$ and the set of subsheaves of $E$ is a subset of the set of subsheaves of $f^*E$, it is clear that if $f^*E$ is semi-stable, then also $E$ is.
Conversely, let $E$ be a semi-stable vector bundle on $X$.
Thanks to the "if" part of the lemma, we may pass to the Galois closure of the finite morphism and assume that $f$ itself is Galois with Galois group $G$.
It is clear that $f^*E$ is fixed by the action of $G$ and that, by uniqueness, its maximal destabilizing subsheaf of it, call it $\overline{F}$, is fixed by the action of $G$ too.
In particular there exists a sheaf $F$ on $C$ with $f^*F=\overline{F}$ and, by the semi-stability assumption we have that $F=E$, from which the assertion follows. 
\end{proof}

\begin{remark}
The "only if" part of this Lemma is not true for every finite morphism in positive characteristic and it is known that every curve of genus greater than or equal to $2$ has at least one counterexample (cf. \cite{langeart} Theorem 1).

If $C$ is the projective line, the statement of Lemma \ref{lemma_finitestable} holds true for every finite morphism, even if it is not separable.
Indeed, in this case we can split $f$ as $g\circ F_k^n$ where $g$ is separable and $F_k^n$ is the $n^{th}$ power of the $k$-linear Frobenius morphism\index{Frobenius morphism!k-linear @ $k$-linear}.
Hence it is enough to prove that $E$ is semi-stable if and only if its pull-back via the $k$-linear Frobenius morphism is.
But this is obvious because over the projective line every vector bundle splits to a sum of line bundles and every line bundle is uniquely determined by its degree.

In the general case in positive characteristic, one can always reduce an inseparable finite morphism to a composition of a power of the $k$-linear Frobenius with a separable morphism (cf. \cite{liu} Proposition 7.4.21 and the discussion preceding it).
Hence Lemma \ref{lemma_finitestable} holds true also in positive characteristic, provided that the pull-back of a semi-stable vector bundle via the $k$-Frobenius morphism is semi-stable: this turns out to be true for varieties with globally generated tangent bundle (cf. \cite{mehta} Theorem 2.1 and Remark 2.2) e.g. rational and elliptic curves.
\label{rem_finitestable}
\end{remark}

\begin{definition}
Let $E$ be a vector bundle on a curve $C$. We say that $E$ is strongly semi-stable\index{vector bundle!strongly semi-stable}, if it is semi-stable after pulling it back via finite morphisms.
\label{def_stronglyss}
\end{definition}

Thanks to Lemma \ref{lemma_finitestable}, in characteristic zero every semi-stable vector bundle is also strongly semi-stable.
As pointed out in Remark \ref{rem_finitestable}, on the projective line every semi-stable vector bundle is strongly semi-stable even in positive characteristic.

In this setting, we can rephrase the second part of Lemma \ref{lemma_finitestable} as follows.

\begin{lemma}
\label{lemma_HNsss}
Let $f\colon C'\to C$ be a finite morphism of smooth curves. 
Then if every quotient of the Harder-Narasimhan filtration of $E$ is strongly semi-stable, we have that the pull-back of the Harder-Narasimhan filtration of $E$ is the Harder-Narasimhan filtration of $f^*E$.
In particular this holds if the characteristic is $0$, if $C=\P^1$ or if $C$ is an elliptic curve.
\end{lemma}

%

\begin{proposition}[\cite{langer} Theorem 6.1]
\label{propo_tensorss}
The tensor product of two strongly semi-stable vector bundles on a smooth curve $C$ is strongly semi-stable. In particular the tensor product of two semi-stable vector bundles is semi-stable if the characteristic of the ground field is zero, if $C=\P^1$ or if $C$ is elliptic.
\end{proposition}

\begin{corollary}
\label{coro_tensorss}
Let $E$ and $F$ be two vector bundles on a smooth curve $C$.
If all the successive quotients of their Harder-Narasimhan filtrations are strongly semi-stable, then   $\mu_{min}(E\otimes F)=\mu_{min}(E)+\mu_{min}(F)$.
\end{corollary}

\begin{proof}
First of all, we notice that the slope is additive on tensor products, indeed 
\begin{equation}
\begin{split}
\mu(A\otimes B)=\deg(A&\otimes B)/\rk(A\otimes B)=\\
=(\rk(B)\deg(A)+\rk(A)\deg(B&))/(\rk(A)\rk(B))=\mu(A)+\mu(B).
\end{split}
\end{equation}

Let $E'$ and $F'$ be the minimal destabilizing quotients of $E$ and $F$: we know that (cf. Remark \ref{rem_minmaxss}) $\mu_{min}(E)=\mu(E')$ and $\mu_{min}(F)=\mu(F')$ and that (cf. Proposition \ref{propo_tensorss}) $E'\otimes F'$ is semi-stable.
Because $E'\otimes F'$ is a quotient of $E\otimes F$ we have that
\begin{equation}
\mu_{min}(E\otimes F)\leq\mu(E'\otimes F')=\mu_{min}(E)+\mu_{min}(F).
\label{eq_mumin}
\end{equation}

Now suppose that $E$ is semi-stable and that $0\subseteq F_1\subseteq F_2\subseteq\ldots\subseteq F_{n-1}\subseteq F$ is the Harder-Narasimhan filtration for $F$ and consider the induced filtration 
\begin{equation}
0\subseteq E\otimes F_1\subseteq E\otimes F_2\subseteq\ldots\subseteq E\otimes F_{n-1}\subseteq E\otimes F.
\label{eq_otimeshn}
\end{equation}
The successive quotients of this filtration are $(E\otimes F_i)/(E\otimes F_{i-1})=E\otimes F_i/F_{i-1}$: they are semi-stable by Proposition \ref{propo_tensorss} and their slope is $\mu(E\otimes F_i/F_{i-1})=\mu(E)+\mu_i(F)$.
Let $A$ be a semi-stable quotient vector bundle of $E\otimes F$ and let $i$ be such that the induced morphism $E\otimes F_{i-1}\to A$ is trivial while $E\otimes F_{i}\to A$ is not: by Corollary \ref{coro_hom0} we have that $\mu(A)\geq\mu(E)+\mu_i(F)\geq\mu(E)+\mu_{min}(F)$.
In particular, because the minimal destabilizing quotient has to be semi-stable, $\mu_{min}(E\otimes F)\geq \mu(E)+\mu_{min}(F)$, which, combined with Equation \ref{eq_mumin}, concludes the proof in this case.

Now suppose that both $E$ and $F$ are not semi-stable and consider, as above, the Harder-Narasimhan filtration $0\subseteq F_1\subseteq F_2\subseteq\ldots\subseteq F_{n-1}\subseteq F$ of $F$.
Let $A$ be a semi-stable quotient of $E\otimes F$ and let, as above,  $i$ be an index such that the induced morphism $E\otimes F_{i-1}\to A$ is trivial while $E\otimes F_i\to A$ is not: by Corollary \ref{coro_hom0} we have that $\mu(A)\geq\mu_{min}(E\otimes F_i/F_{i-1})=\mu_{min}(E)+\mu_i(F)\geq\mu_{min}(E)+\mu_{min}(F)$, where the equality follows by the previous case.
Combining this with Equation \ref{eq_mumin}, we finish the proof.
\end{proof}


\section{Projective bundles}
\label{sec_projbundle}

In this section $X$ is a smooth variety over an algebraically closed field.

\begin{definition}
\label{def_projbundle}
Let $E$ be a vector bundle on $X$, we define the projective bundle associated with $E$ to be $\P(E)=\proju(\bigoplus_{i=1}^\infty S^iE)$\index{projective bundle} \nomenclature{$\P(E)$}{the projective bundle associated with a vector bundle $E$} where $\bigoplus_{i=1}^\infty S^iE$ is the symmetric algebra of $E$ and $\proju$\nomenclature{$\proju$}{the relative $\proj$ functor} denotes the relative $\proj$ functor.
\end{definition}

We recall some properties of projective bundles (cf. \cite{hart} II.7). 
If $E$ is a vector bundle of rank $1$ then $\P(E)=X$.
If the rank $r$ of $E$ is greater than $1$ we have that $\P(E)$ has a natural proper projection $\pi$ to $X$ which makes it into a fibre bundle with fibre isomorphic to $\P^{r-1}$ and is naturally endowed with a line bundle $\O_{\P(E)}(1)$ which restricted to the fibre is the usual ample line bundle of a projective space. 
Moreover we have a natural surjection of sheaves $\pi^*E\to \O_{\P(E)}(1)$ and a canonical isomorphism of graded $\O_X$ algebras $\bigoplus_{i\in\Z}\pi_*\O_{\P(E)}(i)=\bigoplus_{i=1}^\infty S^iE$.
The canonical bundle of $\P(E)$\index{projective bundle!canonical divisor of}\index{canonical divisor!of a projective bundle|see{projective bundle}} is given by the following formula (cf. ibidem Exercise III.8.4(b));
\begin{equation}
\label{eq_canprojbundle}
K_{\P(E)}=\pi^*(K_X\otimes_{\O_X}\det(E))\otimes_{\O_{\P(E)}}\O_{\P(E)}(-r).
\end{equation}

Recall also that, given two different vector bundles $E$ and $E'$, we have that $\P(E)\simeq\P(E')$ if and only if there exists a line bundle $L$ such that $E'\simeq E\otimes L$ (cf. ibidem Lemma II.7.9).
Moreover, if this happens,  we have
\begin{equation}
\O_{\P(E')}(1)=\phi^*\O_{\P(E)}(1)\otimes\pi'^*L
\label{eq_diffprojbundle}
\end{equation}
where $\phi\colon\P(E')\to\P(E)$ is the natural isomorphism and $\pi'\colon\P(E')\to X$ is the natural projection.
In addition, the Picard group $\Pic(\P(E))$ of $\P(E)$ is isomorphic to
\begin{equation}
\Pic(\P(E))=\Pic(X)\oplus \Z
\label{eq_picprojbund}
\end{equation}
where $\Z$ is generated by $\O_{\P(E)}(1)$ and the inclusion $\Pic(X)\hookrightarrow\Pic(\P(E))$ is given by pull-back via $\pi$ (cf. ibidem Exercise II.7.9).
Furthermore we have the following formula for the higher pushforward of multiples of $\O_{\P(E)}(1)$ (cf. ibidem Exercise III.8.4(a)(c)) :
\begin{equation}
\begin{split}
&\pi_*\O_{\P(E)}(l)=S^lE\quad\text{for $l\geq 0$}\\
&\pi_*\O_{\P(E)}(l)=0\quad\text{for $l<0$}\\
&R^i\pi_*\O_{\P(E)}(l)=0\quad \text{for every $l\in\Z$ and $0<i<r-1$}\\
&R^{r-1}\pi_*\O_{\P(E)}(l)=\pi_*(\O_{\P(E)}(-l-r))^\vee\otimes_{\O_X}\det(E)^{-1} \quad \text{for every $l\in\Z$.}\hspace{-20pt}
\end{split}
\label{eq_higherpushprojbund}
\end{equation}

\begin{proposition}[cf. \cite{hart} Proposition II.7.12]
\label{propo_morphtoproj}
Let $X$, $E$ and $\P(E)$ as above and $g\colon Y\to X$ any morphism. 
There is a bijective correspondence between rational maps $f\colon Y\to \P(E)$ such that $\pi\circ f=g$ and morphisms of sheaves $g^*E\to L$  surjective in codimension $1$ where $L$ is a line bundle.
Moreover $f$ is an actual morphism if and only if $g^*E\to L$ is surjective.
\end{proposition}

\begin{proof}[Sketch of proof]
The equivalence is given as follows.
For a rational map $f\colon Y\to\P(E)$ we let $L=f^*\O_{\P(E)}(1)$ and clearly $f^*\pi^*E=g^*E\to f^*\O_{\P(E)}(1)$ is surjective where $f$ is defined  because the pull-back of sheaves is right exact.
Recall that a rational map is always defined in codimension one, hence the morphism of sheaves is surjective in codimension one.

Conversely, let $g^*E\to L$ be a morphism   surjective in codimension $1$. The morphism $f\colon Y\to\P(E)$ is constructed locally where the vector bundle $g^*E$ trivializes using the well known equivalence between rational map to projective spaces and linear systems without fixed components.
\end{proof}

The next Theorem is a generalization of Theorem 3.1 of \cite{miyaoka} in positive characteristic.

\begin{theorem}
\label{teo_nefantican}
Let $\pi\colon\P(E)\to C$ a projective bundle over a smooth curve $C$. 
Suppose that all the quotients of the Harder-Narasimhan filtration are strongly semi-stable (e.g. $\cha(k)=0$ or $C=\P^1_k$). 
Then the rational divisor $\O_{\P(E)}(1)-\mu_{min}(E)F$ is nef\index{divisor!nef}, i.e. it intersects non-negatively every curve of $\P(E)$, where $F$ is a general fibre of $\pi$.
\end{theorem}

\begin{proof}
If $\rk(E)=1$ we have that $\P(E)=C$, $\O_{\P(E)}(1)=E$ and $\O_{\P(E)}(1)-\mu_{min}(E)F$ is a divisor of degree $\deg(E)-\mu_{min}(E)=\mu(E)-\mu_{min}(E)= 0$.

Now let $\rk(E)\geq 2$.
Suppose that the thesis is not true: then there exists an irreducible curve $C'\subseteq\P(E)$ such that $C'.(\O_{\P(E)}(1)-\mu_{min}(E)F)<0$.
We first show that we may assume that $\restr{\pi}{C'}$ is an isomorphism.

Indeed, clearly $C'\to C$ has to be finite (otherwise $C'$ would be contained in a fibre and $C'.(\O_{\P(E)}(1)-\mu_{min}(E)F)$ would be non-negative).
Moreover let $C''$ be the normalization of $C'$, $f\colon C''\to C$ the induced morphism and consider the following Cartesian diagram:
\begin{equation}
\begin{tikzcd}
C''\times_CC'\arrow{d}{i''}\arrow{r}{} & C'\arrow{d}{i}\\
\P(f^*E)\arrow{d}{\pi''}\arrow{r}{\P(f)} & \P(E)\arrow{d}{\pi}\\
C''\arrow{r}{f} & C.
\end{tikzcd}
\label{eq_basechangess}
\end{equation}
It is clear that $C''\times_CC'$ has at least one component isomorphic to $C''$ and the projection restricted to all the other components has degree strictly smaller than  the degree of $f$.
Combining this fact with  Lemma \ref{lemma_HNsss}  and proceeding by induction, it is clear that we may assume that $C'$, and hence $C''$, is isomorphic to $C$.

Proposition \ref{propo_morphtoproj} ensures that the inclusion $i\colon C\to\P(E)$ is uniquely determined by a line bundle $L=i^*\O_{\P(E)}(1)$ and a surjective morphism $E\to L$.
Because the intersection $C.(\O_{\P(E)}(1)-\mu_{min}(E)F)$ is negative, we have that the degree of $L$, which is equal to $C.\O_{\P(E)}(1)$, is strictly smaller than $\mu_{min}(E)$.
On the other hand Corollary \ref{coro_hom0} implies that the degree of $L$ cannot be smaller than $\mu_{min}(E)$, which is a contradiction.
\end{proof}

We also recall here the notion of nefness for vector bundles, which we will need later. 

\begin{definition}[\cite{lazarsfeld} 6.1.1]
\label{def_nefvb}
We say that a vector bundle $E$ is a nef\index{vector bundle!nef} if the line bundle $\O_{\P(E)}(1)$ is nef line bundle on $\P(E)$.
\end{definition}

Notice that, if the rank of $E$ is $1$, this definition is just saying that $E$ is nef as a vector bundle if $E$ is nef as a line bundle.

\begin{proposition}[\cite{lazarsfeld} 6.1.2]
\label{propo_quotnef}
Let $E$ be a vector bundle on a projective scheme $X$.
\begin{enumerate}
	\item If $E$ is nef, then any quotient bundle $G$ of $E$ is nef too.
	\item The pull-back of a nef vector bundle via a finite morphism  is nef.
\end{enumerate}
\end{proposition}

\begin{proof}
Let $E\to G$ a surjective morphism of vector bundles: this induces a closed immersion $\P(G)\hookrightarrow\P(E)$ which satisfies $\restr{\O_{\P(E)}(1)}{\P(G)}=\O_{\P(G)}(1)$. 
Then $1$ follows by the fact that the restriction of a nef line bundle is nef.

A finite morphism $f\colon X\to Y$ functorially gives a morphism $\P(f)\colon \P(f^*E)\to\P(E)$ such that $\O_{\P(f^*E)}(1)=\P(f)^*\O_{\P(E)}(1)$.
Then 2 follows by the same property for line bundles.
\end{proof}

\begin{proposition}
\label{propo_nef>}
If $E$ is a nef vector bundle on a curve $C$, then $\mu_{min}(E)\geq0$.
\end{proposition}

\begin{proof}
We have that $E_n/E_{n-1}$ is nef thanks to Proposition \ref{propo_quotnef}.
Moreover any exterior power of a nef vector bundle is nef (cf. \cite{lazarsfeld} Theorem 6.2.12(IV) and Remark 6.2.16).
Clearly the degree of a nef line bundle is non-negative.
In particular we have that
\begin{equation*}
0\leq\deg(\det(E_n/E_{n-1}))=\deg(E_n/E_{n-1})=\mu(E_n/E_{n-1})\rk(E_n/E_{n-1})
\end{equation*}
from which the thesis easily follows.
\end{proof}

\section{Algebraic fundamental group}
\label{sec_fundgroup}
In this Section we recall some general facts about the algebraic fundamental group\index{algebraic fundamental group}.
All the schemes in this section are defined over an algebraically closed field $k$ for simplicity, even though many of the results here reported are known to be true in more generality.

Let $S$ be a scheme over an algebraically closed field $k$: we denote by $\Fet_S$\nomenclature{$\Fet_S$}{the category of finite \'etale covers of $S$} the category whose objects are finite \'etale covers of $S$ and whose morphisms are  morphisms of $S$-schemes. 
We fix a point $s\colon\spec(k)\to S$ and, for every $X\in\Fet_S$, we consider the fibre $\spec(k)\times_SX$ and its underlying set $\Fib_s(X)$\nomenclature{$\Fib_s(X)$}{the underlying set of the fibre of a morphism $X\to S$ over a closed point $s\in S$}.
Given a morphism $f\colon X\to Y$ in $\Fet_S$, there is a naturally induced morphism $\spec(k)\times_SX\to\spec(k)\times_SY$, hence a set-theoretic map $\Fib_s(X)\to\Fib_s(Y)$ so that $\Fib_s$ is a set valued functor on $\Fet_S$.

\begin{definition}
Given a scheme $S$ and a closed point $s\colon \spec(k)\to S$, we define the algebraic fundamental group $\pi_1(S,s)$ with base point $s$\nomenclature{$\pi_1(S,s)$}{the algebraic fundamental group of a scheme $S$ with base point $s$} as the group of automorphism of the functor $\Fib_s$.
\label{def_fundgroup}
\end{definition}

It can be proved (cf. \cite{szamuely} Proposition 5.5.1) that the algebraic fundamental group is independent of the choice of  base point.
Notice that the algebraic fundamental group $\pi_1(S,s)$ has, by definition, a natural action on $\Fib_s(X)$ for every $X\in\Fet_S$; in particular we have that the functor $\Fib_s$ actually takes value inside the category of $\pi_1(S,s)$-sets, i.e. sets with an action of $\pi_1(S,s)$ on them.

\begin{theorem}[cf. \cite{szamuely} Theorem 5.4.2]
Let $S$ be a connected scheme and $s\colon \spec(k)\to S$ a closed point. 
Then the group $\pi_1(S,s)$ is profinite and its action on $\Fib_s(X)$ is continuous for every $X\in\Fet_S$.
Moreover the functor $\Fib_s$ induces an equivalence of $\Fet_S$ with the category of finite continuous $\pi_1(S,s)$-sets where connected covers correspond to sets with transitive $\pi_1(S,s)$- action and Galois covers to finite quotients of $\pi_1(S,s)$.
\end{theorem}

Notice in particular that the above Theorem implies that two schemes have the same algebraic fundamental group if and only if they have a bijective correspondence between finite \'etale covers.

Now consider a morphism $\phi\colon S'\to S$ of schemes and let $s'\colon \spec(k)\to S'$ be a closed point of $S'$ and $s=\phi\circ s'$.
Then, it is naturally defined a functor $\phi^*\colon \Fet_{S}\to\Fet_{S'}$ which sends an object $X\in\Fet_S$ to $\phi^*(X)=X\times_SS'\in\Fet_{S'}$ and a morphism $f\colon X\to Y$ to $\phi^*(f)\colon X\times_SS'\to Y\times_SS'$.
Moreover $\phi\circ s'=s$ implies that $\Fib_{s}=\Fib_{s'}\circ \phi^*$: consequently, every automorphism of the functor $\Fib_{s'}$ induces and automorphism of the functor $\Fib_{s}$ via post-composition with $\phi^*$ which allow us to define a continuous morphism of profinite groups $\phi_*\colon \pi_1(S',s')\to\pi_1(S,s)$.

\begin{theorem}[\cite{sga1} Corollaire 6.10 expos\'e V]
The morphism $\phi_*$ defined above is an isomorphism if and only $\phi^*\colon\Fet_S\to\Fet_{S'}$ is an equivalence of categories.
Equivalently, it is an isomorphism if and only if the two following conditions are satisfied:
\begin{enumerate}
	\item for all $X\in\Fet_S$ connected,  $\phi^*X$ is connected too;
	\item for all $X'\in\Fet_{S'}$ there exists an $X\in\Fet_S$ such that $\phi^*X\simeq X'$.
\end{enumerate}
\label{teo_fundpushfor}
\end{theorem} 

We have the following version of the Lefschetz hyperplane Theorem for the algebraic fundamental group.

\begin{theorem}[cf. \cite{sga2} Corollaire 3.5 expos\'e XII]
Let $S$ be a smooth projective scheme of dimension at least $3$ with an ample line bundle $L$ on it, $t$ be a section of $L$, $S_0$ be the divisor associated with $t$, $s\in S_0\subseteq S$ be a closed point of $S$ and denote by $i\colon S_0\to S$ the inclusion.
The the induced morphism $i_*\colon\pi_1(S_0,s)\to\pi_1(S,s)$ is an isomorphism.
\label{teo_lefhypsec}
\end{theorem}

We also recall the homotopy exact sequence for flat morphisms.

\begin{proposition}[cf. \cite{szamuely} Proposition 5.6.4 and \cite{sga1} Corollaire 1.4 expos\'e X]
Let $S$ be a Noetherian connected scheme and let $f\colon X\to S$ be a proper flat morphism with integral fibres.
Let $x\colon\spec(k)\to X$ be a closed point of $X$ and let $s=f\circ x$.
then the sequence
\begin{equation}
\pi_1(X_s,x)\to\pi_1(X,x)\to\pi_1(S,s)\to0
\label{eq_flatfundgroup}
\end{equation}
is exact, where $X_s$ is the fibre of $f$ over $s$.
\label{propo_flatfundgroup}
\end{proposition}

As a consequence we have the following Corollary.

\begin{corollary}
Let $f\colon X\to S$ be as in Proposition \ref{propo_flatfundgroup} and suppose that all the fibres are isomorphic to $\P^n$ (i.e. $f\colon X\to S$ is a projective bundle).
Then $f$ induces an isomorphism on the algebraic fundamental groups\index{algebraic fundamental group!of a projective bundle}\index{projective bundle!algebraic fundamental group of |see{algebraic fundamental group}} of $X$ and $S$.
\label{coro_projbundlefund}
\end{corollary} 
\begin{proof}
Let $x\colon\spec(k)\to X$ be a closed point and $s=f\circ x$.
By Proposition \ref{propo_flatfundgroup} we have the following exact sequence:
\begin{equation}
\pi_1(\P^n,x)\to \pi_1(X,x)\to\pi_1(S,s)\to0.
\label{eq_exseqprojbundfund}
\end{equation}
Because the algebraic fundamental group of a projective space is trivial (cf. \cite{szamuely} Exercise 9 Chapter 5), the thesis then easily follows. 
\end{proof}

Before the next Proposition, we recall what is a universal homeomorphism\index{universal homeomorphism}.

\begin{definition}
Let $f\colon X\to Y$ be a morphism of schemes.
Then $f$ is said to be a universal homeomorphism if for all $Y'\to Y$ the induced morphism $f_{Y'}\colon X\times_YY'\to Y'$ is a homeomorphism.
\label{def_homeouniv}
\end{definition}

\begin{remark}
\label{rem_homeouniv}
Examples of universal homeomorphisms include the canonical morphism $X_{red}\to X$ for every scheme $X$ (cf. \cite{egaiv2} Remarques 2.4.8(ii)) and factors of the $k$-linear Frobenius morphism (cf. \cite{egaiv4} Corollaire 18.12.11).
\end{remark}

\begin{proposition}
\label{propo_fundgrouphomeo}
Let $f\colon X\to Y$ be a universal homeomorphism between connected schemes, $x\colon \spec(k)\to X$ be a closed point of $X$ and $y=f\circ x$.
Then $f_*\colon\pi_1(X,x)\to\pi_1(Y,y)$ is an isomorphism.
\end{proposition}
\begin{proof}
By \cite{sga1} Th\'eor\`eme 4.10 Expos\'e IX, we know that $f^*\colon \Fet_Y\to\Fet_X$ is an equivalence of categories: in particular, by Theorem \ref{teo_fundpushfor}, we have the thesis.
\end{proof}

\section{Generalities on surfaces}
\label{sec_surfaces}
In this section we collect some known results for surfaces, where by surface, we mean a smooth projective scheme of dimension two.
The main references are \cite{beauville_1983} for complex surfaces, \cite{badescu} for the general case and also \cite{hart} chapter V.
Another good reference for the classification of surfaces is \cite{liedtke+}.

Let $X$ be a surface\index{surface} and denote by $K_X$ its canonical divisor. 
Recall that on a smooth surface it is defined a bilinear product on the Picard group $\Pic(X)$\nomenclature{$\Pic(X)$}{the Picard group of $X$}\index{Picard! group} which, if $L$ and $L'$ are line bundles and $D$ and $D'$ their corresponding divisors, we denote interchangeably by $L.L'$ or $D.D'$.
Recall that for a line bundle $L$ on $X$, the Euler characteristic\index{sheaf!Euler characteristic of} $\chi(L)$\nomenclature{$\chi(F)$}{the Euler characteristic of a coherent sheaf $F$} of $L$ is defined as the alternating sum of the dimensions of the cohomology groups of $L$, i.e.
\begin{equation}
\chi(L)=h^0(X,L)-h^1(X,L)+h^2(X,L)
\label{eq_eulcharl}
\end{equation}
\nomenclature{$H^i(X,F)$}{the $i$-th cohomology group of the sheaf $F$ on $X$}
\nomenclature{$h^i(X,F)$}{the dimension as a $k$-vector space of $H^i(X,F)$}
in particular $\chi(\O_X)=1-q'+p_g$ where $q'=h^1(X,\O_X)$ is the irregularity\index{surface!irregularity of} and $p_g=h^2(X,\O_X)=h^0(X,K_X)$ is the geometric genus\index{surface!geometric genus of}\nomenclature{$p_g(X)$}{the geometric genus of a surface $X$, i.e. $p_g(X)=h^2(X,\O_X)$} and the last equality follows by Serre duality\index{Serre duality}.
We have the Riemann-Roch Theorem for surfaces\index{surface!Riemann-Roch Theorem for}\index{Riemann-Roch Theorem!for surfaces|see{surface}}

\begin{theorem}[\cite{hart} Theorem V.1.6]
Let $X$ be a surface. 
Then, for every line bundle $L$ on $X$, we have
\begin{equation}
\chi(L)=\chi(\O_X)+\frac{1}{2}L.(L-K_X).
\label{eq_rrss}
\end{equation} 
\label{teo_rrss}
\end{theorem}

Recall that the Betti numbers\index{Betti number} $b^i(X)$\nomenclature{$b^i(X)$}{the $i$-th (singular or $l$-adic) Betti number of $X$} are defined as the dimensions of the $i$-th $l$-adic cohomology groups of $X$ where $l$ is a prime different from the characteristic of the ground field and it is a topological invariant.
If $X$ is a complex variety, it can be shown that $b^i(X)$ is equal to dimension of the $i$-th singular cohomology group.
The topological Euler characteristic\index{topological Euler characteristic} $e(X)$\nomenclature{$e(X)$}{the topological Euler characteristic of a variety $X$} is the alternating sum of the Betti numbers: because Poincar\`e duality holds (i.e. $b^i(X)=b^{4-i}(X)$) we have
\begin{equation}
e(X)=2-2b^1(X)+b^2(X).
\label{eq_topchar}
\end{equation}

There is a relation between the Euler characteristic of the structure sheaf $\O_X$ of $X$ and the topological Euler characteristic $e(X)$, the Noether's formula\index{Noether!formula for surfaces}, which states
\begin{equation}
\chi(\O_X)=\frac{K_X^2+e(X)}{12}.
\label{eq_noetform}
\end{equation}

Over the complex numbers, through Hodge theory, one can shows that $q'(X)=\frac{1}{2}b^1(X)$.
This does not hold in general in positive characteristic, but there is still a relation between these invariants which can be obtained, for example, by the theory of the Picard scheme and the Albanese variety.

\begin{definition}
Given a surface $X$ we define the Picard scheme\index{Picard! scheme} $\Pic_X$\nomenclature{$\Pic_X$}{the Picard scheme of $X$} to be the group scheme which represents the functor $\Pic_X$ (with an abuse of notation, we use the same symbol for the scheme and the functor), which to a scheme $S$ associates the group $\Pic(X\times S)/\Pic(S)$ of all $S$-isomorphism classes of line bundles on $X\times S$: two line bundles $L$ and $L'$ on $X\times S$ are $S$-isomorphic if there exists a line bundle $M$ on $S$ such that $L\simeq L'\otimes_{\O_{X\times S}}\pi_S^*M$ where $\pi_S\colon X\times S\to S$ is the natural projection.

Consider the closed subfunctor  of $\Pic_X$ of line bundles which are algebraically trivial when restricted to fibres of $\pi_S$.
It is represented by a closed subscheme of $\Pic_X$ whose reduced subscheme $\Pic^0_X$\nomenclature{$\Pic^0_X$}{the Picard variety of $X$} is called the Picard variety\index{Picard!variety} of $X$ and  is an Abelian variety\index{Abelian variety}.
\end{definition}

\begin{remark}
\label{rem_picrepresent}
It is not obvious a priori that $\Pic_X$ is a representable functor: a proof of this can be found in \cite{fgav} Section 6.
$\Pic_X$ is a disjoint union of an infinite family of proper schemes and $\Pic_X^0$ is the reduced subscheme of the connected component of the origin.
The tangent space\index{Picard!variety, tangent space of} at the origin is canonically isomorphic to $H^1(X,\O_X)$ while the dimension of the Picard variety\index{Picard!variety, dimension of} is equal to $\frac{1}{2}b_1(X)$ (cf. \cite{badescu} proof of Theorem 5.1): in particular, it is a topological invariant. 

Over a field of characteristic zero, every group scheme is reduced (cf. \cite{mumford} Theorem on page 101), from which we derive $\dim(\Pic^0_X)=q'$; while, in general, we just have $0\leq q'-\dim(\Pic_X^0)\leq p_g(X)$ (cf. \cite{badescu} Theorem 5.1) and there exist examples with $q'-\dim(\Pic_X^0)$ arbitrarily large (cf. \cite{liedtkeuni} Theorem 8.1).
\end{remark}

As we have already noticed, $\Pic_X^0$ is an Abelian variety: hence we may consider its dual Abelian variety\index{Abelian variety!dual} (cf. \cite{mumford} III.13) which turns out to satisfy an important universal property.

\begin{proposition}[\cite{badescu} Theorem 5.3]
\label{propo_albuniv}
Let $X$ be a surface, $x\in X$ be a point, $\Pic^0_X$ be its Picard variety and denote by $\Alb(X)$ the dual Abelian variety of $\Pic^0_X$.
Then there exists a unique morphism $\alb_X\colon X\to\Alb(X)$ which sends $x$ to $0$ and satisfies the following universal property: for any Abelian variety $A$ and any morphism $g\colon X\to A$ that sends $x$ to $0$, there exists a unique morphism $\alpha\colon\Alb(X)\to A$ of Abelian varieties such that $g=\alpha\circ\alb_X$.
In particular, for every morphism $f\colon X\to Y$ that sends $x\in X$ to $y\in Y$, there exists a unique morphism $\alb_f\colon\Alb(X)\to\Alb(Y)$\index{Albanese!morphism}\nomenclature{$\alb_f$}{the morphism $\Alb(X)\to\Alb(Y)$ induced by a morphism $f\colon X\to Y$}  of Abelian varieties such that the following diagram commutes:
\begin{equation}
\begin{tikzcd}
X\arrow{r}{f}\arrow{d}{\alb_X} & Y \arrow{d}{\alb_Y}\\
\Alb(X)\arrow{r}{\alb_f} & \Alb(Y).
\end{tikzcd}
\label{eq_univalbcomm}
\end{equation}
\end{proposition} 

\begin{definition}
Let $X$ be a surface.
We call $\Alb(X)$, as defined in Proposition \ref{propo_albuniv}, the Albanese variety\index{Albanese!variety} of $X$.

Moreover we denote by $q(X)$ the dimension\index{Albanese!variety, dimension of} of $\Alb(X)$.
\label{def_albdef}
\end{definition}

\begin{remark}
It is clear by definition that we have $q'(X)\geq q(X)$.
\label{rem_qq'}
\end{remark}

There is still one numerical invariant of surfaces that we need before giving a brief exposition of the classification of surfaces, i.e. the Kodaira dimension\index{Kodaira dimension}.
\begin{definition}
\label{def_koddim}
Let $X$ be a surface, $K_X$ be its canonical divisor and $\phi_{nK_X}$ the rational map associated with $nK_X$.
We define the Kodaira dimension $k(X)$ of $X$ to be the dimension of the image of $\phi_{nK_X}\colon X\to\P^{h^0(X,nK_X)-1}$ for $n\gg 0$.
We define the Kodaira dimension $k(X)$ of $X$ to be $-\infty$ if $h^0(X,nK_X)=0$ for all integers $n$.
\end{definition}

\begin{remark}
The Kodaira dimension can also be defined as the integer $k(X)$ such that $h^0(nK_X)$ grows like $n^{k(X)}$ up to a scalar.

In particular we have that the Kodaira dimension is a birational invariant because $h^0(nK_X)$ is.
\end{remark}

Recall that a surface $X$ is said to be minimal\index{surface!minimal} if, for every birational morphism $f\colon X\to X'$, we have that $f$ is an isomorphism.
This is equivalent to ask that in $X$ there is no $-1$-curve\index{curve!$-1$-curve} (cf. \cite{badescu} Theorem 3.30 or \cite{beauville_1983} Theorem II.17), i.e smooth rational curve $C$ with $C^2=K_X.C=-1$.

It is also known that every surface $X$ dominates a minimal surface $X'$ (this easily follows from the Castelnuovo's contractibility criterion, cf. \cite{badescu} Theorem 3.30, and from  the finiteness of $h^1(X,\Omega^1_X)$ or of the rank of the N\'eron-Severi group, cf. \cite{badescu} Corollary 6.3 or \cite{beauville_1983} Theorem II.16), i.e. there exists a birational morphism $f\colon X\to X'$ and this $f$ is known to be a composition of blow-ups.
Such a surface $X'$ is called a minimal model of $X$.
For this reason the classification of surfaces reduces to the classification of surfaces up to birational equivalence.
Moreover all the numerical invariants considered for surfaces are birational, with the single exception of $K_X^2$ whose behaviour under blow-ups is well understood.
Here is a brief recap of the classification\index{surface!classification of} of  surfaces based on the Kodaira dimension.

If the Kodaira dimension is $-\infty$ we have that every minimal model of $X$ is isomorphic to a $\P^1$-bundle over a smooth curve $C$ or to the projective plane: this is the only case where a surface $X$ has more than one minimal model.

When the Kodaira dimension of $X$ is greater than or equal to $0$, by the following Theorem, we have the uniqueness of the minimal model of $X$.

\begin{theorem}[\cite{badescu} Theorem 10.21 or \cite{beauville_1983} Theorem V.19]
Let $X$ and $X'$ be two minimal surfaces with non-negative Kodaira dimension.
Then every birational map from $X$ to $X'$ is an isomorphism.
In particular, for every surface $\overline{X}$ birational to $X$, we have $\bir(\overline{X})=\aut(X)$ where $\bir(X)$\nomenclature{$\bir(X)$}{the group of birational maps of a variety $X$} denotes the group of birational maps of $X$ and $\aut(X)$\nomenclature{$\aut(X)$}{the group of automorphism  of a variety $X$} denotes the group of automorphism of $X$.
\label{teo_kxgeq0}
\end{theorem}

The minimal models of surfaces of Kodaira dimension $0$  split into four classes and satisfy $K_X^2=0$.
Actually, if the characteristic of the ground field is $2$ or $3$, there are two more subclasses.
\begin{itemize}
\item
We call such a surface Abelian\index{surface!Abelian}, if it has $q'=q=2$, $p_g=1$ and trivial canonical bundle: it turns out to be an Abelian variety.
\item
If $q'=q=0$, $p_g=1$ and $\omega_X=\O_X$ we talk about K3 surfaces\index{surface!K3}.
\item
When $q'=q=0$, $p_g=0$ and the order of $K_X$ is two, we say that $X$ is an Enriques surface\index{surface!(non-classical) Enriques}: considering the double cover $Y\to X$ given by $\omega_X^2=\O_X$ (which is inseparable if the characteristic is $2$, cf. \cite{cossec} Proposition 0.1.8), we have that $Y$ is a K3 surface (\cite{cossec} Proposition 1.3.1).
If the characteristic is $2$ we have a subclass of the Enriques surfaces: these are called non-classical Enriques surfaces and satisfy $q=p_g=1$, $q'=0$.
Also in this case there is canonically defined a double cover $Y\to X$ where $Y$, which may be singular, is a "K3-type" surface, by which we mean that $h^1(Y,\O_Y)=0$ and $\omega_Y=\O_Y$: it can be shown that if $Y$ is not smooth then it is rational (ibidem).
\item
Finally if $q=q'=1$, $p_g=0$ we say that $X$ is hyperelliptic\index{surface!(quasi-)hyperelliptic}: in this case we have that $X=(E\times F)/G$ where $G$ acts by translation on $E$ and faithfully on $F$ where both $E$ and $F$ are elliptic curves: moreover $X\to E/G$ is an elliptic fibration\index{fibration!elliptic} which coincides with the Albanese morphism, and $X\to F/G=\P^1$ is another elliptic fibration.
If the characteristic is $2$ or $3$ we have another subcase which consists of quasi-hyperelliptic surfaces: in this case we have $q'=2$, $q=p_g=1$ and $X=(E\times F)/G$ where $G$ is a finite group scheme acting by translation on $E$, and on a curve $F$ of arithmetic genus equal to one and, as above, $X\to E/G$ is a quasi-elliptic fibration\index{fibration!quasi-elliptic} (i.e. a genus-one fibration whose general fibre is singular) with fibre $F$ and coincides with the Albanese morphism and $X\to F/G=\P^1$ is an elliptic fibration.
\end{itemize}

A surface $X$ whose Kodaira dimension is $1$ always admits a genus-one fibration\index{fibration!genus-one} $X\to C$: if the characteristic is different from $2$ and $3$ it is an elliptic fibration\index{surface!elliptic|see{elliptic surface}}\index{elliptic surface}, otherwise it may happen that the general fibre is a rational curve with a single cusp (cf. \cite{badescu} Chapter 7) i.e. it is quasi-elliptic. 

If $k(X)=2$ we say that $X$ is of general type\index{surface of general type} (as opposed to surface of special type, i.e. surfaces with $k(X)\leq 1$). 
This is the widest and more mysterious  class of surfaces.

\begin{theorem}[\cite{mireille} Corollaire 4]
\label{teo_autogeneral}
Let $X$ be a surface of general type: then $|\bir(X)|<\infty$ where $\bir(X)$ denotes the group of birational maps of $X$.
In particular (cf. Theorem \ref{teo_kxgeq0}), if $X$ is minimal  $|\aut(X)|<\infty$ where $\aut(X)$ denotes the group of automorphism of $X$.
\end{theorem}

In this thesis we mainly focus on surfaces such that the Albanese morphism is generically finite.

\begin{definition}
We say that $X$ is of maximal Albanese dimension\index{surface!of maximal Albanese dimension} if the image $\im(\alb_X(X))$ of the Albanese morphism $\alb_X$ has dimension $2$.
\label{def_maxalbdim}
\end{definition}

It is clear that no surfaces  with $k(X)=-\infty$ have maximal Albanese dimension (the fibres of the $\P^1$-bundle are clearly contracted by the Albanese morphism).
Abelian surfaces are clearly the only surfaces of Kodaira dimension $0$ of maximal Albanese dimension.
It is also clear that there do not exists quasi-elliptic surfaces of maximal Albanese dimension (the fibres of the rational genus-one fibration are contracted by the Albanese morphism).
In particular we have that for a surface of maximal Albanese dimension 
\begin{equation}
K_X^2\geq 0
\label{eq_kx2alb2}
\end{equation}
holds.

\begin{theorem}
Let $X$ be a surface of maximal Albanese dimension, then $\chi(\O_X)\geq 0$.
Moreover, if $X$ is of general type we have $\chi(\O_X)>0$.
\label{teo_shepbarr}
\end{theorem}

\begin{proof}
By \cite{shepherdgeo} Theorem 2, we know that $e(X)\geq 0$: then the thesis easily follows by Noether's formula \ref{eq_noetform} and Equation \ref{eq_kx2alb2}.
\end{proof}

\section{Curves on surfaces and Gorenstein curves}
\label{sec_curves}
We collect here some known results for curves on surfaces\index{surface!curve on}: the main references are \cite{barth} for $\cha(k)=0$ and \cite{badescu} for the general case. 
We also recall some classical results for smooth curves\index{curve} and generalize them to the case of integral Gorenstein curves\index{curve!Gorenstein}.
For the first part of this section we consider curves on surfaces, i.e.  divisors $C=\sum n_iC_i$ on a smooth projective surface $X$.
Then we will focus on the case where $C$ is a projective integral Gorenstein scheme of dimension $1$ (which includes the case of integral curves on a surface): the Gorenstein condition ensures that the dualizing sheaf is a line bundle.

Now suppose that $C$ is a divisor on a surface.
Its structure sheaf is the quotient of the structure sheaf of the surface as in the following short exact sequence
\begin{equation}
0\to \O_X(-C)\to\O_X\to\O_C\to 0.
\label{eq_structuresheaf}
\end{equation} 

If we can split a curve $C$ as a sum of two curves $A$ and $B$, by the inclusion $\O_X(-C)\subseteq\O_X(-B)$, we obtain the following diagram
\begin{equation}
\begin{tikzcd}
0\arrow{r}{} & \O_X(-C)\arrow[hook]{d}{}\arrow{r}{} &\O_X\arrow[equals]{d}{}\arrow{r}{}&\O_C\arrow[two heads]{d}{}\arrow{r}{}&0\\
0\arrow{r}{} & \O_X(-B)\arrow{r}{} &\O_X\arrow{r}{}&\O_B{}\arrow{r}{}&0
\end{tikzcd}
\label{eq_sumcurves}
\end{equation}
and applying the snake Lemma we get
\begin{equation}
\begin{split}
\ker(\O_C\to\O_B)=\O_X(-B)/\O_X(-C)=\\
=\coker(\O_X(-C)\to\O_X(-B))=\O_A(-B),
\end{split}
\label{eq_snake}
\end{equation}
which leads to the short exact sequence
\begin{equation}
0\to \O_A(-B)\to \O_C\to\O_B\to 0.
\label{eq_decseq}
\end{equation}

For a general curve $C$ we have the adjunction formula for the dualizing sheaf\index{adjunction formula} (cf. \cite{barth} II.1 or \cite{badescu} 3.11.5)
\begin{equation}
\omega_C=\restr{\bigl(\omega_X\otimes_{\O_X}\O_X(C)\bigr)}{C}.
\label{eq_dual}
\end{equation}
Recall that the dualizing sheaf gives the Serre duality\index{Serre duality} (cf. \cite{badescu} 3.11.3), i.e. for every coherent sheaf $F$ on $C$ we have 
\begin{equation}
\label{eq_serredual}
H^i(C,F)=H^{1-i}(C,\omega_C\otimes_{\O_C}F^\vee)^\vee,
\end{equation}
\nomenclature{$F^\vee$}{the dual of a coherent sheaf $F$}
\nomenclature{$H^i(X,F)^\vee$}{the dual of the $i$-th cohomology vector space}
and that it is a line bundle provided that $C$ is Gorenstein.

\begin{theorem}[Riemann-Roch]
\label{teo_RR}
\index{Riemann-Roch Theorem! curves on surfaces}
\index{Riemann-Roch Theorem! Gorenstein curves}
Let $C$ be a curve on a surface $X$ or a projective integral Gorenstein curve. 
Then for every vector bundle $E$ of fixed rank $r$ on every component of $C$ we have
\begin{equation}
\chi(E)=\deg(E)+r\chi(\O_C).
\label{eq_RR}
\end{equation}

If $C$ is a smooth curve, we have the same result for every coherent sheaf $F$  of rank $r$.
\end{theorem}

Before the proof we need to define what is the degree of a vector bundle for singular projective curves.
\begin{definition}
Let $C$ be an irreducible curve and $E$ be a vector bundle on it.
We define the degree of $E$ as $\deg(E)=\deg(\det(\nu^*E))$\index{degree! of a vector bundle on a curve on a suruface}\index{degree! of a vector bundle on a Gorenstein curve} where $\nu\colon\widetilde{C}\to C$ is the normalization morphism of $C$.




Let $C=\sum n_iC_i$ be a curve and $E$ a coherent sheaf on $C$, we define the degree of $E$ as $\sum n_i\deg(\restr{E}{C_i})$.
\label{def_degcurve}
\end{definition}

\begin{proof}[Proof of Theorem \ref{teo_RR}]
If the curve is smooth, this is a well known result: cf. \cite{hart} Theorem IV.1.3 for the case of line bundles, cf. \cite{fulton} Corollary 15.2.1  for vector bundles.
If $F$ is not locally free one can easily reduce to the locally free case by taking a locally free resolution and noticing that  $\chi$, $\deg$ and $\rk$ are additive on exact sequences (recall that on a smooth variety, there exists a finite locally free resolution, cf. \cite{eisenbud} Corollary 19.6).

If the curve is not smooth, the proof is completely identical as over the complex numbers and we refer to \cite{barth} Theorem II.3.1.
\end{proof}

\begin{definition}
\index{curve!arithemtic genus of}
Let $C$ be a curve on a surface or a projective integral Gorenstein curve: we define the arithmetic genus of $C$ to be $g(C)=1-\chi(\O_C)$.
If $C$ is reduced, we define the geometric genus of $C$ to be $p_g=h^1(\O_{\widetilde{C}})$ where $\widetilde{C}$ is the normalization of $C$.\index{curve!geometric genus of} 
\end{definition}

Notice that if the curve is integral, we have that $h^0(C,\O_C)=1$, hence $g(C)=h^1(C,\O_C)=h^0(C,\omega_C)$.
Actually this formula remains true for all curves $C$ that are 1-connected (cf. \cite{barth} II.12).

\begin{lemma}
\label{lemma_gencurves}
Let $C$ be a curve on a surface; then we have the following two formulae:
\begin{enumerate}
	\item $g(C)=1+\frac{1}{2}(K_X+C).C$;
	\item if $C=A+B$, then $g(C)=g(A)+g(B)+A.B-1$.
\end{enumerate}
\end{lemma}

\begin{proof}
By the defining exact sequence \ref{eq_structuresheaf} of $\O_C$ we get $\chi(\O_C)=\chi(\O_X)-\chi(\O_X(-C))$ from which we obtain, using the Riemann-Roch Theorem for smooth surfaces \ref{teo_rrss},
\begin{equation}
g(C)=1-\chi(\O_C)=1+\chi(\O_X(-C))-\chi(\O_X)=1+\frac{1}{2}(K_X+C).C.
\label{eq_genc}
\end{equation}

Using 1, we see that
\begin{equation}
\begin{split}
g(C)=1+\frac{1}{2}(K_X+C).C=&1+\frac{1}{2}(K_X+A+B).(A+B)=\\
=1+\frac{1}{2}(K_X+A).A+1+&\frac{1}{2}(K_X+B).B+A.B-1=\\
=g(A)+&g(B)+A.B-1.\qedhere
\end{split}
\label{eq_genusab}
\end{equation}
\end{proof}

From now on we restrict to study projective integral Gorenstein curves (e.g. integral curves on surfaces).
\begin{lemma}[\cite{hartgen} Theorem 1.6]
\label{lemma_gg}
Let $C$ be a  projective integral Gorenstein curve of genus $g(C)\geq 1$; then the canonical bundle\index{canonical divisor!of an integral Gorenstein curve} $\omega_C$ of $C$ is globally generated.
\end{lemma}


\begin{theorem}[Clifford]
Let $C$ be a projective integral Gorenstein curve of genus $g$ and let $L$ be a line bundle such that $h^0(C,L)\neq0\neq h^1(C,L)$.\index{Clifford Theorem for integral Gorenstein curves}
Then $2h^0(C,L)\leq \deg(L)+2$.
\label{teo_cliff}
\end{theorem}

\begin{proof}
The proof for smooth curves (cf. \cite{arbarello} Clifford's Theorem in Section III.1) relies on Riemann-Roch, which we have already shown to hold also in this setting, and the fact that the dualizing sheaf is a line bundle, which is assured by the Gorenstein hypothesis: in particular it works also in this case. 
%
\end{proof}

\begin{corollary}
\label{coro_luzuo}
Let $L$ be a line bundle on a projective integral Gorenstein curve $C$ and let $V\subseteq H^0(C,L)$ be a globally generated linear subseries.
Denote by $r=\dim(V)$, $d=\deg(L)$, by $\iota\colon C\to\P^{r-1}$\nomenclature{$\P^r$}{the projective space of dimension $r$} the morphism induced by $V$ and by $g$ the geometric genus of the image $\widetilde{C}$ of $\iota$.
Then we have 
\begin{equation}
d\geq \deg(\iota)\min\{2(r-1),r+g-1\}\geq \deg(\iota)(r-1).
\label{eq_corluzuo}
\end{equation}
\end{corollary}

\begin{proof}
Let $\nu\colon C'\to\widetilde{C}$ be the normalization morphism and let $L'=\nu^*\O_{\widetilde{C}}(1)$.
Clearly we have $d=\deg(\iota)\deg(L')$ and, by Riemann-Roch, $\deg(L')=h^0(C',L')-h^1(C',L')+g-1\geq r-h^1(C',L')+g-1$.
If $h^1(C',L')=0$ we are done, otherwise, suppose that $h^1(C',L')\neq 0$ and apply Clifford's Theorem: we get $\deg(L')\geq 2(h^0(C',L')-1)\geq 2(r-1)$.
Summing up we have $d\geq \deg(\iota)\min\{2(r-1),r+g-1\}$.
\end{proof}

\begin{remark}
\label{rem_luzuo}
It is straightforward that in the Corollary above we may substitute $C'$ with any Gorenstein curve $\overline{C}$ which is intermediate between $C'$ and $\widetilde{C}$ and the geometric genus of $\widetilde{C}$ by the arithmetic genus of $\overline{C}$ (indeed we just need to apply Theorems \ref{teo_RR} and \ref{teo_cliff}).
In particular, if $\iota$ is birational, we have 
\begin{equation}
d\geq \min\{2(r-1),r+g(C)-1\},
\end{equation}
where $g(C)$ is the arithmetic genus of $C$.

If, moreover, $L$ is a subsheaf of the canonical sheaf we have that $d\geq 2(r-1)$: indeed we know that $h^1(C,L)=h^0(C,K_C-L)\neq 0$ and Clifford's Theorem applies.
\end{remark}


\begin{corollary}
\label{coro_luzuo2}
Let $L$ be a line bundle on a projective integral Gorenstein curve $C$ and let $V\subseteq H^0(C,L)$ be a globally generated linear subseries.
Denote by $r=\dim(V)$, $d=\deg(L)$, by $\iota\colon C\to\P^{r-1}$ the morphism induced by $V$, by $g$ the geometric genus of the image $\widetilde{C}$ of $\iota$ and by
\begin{equation}
\rho\colon V\otimes V\to H^0(C,L^2)
\end{equation}
the natural multiplication map.
Then we have that 
\begin{equation}
\dim(\im(\rho))\geq\min\{3(r-1),2r+g-1\}\geq 2r-1.
\end{equation}
\end{corollary}

\begin{proof}
As above, let $\nu\colon C'\to\widetilde{C}$ be the normalization morphism and let $L'=\nu^*\O_{\widetilde{C}}(1)$, $V'=\nu^*H^0(\widetilde{C},\O_{\widetilde{C}}(1))$ and $\rho'\colon V'\otimes V'\to H^0(C',(L')^2)$.
It is clear that $r=\dim(V')$ and $\dim(\im(\rho))=\dim(\im(\rho'))$, hence we may assume that $\iota$ is a birational map.

By construction, we know that $L=\iota^*\O_{\P^{r-1}}(1)$, $V=H^0(\P^{r-1},\O_{\P^{r-1}}(1))$ and the linear subseries of $|L^2|$ induced by $\rho$ is the pull-back of the linear series induced by polynomials of degree $2$ on $\P^{r-1}$.
Following the construction of \cite{arbarello} section III.2, we know that 
\begin{equation}
\dim(\im(\rho))\geq 1+\min\{d,r-1\}+\min\{d,2r-3\}.
\label{eq_arba}
\end{equation}

Combining Equation \ref{eq_arba} with Corollary \ref{coro_luzuo} the claim is proven.
\end{proof}

\begin{remark}
\label{rem_luzuo2}
As in Remark \ref{rem_luzuo}, it is straightforward that in the Corollary above we may substitute $C'$ with every Gorenstein curve $\overline{C}$ which is intermediate between $C'$ and $\widetilde{C}$ and the geometric genus of $\widetilde{C}$ by the arithmetic genus of $\overline{C}$.
In particular, if $\iota$ is birational, we have 
\begin{equation}
\dim(\im(\rho))\geq\min\{3(r-1),2r+g(C)-1\}\geq 2r-1,
\end{equation}
where $g(C)$ is the arithmetic genus of $C$.

If, moreover, $L$ is a subsheaf of the canonical sheaf we have that $\dim(\im(\rho))\geq 3(r-1)$.
\end{remark}

\begin{remark}
\label{rem_lineargg}
Notice that in Corollaries \ref{coro_luzuo} and \ref{coro_luzuo2} we can get rid of the globally generated hypothesis.
If $V$ has a fixed part we can consider the base locus $P$ of $V$: because $C$ can be singular, it may happen that $P$ is not a Cartier divisor.
However, because the divisors in the linear series corresponding to $V$ are Cartier, it make sense to consider the linear series $V-P$ (for details see \cite{hartgen} Section 1), whose elements are Cartier divisor.
In particular, as for smooth curve, for any linear series $V$ there exists an associated base-point-free linear system $V'$ whose elements are Cartier divisors and we define the morphism associated with $V$ as the morphism associated with $V'$.
It is clear by definition that the dimension of $V$ is equal to the dimension of $V'$, the line bundle associated to $V'$ is $L'=L(-P)$ where $L$ is the line bundle associated with $V$ (in particular the degree $d$ of $L$ is greater than or equal to the degree of $L'$) and the dimension of the image of $\rho\colon V\otimes V\to H^0(C,L^2)$ is equal to the dimension of the image of $\rho'\colon V'\otimes V'\to H^0(C,L'^2)$.
It is then straightforward to extend Corollaries \ref{coro_luzuo} and \ref{coro_luzuo2} to a general $V$.
\end{remark}

\begin{definition}
A hyperelliptic curve\index{curve!hyperelliptic} $C$ of arithmetic genus $g\geq 2$ is an integral flat double cover of the projective line $\P^1$.
\end{definition}

Notice that, for us, a hyperelliptic curve may be singular, but its singularities are all Gorenstein (cf. Remark \ref{rem_canonical}).
Here we give a generalization to integral curves of the well-known result of Max Noether.

\begin{theorem}[\cite{contiero} Theorem 1]
\label{teo_max}
Let $C$ be a non-hyperelliptic Gorenstein projective integral curve\index{curve!non-hyperelliptic, Max Noether's Theorem for}. 
Then we have that 
\begin{equation}
S^nH^0(C,K_C)\to H^0(C,nK_C)
\label{eq_...} 
\end{equation}
is surjective for every $n\in\N$.
\end{theorem}

We recall also the following result.

\begin{lemma}[\cite{hartgen} Theorem 1.6]
\label{lemma_canbirnonhyp}
Let $C$ be a non-hyperelliptic Gorenstein projective integral curve. 
Then we have that the canonical morphism $\phi_{K_C}\colon C\to \P^{h^0(C,K_C)-1}$ is an embedding.
\end{lemma}

Recall also that a hyperelliptic curve admits exactly one hyperelliptic structure, i.e. we have the following Proposition.

\begin{proposition}
\label{teo_hyperdouble}
Let $C$ be a hyperelliptic curve of arithmetic genus $g\geq2$.
Then there exists one and only one double cover $C\to \P^1$ and it coincides with the canonical morphism.
\end{proposition}

\begin{proof}
Again, the proof is completely identical to the case of a smooth curve and we refer to \cite{hart} Proposition IV.5.3.
%
%
\end{proof}

\section{Elliptic surfaces}
\label{sec_ell}

As\index{elliptic surface} we have already noticed, every surface $X$ of Kodaira dimension $1$ admits a surface fibration $f\colon X\to C$ with fibres of arithmetic genus equal to $1$.
In characteristic zero we have that the general fibre is smooth thanks to  generic smoothness Theorem (\cite{hart} III.10.7), while in characteristic different from zero it may happen that the general fibre is integral but not smooth (cf. \cite{badescu} Corollary 7.3 and Definition 7.6 and the discussion following to it): if the general fibre is smooth we say that $X$ is an elliptic surface, otherwise it is a quasi-elliptic surface.
Notice that the general fibre of a quasi-elliptic surface is rational, in particular quasi-elliptic surfaces can not have maximal Albanese dimension.
In this thesis we are mainly interested in surfaces of maximal Albanese dimension, that's why  here we focus on elliptic surfaces.
In particular we will show that a smooth elliptic surface of maximal Albanese dimension is a product $C\times E$.
This result is known and over the complex numbers can be derived from Theorem 2.1 of \cite{serrano}: the argument we present here works over algebraically closed fields of any characteristic.

Denote by $q(X)$ the dimension of the Albanese variety of $X$, by $q'(X)=h^1(X,\O_X)$ and by $g(C)$ the genus of $C$.

Consider the first higher pushforward\nomenclature{$R^if_*F$}{the $i$-th derived pushforward of a coherent sheaf $F$ via a morphism $f$} $R^1f_*\O_X=\L\oplus T$ of $\O_X$ where $T=T(R^1f_*\O_X)$ is the torsion subsheaf of $R^1f_*\O_X$ and $\L=R^1f_*\O_X/T$ is a line bundle: it is known that a point $c\in C$ is in the support of $T$ if and only if $h^0(F_c,\O_{F_c})\geq 2$ where $F_c$ is the fibre of $f$ over $c$\index{fibre!multiple}, in particular we have that $F_c$ has to be a multiple fibre (cf. \cite{badescu} page 100).

\begin{definition}[ibidem Definition 7.14]
We call such a fibre $F_c$ an exceptional fibre\index{fibre!exceptional}.
\label{def_excfib}
\end{definition}

It is known that over the complex numbers (cf \cite{barth} Corollary III.11.2) $T=0$ i.e. there are no exceptional fibres.

\begin{theorem}[\cite{badescu} Theorem 7.15]
\label{teo_canellsur}
Let $f\colon X\to C$ be an elliptic fibration and let $R^1f_*\O_X=\L\oplus T$ as above.
Then we have that
\begin{equation}
K_X=f^*(\L^{-1}\otimes\omega_C)\otimes\O_X\Bigl(\sum_{i=1}^ma_iF_i\Bigr)
\label{eq_canellsur}
\end{equation}
where
\begin{itemize}
\item $n_iF_i=F_{c_i}$ are all the multiple fibres of $f$ for $i=1,\ldots m$ with $F_i$ reduced;
\item $0\leq a_i<n_i$;
\item $a_i=n_i-1$ if $F_{b_i}$ is not an exceptional fibre.
\end{itemize}

In particular the canonical divisor $K_X$\index{canonical divisor!of an elliptic surface|see{elliptic surface}}\index{elliptic surface!canonical divisor of} is contracted by $f$ and satisfies $K_X^2=0$.
\label{teo_canellsur}
\end{theorem}

\begin{remark}
\label{rem_topchar+}
Recall that we have the following formula (cf. \cite{cossec} Proposition 5.1.6) for the topological Euler characteristic\index{elliptic surface!topological Euler characteristic of|see{topological Euler characteristic}}\index{topological Euler characteristic! of an elliptic surface} of an elliptic surface $f\colon X\to C$:
\begin{equation}
e(X)=e(F)e(C)+\sum_{i=1}^m(e(F_{c_i})-e(F)+\delta_i)=\sum_{i=1}^m(e(F_{c_i})+\delta_i)\geq 0
\label{eq_eultopchar+}
\end{equation}
where $F$ is the generic fibre of $f$,  $c_i$ is as in Theorem \ref{teo_canellsur} and $\delta_i$ is a non negative integer which is zero if, for example, $F_i$ is smooth.
In particular we have that $e(X)=0$ if and only if  $f$ has only smooth, but possibly multiple, fibres (ibidem or \cite{liedtke+} Section 5.3).

Moreover, by Noether's formula \ref{eq_noetform} and the formula for the canonical bundle of Theorem \ref{teo_canellsur}, we have
\begin{equation}
\chi(\O_X)=\frac{e(X)}{12}\geq0
\label{eq_eulercharellsufc}
\end{equation}
and equality occurs if and only if all the fibres of $f$ are smooth (but posiibly multiple).
\end{remark}

\begin{remark}
\label{rem_spectralsfib}
Consider the Leray spectral sequence
\begin{equation}
E^{pq}_2=H^p(C,R^qf_*\O_X)\Rightarrow H^{p+q}(X,\O_X).
\label{eq_ssellfib}
\end{equation}
The exact sequence on low terms (cf. \cite{cartan} Theorem XV.5.11) tells us
\begin{equation}
\begin{split}
0\to H^1(C,\O_C)\to H^1(X,\O_X)\to H^0(C,R^1f_*\O_X)\to\\
\to H^2(C,\O_C)=0\to H^2(X,\O_X)\to H^1(R^1f_*\O_X)\to 0.
\end{split}
\label{eq_lowdegreeterm}
\end{equation}

In particular from this, using Riemann-Roch for curves \ref{teo_RR}, we derive that 
\begin{equation}
\begin{split}
\chi(\O_X)=\chi(&\O_C)-\chi(R^1f_*\O_X)=\\
=-\deg(R^1f_*\O_X)&=-\deg(\L)-h^0(C,T)
\end{split}
\label{eq:eulerellfib}
\end{equation}
and
\begin{equation}
q'(X)=g(C)+h^0(C,R^1f_*\O_X)=g(C)+h^0(C,\L)+h^0(C,T).
\label{eq_irralb}
\end{equation}

By Remark \ref{rem_topchar+}, we have that 
\begin{equation}
-\deg(\L)\geq h^0(C,T)\geq 0
\label{eq_bohdeg}
\end{equation}
and $\deg(\L)=0$ implies that there are no exceptional fibres and $q'(X)=g(C)$ or $q'(X)=g(C)+1$ depending on  whether $\L$ is trivial or not.
\end{remark}

\begin{lemma}
\label{lemma_ellalbjac}
Let $f\colon X\to C$ be an elliptic surface, $\alb_X\colon X\to\Alb(X)$ be it Albanese morphism\index{elliptic surface!Albanese variety of}\index{Albanese!variety of an elliptic surface|see{elliptic surface}}, $AJ_C\colon C\to\jac(C)$\nomenclature{$\jac(C)$}{the Jacobian of a smooth curve $C$}\nomenclature{$AJ_C$}{the Abel-Jacobi map from a curve $C$ to its Jacobian} be the Abel-Jacobi map of $C$ and $\alb_f\colon\Alb(X)\to\jac(C)$ be the induced morphism.
Then we have that $\Alb(X)=\jac(C)$ if and only if every fibre of $f$ (or, equivalently, there exists a fibre which) is contracted by $\alb_X$, otherwise we have that $q(X)=g(C)+1$ and all the fibres of $X$ are isogenous to the kernel of $\alb_f$.
\end{lemma}

\begin{proof}
We have the following Cartesian diagram
\begin{equation}
\begin{tikzcd}
X\arrow{r}{f}\arrow{d}{\alb_X} & C\arrow{d}{AJ_C}\\
\Alb(X)\arrow{r}{\alb_f} & \jac(C).
\end{tikzcd}
\label{eq_albellipticqg}
\end{equation}
By the universal property of the Albanese morphism and the surjectivity of $f$, we have that $\alb_f$ is surjective.
If $\alb_f$ is an isomorphism, clearly, by the commutativity of the diagram \ref{eq_albellipticqg}, it follows that $\alb_X$ contracts the fibres of $f$. 

Now assume that there exists a fibre $F$, such that $\alb_X(F)$ is a point inside $\Alb(X)$.
If $q(X)=0$, clearly $\alb_X$ contracts the fibres of $f$.
Suppose there exists a fibre $F'$ of $f$, such that $\alb_X(F')$ is a curve.
Then, there is a hyperplane section $H$ which intersects $\alb_X(F')$ but does not intersect $\alb_X(F)$: hence $\alb_X^*H$ intersects $F'$ but does not intersect $F$, which is not possible because all the fibres are numerically equivalent.
In particular $\alb_X$ contracts every fibre of $f$ and there exists a morphism $h\colon C\to \Alb(X)$ such that
\begin{equation}
\begin{tikzcd}
X\arrow{r}{f}\arrow{d}{\alb_X} & C\arrow{dl}{h}\arrow{d}{AJ_C}\\
\Alb(X)\arrow{r}{\alb_f} & \jac(C).
\end{tikzcd}
\label{eq_albellipticqg2}
\end{equation}
commutes.
Again, by the universal property of the Abel-Jacobi morphism, we have that $\Alb(X)=\jac(C)$.

Now suppose that $\alb_X$ does not contract any fibre of $f$.
It is clear, by the commutativity of the diagram \ref{eq_albellipticqg}, that $\alb_X(F)$ is an elliptic curve contracted by $\alb_f$ for every fibre $F$ of $f$.
Consider
\begin{equation}
\begin{tikzcd}
\Alb(X)\arrow{r}{\pi}\arrow{d}{\alb_f} & \Alb(X)/\alb_X(F_0)\arrow{dl}{\pi'}\\
\jac(C) &
\end{tikzcd}
\label{eq_variaell}
\end{equation}
where $F_0$ is an elliptic fibre such that $0\in\alb_X(F_0)$: if $q(X)=1$ then $\dim(\Alb(X)/\alb_X(F_0))=g(C)=0$ and we are done.
If $\dim(\Alb(X)/\alb_X(F_0))\geq 1$, arguing as before, the image of $\pi\circ\alb_X$ is a curve: in particular there exists a morphism $h\colon C\to \Alb(X)/\alb_X(F_0)$ such that $h\circ f=\pi\circ \alb_X$.
Hence, by the universal property of the Jacobian variety, $\jac(C)=\Alb(X)/\alb_X(F_0)$ and every fibre of $f$ is isogenous to $\alb_X(F_0)$.
\end{proof}

\begin{remark}
\label{rem_ellsurf1}
Let $f\colon X\to C$ be an elliptic fibration such that $\alb_X$ does not contract the fibres: in this case we can prove, with the help of Remark \ref{rem_topchar+} that the  topological Euler characteristic $e(X)$ and the Euler characteristic $\chi(\O_X)$ of the structure sheaf have to be zero.
Indeed, we have that the fibres are smooth, but possibly multiple, otherwise they would be rational and hence contracted by the Albanese morphism (cf. \cite{cossec} Proposition 5.1.4).

In particular we conclude that $e(X)=\chi(\O_X)=0$.
If we further assume that $f\colon X\to C$ has no multiple fibres, we can conclude that $T=0$ and, consequently (using Remark \ref{rem_spectralsfib} and Lemma \ref{lemma_ellalbjac}), that $q(X)=q'(X)$. 
\label{rem_notcontractsfibre}
\end{remark}

\begin{lemma}
\label{lemma_isoell}
Let $f\colon X\to C$ be a smooth elliptic surface fibration, i.e. all the fibres of $f$ are smooth curves of genus $1$.
Then $f$ is isotrivial, i.e. all the fibres are mutually isomorphic.
Moreover the fibration is \'etale locally trivial, i.e. for every point $c\in C$ there exists an \'etale neighbourhood $U$ (by which we mean an \'etale morphism $U\to C$ whose image contains $x$) such that $X_U\simeq F\times U$ as a scheme over $U$, where  $X_U=X\times_C U$.
\end{lemma}

\begin{proof}
We know that an elliptic curve is uniquely determined by its $j$-invariant which takes value in the ground field $k$ and consequently that the coarse moduli space of pointed smooth curves of genus $1$ is isomorphic to the affine line $\A^1_k$ (cf. \cite{olsson} Theorem 13.1.15).
By this we get a morphism $C\to \A^1_k$\nomenclature{$\A_k^i$}{the affine space of dimension $i$ over the field $k$} which has to be constant, in particular all the fibres of $f$ are mutually isomorphic.

We know that every smooth morphism locally admits an \'etale multisection $U'\to X$ (cf. \cite{egaiv4} Corollaire 17.16.3(ii)): in particular we have a section of the fibration $X_{U'}\to U'$  which gives to its fibres a natural structure of elliptic curves.
It is a standard result (cf. \cite{hartdef} Corollary 26.5) that a family of mutually isomorphic elliptic curves (i.e. a family of mutually isomorphic smooth genus-$1$ curve with a section) becomes trivial after an \'etale cover.
\end{proof}

\begin{proposition}
\label{propo_isoell}
Let $f\colon X\to C$ be a smooth elliptic surface fibration.
Then there exists a curve $C'$ and a finite group scheme $G$  (both of which may not be smooth) acting freely on $C'$ and not necessarily freely on $F$ such that $C'/G=C$ and
\begin{equation}
\begin{tikzcd}
(C'\times F)/G\arrow{d}{\pi_1}\arrow{r}{\simeq} & X\arrow{d}{f}\\
C'/G\arrow{r}{\simeq} & C
\end{tikzcd}
\label{eq_ellisotriv}
\end{equation} 
where the action of $G$ on $C'\times F$ is the diagonal one.
We can always take $G=H\ltimes F_d$ where $H$ is the finite group of automorphism of $F$ as an elliptic curve and $F_d$ is the $d$-torsion of the general fibre $F$ of $f\colon X\to C$.
\end{proposition}

\begin{proof}
By Lemma \ref{lemma_isoell}, we know that all the fibres of $f$ are mutually isomorphic. 
Denote by $F$ one of its fibres and by $A$ the group scheme of automorphism of $F$.
Choose a base point $0$ of $F$ so that $F$ has the structure of an elliptic curve: it is known that $A=H\ltimes F$.
Denote by $F_d$\nomenclature{$E_d$}{the finite subgroup scheme of $d$-torsion points of an elliptic curve $E$} the finite subgroup scheme of $d$-torsion points (i.e. the kernel of the multiplication by $d$ morphism $\mu_d\colon F\to F$ cf. \cite{mumford} Section 4) and denote by $A_d=H\ltimes F_d$.
Notice that, if $d$ is coprime with the characteristic of the ground field, $A_d$ is an \'etale finite group scheme (cf. \cite{mumford} Corollary 1 Section 7).

Consider the functor $\isou_C(X,F\times C)$\nomenclature{$\isou_S(X,Y)$}{the functor of isomorphisms of the $S$-schemes $X$ and $Y$} which associates to every scheme $S$ over $C$ the set $\iso_S(X_S,F\times S)$\nomenclature{$\iso_S(X,Y)$}{the group of isomorphisms of the $S$-schemes $X$ and $Y$} of isomorphisms between $X_S$ and $F\times S$ where $X_S=X\times_C S$.
Because $f$ is a flat morphism and $X$ is projective, this functor is represented by a scheme (cf. \cite{fgaiv} 4.c) which we call $P$.
Notice that there exists a natural action of $A$ on $P$ such that the natural morphism $p\colon P\to C$ is $A$-invariant.
By the \'etale local triviality of $f\colon X\to C$ (cf. Lemma \ref{lemma_isoell}) we easily derive that $p\colon P\to C$ is an $A$-torsor, i.e. there exists  a covering $\{U_i\to C\}$ for the flat topology such that $P\times_C U_i\simeq A\times U_i$, (cf. \cite{milne} Proposition III.4.1) in particular $P/A=C$.

Let $X_d=\Pic_{X/C}^d\to C$ be the relative Picard scheme parametrizing line bundles of degree $d$ along the fibres: this is the scheme representing the functor $\Pic_{X/C}^d$ which associates to a scheme $S$ the set of $S$-isomorphism classes of line bundles on $X_S$ whose restriction to the fibres of $f_S$ is a degree-$d$ line bundle; its existence is ensured by the flatness of $f$ (cf. \cite{bosch} Theorem 2 Section 8.2).
It is clear by definition that $X_1=X$.
Observe that the natural morphism $X_d\to C$ is \'etale locally trivial because $f\colon X\to C$ is, which means that there exists a covering in the \'etale topology $\{U_i\to C\}$ such that $(X_d)_{U_i}\simeq U_i\times\Pic^d(F)$.
Notice 
that there is a natural morphism $\pi_d\colon P\to X_d$ functorially defined as follows: for every scheme $S$, to an isomorphism $X_S=F\times S$ it associates the line bundle corresponding to the divisor $dS_0$ which clearly has degree $d$ when restricted to a fibre of $f_S$ ($S_0$ denotes the fibre over $0\in F$ of the natural projection $\pi_F\colon F\times S\to F$).
It is clear that $\pi_d$ is $A_d$-invariant (the action of $A_d$ on $P$ is induced by the one of $A$) and, by the \'etale local triviality of $X_d\to C$, is an $A_d$-torsor, in particular we have that $X_d=P/A_d$.

From now on we fix $d$ as the degree of a multisection $\overline{C}$ of $f\colon X\to C$: this clearly defines a section of $X_d\to C$.
Denote by $G=A_d$ and by $C'=\pi_d^*C$: then  we have that the free action of $G$ on $P$ restricts to a free action of $G$ on $C'$ such that $C'/G=C$ and $C'\to C$ is a $G$-torsor.
Observe that $C'\to C$ is a reduction of structure group of the $A$-torsor $P\to C$ to $G$ (cf. \cite{balaji} Remark 2.5): it is clear by definition of reduction of structure group (ibidem Definition 2.4), that an $A$-torsor is trivial if and only if it admits a trivial reduction of structure group.
In particular, because we have $C'\times_C C'=C'\times G$ (cf. \cite{milne} Proposition III.4.1(a)), that is the $G$-torsor $C'\to C$ becomes the trivial torsor after base change to $C'$, the same holds for the $A$-torsor $P\to C$, i.e. $P\times_C C'\simeq A\times C'$.
From this we derive that $X\times_C C'\simeq F\times C'$ and the Proposition is proven. 
\end{proof}

\begin{remark}
\label{rem_notconn}
Notice that in the construction of Proposition \ref{propo_isoell}, nothing assures that $C'$ is a connected curve: as we are going to see in the next Lemma, when $X$ has maximal Albanese dimension, $C'$ can not be connected.
In particular we don't really need an \'etale cover with Galois group $G$ in order to trivialize the elliptic surface, but it suffices to take an Abelian subgroup of $F_d$ in order to do so. 
\end{remark}

\begin{lemma}
Let $E$ be an elliptic curve and $G$ a subgroup of $H\ltimes E_d$ where $H$ is the automorphism group of $E$ as an elliptic curve and $E_d$ is the $d$-torsion subgroup of $E$ for an integer $d$ coprime with the characteristic of the ground field such that $G$ is not contained in $E_d$. 
Assume that $G$ acts freely on a smooth curve $C$, consider the diagonal action on $C\times E$ (i.e. $g\cdot(c,e)=(g\cdot c,g\cdot e)$) and its naturally defined elliptic fibration $(C\times E)/G\to C/G$.
Then $q\bigl((C\times E)/G\bigr)=g(C/G)$, in particular $(C\times E)/G$ is not of maximal Albanese dimension\index{elliptic surface!of maximal Albanese dimension}.
\label{lemma_pic0maxalb}
\end{lemma}

\begin{proof}
By Riemann-Hurwitz formula (cf. \cite{hart} Corollary IV.2.4), it is clear that $E/G\simeq\P^1$.
Observe that for every curve $C$ and finite group $G$ acting on it, we have $\Pic^0(C)^G=\sum_{\eta\in\Pic^0(C)[l]^G}\bigl(f^*\Pic^0(C/G)+\eta\bigr)$ where $l=|G|$, by $\Pic^0(C)^G$\nomenclature{$\Pic^0(X)$}{the subgroup of $\Pic(X)$ of line bundles algebraically equivalent to zero on $X$}\nomenclature{$\Pic^0(X)^G$}{the subgroup of $\Pic^0(X)$ of line bundles fixed by the action of a group $G$ on $X$}\nomenclature{$\Pic(X)^G$}{the subgroup of $\Pic(X)$ of line bundles fixed by the action of a group $G$ on $X$} we mean the subgroup of line bundles of degree zero fixed by the action of $G$, by $\Pic^0(C)[l]^G$ we mean the torsion points of $\Pic^0(C)$ of order the cardinality of $G$ fixed by the action of $G$ and $f\colon C\to C/G$ is the natural projection.
Indeed, for every Cartier divisor $D$ on $C$ fixed by the action of $G$ (here we are using the natural identification between line bundles and Cartier divisor), we have that $f^*f_*D=\sum_{g\in G}g^*D\equiv lD$ (cf. \cite{hart} Exercise IV.2.6 for the definition of $f_*$ for Cartier divisors).
If moreover $D$ has degree zero, the same is true for  $f_*D$: hence, because $\Pic^0(C/G)$ is a divisible group, there exists a Cartier divisor $\overline{D}$ on $C/G$ such that $l\overline{D}\equiv f_*D$, i.e. $D\equiv f^*\overline{D}+\eta$ for a suitable $\eta\in\Pic^0(C)^G[l]$.

Now consider the following commutative diagram
\begin{equation}
\begin{tikzcd}
C\times E\arrow{r}{\alpha}\arrow{d}{\pi_C} & (C\times E)/G \arrow{r}{\beta}\arrow{d}{\pi_{C/G}} & C/G\times E/G \arrow{d}{p_{C/G}} \\
C \arrow{r}{f} & C/G\arrow{r}{Id} & C/G.
\end{tikzcd}
\label{eq_albpicmax}
\end{equation}

By universal property of the Albanese variety, it easily follows that $\alb_{\alpha}$, $\alb_{\beta}$ and $\alb_{\beta\circ\alpha}$ are surjective morphisms of Abelian varieties (cf. \cite{beauville_1983} Remark V.14(3)).
By the theory of dual Abelian varieties (cf. \cite{polishchuk} Theorem 10.1 and Exercise 10.1), this shows that their dual morphisms $\alpha^*\colon\Pic^0_{(C\times E)/G}\to\Pic^0_{C\times E}$, $\beta^*\colon\Pic^0_{C/G\times E/G}\to\Pic^0_{(C\times E)/G}$ and $(\beta\circ\alpha)^*\colon\Pic^0_{C/G\times E/G}\to\Pic^0_{C\times E}$ are isogenies on their images.

We know that
\begin{equation}
\begin{split}
\alpha^*\Pic^0((C\times E)/G)\subseteq \Pic^0(C\times E)^G=\Pic^0(C)^G\times \Pic^0(E)^G=\\
=\sum_{\eta\in\Pic^0(C)[l]^G}\bigl(f^*\Pic^0(C/G)+\eta\bigr)\times \sum_{\eta\in\Pic^0(E)[l]^G}\bigl(g^*\Pic^0(E/G)+\eta\bigr)=\\
=\sum_{\eta\in\Pic^0(C\times E)[l]^G}\bigl(\alpha^*\circ\beta^*\Pic^0(C/G\times E/G)+\eta\bigr) \quad \quad \quad \quad
\end{split}
\label{eq_naltra}
\end{equation}
where $g\colon E\to E/G$ is the natural projection.
%
In particular, because $\alpha^*\Pic^0((C\times E)/G)$ (respectively $\sum_{\eta\in\Pic^0(C)[l]^G}\bigl(\alpha^*\circ\beta^*\Pic^0(C/G\times E/G)+\eta\bigr)$) is the set of rational points of the image of $\alpha^*\colon\Pic^0_{(C\times E)/G}\to\Pic^0_{C\times E}$ (of finite copies of $(\beta\circ\alpha)^*\colon \Pic^0_{C/G\times E/G}\to \Pic^0_{C\times E}$ respectively), we have proved that the dimension of the Picard variety of $(C\times E)/G$ is smaller than or equal to the dimension of the Picard variety of $C/G\times E/G$.
This means that $q((C\times E)/G)\leq g(C/G)$ and equality follows by surjectivity of $\alb_{\beta}$, i.e. $\Alb((C\times E)/G)\to\jac(C/G)$ is an isogeny.
In particular, because the composition $\alb_{\pi_{C/G}}\circ\alb_{(C\times E)/G}\colon(C\times E)/G\to\jac(C/G)$ factors through $\pi_{C/G}\colon (C\times E)/G\to C/G$, we have that $(C\times E)/G$ has not maximal Albanese dimension.
\end{proof}

\begin{corollary}
\label{coro_isoell}
Let $f\colon X\to C$ be an elliptic fibration as in Lemma \ref{lemma_isoell} with maximal Albanese dimension and suppose that there exists a line bundle $L$ such that $L.F=d>0$ is coprime with the characteristic of the ground field.
Then $f\colon X\to C$ becomes trivial after a suitable \'etale base change.
More precisely, there exists a subgroup $G$ of the $d$-torsion points of $F$ acting by translation on $F$, freely on a smooth curve $C'$ and diagonally on $C'\times F$ such that 
\begin{equation}
\begin{tikzcd}
C'\times F\arrow{d}{}\arrow{r}{} & X \arrow{d}{f}\\
C'\arrow{r}{} & C
\end{tikzcd}
\label{eq_isoell+0}
\end{equation}
is commutative with $C'/{G}=C$, $(C'\times F)/{G}=X$ and the horizontal arrows are the quotients by ${G}$.

In particular, if $q(X)=2$, then $X$ is an Abelian variety.
\end{corollary}

\begin{proof}
Let $n\gg 0$ be a sufficiently big integer and consider $L+nF$.
Because the canonical bundle of an elliptic surface is numerically equivalent to a rational multiple $mF$ of a fibre (cf. Theorem \ref{teo_canellsur}) we have that $K_X-L-nF\sim_{num}-L+(m-n)F$ is not effective: hence, by Serre duality, $h^2(X,L+nF)=h^0(X,K_X-L-nF)=0$.
In particular,  by Riemann-Roch \ref{teo_rrss}, we obtain
\begin{equation}
\begin{split}
h^0(X,L+nF)\geq\chi(L+nF)=&\chi(\O_X)+\frac{1}{2}(L+nF)(L+(n-m)F)=\\
=\chi(\O_X)+&\frac{1}{2}(L^2+(2n-m)d)>0
\end{split}
\label{eq_forseserve}
\end{equation}
for a sufficiently big $n$.
In particular we have proved that $f\colon X\to C$ admits a multisection of degree $d$ coprime with the characteristic, hence, by Proposition \ref{propo_isoell}, we may assume that $X\to C$ becomes trivial after an \'etale Galois cover whose Galois group is $G'=H\ltimes F_d$ where $H$ is the automorphism group of $F$ as an elliptic curve and $F_d$ is the subgroup of $d$-torsion points of $F$.

In particular there exists an \'etale cover $C'\to C$ with Galois group $G'$ such that $X\to C$ becomes trivial after base change to $C'$.
Now let $\widetilde{C}$ be a connected component of $C'$ and consider the induced morphism $\widetilde{C}\to C$: this is a Galois morphism with Galois group ${G}=\{g\in G'\ |\ g(\widetilde{C})=\widetilde{C}\}$.

By Lemma \ref{lemma_pic0maxalb}, we get that ${G}\subseteq F_d$ otherwise $X$ would not have maximal Albanese dimension: in particular ${G}$ is Abelian and acts freely on $F$.

Observe that if $q(X)=2$ then $g(C)=1$ by Lemma \ref{lemma_ellalbjac} and in particular $X$, after an \'etale base change, becomes an Abelian variety, hence it is Abelian itself.
\end{proof}

An immediate Corollary is the following.

\begin{corollary}
\label{coro_isotrivell}
Let $f\colon Y\to C$ be a smooth elliptic fibration of maximal Albanese dimension for which there exists a line bundle $L$ with $L.F=1$ where $F$ is a fibre of $f$.
Then $Y\simeq C\times F$.
\end{corollary}

At the end of this Section we recall the following result on elliptic surfaces.

\begin{lemma}[\cite{Mukai} Proposition 3.3]
\label{lemma_kodvan}
Let $f\colon X\to C$ be an elliptic surface.
Then Kodaira vanishing Theorem holds, i.e. for every ample line bundle $L$ we have that $h^1(X,L^{-1})=0$.
\end{lemma}

\section{Picard group of $C\times E$}
\label{sec_piccxe}
\index{elliptic surface!product}
The Picard group of the product of two curves is known\index{Picard!group of a product elliptic surface}. Here we give a proof of how to compute it in the case one of the two curves is elliptic. We also consider the equivariant setting: that is we see how the formula for the Picard group of the product behaves with respect to the group action, in the case there is a finite Abelian group acting freely on both curves.

Denote by $\Hom_{c_0}(C,E)$\nomenclature{$\Hom_{C_0}(C,E)$}{the group of morphisms between a curve $C$ and an elliptic curve $E$ which sends $c_0$ to the origin of $E$} the group of morphisms between $C$ and $E$ for which the image of $c_0\in C$ is the origin of the elliptic curve (the group structure is given by that of $E$). 
Denote by $i_c\colon E\to C\times E$ the map defined by $e\mapsto (c,e)$ and by $E_c$ its image.
Similarly denote by $j_e\colon C\to C\times E$ the map defined by $c\mapsto (c,e)$ and by $C_e$ its image.

\begin{proposition}
\label{propo_piccxe}
In the above settings we have the following split exact sequence of groups:
\begin{equation}
0\rightarrow \Pic (C)\times \Pic(E) \xrightarrow{\alpha} \Pic(C\times E) \xrightarrow{\beta} \Hom_{c_0}(C,E) \rightarrow 0,
\label{eq_piccxe}
\end{equation}
where $\beta$ is defined by $\beta(D)(c)=i_c^*(D)-i_{c_0}^*(D)$ (the sum is that of the elliptic curve) and the section $s$ of $\beta$ is given by
\[s\colon \Hom_{c_0}(C,E) \to  \Pic(C\times E)\quad s(f)=\Gamma_f-C_0-\sum_{c\in f^{-1}(0)} a_cE_c,\]
\nomenclature{$\Gamma_f$}{the graph of a morphism $f$} where $a_c$ is  the multiplicity of $f$ at $c$.
\end{proposition}

Before the proof we recall the See-saw Principle (cf. \cite{mumford} Corollary II.5.6) which we state here in the form we are going to use.

\begin{theorem}[See-saw Principle\index{see-saw principle}]
\label{saw}
Let $A$ and $B$ be two smooth curves, $\pi_2\colon A\times B\to B$ be the second projection and let $L\in\Pic(A\times B)$ be a line bundle such that $\restr{L}{A\times\{b\}}$ is trivial for every $b\in B$ then there exists $L_B\in\Pic(B)$ such that $L=\pi_2^*L_B$. 
\end{theorem}

\begin{proof}[Proof of Proposition \ref{propo_piccxe}]
Let $\pi_C\colon C\times E\to C$ and $\pi_E\colon C\times E\to E$ be the two projections. The morphism $\alpha:=\pi_C^*\times\pi_E^*$ is defined as $\alpha(L_C,L_E)=\pi_C^*L_C\otimes_{\O_{C\times E}} \pi_E^*L_E$. 
The pull-back $\pi_C^*$, respectively $\pi_E^*$, is injective (the composition with the restriction to $C\times\{0\}$, respectively to $\{c_0\}\times E$, is the identity), and $\im(\pi^*_C)\cap\im(\pi^*_E)=\{0\}\subseteq\Pic(C\times E)$, hence also $\alpha$ is injective. 
To define $\beta$ we first need to notice that 
\[\Hom_{c_0}(C,E)\simeq\Hom(\jac(C),E),\]
where the isomorphism, which we will call $\phi$, is induced  by the Abel-Jacobi map $AJ_C\colon C\to\jac(C)$ and the canonical isomorphism of $E$ with its Jacobian variety:
\begin{equation}
\label{eq_cxejac}
\begin{tikzcd}
C\arrow{d}{AJ_C} \arrow{dr}{f} &  \\
\jac(C) \arrow[swap]{r}{\phi(f)} & E\simeq\jac(E).
\end{tikzcd}
\end{equation}

 By this, it is possible to define 
\begin{equation}
\label{eq_beta}
\beta\colon\Pic(C\times E)\to \Hom_{c_0}(C,E)
\end{equation}
as
\begin{equation}
\label{eq_betadef}
\O(D)\mapsto \phi^{-1}(f_D)
\end{equation}
where $f_D$ is the morphism in $\Hom(\jac(C),E)$ defined as 
\begin{equation}
f_D(P)=(\pi_E)_*(\pi^*_CP\cdot D):
\label{eq_betadef2}
\end{equation}
we refer to \cite{fulton} for pushforward and pull-back of cycles of a proper flat morphism (chapter 1), for intersection product of cycles of a smooth variety (chapter 8) and their behaviour with respect to linear equivalence.

Actually $f_D$ is a morphism from $\Pic(C)$ to $\Pic(E)$ and one has to show that its restriction to $\Pic^0(C)$ has image inside $\Pic^0(E)\simeq E$. 
This goes as follows: let $c-c_0\in \Pic^0(C)$ where $c$ is a point of $C$. Then $\pi_C^*(c-c_0)=E_c-E_{c_0}$ is a divisor inside $C\times E$. 
Because all the fibres are algebraically equivalent the intersection numbers $E_c.D$ and $E_{c_0}.D$ coincide. 
From this it follows that $f_D(c-c_0)\in\Pic^0(E)$. 
Moreover, it is immediate by definition, that $\beta$ is a morphism of groups. 
Let $i_c\colon E\to C\times E$ be as in the statement, we see that $i_c^*(D)=(\pi_E)_*(D.E_c)$ where $E_c$ is the fibre of $\pi_C$ over $c\in C$. 
By this 
\[\beta(D)(c)=\phi^{-1}(f_D(c-c_0))=\phi^{-1}((\pi_E)_*((E_c-E_{c_0})\cdot D))=i_c^*(D)-i_{c_0}^*(D).\]

Now let $E_c$ be  the fibre of $\pi_C$ above $c$. Because the intersection of two fibres is zero, we have that $f_{E_c}(P)=(\pi_E)_*(\pi_C^*P.E_c)=0$ for every $P\in\Pic^0(C)$. 
Let $C_e$ be the fibre of $\pi_E$ above $e$, we have that $\restr{\pi^*_CP}{E_c}\in\Pic^0(E_c)$ for every $P\in\Pic^0(C)$ and hence, when pushing forward to $E$, we get $f_{E_c}(P)=0$.
By additivity,  $f_{D}(P)=0$ for every  $D=\pi_C^*(D_C)$ or $D=\pi_E^*(D_E)$.
In particular we have shown that $\beta\circ\alpha=0$.

Now let $f\colon C\to E$ be a morphism which sends $c_0$ to $0$ and let $\Gamma_{f+e}$ be the graph of the morphism $(f+e)(c)=f(c)+e$ inside the product $C\times E$, where $e$ is a point of $E$. We see that 
\begin{equation}
\label{eq_invtranslation}
\beta(\Gamma_{f+e})(c)=i_c^*(\Gamma_{f+e})-i_{c_0}^*(\Gamma_{f+e})=(f(c)+e)-(f(c_0)+e)=f(c);
\end{equation}
thus the morphism $\beta$ is surjective. Notice also that, by a similar argument, $\beta$ is invariant under translation by $e\in E$ for a general divisor $D\in\Pic(C\times E)$.

To show the exactness in the middle, suppose that $D$ is a divisor for which $\beta(D)(c)=0$ for all $c\in C$: we want to show that it is linearly equivalent to a sum of vertical and horizontal fibres.  $\beta(D)(c)=0$ means that $i_c^*(D)-i_{c_0}^*(D)$ is equivalent to zero in the Picard group.  Let $e_1,\ldots,e_n$ be the points in $E$ such that $D\cap E_{c_0}=\{(e_i,c_0)\ |\ i=1,\ldots,n\}$ with multiplicity $a_i$. Then $i_c^*(D-\sum_{i=1}^n a_iC_{e_i})=0$  for all $c \in C$. This means that the restriction of $D-\sum_{i=1}^n a_iC_{e_i}$ is trivial in all the fibres of $\pi_C$. Hence, using the See-saw Principle, $D-\sum_{i=1}^n a_i C_{e_i}\sim_{lin}\sum_{j=1}^m b_jE_{c_j}$ which proves the exactness in the middle of Equation \ref{eq_piccxe}.

It remains to show that the sequence split.  
We consider the map
\[s\colon \Hom_{c_0}(C,E)\to \Pic(C\times E)\]
defined by
\[f\mapsto \Gamma_f-C_0-\quad\sum_{\mathclap{c\in f^{-1}(0)}}a_cE_c,\]
where $\Gamma_f$ is the graph of $f$ and $a_c$ is the multiplicity of $c$ in the fibre of $0$. 
It is clear that $s$ is a section of $\beta$; we show that $s$ is a morphism of groups, i.e. that
\[\Gamma_{f+g}-C_0-\qquad\sum_{\mathclap{c\in (f+g)^{-1}(0)}}a_cE_c\sim_{lin} \Gamma_f-C_0-\quad\sum_{\mathclap{c'\in f^{-1}(0)}}a_{c'}E_{c'}+\Gamma_g-C_0-\quad\sum_{\mathclap{c''\in g^{-1}(0)}} a_{c''}E_{c''} \]
or, equivalently,
\[\Gamma_{f+g}-\Gamma_f-\Gamma_g\sim_{lin} -C_0+\quad\sum_{\mathclap{c\in (f+g)^{-1}(0)}}a_cE_c-\quad\sum_{\mathclap{c'\in f^{-1}(0)}}a_{c'}E_{c'}-\quad\sum_{\mathclap{c''\in g^{-1}(0)}}a_{c''}E_{c''}.\]  

To prove this, we use once again the See-saw Principle. Consider, for every $c\in C$,  $i_c^*(\Gamma_{f+g}-\Gamma_f-\Gamma_g+C_0)=(f+g)(c)-f(c)-g(c)+0=\Bigl((f+g)(c)-0\Bigr)-\Bigl(f(c)-0\Bigr)-\Bigl(g(c)-0\Bigr)$. Using the Abel-Jacobi isomorphism of $E$ with its Jacobian, we obtain that $i_c^*(\Gamma_{f+g}-\Gamma_f-\Gamma_g+C_0)= 0$, in particular, by the See-saw Principle, we have that
\begin{equation}
\Gamma_{f+g}-\Gamma_f-\Gamma_g+C_0\sim_{lin}a_1E_1+\ldots+a_kE_k.
\label{eqsee-saw1}
\end{equation}

Similarly as before, we consider, for every $e\in E$, $j_e\colon C\to C\times E$ defined by $c\mapsto (c,e)$. Then
\[j_e^*\Biggl(\Gamma_{f+g}-\Gamma_f-\Gamma_g-\qquad\sum_{\mathclap{c\in (f+g)^{-1}(0)}}a_cE_c+\quad\sum_{\mathclap{c'\in f^{-1}(0)}}a_{c'}E_{c'}+\quad\sum_{\mathclap{c''\in g^{-1}(0)}}a_{c''}E_{c''}\Biggr)=\]
\[=(f+g)^*(e)-(f+g)^*(0)-f^*(e)+f^*(0)-g^*(e)+g^*(0)=\]
\[=(f+g)^*(e-0)-f^*(e-0)-g^*(e-0)\sim_{lin} 0.\]
Here the last equality follows from the fact that $f^*\colon\Pic^0(E)\to\Pic^0(C)$ corresponds to $f$ in the isomorphism of groups $\Hom_{c_0}(C,E)\xrightarrow{\phi}\Hom(\jac(C),\jac(E))\simeq\Hom(\jac(E),\jac(C))$: in particular $f^*+g^*=(f+g)^*\colon\Pic^0(E)\to\Pic^0(C)$.
From this we get
\begin{equation}
\begin{split}
\Gamma_{f+g}-\Gamma_f-\Gamma_g-\sum_{\mathclap{c\in (f+g)^{-1}(0)}}a_cE_c+\quad &\sum_{\mathclap{c'\in f^{-1}(0)}}a_{c'}E_{c'}+\quad\sum_{\mathclap{c''\in f^{-1}(0)}}a_{c''}E_{c''}\sim_{lin}\\
\sim_{lin} b_1C_1&+\ldots+b_lC_l.
\end{split}
\label{eqsee-saw2}
\end{equation}

Comparing Equations \ref{eqsee-saw1} and \ref{eqsee-saw2} and using the injectivity of $\alpha$, we obtain
\[\Gamma_{f+g}-\Gamma_f-\Gamma_g\sim_{lin} -C_0+\qquad\sum_{\mathclap{c\in (f+g)^{-1}(0)}}a_cE_c-\quad\sum_{\mathclap{c'\in f^{-1}(0)}}a_{c'}E_{c'}-\quad\sum_{\mathclap{c''\in g^{-1}(0)}}a_{c''}E_{c''}.\qedhere\]
\end{proof}

Now we give an equivariant version of Proposition \ref{propo_piccxe}.

\begin{proposition}
\label{propo_piccxeg}
Suppose there exists a finite Abelian group $G$ acting freely on $C$, $E$ and diagonally on $C\times E$ (i.e. $g\cdot(c,e)=(g\cdot c,g\cdot e)$). Then it is possible to give to $\Pic(C)$, $\Pic(E)$, $\Pic(C\times E)$ and $\Hom_{c_0}(C,E)$ a $G$-module structure such that
\begin{equation}
0\rightarrow \Pic (C)\times \Pic(E) \xrightarrow{\alpha} \Pic(C\times E) \xrightarrow{\beta} \Hom_{c_0}(C,E) \rightarrow 0,
\end{equation}
is an exact sequence of $G$-modules, where $\alpha$ and $\beta$ are the same morphisms defined in Proposition \ref{propo_piccxe}.
\end{proposition}
\begin{proof}
Recall that, if $X$ is a variety and $G$ a group acting on $X$, we can naturally define an action of $G$ on $\Pic(X)$ given by $g\cdot L=((g)^{-1})^*L$ where, with an abuse of notation, we are identifying $g$ with the corresponding automorphism of $X$; equivalently we easily see that this action is induced by $g\cdot D=(g^{-1})^*D=\{x\in X\ |\ \exists\ d\in D\ \text{with}\ g\cdot d=x\}$ at the level of divisors. 
Hence, it is clear that the morphism $\alpha$ in the statement is a morphism of $G$-modules. 
Moreover, we can consider $G$ as a finite subgroup of $E$ when considering its action on $E$ (hence we use the additive notation for this factor), while we use the multiplicative notation when considering $G$ acting on $C$. 
Clearly, the action of $G$ on the Picard groups of $C$, $E$ and $C\times E$ is faithful (because two points on a curve of genus greater than or equal to $1$ are never linearly equivalent) but not free (every divisor of the form $\sum_{g\in G}g\cdot D$ is fixed by every $g\in G$).

Now, we would like to give to $\Hom_{c_0}(C,E)$ the structure of a $G$-module such that the morphism $\beta$ preserves the $G$-module structure. 
In order to do this we look at how $G$ acts on divisors $\Gamma_f$ which are graphs  of functions $f$ in $\Hom_{c_0}(C,E)$. 
We see that
\[g\cdot \Gamma_f=\{(g c,f(c)+g)\ | \ c\in C\}=\{(c,f(g^{-1} c)+g)\ | \ c\in C\}=\Gamma_{f_g},\]
where $f_g$ is defined by $c\mapsto f(g^{-1}c)+g$. However, in general $f_g\notin \Hom_{c_0}(C,E)$, but if we take $g\cdot f(c)=f_g(c)-f_g(c_0)=f(g^{-1}c)-f(g^{-1}c_0)$, then we see that $g\cdot f\in\Hom_{c_0}(C,E)$. 
It is clear that $f\mapsto g\cdot f$ is a well defined action of $G$ on $\Hom_{c_0}(C,E)$: indeed the axioms are easily verified. 
We have already noticed that $\beta$ is invariant under translation by $e\in E$ (Equation \ref{eq_invtranslation}); moreover, by the splitting exact sequence \ref{eq_piccxe}, we know that every divisor $D$ on $C\times E$ is linearly equivalent to
\[\sum_{i}E_{c_i}+\sum_{j}C_{e_j}+\Gamma_{f}\]
with $f=\beta(\O(D))$ and suitable $c_i$ and $e_j$. These two facts imply that $\beta$ preserves the $G$-module structure.
\end{proof}


The following Lemma gives a necessary and sufficient condition for the quotient of a product elliptic surface to be trivial and we will use it in the construction of Example \ref{es4}.

\begin{lemma}
\label{lem2}
In the same settings as Proposition \ref{propo_piccxeg}, the elliptic fibration $(C\times E)/G\to C/G$  with general fibre $E$ is trivial if and only if there is a line bundle $L$ on $C\times E$ which is fixed by the action of $G$ for which $L.E=1$ where $E$ (by an abuse of notation) is a general fibre of the first projection.
\end{lemma}

\begin{proof}
Denote by $\pi_{C/G}\colon(C\times E)/G\to C/G$, $\pi_C\colon C\times E\to C$ the two elliptic fibrations and by $\phi\colon C\times E\to (C\times E)/G$ and $\psi\colon C\to C/G$ the two quotients by $G$. 
Suppose that  $(C\times E)/G\to C/G$ is a product elliptic fibration: in particular there exists a section $i_{C/G}\colon C/G\to (C\times E)/G$ of $\pi_{C/G}$.
Then, the pull-back $i_C$ of $i_{C/G}$ to $C\times E$ is a section of $\pi_C$, i.e.
\begin{equation}
\begin{tikzcd}
C\arrow[dr, phantom, "\square"]\arrow[bend right=80, swap]{dd}{Id}\arrow{r}{\psi}\arrow{d}{i_C} & C/G\arrow{d}{i_{C/G}}\arrow[bend left=80]{dd}{Id}\\
C\times E \arrow[dr, phantom, "\square"] \arrow{r}{\phi} \arrow{d}{\pi_C} & (C\times E)/G \arrow{d}{\pi_{C/G}}\\
C\arrow{r}{\psi} & C/G
\end{tikzcd}
\label{eq_variabella}
\end{equation}
commutes.
This means that  $\phi^*\O_{(C\times E)/G}{(\im(i_{C/G}))}$ is a line bundle on $C\times E$ fixed by the action of $G$ such that $\phi^*\O_{(C\times E)/G}{(\im(i_{C/G}))}.E=1$.

 Now let $L$ be a line bundle on $C\times E$ fixed by the action of $G$ with $L.E=1$  and let $f=\beta(L)$: because $L$ is fixed by the action of $G$, thanks to Proposition \ref{propo_piccxeg}, we can say that also $f$ is. In particular we obtain 
\begin{equation*}
f(g^{-1}c)-f(g^{-1}c_0)=f(c)-f(c_0)=f(c),
\end{equation*}
for every $g\in G$ or, equivalently, if we denote by $e_g=f(gc_0)$,
\begin{equation}
\label{eqG}
f(gc)-f(c)=e_g\in E.
\end{equation}
Observe that
\begin{equation}
e_{g+h}=f((g+h)c_0)=f(g(hc_0))=f(gc_0)+f(hc_0)=e_g+e_h
\label{eq_eqG}
\end{equation}
(notice that for $f(g(hc_0)=f(gc_0)+f(hc_0)$ we are using Equation \ref{eqG} with $c=hc_0$ and $e_g=f(gc_0)$).

Let $\alpha_f\in \aut(C\times E)$ defined by  
\[(c,e)\mapsto(c,e-f(c))\]
and notice that if fixes the fibres of the natural projection $\pi_C\colon C\times E\to C$.
This gives the following diagram
\[
\begin{tikzcd}
\Pic(C\times E) \arrow{r}{\beta}\arrow{d}{(\alpha_f^{-1})^*} & \Hom_{c_0}(C,E) \arrow{r}{}\arrow{d}{h} & 0\\
\Pic(C\times E) \arrow{r}{\beta'} & \Hom_{c_0}(C,E) \arrow{r}{} & 0,
\end{tikzcd}
\]
which is commutative for line bundles with $L.E=1$ (indeed such line bundles, by Proposition \ref{propo_piccxe}, are linearly equivalent to $\Gamma_f+\pi^*_CB$ and the isomorphism $\alpha_f$ turns the translated graphs $\Gamma_f+e$ into fibres of the second projection) and where the map $h$ is bijective and defined by $h(\psi)=\psi-f$. Notice that $\beta'\circ(\alpha_f^{-1})^*(L)=h\circ\beta(L)=0$.    We would like to know how $G$ acts on $C\times E$ after this change of coordinates, i.e. what is $\alpha_f\circ g\circ \alpha^{-1}_f(c,e)$: we see that
\begin{equation}
\begin{split}
\alpha_f\circ g\circ \alpha_f^{-1}(c,e)=&\alpha_f\circ g(c,e+f(c))=\alpha_f(gc,e+f(c)+g)=\\
=(gc,e+&f(c)+g-f(gc))=(gc,e+g-e_g),
\end{split}
\label{eq_alfafinv}
\end{equation}
where the last equality follows by Equation \ref{eqG}. Notice that $(\alpha_f^{-1})^*(L).E=(\alpha_f^{-1})^*(L).(\alpha_f^{-1})^*(E)=1$ and that $(\alpha_f^{-1})^*(L)$ is still $G$-invariant after conjugating the action of $G$ with $\alpha_f$, i.e. $(\alpha_f^{-1})^*g(\alpha_f)^*(\alpha_f^{-1})^*L=(\alpha_f^{-1})^*L$. 

The set of elements $G'=\{g\in G\ |\ g=e_g\}$ is a subgroup of $G$: indeed, by Equation \ref{eq_eqG}, we have that if $g,h\in G'$, then 
\begin{equation}
e_{g+h}=e_g+e_h=g+h
\label{eq_ambnonso}
\end{equation}
and, in particular, $g+h\in G'$ (this is enough to show that $G'$ is a subgroup because $G$ is a finite group, indeed $g^{-1}=(n-1)g$ where $n$ is the order of $g$ in $G$).

It is immediate by Equation \ref{eq_alfafinv}, that $(C\times E)/G'=C/G'\times E$. Hence, suppose by contradiction that there exists $\gamma\in G$ such that $\gamma-e_\gamma\neq 0$ or, equivalently that $G'\neq G$: after quotienting by $G'$ we may assume that $\gamma-e_\gamma\neq 0$ for all $\gamma\in G$ and, in particular, the action of $G$ on $E$ is still free after the change of coordinates given by $\alpha_f$. Because $\beta'((\alpha_f^{-1})^*L)=0$ and $(\alpha_f^{-1})^*L.E=1$, we have that 
\[(\alpha_f^{-1})^*L=C_e+\pi_C^*(B),\]
where $C_e$ is the fibre of the second projection over $e\in E$ and $B\in \Pic(C)$.  Because $(\alpha_f^{-1})^*(L)$ is $G$-invariant we can conclude that
\[
\begin{split}
(\alpha_f^{-1})^*L=(\alpha_f^{-1})^*(\gamma^{-1})^*(\alpha_f)^*(\alpha_f^{-1})^*L=C_{e+\gamma-e_\gamma}&+\pi_C^*(B');
\end{split}
\] 
in particular it follows $\gamma-e_\gamma=0$, a contradiction.
\end{proof}

\begin{remark}
The "if" part of of Lemma \ref{lem2} could also be derived more directly by Corollary \ref{coro_isoell} using the fact that a line bundle fixed by the action of an Abelian group passes to the quotient.
We decided to write this rather than the other proof because this is how we first proved this result.
Moreover this proof is more elementary and is an application of Proposition \ref{propo_piccxeg}, which we think is interesting on its own and we don't use anywhere else in this Thesis. 
\label{rem_lem2}
\end{remark}

\section{Surfaces of general type}
\label{sec_gentype}
In\index{surface!of general type|see{surface of general type}} this section we recall all the inequalities that the numerical invariants $K_X^2$ and $\chi(\O_X)$ of minimal surfaces of general type satisfy.
As usual, the situation in positive characteristic is more complicated but the further assumption of maximal Albanese dimension makes things easier:
indeed, most of the counterexamples that arise in positive characteristic are uniruled surfaces.

Almost immediately by definition, we know that for a minimal surface of general type we have that $K_X^2>0$ (\cite{beauville_1983} Proposition X.1). 
If the characteristic of the ground field is zero it is known that both the Euler characteristic of the structure sheaf $\chi(\O_X)$ of $X$ and the topological Euler characteristic $e(X)$ are strictly greater than $0$ (cf. \cite{beauville_1983}  Theorem X.4 and \cite{barth} Proposition VII.2.4).
If the characteristic of the ground field is positive, it is known that this is no longer true: there exist examples of surfaces with $e(X)\leq 0$ for every characteristic $p\geq 5$ (cf. \cite{raynaud}).
On the other hand, it has been recently proved that $\chi(\O_X)$ is still positive over algebraically closed fields of any characteristic (cf. \cite{gu2019slope} Corollary 3.4).
Moreover, as we have already observed in the proof of Theorem \ref{teo_shepbarr}, the dimension of the image of the Albanese morphism of a minimal surface $X$ with $e(X)<0$ is a curve and the fibres are rational (cf. \cite{shepherdgeo} Theorem 2 and 7).

Recall that for a surface of general type the Noether's inequality\index{Noether!inequality for surfaces of general type} \index{surface of general type!Noether' inequality|see{Noether}} holds (cf. \cite{barth} Theorem VII.3.1 in characteristic 0 and \cite{liedtkec2} Theorem 2.1 for positive characteristic):
\begin{equation}
K_X^2\geq 2p_g(X)-4.
\label{eq_noetine}
\end{equation}
From this we easily derive 
\begin{equation}
K_X^2\geq 2\chi(\O_X)+2q'(X)-6\geq 2\chi(\O_X)-6
\label{eq_noetin2}
\end{equation}
which we still call Noether's inequality.
Surfaces which satisfy equality in \ref{eq_noetine} are called Horikawa surfaces\index{surface!Horikawa}.
These surfaces also satisfy the equality in \ref{eq_noetin2}, in particular $q'(X)=0$, have trivial algebraic fundamental group and the Picard scheme is reduced (cf. \cite{horikawa} Theorem 3.4 and \cite{barth} Corollary VII.3.3 for characteristic $0$ and \cite{liedtkec2} Proposition 3.7 for positive characteristic).

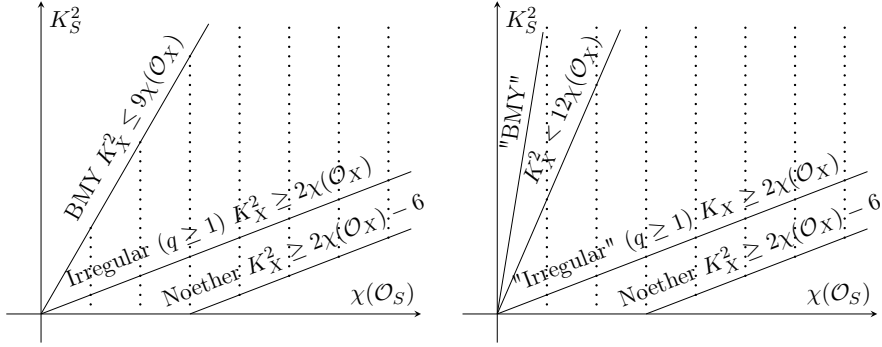
\begin{figure}%
\centering
\begin{tikzpicture}[shorten >=1pt,node distance=2cm,on grid,auto, xscale=0.8, yscale=0.8, every node/.style={scale=1}]
    \begin{scope}
        \begin{axis}[ticks=none, xmin=0, xmax=7, ymin=0, ymax=30, axis lines=middle, xlabel=$\chi(\O_S)$, ylabel=$K^2_S$, enlargelimits]

	\draw (0,0) -- (3.4,30.6) node [sloped, midway, above] {\quad \quad \quad BMY $K_X^2\leq 9\chi(\O_X)$};
	
	
	\draw (3,0) -- (7.5,9) node [sloped, midway, above] {Noether $K_X^2\geq 2\chi(\O_X)-6$};
	
	
	\draw (0,0) -- (7.5,15) node [sloped,midway, above] {Irregular ($q\geq1$) $K_X^2\geq2\chi(\O_X)$};
	
	
	\filldraw  (1,1) circle (0.3pt);
	\filldraw  (1,2) circle (0.3pt);
	\filldraw  (1,3) circle (0.3pt);
	\filldraw  (1,4) circle (0.3pt);
	\filldraw  (1,5) circle (0.3pt);
	\filldraw  (1,6) circle (0.3pt);
	\filldraw  (1,7) circle (0.3pt);
	\filldraw  (1,8) circle (0.3pt);
	\filldraw  (1,9) circle (0.3pt);
	\filldraw  (2,1) circle (0.3pt);
	\filldraw  (2,2) circle (0.3pt);
	\filldraw  (2,3) circle (0.3pt);
	\filldraw  (2,4) circle (0.3pt);
	\filldraw  (2,5) circle (0.3pt);
	\filldraw  (2,6) circle (0.3pt);
	\filldraw  (2,7) circle (0.3pt);
	\filldraw  (2,8) circle (0.3pt);
	\filldraw  (2,9) circle (0.3pt);
	\filldraw  (2,10) circle (0.3pt);
	\filldraw  (2,11) circle (0.3pt);
	\filldraw  (2,12) circle (0.3pt);
	\filldraw  (2,13) circle (0.3pt);
	\filldraw  (2,14) circle (0.3pt);
	\filldraw  (2,15) circle (0.3pt);
	\filldraw  (2,16) circle (0.3pt);
	\filldraw  (2,17) circle (0.3pt);
	\filldraw  (2,18) circle (0.3pt);
	\filldraw  (3,1) circle (0.3pt);
	\filldraw  (3,2) circle (0.3pt);
	\filldraw  (3,3) circle (0.3pt);
	\filldraw  (3,4) circle (0.3pt);
	\filldraw  (3,5) circle (0.3pt);
	\filldraw  (3,6) circle (0.3pt);
	\filldraw  (3,7) circle (0.3pt);
	\filldraw  (3,8) circle (0.3pt);
	\filldraw  (3,9) circle (0.3pt);
	\filldraw  (3,10) circle (0.3pt);
	\filldraw  (3,11) circle (0.3pt);
	\filldraw  (3,12) circle (0.3pt);
	\filldraw  (3,13) circle (0.3pt);
	\filldraw  (3,14) circle (0.3pt);
	\filldraw  (3,15) circle (0.3pt);
	\filldraw  (3,16) circle (0.3pt);
	\filldraw  (3,17) circle (0.3pt);
	\filldraw  (3,18) circle (0.3pt);
	\filldraw  (3,19) circle (0.3pt);
	\filldraw  (3,20) circle (0.3pt);
	\filldraw  (3,21) circle (0.3pt);
	\filldraw  (3,22) circle (0.3pt);
	\filldraw  (3,23) circle (0.3pt);
	\filldraw  (3,24) circle (0.3pt);
	\filldraw  (3,25) circle (0.3pt);
	\filldraw  (3,26) circle (0.3pt);
	\filldraw  (3,27) circle (0.3pt);
	\filldraw  (4,2) circle (0.3pt);
	\filldraw  (4,3) circle (0.3pt);
	\filldraw  (4,4) circle (0.3pt);
	\filldraw  (4,5) circle (0.3pt);
	\filldraw  (4,6) circle (0.3pt);
	\filldraw  (4,7) circle (0.3pt);
	\filldraw  (4,8) circle (0.3pt);
	\filldraw  (4,9) circle (0.3pt);
	\filldraw  (4,10) circle (0.3pt);
	\filldraw  (4,11) circle (0.3pt);
	\filldraw  (4,12) circle (0.3pt);
	\filldraw  (4,13) circle (0.3pt);
	\filldraw  (4,14) circle (0.3pt);
	\filldraw  (4,15) circle (0.3pt);
	\filldraw  (4,16) circle (0.3pt);
	\filldraw  (4,17) circle (0.3pt);
	\filldraw  (4,18) circle (0.3pt);
	\filldraw  (4,19) circle (0.3pt);
	\filldraw  (4,20) circle (0.3pt);
	\filldraw  (4,21) circle (0.3pt);
	\filldraw  (4,22) circle (0.3pt);
	\filldraw  (4,23) circle (0.3pt);
	\filldraw  (4,24) circle (0.3pt);
	\filldraw  (4,25) circle (0.3pt);
	\filldraw  (4,26) circle (0.3pt);
	\filldraw  (4,27) circle (0.3pt);
	\filldraw  (4,28) circle (0.3pt);
	\filldraw  (4,29) circle (0.3pt);
	\filldraw  (4,30) circle (0.3pt);
	\filldraw  (5,4) circle (0.3pt);
	\filldraw  (5,5) circle (0.3pt);
	\filldraw  (5,6) circle (0.3pt);
	\filldraw  (5,7) circle (0.3pt);
	\filldraw  (5,8) circle (0.3pt);
	\filldraw  (5,9) circle (0.3pt);
	\filldraw  (5,10) circle (0.3pt);
	\filldraw  (5,11) circle (0.3pt);
	\filldraw  (5,12) circle (0.3pt);
	\filldraw  (5,13) circle (0.3pt);
	\filldraw  (5,14) circle (0.3pt);
	\filldraw  (5,15) circle (0.3pt);
	\filldraw  (5,16) circle (0.3pt);
	\filldraw  (5,17) circle (0.3pt);
	\filldraw  (5,18) circle (0.3pt);
	\filldraw  (5,19) circle (0.3pt);
	\filldraw  (5,20) circle (0.3pt);
	\filldraw  (5,21) circle (0.3pt);
	\filldraw  (5,22) circle (0.3pt);
	\filldraw  (5,23) circle (0.3pt);
	\filldraw  (5,24) circle (0.3pt);
	\filldraw  (5,25) circle (0.3pt);
	\filldraw  (5,26) circle (0.3pt);
	\filldraw  (5,27) circle (0.3pt);
	\filldraw  (5,28) circle (0.3pt);
	\filldraw  (5,29) circle (0.3pt);
	\filldraw  (5,30) circle (0.3pt);
	\filldraw  (6,6) circle (0.3pt);
	\filldraw  (6,7) circle (0.3pt);
	\filldraw  (6,8) circle (0.3pt);
	\filldraw  (6,9) circle (0.3pt);
	\filldraw  (6,10) circle (0.3pt);
	\filldraw  (6,11) circle (0.3pt);
	\filldraw  (6,12) circle (0.3pt);
	\filldraw  (6,13) circle (0.3pt);
	\filldraw  (6,14) circle (0.3pt);
	\filldraw  (6,15) circle (0.3pt);
	\filldraw  (6,16) circle (0.3pt);
	\filldraw  (6,17) circle (0.3pt);
	\filldraw  (6,18) circle (0.3pt);
	\filldraw  (6,19) circle (0.3pt);
	\filldraw  (6,20) circle (0.3pt);
	\filldraw  (6,21) circle (0.3pt);
	\filldraw  (6,22) circle (0.3pt);
	\filldraw  (6,23) circle (0.3pt);
	\filldraw  (6,24) circle (0.3pt);
	\filldraw  (6,25) circle (0.3pt);
	\filldraw  (6,26) circle (0.3pt);
	\filldraw  (6,27) circle (0.3pt);
	\filldraw  (6,28) circle (0.3pt);
	\filldraw  (6,29) circle (0.3pt);
	\filldraw  (6,30) circle (0.3pt);
	\filldraw  (7,8) circle (0.3pt);
	\filldraw  (7,9) circle (0.3pt);
	\filldraw  (7,10) circle (0.3pt);
	\filldraw  (7,11) circle (0.3pt);
	\filldraw  (7,12) circle (0.3pt);
	\filldraw  (7,13) circle (0.3pt);
	\filldraw  (7,14) circle (0.3pt);
	\filldraw  (7,15) circle (0.3pt);
	\filldraw  (7,16) circle (0.3pt);
	\filldraw  (7,17) circle (0.3pt);
	\filldraw  (7,18) circle (0.3pt);
	\filldraw  (7,19) circle (0.3pt);
	\filldraw  (7,20) circle (0.3pt);
	\filldraw  (7,21) circle (0.3pt);
	\filldraw  (7,22) circle (0.3pt);
	\filldraw  (7,23) circle (0.3pt);
	\filldraw  (7,24) circle (0.3pt);
	\filldraw  (7,25) circle (0.3pt);
	\filldraw  (7,26) circle (0.3pt);
	\filldraw  (7,27) circle (0.3pt);
	\filldraw  (7,28) circle (0.3pt);
	\filldraw  (7,29) circle (0.3pt);
	\filldraw  (7,30) circle (0.3pt);

	\end{axis}
    \end{scope}

    \begin{scope}[xshift=7.5cm]
        \begin{axis}[ticks=none, xmin=0, xmax=7, ymin=0, ymax=30, axis lines=middle, xlabel=$\chi(\O_S)$, ylabel=$K^2_S$, enlargelimits]
	
	\draw (0,0) -- (0.93,29.76) node [sloped, midway, above] {\quad \quad \quad \quad \quad "BMY"};
	
	\draw (0,0) -- (2.5,30) node [sloped, midway, above] {\quad \quad \quad \quad \quad $K_X^2<12\chi(\O_X)$};
	
	
	\draw (3,0) -- (7.5,9) node [sloped, midway, above] {Noether $K_X^2\geq2\chi(\O_X)-6$};
	
	
	\draw (0,0) -- (7.5,15) node [sloped,midway, above] {"Irregular" ($q\geq 1$)  $K_X\geq2\chi(\O_X)$};
	

  \filldraw  (1,1) circle (0.3pt);
	\filldraw  (1,2) circle (0.3pt);
	\filldraw  (1,3) circle (0.3pt);
	\filldraw  (1,4) circle (0.3pt);
	\filldraw  (1,5) circle (0.3pt);
	\filldraw  (1,6) circle (0.3pt);
	\filldraw  (1,7) circle (0.3pt);
	\filldraw  (1,8) circle (0.3pt);
	\filldraw  (1,9) circle (0.3pt);
	\filldraw  (1,10) circle (0.3pt);
	\filldraw  (1,11) circle (0.3pt);
	\filldraw  (1,12) circle (0.3pt);
	\filldraw  (1,13) circle (0.3pt);
	\filldraw  (1,14) circle (0.3pt);
	\filldraw  (1,15) circle (0.3pt);
	\filldraw  (1,16) circle (0.3pt);
	\filldraw  (1,17) circle (0.3pt);
	\filldraw  (1,18) circle (0.3pt);
	\filldraw  (1,19) circle (0.3pt);
	\filldraw  (1,20) circle (0.3pt);
	\filldraw  (1,21) circle (0.3pt);
	\filldraw  (1,22) circle (0.3pt);
	\filldraw  (1,23) circle (0.3pt);
	\filldraw  (1,24) circle (0.3pt);
	\filldraw  (1,25) circle (0.3pt);
	\filldraw  (1,26) circle (0.3pt);
	\filldraw  (1,27) circle (0.3pt);
	\filldraw  (1,28) circle (0.3pt);
	\filldraw  (1,29) circle (0.3pt);
	\filldraw  (1,30) circle (0.3pt);
	\filldraw  (2,1) circle (0.3pt);
	\filldraw  (2,2) circle (0.3pt);
	\filldraw  (2,3) circle (0.3pt);
	\filldraw  (2,4) circle (0.3pt);
	\filldraw  (2,5) circle (0.3pt);
	\filldraw  (2,6) circle (0.3pt);
	\filldraw  (2,7) circle (0.3pt);
	\filldraw  (2,8) circle (0.3pt);
	\filldraw  (2,9) circle (0.3pt);
	\filldraw  (2,10) circle (0.3pt);
	\filldraw  (2,11) circle (0.3pt);
	\filldraw  (2,12) circle (0.3pt);
	\filldraw  (2,13) circle (0.3pt);
	\filldraw  (2,14) circle (0.3pt);
	\filldraw  (2,15) circle (0.3pt);
	\filldraw  (2,16) circle (0.3pt);
	\filldraw  (2,17) circle (0.3pt);
	\filldraw  (2,18) circle (0.3pt);
	\filldraw  (2,19) circle (0.3pt);
	\filldraw  (2,20) circle (0.3pt);
	\filldraw  (2,21) circle (0.3pt);
	\filldraw  (2,22) circle (0.3pt);
	\filldraw  (2,23) circle (0.3pt);
	\filldraw  (2,24) circle (0.3pt);
	\filldraw  (2,25) circle (0.3pt);
	\filldraw  (2,26) circle (0.3pt);
	\filldraw  (2,27) circle (0.3pt);
	\filldraw  (2,28) circle (0.3pt);
	\filldraw  (2,29) circle (0.3pt);
	\filldraw  (2,30) circle (0.3pt);
	\filldraw  (3,1) circle (0.3pt);
	\filldraw  (3,2) circle (0.3pt);
	\filldraw  (3,3) circle (0.3pt);
	\filldraw  (3,4) circle (0.3pt);
	\filldraw  (3,5) circle (0.3pt);
	\filldraw  (3,6) circle (0.3pt);
	\filldraw  (3,7) circle (0.3pt);
	\filldraw  (3,8) circle (0.3pt);
	\filldraw  (3,9) circle (0.3pt);
	\filldraw  (3,10) circle (0.3pt);
	\filldraw  (3,11) circle (0.3pt);
	\filldraw  (3,12) circle (0.3pt);
	\filldraw  (3,13) circle (0.3pt);
	\filldraw  (3,14) circle (0.3pt);
	\filldraw  (3,15) circle (0.3pt);
	\filldraw  (3,16) circle (0.3pt);
	\filldraw  (3,17) circle (0.3pt);
	\filldraw  (3,18) circle (0.3pt);
	\filldraw  (3,19) circle (0.3pt);
	\filldraw  (3,20) circle (0.3pt);
	\filldraw  (3,21) circle (0.3pt);
	\filldraw  (3,22) circle (0.3pt);
	\filldraw  (3,23) circle (0.3pt);
	\filldraw  (3,24) circle (0.3pt);
	\filldraw  (3,25) circle (0.3pt);
	\filldraw  (3,26) circle (0.3pt);
	\filldraw  (3,27) circle (0.3pt);
	\filldraw  (3,28) circle (0.3pt);
	\filldraw  (3,29) circle (0.3pt);
	\filldraw  (3,30) circle (0.3pt);
	\filldraw  (4,2) circle (0.3pt);
	\filldraw  (4,3) circle (0.3pt);
	\filldraw  (4,4) circle (0.3pt);
	\filldraw  (4,5) circle (0.3pt);
	\filldraw  (4,6) circle (0.3pt);
	\filldraw  (4,7) circle (0.3pt);
	\filldraw  (4,8) circle (0.3pt);
	\filldraw  (4,9) circle (0.3pt);
	\filldraw  (4,10) circle (0.3pt);
	\filldraw  (4,11) circle (0.3pt);
	\filldraw  (4,12) circle (0.3pt);
	\filldraw  (4,13) circle (0.3pt);
	\filldraw  (4,14) circle (0.3pt);
	\filldraw  (4,15) circle (0.3pt);
	\filldraw  (4,16) circle (0.3pt);
	\filldraw  (4,17) circle (0.3pt);
	\filldraw  (4,18) circle (0.3pt);
	\filldraw  (4,19) circle (0.3pt);
	\filldraw  (4,20) circle (0.3pt);
	\filldraw  (4,21) circle (0.3pt);
	\filldraw  (4,22) circle (0.3pt);
	\filldraw  (4,23) circle (0.3pt);
	\filldraw  (4,24) circle (0.3pt);
	\filldraw  (4,25) circle (0.3pt);
	\filldraw  (4,26) circle (0.3pt);
	\filldraw  (4,27) circle (0.3pt);
	\filldraw  (4,28) circle (0.3pt);
	\filldraw  (4,29) circle (0.3pt);
	\filldraw  (4,30) circle (0.3pt);
	\filldraw  (5,4) circle (0.3pt);
	\filldraw  (5,5) circle (0.3pt);
	\filldraw  (5,6) circle (0.3pt);
	\filldraw  (5,7) circle (0.3pt);
	\filldraw  (5,8) circle (0.3pt);
	\filldraw  (5,9) circle (0.3pt);
	\filldraw  (5,10) circle (0.3pt);
	\filldraw  (5,11) circle (0.3pt);
	\filldraw  (5,12) circle (0.3pt);
	\filldraw  (5,13) circle (0.3pt);
	\filldraw  (5,14) circle (0.3pt);
	\filldraw  (5,15) circle (0.3pt);
	\filldraw  (5,16) circle (0.3pt);
	\filldraw  (5,17) circle (0.3pt);
	\filldraw  (5,18) circle (0.3pt);
	\filldraw  (5,19) circle (0.3pt);
	\filldraw  (5,20) circle (0.3pt);
	\filldraw  (5,21) circle (0.3pt);
	\filldraw  (5,22) circle (0.3pt);
	\filldraw  (5,23) circle (0.3pt);
	\filldraw  (5,24) circle (0.3pt);
	\filldraw  (5,25) circle (0.3pt);
	\filldraw  (5,26) circle (0.3pt);
	\filldraw  (5,27) circle (0.3pt);
	\filldraw  (5,28) circle (0.3pt);
	\filldraw  (5,29) circle (0.3pt);
	\filldraw  (5,30) circle (0.3pt);
	\filldraw  (6,6) circle (0.3pt);
	\filldraw  (6,7) circle (0.3pt);
	\filldraw  (6,8) circle (0.3pt);
	\filldraw  (6,9) circle (0.3pt);
	\filldraw  (6,10) circle (0.3pt);
	\filldraw  (6,11) circle (0.3pt);
	\filldraw  (6,12) circle (0.3pt);
	\filldraw  (6,13) circle (0.3pt);
	\filldraw  (6,14) circle (0.3pt);
	\filldraw  (6,15) circle (0.3pt);
	\filldraw  (6,16) circle (0.3pt);
	\filldraw  (6,17) circle (0.3pt);
	\filldraw  (6,18) circle (0.3pt);
	\filldraw  (6,19) circle (0.3pt);
	\filldraw  (6,20) circle (0.3pt);
	\filldraw  (6,21) circle (0.3pt);
	\filldraw  (6,22) circle (0.3pt);
	\filldraw  (6,23) circle (0.3pt);
	\filldraw  (6,24) circle (0.3pt);
	\filldraw  (6,25) circle (0.3pt);
	\filldraw  (6,26) circle (0.3pt);
	\filldraw  (6,27) circle (0.3pt);
	\filldraw  (6,28) circle (0.3pt);
	\filldraw  (6,29) circle (0.3pt);
	\filldraw  (6,30) circle (0.3pt);
	\filldraw  (7,8) circle (0.3pt);
	\filldraw  (7,9) circle (0.3pt);
	\filldraw  (7,10) circle (0.3pt);
	\filldraw  (7,11) circle (0.3pt);
	\filldraw  (7,12) circle (0.3pt);
	\filldraw  (7,13) circle (0.3pt);
	\filldraw  (7,14) circle (0.3pt);
	\filldraw  (7,15) circle (0.3pt);
	\filldraw  (7,16) circle (0.3pt);
	\filldraw  (7,17) circle (0.3pt);
	\filldraw  (7,18) circle (0.3pt);
	\filldraw  (7,19) circle (0.3pt);
	\filldraw  (7,20) circle (0.3pt);
	\filldraw  (7,21) circle (0.3pt);
	\filldraw  (7,22) circle (0.3pt);
	\filldraw  (7,23) circle (0.3pt);
	\filldraw  (7,24) circle (0.3pt);
	\filldraw  (7,25) circle (0.3pt);
	\filldraw  (7,26) circle (0.3pt);
	\filldraw  (7,27) circle (0.3pt);
	\filldraw  (7,28) circle (0.3pt);
	\filldraw  (7,29) circle (0.3pt);
	\filldraw  (7,30) circle (0.3pt);
	
	\end{axis}
    \end{scope}
\end{tikzpicture}\caption{A comparison for the geography of Surfaces of general type for $\cha(k)=0$ (on the left) and $\cha(k)>0$ (on the right). Every dot symbolizes a couple $(a,b)$ in the plane for which it could be possible to find a surface of general type $X$ with invariants $K_X^2=a$ and $\chi(\O_X)=b$. Here by "irregular" we mean with non-trivial Albanese morphism and by "BMY" we mean the Bogomolo-Miyaoka-Yau type inequality $K_X^2\leq32\chi(\O_X)$ which holds over algebraically closed fields of positive characteristic. Notice that, over fields of positive characteristic, it is possible to find surfaces of general type with $K_X\geq12\chi(\O_X)$ only if they are not of maximal Albanese dimension.}
\label{geogen}%
\end{figure}

Observe that if the Albanese morphism $\alb_X$ of $X$ is non-trivial, there exists an \'etale morphism of degree $d^{2q}$ for every integer $d$ coprime with the characteristic of the ground field which is the pull-back of the multiplication by $d$ on the Albanese variety $\Alb(X)$:
\begin{equation}
\label{eq_cartalbbis}
\begin{tikzcd}
X_d\arrow{r}{\nu_d} \arrow[dr, phantom, "\square"] \arrow{d}{a_d} & X\arrow{d}{\alb_X}\\
\Alb(X)\arrow{r}{\mu_d} & \Alb(X).
\end{tikzcd}
\end{equation}
 In this case we have that $K_{X_d}^2=d^{2q}K_X^2$ and $\chi(\O_{X_d})=d^{2q}\chi(\O_X)$ and,  if $K_X=2\chi(\O_X)-\epsilon$ for a positive number $\epsilon$,  we would have that  $K_{X_d}^2<2\chi(\O_{X_d})-6$ for $d\gg 0$ which contradicts the Noether's inequality, i.e. we have
\begin{equation}
K_X^2\geq2\chi(\O_X).
\label{eq_noetin4}
\end{equation}
In particular we have shown that surfaces which lie near the Noether line have trivial Albanese variety: in positive characteristic it is not true that they are regular (i.e. $q'(X)=0$) and there are known example of surfaces satisfying $K_X^2=q'(X)=p_g(X)=1$ in characteristic at most $5$ which are called non classical Godeaux\index{surface!non-classical Godeaux} surfaces (cf. \cite{liedtkegodeaux}).

The last inequality we are going to recall is valid only in characteristic zero and is the so called Bogomolov-Miyaoka-Yau inequality\index{Bogomolov-Miyaoka-Yau inequality}, which states that for a minimal surface of general type $K_X^2\leq 9\chi(\O_X)$ holds.
In positive characteristic this inequality can not be true: as we have already remarked, there exist surfaces $X$ with $e(X)<0$ and if this happens, by Noether's formula \ref{eq_noetform}, we have $K_X^2>12\chi(\O_X)$. 
Nevertheless, there exists a weaker Bogomolov-Miyaoka-Yau type inequality (cf. \cite{gu2019slope} Corollary 3.4) which states that for a minimal surface of general type $K_X^2\leq 32\chi(\O_X)$ holds.
Observe that this inequality is sharp over algebraically closed fields of characteristics $2$ and $5$ (ibidem Paragraph 4.3.2 and Proposition 4.6) while it is not sharp over algebraically closed fields of characteristic $2$ and of characteristic $p\gg 0$ (ibidem Paragraph 4.3.4 and \cite{gu_2016} Theorem 1.3(3)).
Again, if we assume $X$ to be of maximal Albanese dimension, by Theorem \ref{teo_shepbarr} and Noether's formula we get 
\begin{equation}
K_X^2<12\chi(\O_X).
\label{eq_bmyp}
\end{equation}

\chapter{Double covers of surfaces}
\label{chap_double}

In this chapter we are going to study morphisms of degree two between surfaces. 
Usually, for a morphism of degree two $f\colon X\to Y$, we assume $Y$ to be smooth and $X$ normal, so that f is flat: this case is studied in Section \ref{sec_double}.
In characteristic two, we see that it may be useful also to consider morphisms $g\colon Y\to X$ of degree two again with $X$ normal and $Y$ smooth, but this time $X$ is the target and $Y$ the source of  the morphism.
These morphisms are studied through the theory of foliations and rational vector fields in Section \ref{sec_foliations}.
In Section \ref{sec_insdouble} we see how these two different types of morphisms of degree two are related to one another, while in Section \ref{sec_canres} we study a process, the canonical resolution, that allows us to reduce to the study of morphisms of degree two where both $X$ and $Y$ are smooth.
In section \ref{sec_doublefundgroup} we prove that a double cover branched over an ample divisor induces an isomorphism of algebraic fundamental groups. 
In Section \ref{sec_insdoubleprod} we give a general construction for surfaces of general type as inseparable double cover of product surfaces such that the canonical resolution is a non-splittable double cover.

As we will see in Chapters \ref{chap_severi0} and \ref{chap_severi+}, we will need double covers in the study of minimal surfaces of general type with maximal Albanese dimension which lie on or near the Severi lines: indeed they turn out to be double covers of smooth elliptic surface fibrations. 

\section{Double covers}
\label{sec_double}
In this section we study morphisms of degree two from a normal surface to a smooth one. While, in the case the characteristic of the ground field is different from two, this is a well known topic, when $2$ is not invertible many pathologies may arise. The material of this part is taken from \cite{barth} section I.17 when the characteristic is different from two (actually there, only the case of complex varieties is treated, but everything works in the same way when the characteristic is not two). In case the characteristic is two, the theory of double covers is more fragmentary, we try to collect all the results concerning them and explain how they are related one another: the main references for what we need here are \cite{cossec} section 0.1 and \cite{gusunzhou} section 6.

\begin{definition}
\label{def_double}
A double cover\index{double cover} $f\colon X\to Y$ is a finite morphism of degree two, i.e. the extension of fields $k(Y)\subseteq k(X)$\nomenclature{$k(X)$}{the field of rational functions of an irreducible scheme $X$} has degree two.
The morphism is said to be separable (respectively inseparable) if the corresponding morphism on rational functions is. 
\end{definition}

We will always consider double covers of surfaces, one of which is smooth and the other at least normal.
In this section we will always consider $X$ to be normal and $Y$ to be smooth.
In this case we have that there exists an affine covering $\{V_i\}$ of $Y$, for which $U_i=f^{-1}(V_i)$ is still affine and the induced morphism on regular functions makes $B_i=\O_X(U_i)$ an $A_i=\O_Y(V_i)$-module generated by two elements.

By the miracle flatness Theorem (\cite{matsumura} Theorem  23.1), such a double cover $f\colon X\to Y$ is flat: in particular $B=f_*\O_X$ is a locally free $\O_Y$-algebra of rank two and $X\simeq\specu(B)$, where $\specu(B)$\nomenclature{$\specu(F)$}{the relative spectrum of the sheaf of algebras $F$} is the relative spectrum of the sheaf of algebras $B$ (cf. \cite{hart} Exercise II.5.17).

Up to a refinement of the affine covering $\{V_i\}$, we have that
\begin{equation}
B_i=A_i[t_i]/(t_i^2+a_it_i+b_i),
\label{eq_doublecov}
\end{equation}
where $a_i, b_i\in A_i$.
Clearly two generators of $B_i$ as an $A_i$-module are $1$ and $t_i$. 
Changing coordinates, the glueing conditions give us that 
\begin{equation}
t_i=g_{ij}t_j+c_{ij},
\label{eq_ti}
\end{equation}
with $g_{ij}\in A_{ij}^*=\O_{Y}(V_{ij})^*$ and $c_{ij}\in A_{ij}=\O_Y(V_{ij})$.
This shows that $f_*\O_X$ fits into the following short exact sequence
\begin{equation}
0\to\O_Y\to f_*\O_X\to L^{-1}\to 0
\label{eq_seqdouble}
\end{equation}
where the transition functions of the sections of $L$ are $g_{ij}$ and the transition functions of the sections $f_*\O_X$ are the matrices
\begin{equation}
\begin{pmatrix}
	1 & -g_{ij}^{-1}c_{ij}\\
	0 & g_{ij}^{-1}
\end{pmatrix}
\label{eq_transitions}
\end{equation}
(indeed if the basis of a free module transforms via a matrix $A$, then the sections expressed in that basis transform via the transpose of the inverse matrix $^tA^{-1}$).

\begin{definition}
The line bundle $L$ appearing in the exact sequence \ref{eq_seqdouble} is called the line bundle associated with the double cover\index{double cover!line bundle associated with} $f\colon X\to Y$.
\end{definition}

\begin{remark}
\label{rem_doubleeuler}
From equation \ref{eq_seqdouble} one easily derives
\begin{equation}
\chi(\O_X)=\chi(f_*\O_X)=\chi(\O_Y)+\chi(L^{-1})=2\chi(\O_Y)+\frac{1}{2}(L^2+K_Y.L).
\label{eq_doubleeuler}
\end{equation}
\end{remark}

\begin{remark}
\label{rem_split}
Notice that $c_{ij}$ defines a cocycle, which can be considered as an element of $H^1(Y,L)=\ext(L^{-1},\O_Y)$ which defines the short exact sequence \ref{eq_seqdouble}.
In particular the sequence splits if and only if $c_{ij}=0\in H^1(Y,L)$.   
\end{remark}

\begin{remark}
\label{rem_transitions}
On $A_i$ we have $t_i^2+a_it_i+b_i=0$ and a similar relation on $A_j$: comparing the two equations on $A_{ij}$, thanks to Equation \ref{eq_ti} we get
\begin{equation}
\begin{split}
t_i^2+a_it_i+b_i=(g_{ij}t_j+c_{ij})^2+a_i(g_{ij}t_j+c_{ij})+b_i=\\
=g_{ij}^2\Bigl(t_j^2+(2c_{ij}+a_i)g_{ij}^{-1}t_j+(b_i+a_ic_{ij}+c_{ij}^2)g_{ij}^{-2}\Bigr)
\end{split}
\label{eq_changecoordinatestab}
\end{equation}
which leads to
\begin{equation}
\label{eq_ai}
a_i=g_{ij}a_j-2c_{ij}
\end{equation}
and
\begin{equation}
\label{eq_bi}
b_i=g_{ij}^2b_j-a_jg_{ij}c_{ij}+c_{ij}^2.
\end{equation}
\end{remark}

\begin{remark}
\label{rem_canonical}
Combining Equations \ref{eq_changecoordinatestab}, \ref{eq_ai} and \ref{eq_bi} we get
\begin{equation}
t_i^2+a_it_i+b_i=g_{ij}^2(t_j^2+t_ja_j+b_j)
\label{eq_eqdouble}
\end{equation} which shows that $X$ can be seen as the zero divisor of a section of the pull-back of $L^2$ in the affine bundle $W$ of rank one (which need not to be an actual line bundle, i.e. may have no sections) defined by the glueing $t_i=g_{ij}t_j+c_{ij}$.
In particular, because it is a hypersurface in a smooth variety, $X$ is Gorenstein, i.e. the dualizing sheaf is a line bundle.

This can be used (cf. \cite{cossec} Proposition 0.1.3) to determine the canonical bundle of $X$. We have
\begin{equation}
\label{eq_candouble}
K_X=f^*(K_Y+L),
\end{equation}
which implies
\begin{equation}
K_X^2=2K_Y^2+4K_Y.L+2L^2.
\label{eq_doublevolume}
\end{equation}
\end{remark}

 It is now time to distinguish whether the characteristic of the ground field is two or not. First we treat the case $\cha(k)\neq 2$, which is easier.

\paragraph{$\mathbf{\bcha(k)\neq 2}$}
Replacing $t_i$ by $t_i+\frac{1}{2}a_i$ we may assume that $a_i=0$, in particular (cf. Equation \ref{eq_ai}) $c_{ij}=0$ and the sequence \ref{eq_seqdouble} splits.
This means that 
\begin{equation}
f_*\O_X=\O_Y\oplus L^{-1}.
\label{eq_splitdouble}
\end{equation}

Equation \ref{eq_bi} becomes
\begin{equation}
b_i=g_{ij}^2b_j,
\label{eq_bisplit}
\end{equation}
in particular the $b_i$'s glue to a section $b$ of $L^2$ which defines a divisor $R$.

\begin{definition}
\label{def_branch}
The divisor $R$ is called the branch divisor\index{double cover!branch divisor of} of the double cover $f\colon X\to Y$.
\end{definition}

\begin{remark}
The condition $X$ normal is equivalent to the condition $R$ reduced (cf. \cite{barth} I.17).

Moreover, the singularities of $X$ lie above the singularities of the branch divisor $R$ (ibidem).
\label{rem_normal}
\end{remark}

Summing up, a double cover in characteristic different from two, is uniquely determined by the line bundle $L$ and the branch divisor $R$ which is reduced and satisfies 
\begin{equation}
L^{\otimes 2}=\O_Y(R).
\label{eq_doubleramification}
\end{equation}

\paragraph{$\mathbf{\bcha(k)=2}$}

In this case Equations \ref{eq_ai} and equations \ref{eq_bi} become
\begin{equation}
a_i=g_{ij}a_j
\label{eq_ai2}
\end{equation}
and
\begin{equation}
b_i=g_{ij}^2b_j+a_jg_{ij}c_{ij}+c_{ij}^2.
\label{eq_bi2}
\end{equation}

It is now clear, that in this case the $a_i$'s glue to a section $a$ of $L$ and $(1,a_i,b_i)$ glue to a section $\sigma$ of a vector bundle $V$ of rank $3$ with transition functions the matrices
\begin{equation}
\begin{pmatrix}
	1 & 0 & 0\\
	0 & g_{ij} & 0\\
	c_{ij}^2 & g_{ij}c_{ij} & g_{ij}^2
\end{pmatrix}
\label{eq_transv}
\end{equation}
which fits in the following short exact sequence
\begin{equation}
0\to L^2\to V\to \O_Y\oplus L\to 0.
\label{eq_exactv}
\end{equation}

Hence a double cover in characteristic two is uniquely determined by an exact sequence as \ref{eq_exactv} and a global section $\sigma=(1,a,b)$ of $V$ (or, equivalently, by an exact sequence as \ref{eq_seqdouble} and local sections $a_i$ and $b_i$ satisfying equations \ref{eq_ai2} and $\ref{eq_bi2}$). 

\begin{definition}
\label{def_split2}
A double cover in characteristic $2$ is said to be splittable\index{double cover!splittable} if the exact sequence \ref{eq_seqdouble} splits.
\end{definition}

\begin{remark}
Notice that the exact sequence \ref{eq_seqdouble} splits if and only if  sequence \ref{eq_exactv} does: indeed, looking at the transition functions \ref{eq_transitions}  and \ref{eq_transv} we see that both sequences split if and only if it is possible to take $c_{ij}=0$.
\label{oss_splitseqdoubcov}
\end{remark}

\begin{proposition}[\cite{cossec} Proposition 0.1.2]
\label{propo_inssep}
A double cover in characteristic $2$ is separable if and only if the section $a$ of $L$ defined by equation \ref{eq_ai2} is not identically zero.
\end{proposition}

We also have a notion of branch divisor for double covers in characteristic $2$ except when the double cover is inseparable and non-splittable. Notice, indeed, that if the double cover is inseparable and splittable, then the $b_i$'s glue to a global section $b$ of $L^2$ (cf. Equation \ref{eq_bi2}).

\begin{definition}
\label{def_branch2}
If $f\colon X\to Y$ is a separable\index{double cover!separable} (splittable or not) double cover, we define the branch divisor\index{double cover!branch divisor of} $R$ to be the zero divisor of the section $a$ of $L$.
In this case we have $\O_Y(R)=L$.

If $f\colon X\to Y$ is an inseparable\index{double cover!inseparable} splittable double cover, we define the branch divisor $R$ to be the zero divisor of $b$. In this case we have
$L^{ 2}=\O_Y(R)$ as when the characteristic is different from two.
\end{definition}

\begin{remark}
Notice that for a separable double cover in characteristic two, the branch divisor is a divisor in the linear series of $L$ instead of a divisor in the linear series of $L^{\otimes 2}$ as it is over fields of characteristic different from two: this is related to the fact that every double cover in characteristic two is wildly ramified (cf. \cite{hart} Proposition 2.2(c)).

In particular a splittable double cover in characteristic $2$ is uniquely determined by a line bundle $L$ and two effective divisors $R_1$ and $R_2$ such that
\begin{equation}
L=\O_Y(R_1)
\label{eq_ram12}
\end{equation} 
and
\begin{equation}
L^{\otimes2}=\O_Y(R_2).
\label{eq_ram22}
\end{equation}
Moreover it is inseparable if and only if $R_1=0$.
Clearly $R_1$ coincides with the zero locus of the global section $a$ of $L$ while $R_2$ with the zero locus of the global section $b$ of $L^{\otimes2}$.
\label{rem_splitdoublecover2}
\end{remark}

\begin{proposition}
\label{propo_sing2}
Let $f\colon X\to Y$ be a double cover of surfaces over a field of characteristic $2$ and let $\pi\colon W\to Y$ the affine rank-$1$ bundle on which $X$ is embedded as the zero of $t^2+at+b$ as a section of $\pi^*L^2$. Then there exists a section $\alpha(f)$ of $\restr{\bigr(\pi^*L^2\otimes\Omega^1_W\bigl)}{a=0}$ such that the singular points of $X$ are the points where $\alpha(f)$ vanishes.
\end{proposition}
\begin{proof}
We know that locally on $U_i\times \A^1$, the equation for $X$ is $t_i^2+t_ia_i+b_i$: hence we know that the singularities of $X$ lie where its differential vanishes, that is where $a_idt_i+t_ida_i+db_i$ vanishes. 
Using the natural splitting of the K\"ahler differentials on $U_i\times \A^1$, this is equivalent to asking $a_i=0$ and $t_ida_i+db_i=0$.

Differentiating equations \ref{eq_ti}, \ref{eq_ai2} and \ref{eq_bi2} we get
\begin{equation}
dt_i=g_{ij}dt_j+t_jdg_{ij}+dc_{ij},
\label{eq_dti}
\end{equation}
\begin{equation}
da_i=g_{ij}da_j+a_jdg_{ij}
\label{eq_dai}
\end{equation}
and
\begin{equation}
db_i=g_{ij}^2db_j+a_jg_{ij}dc_{ij}+a_jc_{ij}dg_{ij}+c_{ij}g_{ij}da_j.
\label{eq_dbi}
\end{equation}

Combining these equations where $a_i=0$ we get $t_ida_i+db_i=g_{ij}^2(t_jda_j+b_j)$, that is these local expressions glue to a section $\alpha(f)=tda+db$ of $\restr{\bigl(\pi^*L^2\otimes\Omega^1_W\bigr)}{a=0}$ whose zeroes are clearly the singular points of $X$.
\end{proof}

\begin{remark}
\label{rem_globaldiff}
Let $f\colon X\to Y$ be an inseparable double cover\index{double cover!inseparable} (splittable or not): we have already said (cf. Proposition \ref{propo_inssep}) that  $a$ is identically zero. 
In this case the twisted differential 1-form of Proposition \ref{propo_sing2} is the pull-back of a section on $Y$ and is locally given by $db_i$ which glue to a global section $db$ of $L^2\otimes\Omega^1_Y$ thanks to equation \ref{eq_dbi}.

We also notice that in this case the condition $X$ normal is equivalent to the fact that $db$ has only isolated zeroes. Indeed, because $X$ is a hypersurface, hence Cohen-Macaulay, clearly satisfies Serre's condition $S_2$ (\cite{egaiv2} Definition 5.7.2): so, thanks to Serre's criterion for normality (ibidem Th\'eor\`eme 5.8.6), $X$ is normal exactly when the zeroes of $db$ are isolated.
In particular the quotient $\Omega^1_Y/L^{-2}$ induced by $db$ is free in codimension $1$ and hence torsion-free.

In other words, we get the following exact sequence
\begin{equation}
0\to L^{-2}\to \Omega_Y^1\to I_Z\F^{-1}\to 0:
\label{eq_lfoliation}
\end{equation}
Indeed, for the well-known result that over a principal ideal domain every torsion-free module is free, any torsion-free sheaf $\F'$ is isomorphic in codimension $1$ to its double dual $\Hom(\Hom(\F',\O_Y),\O_Y)$ which is reflexive and hence free, because we are on a smooth surface (cf. Corollary 1.4 \cite{hartref}). In particular $\F'=I_Z\F^{-1}$ where $Z$ is the subscheme of codimension $2$  where $\F'$ and $\F^{-1}$ do not coincide (the $-1$ exponent is just for coherence of notation and will be more clear in section \ref{sec_insdouble}).
Actually the fact that $\F'$ is torsion-free is equivalent to the fact that $L^{-2}$ is a subbundle of $\Omega^1_Y$ in codimension $1$.
\end{remark}

\begin{remark}
\label{rem_normsep}
Suppose $f\colon X\to Y$ is a separable double cover: as we have already said, locally $X$ is given by the equation $t^2+at+b$ where $a$ and $b$ are local functions on $Y$.
Because we are assuming that $X$ is normal, we know it has only isolated singularities.
In particular this implies that $b$ has no factors of multiplicity greater than or equal to $2$ along factors of $a$.
Indeed, if this were not the case, we could rewrite the local equation of $X$ as $t^2+a'ct+c^2b'$: then, by the Jacobian criterion we know that the singular points are those were the differential vanishes. i.e. $a'c=0$ and $ta'dc+tcda'+c^2db'= 0$ and this is true where $c$ vanishes, which would mean that $X$ has a divisorial singularity.
\end{remark}

\begin{lemma}[cf. \cite{cossec} Remark 0.2.2]
\label{lemma_branchan}
Let $f\colon X\to Y$ be a separable double cover\index{double cover!separable} with branch divisor $R$.
Then the singular points of $X$ lying over smooth points of $R$ (if they exist) are rational singularities of type $A_n$.
In particular if $R$ is smooth, $X$ has at most rational double points\index{rational double point} as singularities.
\end{lemma}

\begin{proof}
We already know that the singularities of $X$ lie over the branch divisor $R$ (Proposition \ref{propo_sing2}).
Recall that a normal singularity\index{singularity!rational} of $X$ is called rational if, after taking a resolution $X'\to X$ of the singularity, one has $\chi(\O_{X'})=\chi(\O_X)$.
In particular, as we will see, this implies that for a surface $X$ which is the double cover of a smooth surface $Y$, in order to compute the invariants of a minimal resolution of $X$, one could operate as if $X$ is non-singular  and calculate the invariants directly on $X$.
To see that in this case the singularity is of type $A_n$ is a local computation so, we may restrict to the completion of the local ring, i.e. $k[[x,y,t]]$ where  $x$ and $y$ are local parameters for $p$ and $t$ is a local parameter for the affine bundle where $X$ is defined.
Because $a$ is smooth we can take it as a local coordinate (e.g. $a=x$), which means that we may assume that $X$ is given by the equation
\begin{equation}
t^2+xt+b(x,y),
\label{eq_locsing}
\end{equation}
where $b$ is a polynomial vanishing at the origin (we can always assume $b((0,0))=0$ up to a substitution $t\mapsto t+\sqrt{b((0,0))}$).
The partial derivatives of the equation are
\begin{equation}
\begin{split}
t+b_x&(x,y);\\
b_y(x&,y);\\
x&.
\end{split}
\label{eq_partderiv}
\end{equation}
In particular, to have that the origin is a singular point, we need that every monomial term of $b$ has at least degree $2$.
Moreover we have seen that, in order to have an isolated singularity, $x^2$ does not divide  $b$, which means that we can write $b=x^2f(x,y)+xg(y)+h(y)$ where $g$ and $h$ vanish at zero with multiplicity at least $1$ and $2$ respectively and at least one of them is non-trivial.
We may assume that $h(y)\neq 0$ and $g(y)=0$: indeed  the substitution $t\mapsto t+g(y)$ turns our equation into
\begin{equation}
\label{eq_intermediatecoordinates}
t^2+xt+x^2f(x,y)+g^2(y)+h(y).
\end{equation}
We have just seen that we may reduce our equation to
\begin{equation}
t^2+xt+x^2f(x,y)+h(y).
\label{eq_variecoordinate}
\end{equation}

\begin{lemma}
\label{lemma_quadre}
The equation 
\begin{equation}
\label{eq_lemmaquadre}
\lambda^2+x\lambda=x^2f(x,y),
\end{equation}
 with $f\in k[[x,y]]$, has always a solution $\psi(x,y)\in k[[x,y]]$.
\end{lemma}

\begin{proof}
Consider $\lambda$ and $x^2f(x,y)$ as formal power series in $x$ with coefficients $\lambda_i$, $\phi_i$ respectively, in $k[[y]]$.
With this notation, Equation \ref{eq_lemmaquadre} becomes
\begin{equation}
\begin{split}
\lambda^2_{j}+\lambda_{2j-1}=\phi_{2j};\\
\lambda_{2j}=\phi_{2j+1}.
\end{split}
\label{eq_finalesing}
\end{equation}
If $j=1$ we get $\lambda_1^2+\lambda_1=\phi_2$: this equation, with unknown $\lambda_1$, has a solution in $k[[y]]$ because $k[[y]]$ is complete, hence Henselian, and the reduction modulo $(y)$ of the equation has simple zeroes.
The other equations for $j>1$ are solved by induction and we get a solution $\psi(x,y)\in k[[x,y]]$.
\end{proof}

Now take the substitution $t\mapsto t+\psi(x,y)$, where $\psi(x,y)$ is taken as in Lemma \ref{lemma_quadre}: then we obtain 
\begin{equation}
t^2+xt+h(y).
\label{eq_quasifinale}
\end{equation}

Notice that $h(y)=y^n\alpha$, where $\alpha$ is an invertible element in $k[[y]]$.
As a final step, we take $t\mapsto \alpha t$, $x\mapsto x+t+\alpha t$ and the equation finally becomes
\begin{equation}
t^2+xt+y^n,
\label{eq_finalequadre}
\end{equation}
which is the standard form for an $A_n$ singularity.
\end {proof}

\begin{remark}
\label{rem_normsep2}
We have proved a relation between the singularities of the branch divisor and the singularities of the double covers, in particular we have seen that if the branch divisor is smooth, the singularities of $X$ are as good as possible.
The opposite is not true: it is indeed quite easy to find example of separable double covers with branch divisor as singular as possible, which still have at most rational singularities.

Suppose that the double cover is splittable and denote, as usual, by $a$ and $b$ the sections of $L$ and $L^2$ which defines the double cover.
In addition, assume that the divisor $D(b)$ associated with $b$  passes through every singular point of the divisor $D(a)$ associated with $a$ and that $D(b)$ is smooth in those points.
Then we can conclude that the associated double cover $X$ has at most $A_n$ singularities.
Indeed, outside the singular points of $D(a)$, this is what we have just shown in Lemma \ref{lemma_branchan}.
Restricting to a neighbourhood of a singular point $p$ of $D(a)$
we have the equation
\begin{equation}
t_i^2+a_it_i+b_i
\label{eq_boh}
\end{equation}
where, by hypothesis, $a_i\in m_p^2$ and $b_i\in m_p\setminus m_p^2$. 
Hence the Jacobian criterion applies and shows that $X$ is smooth above $p$.
\end{remark}

\begin{remark}
\label{rem_effdouble}
Notice that we have proved that if the characteristic is different from two or the characteristic is two and $f$ is separable or $f$ is inseparable but splittable then $L$ is effective up to a multiple. 
Indeed in this case the $a_i$'s or he $b_i$'s glue together to form a section of $L$, respectively $L^2$.
On the opposite if $f\colon X \to Y$ is non splittable inseparable, in general $L$ will not be effective even up to a multiple (as we will see in Section \ref{sec_insdoubleprod}).
\end{remark}

\section{Foliations and rational vector fields}
\label{sec_foliations}
Miyaoka (cf. \cite{miyapet}) and Ekedahl (cf. \cite{ekeins} and \cite{ekecan}) introduce the tool of foliations in the study of purely inseparable finite covers: their works are more general but, as we need just for double covers of surfaces, we will assume, as before, that every variety is a surface and the algebraically closed ground field $k$ is of characteristic $2$. 
In this section we will consider double covers from a smooth surface to a normal one: in the case the morphism is inseparable, we will see that it is equivalent to consider double covers from a normal surface to a smooth one or vice versa.
The fact that the quotient of the tangent sheaf by a 1-foliation is a 2-divisible line bundle up to an ideal sheaf is proved in Section \ref{sec_insdouble}: in this section we will take this for granted for consistence of notation.

For the theory of rational vector fields in the study of double covers we refer  to Rudakov and Shafarevich \cite{rudakov} and Cossec and Dolgachev \cite{cossec} chapter 0.

\begin{definition}
Let $Y$ be a surface; a 1-foliation\index{foliation} is a saturated subsheaf $\F$ of the tangent sheaf $\T_{Y}$\nomenclature{$\T_{X}$}{the tangent sheaf of a scheme $X$} which is involutive (i.e. $[\F,\F]\subseteq\F$ where $[\cdot,\cdot]$ denotes the Lie bracket on the tangent space) and 2-closed (i.e. $\F^2\subseteq \F$, which means that the square of any local derivative of $\F$, which is again a local derivative in characteristic two, is again in $\F$). 
It is called smooth if it is a subbundle of $\T_{Y}$. 
\label{def_fol}
\end{definition}

We will consider 1-foliations only on smooth surfaces.
 Observe that in this case, thanks to Remark \ref{rem_folfree}, we have that every 1-foliation is actually a line bundle even if it is not smooth and, arguing as in Remark \ref{rem_globaldiff}, $\F$ is a subbundle of the tangent sheaf $\T_Y$ (hence it is a smooth foliation) in codimension $1$ and fits in the following short exact sequence
\begin{equation}
0\to \F\to \T_Y\to I_Z L^2\to 0
\label{eq_foliationl}
\end{equation}
(as in the discussion for Equation \ref{eq_lfoliation} a torsion-free sheaf of rank one in a surface is isomorphic to $I_ZL^2$ where $Z$ is a subscheme of codimension $2$ and $L^2$ is a line bundle; so far, and throughout all of this section, $L^2$ has to be intended simply as a symbol which denotes a line bundle, the square here does not mean it to be a priori a square in the Picard group: the notation will become clear in section \ref{sec_insdouble}).

\begin{definition}
\label{def_zeroeset}
The subscheme $Z$ is called the zero set of the foliation\index{foliation!zero set of}.
\end{definition}

A simple calculation with the Chern classes of these sheaves lead to
\begin{equation}
K_Y+\F=-2L
\label{eq:}
\end{equation}
and
\begin{equation}
\begin{split}
e(Y)=c_2(\T_Y)=&\deg(Z)+2L.\F=\\
=\deg(Z)-\F.(K_Y+\F)&=\deg(Z)-2L.(2L+K_Y).
\label{eq_eulfol}
\end{split}
\end{equation}

\begin{definition}
\label{def_quotfol}
Given a 1-foliation $\F$ on a surface $Y$ we define the quotient\index{foliation!quotient via} $X=Y/\F$\nomenclature{$Y/\F$}{the quotient of a smooth surface $Y$ via the 1-foliation $\F$} to be the scheme with the same topological space as $Y$ and whose structure sheaf is given by the functions $\ann(\F)\subseteq\O_Y$ on which every local derivative of $\F$ vanishes.
Notice that, given a 1-foliation, it is defined a natural morphism $g\colon Y\to Y/\F$ which is a homeomorphism at topological level, i.e. it is inseparable.
\end{definition}

\begin{remark}
\label{rem_folfrob}
It is quite clear from Definition \ref{def_quotfol} that the morphism $g\colon Y\to Y/\F$ induced by a 1-foliation $\F$ is a factor of the $k$-linear Frobenius morphism.
Indeed, if we take the quotient of $Y$ by the entire tangent sheaf instead of a 1-foliation, we get the $k$-linear Frobenius morphism, i.e. $Y^{(-1)}=Y/\T_Y$.
In particular, for every 1-foliation $\F$, we get (cf. \cite{ekecan} section I.1)
\begin{equation}
\begin{tikzcd}
Y\arrow{r}{g}\arrow[bend left=30]{rr}{F_k} & X=Y/\F\arrow{r}{} & Y^{(-1)}.
\end{tikzcd}
\label{eq_eq_folfrob}
\end{equation}

In the definition of 1-foliation, the 1 means exactly that the quotient morphism $Y\to Y/\F$ is a factor of the $k$-linear Frobenius morphism.
It is possible to give a definition of an n-foliation for every positive integer (cf. \cite{ekeins} Definition 3.1): this is defined as a  suitable subsheaf of the algebra $\diff(Y)$ of differential operators.
In the case of 1-foliations it turns out to be equivalent to consider subsheaves of $\diff(Y)$ and their intersection with $\T_Y\subseteq\diff(Y)$ (ibidem Lemma 4.1).
For every n-foliation $\F$, it is possible to define a quotient $Y/\F$ such that the natural morphism $Y\to Y/\F$ is a factor of the n-th iterated $k$-linear Frobenius morphism (ibidem Proposition 3.2).   
\end{remark}

\begin{remark}
\label{rem_singfol}
It is known (cf. \cite{ekecan} section I.1), that it is always possible to obtain a 1-foliation $\F$ as the unique  extension of a smooth 1-foliation $\restr{\F}{U}$ defined over a dense open subset $U$ whose complement is  a codimension 2 subset.
 
Moreover (ibidem or \cite{miyapet} Proposition 1.9(3) Lecture III) the singularities of $X=Y/\F$ lie exactly on  the image of the subscheme $Z$ where $\F$ is not a subbundle of the tangent space $\T_Y$.

In addition, $X$ is always normal: indeed we have (cf. Remark \ref{rem_folfrob}) $\O_Y\supseteq g^*\O_X=\ann(\F)\supseteq \O_Y^2$. Let $\alpha\in k(X)\subseteq k(Y)$ a rational function on $X$ which is the zero of a polynomial $p$ with coefficients in $\O_X=\ann(\F)$: in particular it is the zero of a polynomial with coefficients in $\O_Y$ and hence $\alpha\in \O_Y$ ($Y$ is smooth and hence normal). Moreover $\F(\alpha)=0$ and hence $\alpha\in\O_X$, which implies that $X$ is normal.
\end{remark}

\begin{theorem}[\cite{miyapet} Proposition 1.9 Lecture III or \cite{ekeins} Lemma 4.1]
\label{teo_doublefol}
There is a bijective correspondence between the set of 1-foliations of rank $1$ on $Y$ and the set of inseparable double covers $g\colon Y\rightarrow X$ where $X$ is normal and its singularities lie in the image of the zero set $Z$ of the foliation.

The correspondence is given as follows: to a double cover $g\colon Y\to X$ we associate the 1-foliation given by the relative tangent space $\T_{Y/X}$\nomenclature{$\T_{X/Y}$}{the relative tangent sheaf for a morphism $f\colon X\to Y$ of schemes}, while to a 1-foliation $\F$ we associate the scheme $Y/\F$. 
Moreover the k-linear Frobenius morphism $F_k$ factors through $g$, that is the following diagram commutes:
\[
\begin{tikzcd}
Y \arrow{r}{g}\arrow{d}{F_k} & X\arrow{dl}{f}\\
Y^{(-1)}. &
\end{tikzcd}
\]
\end{theorem}

\begin{theorem}[\cite{ekeins} Corollary 3.4 or \cite{miyapet} Corollary 1.11 Lecture III]
\label{teo_canfol}
Let $g\colon Y\to X=Y/\F$ be an inseparable double cover corresponding to the 1-foliation $\F$ with $Y$ smooth.
Let $U\subseteq X$ be the smooth locus of $X$ and $V=f^{-1}(U)$.
Then we have the following exact sequence
\begin{equation}
0\to\restr{\F}{V}\to \restr{\T_Y}{V}\to \restr{g^*\T_X}{V}\to\restr{\F^2}{V}\to 0.
\label{eq_tanfol}
\end{equation}

In particular we have the following formula for the first Chern class of the canonical bundle of $Y$
\begin{equation}
c_1(\omega_Y)=c_1(g^*\omega_X)+c_1(\F),
\label{eq_canfol}
\end{equation}
where $\omega_X$ is the dualizing sheaf of $X$.
\end{theorem}

Let us introduce rational vector fields and show how these are related to 1-foliations in the study of inseparable double covers.

\begin{definition}
\label{def_ratvect}
A rational vector field\index{rational vector field} on $Y$ is a rational section $s$ of the tangent space $\T_Y$, that is $s\in H^0(\T_Y\otimes k(Y))$ where $k(Y)$ is the field of rational functions of $Y$.

Two rational vector fields $s_1$ and $s_2$ are said to be equivalent if there exists a rational function $\alpha$ such that $s_1=\alpha s_2$.

A rational vector field is said to be 2-closed if $s^2=\lambda s$ for a constant  $\lambda\in k$ (where the multiplication is the composition of rational derivations): it is called multiplicative if $\lambda=1$, whereas it is called additive if $\lambda=0$.
\end{definition}

Notice that for every rational vector field $s$ we have $[s,s]=0$, in particular, with the notation used for 1-foliations, $s$ is always involutive and that, if $s^2=\lambda s$ with $0\neq\lambda \in k$, then $\sqrt{\lambda^{-1}}s$ is multiplicative.

\begin{proposition}
\label{propo_folratvect}
There exists a bijective correspondence between the set of 1-foliations on a smooth surface $Y$ and the equivalence classes of rational 2-closed vector fields on $Y$.
\end{proposition}
\begin{proof}
Let $s\in H^0(\T_Y\otimes K(Y))$ be a rational 2-closed vector field. 
The tangent sheaf is clearly reflexive: hence combining Theorem 1.9 and Proposition 1.11 of \cite{hartgor} (cf. also \cite{schwede} Theorems 1.17 and 2.10) we know that $s$ is defined up to an effective divisor $D_1$.
In particular there exists a section $l\in H^0(Y,\O_Y(D_1))$ such that $sl\in H^0(Y,\T_Y(D_1))$.
Consider the inclusion $\O_Y\hookrightarrow \T_Y(D_1)$ induced by the global section $sl$: there exists a divisor $D_2$ and a subscheme $Z$ of dimension $0$ such that $\O_Y$ is not a subbundle of $\T_Y(D_1)$ exactly on $D_2\cup Z$.
In particular there exists a section $m\in H^0(Y,\O_Y(D_2))$ such that the morphism $\O_Y\to\T_Y(D_1-D_2)$ induced by $slm^{-1}$ is not a subbundle on $Z$, i.e. it is saturated (cf. Remark \ref{rem_globaldiff}).
Thus, we have the following exact sequence
\begin{equation}
0\to\O_Y\xrightarrow{\cdot slm^{-1}}\T_Y(D_1-D_2)\to Q'\to 0,
\label{eq_ratvecfol}
\end{equation}
where $Q'$ is torsion-free.
Tensoring \ref{eq_ratvecfol} by $\O_Y(D_2-D_1)$ we get
\begin{equation}
0\to\O_Y(D_2-D_1)\to\T_Y\to Q\to 0,
\label{eq_ratvecfol2}
\end{equation}
where $Q$ is torsion-free.
It is clear by the hypothesis on $s$ that $\O_Y(D_2-D_1)$ is a 1-foliation

Conversely, given a 1-foliation $\F$ of $Y$, we already observed that  it is a line bundle: hence there always exists a section $s\in H^0(\F\otimes k(Y))$.
It is straightforward to see that these two construction are one the inverse of the other.
\end{proof}

\section{Inseparable double covers}
\label{sec_insdouble}

We have seen two ways to construct inseparable double covers\index{double cover!inseparable} of surfaces: in Section \ref{sec_double} we have seen how to construct a double cover $f\colon X\to Y$ (separable or inseparable) starting from a surface $Y$ and an exact sequence $0\to O_Y\to f_*\O_X\to L^{-1}\to 0$ and we have given some condition for $f$ to be inseparable; in Section \ref{sec_foliations} we have seen how to construct an inseparable double cover $g\colon Y\to X$ where $X$ is the quotient of $Y$ via a 1-foliation or, equivalently, a rational 2-closed vector field.
In both cases, we assumed $Y$ to be smooth and $X$ to be normal.

We would like to see how these two theories of inseparable double covers are related to one another and the way to do this is to consider the $k$-linear Frobenius morphism $F_k$.
Indeed it is known (cf. Remark \ref{rem_folfrob}) that every inseparable double cover between two schemes is a factor of the $k$-linear Frobenius morphism.
That is, let $f\colon X\to Y$ be an inseparable double cover between two surfaces $X$ and $Y$, then there exists a double cover $g\colon Y\to X^{(-1)}$ such that the following diagram commutes
\begin{equation}
\begin{tikzcd}
X\arrow{r}{f}\arrow[bend left=30]{rr}{F_k} & Y \arrow{r}{g} & X^{(-1)}.
\end{tikzcd}
\label{eq_facfrob}
\end{equation}

\begin{definition}
In the situation above, we define $f$ to be the dual inseparable double cover associated with $g$ and vice-versa. \index{double cover!inseparable dual}
\end{definition}

Now suppose that $X$ is normal (hence also $X^{(-1)}$ is, because it is abstractly isomorphic to $X$) and $Y$ smooth.
In this case we have seen (cf. Section \ref{sec_double}) that $f$ is uniquely determined by a line bundle $L$, an element $[c_{ij}]\in H^1(L)$ which gives an exact sequence as in \ref{eq_seqdouble} and by $b_i$ satisfying the condition 
\begin{equation}
b_i=g_{ij}^2b_j+c_{ij}^2
\label{eq_bi2ins}
\end{equation}
(recall that $a_i=0$ because $f$ is inseparable), while (cf. Section \ref{sec_foliations}) $g$ is uniquely determined by a 1-foliation $\F$ which turns out to be a line bundle.
Ekedahl has proved (cf. \cite{ekecan} pages 105-106) that the relation between $\F$ and $L$ is the following:
\begin{equation}
-2L=K_Y+\F
\label{eq_lf}
\end{equation}
and this clarifies the notation used in the exact sequences \ref{eq_lfoliation} and \ref{eq_foliationl}.

We have seen how to associate to a given line bundle $L$ and local sections $b_i$ satisfying equations \ref{eq_bi2ins} an exact sequence as in \ref{eq_lfoliation}.
Taking the dual of this exact sequence we get the exact sequence \ref{eq_foliationl}. 
So we have that $L^{-2}$ is a saturated subsheaf of $\Omega^1_Y$.
In particular another relation between $\F$ and $L$ (ibidem) is given by the fact that $\F$ (considered as a subsheaf of $\T_Y$) is the annihilator of the image of $L^{-2}$ inside the cotangent sheaf.

\begin{remark}
\label{rem_torsor}
Because the $a_i$'s vanish, we can replace the exact sequence \ref{eq_exactv} by the following one:
\begin{equation}
0\to L^2\to V'\to \O_Y\to 0,
\label{eq_exactv'}
\end{equation}
where $V'$ is a vector bundle whose transition functions are the matrices
\begin{equation}
\begin{pmatrix}
	1 & 0\\
	c_{ij}^2 & g_{ij}^2
\end{pmatrix}
\label{eq_transv'}
\end{equation}
and $(1,b_i)$'s glue to a global section of $V'$.
In particular this proves that the exact sequence \ref{eq_exactv'} splits, i.e. $c_{ij}^2=0\in H^{1}(L^2)$.
Notice that the vanishing of $c_{ij}^2$ as a cohomology class is given by equation \ref{eq_bi2ins}. 

If one looks carefully at the transition functions, immediately  notices that the sequence \ref{eq_exactv'} is nothing but the pull-back via the $k$-linear Frobenius morphism of the dual of the sequence \ref{eq_seqdouble} (here we are using that $Y$ and $Y^{(1)}$ are abstractly isomorphic).

In \cite{ekecan} Proposition 1.11, this result is used to prove that every such double cover $f\colon X\to Y$ is an $\alpha_L$-torsor where $\alpha_L$ is the finite flat group scheme which fits in the following exact sequence
\begin{equation}
0\to\alpha_L\to L\to L^2\to 0
\label{eq_seqtorsor}
\end{equation}
where we are looking at $L$ and $L^2$ as smooth $Y$-group schemes and the morphism is given by the Frobenius morphism.
We will not use this interpretation of double covers in this thesis.
\end{remark}

If we assume that $X$ is smooth too (and we can always reduce to this case after a finite number of blow-ups, cf. Section \ref{sec_canres}), we can consider $f\colon X\to Y$ as a quotient via a 1-foliation $\G$ and the morphism $g^{(1)}\colon Y^{(1)}\to X$ as a double cover as in Section \ref{sec_double} with associated line bundle $M$. 
A simple calculation leads to the following.
\begin{lemma}
\label{lemma_flgm}
We have $f^* L=\G$, $f^*\F=2M$ and $-2M=K_X+\G$.
\end{lemma} 

\begin{proof}
We know that the canonical bundle of $X$ is isomorphic to $f^*K_Y+\G$ (Theorem \ref{teo_canfol}) and to $f^*(K_Y+L)$ (Equation \ref{eq_candouble}) from which we have $f^* L=\G$.
Then, $-2M=K_X+\G$ is simply Equation \ref{eq_lf} rephrased on $X$.
Moreover, pulling back Equation \ref{eq_lf}, we get $f^*\F=-f^*K_Y-2f^*L=-K_X-\G=2M$.
\end{proof}

Let's now consider a smooth surface $X$ and its Albanese morphism $\alb_X\colon X\to\Alb(X)$. 
In addition, suppose that the dimension of the image of the Albanese morphism is two and that $\alb_X$ is inseparable.

\begin{lemma}
\label{lemma_albinsfact}
There is a canonical way to factor $\alb_X$ as
\begin{equation}
X\xrightarrow{f} Y\xrightarrow{h} \Alb(X),
\label{eq_factalbins}
\end{equation}
where $f$ is inseparable of degree two.
\end{lemma}
\begin{proof}
Because $\alb_X$ is inseparable we have that the image of $\alb_X^*\Omega^1_{\Alb(X)}\to\Omega^1_X$ or, equivalently, the relative tangent sheaf $\T_{X/\Alb(X)}$ has rank $1$ (cf. Remark \ref{rem_albfrob}).
We want to show that $\T_{X/\Alb(X)}$ is a 1-foliation: this will give rise to the desired morphism $f\colon X\to Y$.
Using the short exact sequence of the relative tangent bundle
\begin{equation}
0\to\T_{X/\Alb(X)}\to\T_X\to \alb_X^*\T_{\Alb(X)},
\label{eq_shorttan}
\end{equation}
we have that $\T_{X/\Alb(X)}$ is torsion-free and saturated: in particular it is a line bundle (cf. Remark \ref{rem_folfree}).
It is also clear that $\T_{X/\Alb(X)}$ is involutive because it has rank one.
The short exact sequence also shows that $\T_{X/\Alb(X)}$ is the sheaf of local sections of $\T_X$ which vanishes on the pull-back of local functions of $\Alb(X)$: this is clearly a condition which is preserved by taking squares, in particular $\T_{X/\Alb(X)}$ is 2-closed and hence a 1-foliation.
It is then clear then that $\alb_X$ factorizes as $X\xrightarrow{f}Y=X/\T_{X/\Alb(X)}\xrightarrow{g}\Alb(X)$.
\end{proof}

\begin{remark}
\label{rem_albinsfact}
In the same setting as in Lemma \ref{lemma_albinsfact} denote by $\widetilde{M}=\im(\alb_X^*\Omega^1_{\Alb(X)}\to\Omega^1_X)$ and by $\overline{M}$ its saturation inside $\Omega^1_X$.
We have the short exact sequence
\begin{equation}
0\to \overline{M}\to \Omega^1_X \to \Omega^1_{X/\Alb(X)}/T(\Omega^1_{X/\Alb(X)})\to 0.
\label{eq_shortcantanfol}
\end{equation}
Notice that the dual of $\Omega^1_{X/\Alb(X)}/T(\Omega^1_{X/\Alb(X)})$ is equal to the relative tangent sheaf $\T_{X/\Alb(X)}$ (the operation of taking the dual ignores the torsion subsheaf): hence dualizing Equation \ref{eq_shortcantanfol} we get
\begin{equation}
0\to \T_{X/\Alb(X)}\to \T_X \to I_Z\overline{M}^{-1}\to 0,
\label{eq_shortcantanfoldual}
\end{equation}
where $Z$ is the subscheme of zeroes of the foliation $\T_{X/\Alb(X)}$ or, equivalently, the subscheme where $\overline{M}$ is not a subbundle of $\Omega^1_X$.

A Theorem of Igusa \cite{igusa} ensure that the pull-back of global 1-forms via the Albanese morphism is injective, in particular $h^0(\overline{M})\geq q$ where $q$ is the dimension of the Albanese variety.
Let $M=(-K_X-\T_{X/\Alb(X)})/2=-\overline{M}/2$: we have seen that this is the line bundle associated with the double cover $g^{(1)}\colon Y^{(1)}\to X$ (cf. Equation \ref{eq_lf} and Lemma \ref{lemma_albinsfact}) and the square of its inverse has at least $q$ global sections.
This proves that the double cover $g^{(1)}$ is non-splittable (otherwise we would have that $M$ is effective up to a multiple thanks to Remark \ref{rem_effdouble}).
\end{remark}

\section{Canonical resolution of double covers}
\label{sec_canres}
So far we have considered double covers $f\colon X\to Y$ from a normal surface $X$ to a smooth surface $Y$, or, equivalently if the characteristic of the ground field is $2$, from a smooth surface $Y$ to a normal surface $X$ via the theory of foliations.
We would like, somehow, to give a "canonical way" to reduce to the case were both $X$ and $Y$ are smooth.
This is done by taking the canonical resolution\index{double cover!canonical resolution of} and again we need to distinguish the case were the characteristic of the ground field is two or different from two.

\paragraph{$\mathbf{\bcha(k)\neq2}$}
In this case the canonical resolution (cf. \cite{barth} III.7) is classical result  and gives the following diagram:
\begin{equation}
\label{eq_canres}
\begin{tikzcd}
\widetilde{X}=X_t \arrow{d}{f_t}\arrow{r}{\phi_t} & X_{t-1}\arrow{d}{f_{t-1}}\arrow{r}{\phi_{t-1}} & \ldots \arrow{r}{\phi_2} & X_1\arrow{d}{f_1} \arrow{r}{\phi_1} & X \arrow{d}{f}\\
Y_t \arrow{r}{\psi_t} & Y_{t-1}\arrow{r}{\psi_{t-1}} & \ldots \arrow{r}{\psi_2} & Y_1 \arrow{r}{\psi_1} & Y,
\end{tikzcd}
\end{equation}
where the $\psi_i$ are successive blow-ups that resolve the singularities of the branch divisor $R$, the morphism $f_i$ is the double cover branched over $R_i=\psi_i^*R_{i-1}-2m_iE_i$ (equivalently $L_i=\psi_i^*L_{i-1}(-m_iE_i)$), where $E_i$ is the exceptional divisor of $\psi_i$, $m_i=\intinf{d_i/2}$\nomenclature{$\intinf{d}$}{the integral part of a real number $d$} with $d_i$ the multiplicity in $R_{i-1}$ of the blown-up point and $\intinf{d_i/2}$ denotes the integral part of $d_i/2$. 
One has the following relations (cf. \cite{barth} V.22):
\begin{equation}
K_{\widetilde{X}}^2=2K_Y^2+4K_Y.L+2L^2-2\sum_{i=1}^{t}(m_i-1)^2,
\label{cansq}
\end{equation}
and
\begin{equation}
\chi(\O_{\widetilde{X}})=2\chi(\O_Y)+\frac{1}{2}K_Y.L+\frac{1}{2}L^2-\frac{1}{2}\sum_{i=1}^{t}m_i(m_i-1).
\label{eulchar}
\end{equation}

\begin{definition}
\label{def_cansing}
The singularities\index{singularity!negligible} of the branch locus $R$ are said to be negligible if $m_i=1$ (or, equivalently, $d_i=2,3$) for all $i=1,\ldots,t$.
\end{definition}

\begin{remark}
\label{rem_cansing}
If the branch divisor $R$ has at most negligible singularities, then  $K_{\widetilde{X}}=(f_t\circ\psi_t\circ\ldots\circ\psi_1)^*(K_Y+L)$ (\cite{barth} Theorem III.7.2).
If we further assume that $\widetilde{X}$ is of general type and that $Y$ contains no rational curves, we have that $\widetilde{X}$ is minimal (in general it can have exceptional divisors even if $Y$ contains no rational curves) and $X$ is its canonical model (ibidem III.7 table 1).
%
\end{remark}

\begin{remark}
\label{alb}
Suppose that we have a double cover $f\colon X\to Y$ with non-trivial smooth branch divisor $R$ and $\cha(k)=0$. 
Then, if $q(X)=q(Y)$, it follows that $\alb_{f}\colon \Alb(X)\to\Alb(Y)$ is an isomorphism. 
Indeed, because $q(X)=q(Y)$, the morphism $\alb_{f}$ is an isogeny and so is, by duality, $f^*\colon\Pic^0(Y)\to\Pic^0(X)$. 
Suppose that there exists a non-trivial element $\eta\in \ker({f}^*)$. 
This in particular means that $\eta$ is a torsion element ($\ker(f^*)$ is a finite group) and $f^*(\eta)=0$. 
If we consider the \'etale cover $Z\to{Y}$ given by $\eta$ and we complete the diagram as follows
\[
\begin{tikzcd}
\bigsqcup_{i=1}^{ord(\eta)}{X} \arrow[dr, phantom, "\square"] \arrow[r] \arrow[d]  & {X}\arrow[d, "{f}"]\\ 
Z \arrow[r] & {Y},
\end{tikzcd}
\]
we see that ${f}$ factors through $Z$, but this is impossible because it has degree two and has ramification. So $\alb_{f}$ is an isomorphism.




Notice that, if the characteristic is different from zero, the same proof above just shows that $f^*\colon\Pic_Y^0\to \Pic_X^0$ is bijective on closed points, but it could be a purely inseparable isogeny.
Hence we need another proof which works for every characteristic.

As above suppose that $f\colon X\to Y$ is a separable double cover between two smooth surfaces with non-trivial ramification such that $q(X)=q(Y)$.
In particular there exists an involution $\sigma\colon X\to X$ such that $Y=X/\sigma$.
The hypothesis on the ramification tells us that there exists a point $x\in X$ fixed by the action of $\sigma$ which, up to a translation on the Albanese variety, is sent to $0$ by $\alb_X$.
By the universal property of the Albanese variety, we know that $\sigma$ extends to an involution $\alb_\sigma$ on $\Alb(X)$ (which preserves the group structure) such that  the following Cartesian diagram commutes:
\begin{equation}
\begin{tikzcd}
X\arrow{d}{\alb_X}\arrow{r}{\sigma}\arrow[bend left=30]{rr}{f} & X\arrow{d}{\alb_X}\arrow{r}{f} & Y\arrow{d}{\alb_Y}\\
\Alb(X) \arrow{r}{\alb_\sigma}\arrow[bend right=30]{rr}{\alb_f} & \Alb(X)\arrow{r}{\alb_f} & \Alb(Y).
\end{tikzcd}
\label{eq_commalbq=q}
\end{equation}
From $\alb_f\circ\alb_\sigma=\alb_f$ we see that $\alb_f\circ(\alb_\sigma-Id)=0$, from which we derive $\alb_\sigma=Id$ because $\alb_f$ is an isogeny.
The commutative  diagram \ref{eq_commalbq=q}, shows that there exists a morphism $\alpha\colon Y\to \Alb(X)$ such that $\alb_f\circ\alpha=\alb_Y$ from which we can derive $\Alb(X)=\Alb(Y)$.
%
%
\end{remark}

\paragraph{$\mathbf{\bcha(k)=2}$}
As expected, in this case things get more complicated. 
The good thing is that we obtain the same diagram  in \ref{eq_canres} as a canonical resolution and the same equations hold for the square of the canonical divisor and the Euler characteristic of the structure sheaf as in equations \ref{cansq} and \ref{eulchar}; the bad thing is that we have not a characterization of the $m_i$'s.
We just have some partial results, which ensure that in some particular cases we have $m_i=1$ which is particularly nice in order to compute the invariants of $X$.
The results here presented are taken from \cite{gusunzhou} sections 6 and 7.

Let $f\colon X\to Y$ be a double cover with $Y$ smooth and $X$ normal, $q\in Y$ be a point such that the (unique) point $p$ lying over it is singular in $X$ and $\psi\colon Y'\to Y$ be the blow-up of $q$.
It is clear that the short exact sequence which determines the double cover $\overline{X}=Y'\times_Y X\to Y'$ is the pull-back via $\psi$ of the short exact sequence \ref{eq_seqdouble} associated with $f$ and the same is true for the local sections defining the double cover, i.e. $a_i'=\psi^* a_i$ and $b_i'=\psi^*b_i$.
This shows that the singularities of $\overline{X}$ (which are defined by equations depending on $a_i'$ and $b_i'$) lie exactly above those of $X$ (which are defined by equations depending on $a_i$ and $b_i$).
In particular we have that $\overline{X}$ has a divisorial singularity and can not be normal.
Let $\nu\colon X' \to\overline{X}$ the normalization  and consider the following diagram:
\begin{equation}
\begin{tikzcd}
X'\arrow[bend left=20]{drrr}{\phi}\arrow[bend right=20]{dddr}{f'}\arrow{dr}{\nu} & & &\\
 & \overline{X} \arrow[ddrr, phantom, "\square"]\arrow{rr}{\widetilde{\phi}}\arrow{dd}{\widetilde{f}} & & X\arrow{dd}{f}\\
 & & & \\
 & Y'\arrow{rr}{\psi} & & Y.
\end{tikzcd}
\label{eq_doublecvonorm}
\end{equation}  
We have the following relation between the two sequences defining the double covers $f'$ and $\widetilde{f}$:
\begin{equation}
\begin{tikzcd}
0 \arrow{r} & \O_{Y'}\arrow{r}{}\arrow[equal]{d}{} & \widetilde{f}_*\O_{\overline{X}}\arrow[hook]{d}{}\arrow{r}{} & (\psi^*L)^{-1}\arrow{r}{}\arrow[hook]{d}{} & 0\\
0 \arrow{r} & \O_{Y'}\arrow{r}{} & f'_*\O_{X'}\arrow{r}{} & (L')^{-1}\arrow{r}{} & 0,
\end{tikzcd}
\label{eq_definingequations}
\end{equation}
where the vertical injections follow from the fact that $X'$ is the normalization of $\overline{X}$.
In particular this proves that $L'=\psi^*L(-mE)$ for a strictly positive integer $m$. 
Reiterating this process we obtain the canonical resolution: the fact that this process leads to a smooth $\widetilde{X}$ is proven in \cite{gusunzhou}   section 6.2 for the inseparable case and section 6.3 for the separable one.
The invariants of $X$ satisfy Equations \ref{cansq} and \ref{eulchar}.
Notice that in this case it is not clear a priori what are the values of the $m_i$'s.

\begin{remark}
\label{rem_mi2}
Observe that, if all the $m_i$'s are equal to $1$ and $\widetilde{X}$ is a surface of general type, all the singularities of $X$ are rational double points\index{rational double point} (we can derive this from Equation \ref{eulchar}). 
In particular, if $Y$ contains no rational curves, then $\widetilde{X}$ is minimal and $X$ is its canonical model.
\end{remark}

\paragraph{Inseparable double cover}
We want to focus a little bit more on the canonical resolution of inseparable double covers.
As we have already seen (cf. Section \ref{sec_insdouble}), for every inseparable double cover $f\colon X\to Y$ there is a canonically defined inseparable double cover $g\colon Y\to X^{(-1)}$.
Fitting this in the canonical resolution we get:
\begin{equation}
\label{eq_canresins}
\begin{tikzcd}
\widetilde{X}=X_t \arrow{d}{f_t}\arrow{r}{\phi_t} & X_{t-1}\arrow{d}{f_{t-1}}\arrow{r}{\phi_{t-1}} & \ldots \arrow{r}{\phi_2} & X_1\arrow{d}{f_1} \arrow{r}{\phi_1} & X \arrow{d}{f}\\
Y_t \arrow{d}{g_t}\arrow{r}{\psi_t} & Y_{t-1}\arrow{r}{\psi_{t-1}}\arrow{d}{g_{t-1}} & \ldots \arrow{r}{\psi_2} & Y_1 \arrow{r}{\psi_1}\arrow{d}{g_1} & Y\arrow{d}{g}\\
X_t^{(-1)} \arrow{r}{\phi_t^{(-1)}} & X_{t-1}^{(-1)}\arrow{r}{\phi_{t-1}^{(-1)}} & \ldots \arrow{r}{\phi_2^{(-1)}} & X_1^{(-1)} \arrow{r}{\phi_1^{(-1)}} & X^{(-1)},
\end{tikzcd}
\end{equation}
where the $\psi_i$ are successive blow-ups among the points of $Y_{i-1}$ dominating singular points of $X_{i-1}^{(-1)}$, the morphism $f_i$ is the double cover defined by the short exact sequence $0\to \O_{Y_i}\to V_i \to L_i^{-1}\to 0$ with $L_i=\psi_i^*L_{i-1}(-m_iE_i)$), where $E_i$ is the exceptional divisor of $\psi_i$ and $\psi_i^*(f_{i-1})_*\O_{X_{i-1}}\subseteq V_i=(f_i)_*\O_{X_i}$; similarly the $g_i$ are defined by the 1-foliation $\F_i=-2L_i-K_{Y_i}$ (cf. Equation \ref{eq_lf}), in particular $\F_i=\psi_i^*\F_{i-1}\bigl((2m_i-1)E\bigr)$. 

Denoting the by $Z_i$ the zero set of the 1-foliation $\F_i$, thanks to Equation \ref{eq_eulfol}, we get that the degree $d_i$ of $Z_i$ is equal to $e(Y_i)+2L_i.(2L_i+K_{Y_i})$.
In particular one has
\begin{equation}
\begin{split}
d_i=e(Y_i)+&2L_i.(2L_i+K_{Y_i})=\\
=e(Y_{i-1}&)+1+\\
+2\psi_i^*L_{i-1}(-m_iE).\Bigl(2\psi_i^*&L_{i-1}(-m_iE)+\psi_i^*K_{Y_{i-1}}(E)\Bigr)=\\
=e(Y_{i-1})+2L_{i-1}.(2L_{i-1}&+K_{Y_{i-1}})-4m_i^2+2m_i+1=\\
=d_{i-1}-4&m_i^2+2m_i+1.
\end{split}
\label{eq_degzi}
\end{equation}
This equation is used in \cite{gusunzhou} to prove that the canonical resolution process terminates after a finite number of blow-ups, and that's because the degree of $Z_i$ has to be non-negative and $X_i$ is smooth as soon as $d_i=0$.

\begin{remark}
Notice that Equation \ref{eq_degzi} implies that if the zero set of the foliation has only points with multiplicity at most $10$, then all the $m_i$'s appearing in the canonical resolution process have to be equal to $1$. 
Indeed if it were not the case assume, for simplicity, $m_1\geq 2$. 
Then we see that  $d_1=\deg(Z)-4m_1^2+2m_1+1\leq \deg(Z)-11$: observe that $Z_1$ coincides with $\psi_1^*Z$ outside of $E_1$ and we get a contradiction.
As far as we know the opposite may not be true: a priori it may happen that $Z$ has points with multiplicity bigger then $10$ and $X$ has only rational singularities.
\label{rem_atmost10}
\end{remark}

\section{Double covers and fundamental group}
\label{sec_doublefundgroup}
In this Section we are going to study when a double cover induces an isomorphism of algebraic fundamental groups.
Here, for double cover, we mean, as in Section \ref{sec_double}, a finite morphism $f\colon X\to Y$ of degree $2$ with $Y$ smooth and $X$ normal.
If the double cover is inseparable, this is clear because, as we have seen in Section \ref{sec_insdouble}, any such morphism is a factor of the $k$ linear Frobenius morphism and, in particular, is a universal homeomorphism (cf. Remark \ref{rem_homeouniv}) and Proposition \ref{propo_fundgrouphomeo} applies. 
The main Theorem of this Section proves that a separable double cover branched over an ample divisor induces an isomorphism of the algebraic fundamental group.
Before doing this we see a different construction of double covers with respect to the one presented in Section \ref{sec_double}.

Let $f\colon X\to Y$ be a double cover with $Y$ smooth and $X$ normal so that $f$ is flat (cf. miracle flatness Theorem \cite{matsumura} Theorem  23.1).
Recall that such a double cover determines (cf. Section \ref{sec_double}) a sequence
\begin{equation}
0\to\O_Y\to f_*\O_X\to L^{-1}\to 0
\label{eq_seqdoublebis}
\end{equation}
where the transition functions of the sections of $L$ are $g_{ij}$ and the transition functions of the sections $f_*\O_X$ are the matrices
\begin{equation}
\begin{pmatrix}
	1 & -g_{ij}^{-1}c_{ij}\\
	0 & g_{ij}^{-1}
\end{pmatrix}
\label{eq_transitionsbis}
\end{equation} and $L$ is called the line bundle associated with the double cover $f$.
Recall also that in characteristic different from $2$, the exact sequence \ref{eq_seqdoublebis} splits or, equivalently, $c_{ij}$ is cohomologous to zero while in characteristic two it can happen that it does not split.

Now consider the projective bundle $\pi\colon\P=\P(f_*\O_X) \to Y$ with its corresponding line bundle $\O_\P(1)$ and the affine bundle of rank two $\V=\specu(\bigoplus_{i=1}^\infty S^if_*\O_X)$ such that, if we denote by $Z$ the zero section of $\V$, we have a natural surjective morphism $\V\setminus Z\twoheadrightarrow \P$.
We can see $X$ as a closed subscheme of $\P$ such that the restriction $\restr{\pi}{X}$ is equal to $f$ (recall that $\P=\proju(\bigoplus_{i=1}^\infty S^if_*\O_X)$, cf. Definition \ref{def_projbundle}, and $X=\specu(f_*\O_X)$, cf. Section \ref{sec_double}).
Indeed, via the natural surjective morphism of sheaves of algebras $\bigoplus_{i=1}^\infty S^if_*\O_X\twoheadrightarrow f_*\O_X$, we can see $X$ as a closed subscheme of $\V$:  in local coordinates (cf. Equation \ref{eq_doublecov}) we get that this morphism is defined by
\begin{equation}
\begin{split}
A_i[x_i,y_i]&\to A_i[t_i]/(t_i^2+a_it_i+b_i)\\
x_i&\mapsto 1\\
y_i&\mapsto t
\end{split}
\label{eq_doublecovtetra}
\end{equation}
and its kernel is generated by $x_i-1$ and $y_i^2+a_ix_iy_i+b_ix_i^2$ where $x_i$ and $y_i$ are, respectively, $1$ and $t_i$ seen as generators of the free $A_i$-module $A_i[t_i]/(t_i^2+a_it_i+b_i)$, i.e. they are local coordinates for $\V$.
It is then clear that $X$, as a closed subset of $\V$, does not have points in common with $Z$ and that the composition $X\to\V\setminus Z\to\P$ is a closed immersion with the desired property.
Observe that the glueing conditions for $\P$ (or, equivalently, for $\V$) give
\begin{equation}
x_i=x_j
\label{eq_xidc}
\end{equation}
and
\begin{equation}
y_i=g_{ij}y_j+c_{ij}
\label{eq_yidc}
\end{equation}
where $g_{ij}$ and $c_{ij}$ are as in Equation \ref{eq_ti} and that $X$ is defined inside $\P$ by the sheaf of homogeneous ideals generated locally by $y_i^2+a_ix_iy_i+b_ix_i^2$.

\begin{remark}
\label{rem_projbunddoubl}
Observe that the projective bundle $\P$ defined above, is nothing but the projectivization of the affine bundle $W$ described in Remark \ref{rem_canonical}, where the divisor at infinity is in the linear class of $\O_\P(1)$: in particular we have $\restr{\O_\P(1)}{X}=\O_X$.
Let's make this more explicit; 
in local coordinates, we can factor the morphism $A_i[x_i,y_i]\to A_i[t_i]/(t_i^2+a_it_i+b_i)$ described above as $A_i[x_i,y_i]\to A_i[t_i]\to A_i[t_i]/(t_i^2+a_it_i+b_i)$ where the left arrow is defined by $x_i\mapsto 1$, $y_i\mapsto t_i$ and the right arrow is the natural surjection.
Recall that $W$ is given by the glueing of $A_i[t_i]$ and notice that $W$, seen as a closed subscheme of $\V$, does not have points in common with $Z$ and is isomorphic with its image through the composition morphism $X\to\V\setminus Z\to\P$: in particular this shows that the closed inclusion $X\hookrightarrow\P$ is a composition of a closed inclusion $X\hookrightarrow W$ and an open immersion $W\hookrightarrow\P$.
Observe that the complement of $W$ defines a section $i\colon Y\to\P$ of $\pi\colon\P\to Y$ which is defined by the sheaf of homogeneous ideals locally generated by $x_i$.


By Proposition \ref{propo_morphtoproj}, the surjective morphism $f_*\O_X\twoheadrightarrow L^{-1}$ corresponds to a section $i\colon Y\to\P$ of $\pi\colon\P\to Y$ such that $i^*\O_\P(1)=L^{-1}$.
Observe that locally this map of $\O_Y$-modules is defined as 
\begin{equation}
\begin{split}
A_ix_i\oplus A_iy_i&\to A_i\\
x_i&\mapsto 0\\
y_i&\mapsto 1:
\end{split}
\label{eq_doublecovtetrabis}
\end{equation}
in particular the embedding of $Y$ in $\P$ via $i$ coincides with the complement of $W$ and $\restr{\O_\P(Y)}{X}=\O_X$ because $X$ and $Y$ have no points in common.

By adjunction formula, seeing $Y$ as a closed subscheme of $\P$ via $i$, we see that
\begin{equation}
K_Y=i^*K_\P(Y):
\label{eq_adjdoublprojbund}
\end{equation}
by the Equation of the canonical divisor of a projective bundle \ref{eq_canprojbundle} we can conclude that
\begin{equation}
\begin{split}
K_Y=i^*K_\P(Y)=K_Y\otimes_{\O_Y}\det(f_*\O_X)&\otimes_{\O_Y}i^*(\O_\P(-2)\otimes_{\O_\P}\O_\P(Y))=\\
=K_Y\otimes_{\O_Y}L&\otimes_{\O_Y}i^*\O_\P(Y)
\end{split}
\label{eq_adjdoublprojbundbis}
\end{equation}
which means that 
\begin{equation}
i^*\O_\P(Y)=L^{-1}=i^*\O_\P(1).
\label{eq_inftydivproj}
\end{equation}
By the formula of the Picard group $\Pic(\P)$ of $\P$ (cf. Equation \ref{eq_picprojbund}) we can conclude that $\O_\P(Y)=\O_\P(1)$.
In particular $\restr{\O_\P(1)}{X}=\O_X$.
\end{remark}

\begin{lemma}
Let $f\colon X\to Y$ be a double cover and $\P=\P(f_*\O_X)$ the projective bundle associated with $f_*\O_X$.
Then the line bundle $\O_\P(X)$ associated with the divisor $X$ is isomorphic to  $\pi^*L^2\otimes_{\O_\P} \O_\P(2)$.
\label{lemma_doubleprojbund}
\end{lemma}

\begin{proof}
By the formula of the Picard group $\Pic(\P)$ of $\P$ (cf. Equation \ref{eq_picprojbund}), we have that $\O_\P(X)=\pi^*M\otimes_{\O_\P}\O_\P(k)$ for suitable $M\in\Pic(Y)$ and $k\in\Z$.
Because $X$ is a double cover of $Y$, we have that $k=2$ (indeed the intersection between $X$ and a fibre of $\pi$ consists of two points).

Pushing forward with $\pi$ the exact sequence
\begin{equation}
0\to\O_\P(-X)\to\O_\P\to\O_X\to 0
\label{eq_defxp}
\end{equation}
and using the projection formula (cf. \cite{hart} Exercise III.8.3) we get
\begin{equation}
0\to\O_Y\to f_*\O_X\to L\otimes_{\O_Y}M^{-1}=R^1\pi_*\O_\P(-X)\to 0
\label{eq_defxppush}
\end{equation}
(cf. Equation \ref{eq_higherpushprojbund} for the formulae of the derived pushforward of $\O_\P(n)$ for a projective bundle).

Comparing the exact sequence \ref{eq_defxppush} with \ref{eq_seqdoublebis}, we get $M=L^2$: hence $\O_\P(X)=\O_\P(2)\otimes_{\O_\P}\pi^*L^2$ as desired.
%
\end{proof}

We are ready to prove the main result of this Section.

%
%

\begin{theorem}
\label{teo_doubcovamplefundgr}
Let $f\colon X\to Y$ be a separable  double cover branched over an ample divisor.
Then the induced morphism on algebraic fundamental groups is an isomorphism.
\end{theorem}

\begin{proof}
Let $\pi\colon\P=\P(f_*\O_X)\to Y$ be the projective bundle associated with $f_*\O_X$, $i\colon X\hookrightarrow \P$ be the inclusion seen above such that $\pi\circ i=f$.
Let $x\colon\spec(k)\to X$ be a closed point of $X$ and $y=f\circ x$.
By Corollary \ref{coro_projbundlefund}, we know that $\pi_*\colon\pi_1(\P,x)\to\pi_1(Y,y)$ is an isomorphism: we are going to show that $i_*\colon\pi_1(X,x)\to\pi_1(\P,x)$ is an isomorphism using the Lefschetz hyperplane Theorem \ref{teo_lefhypsec}.

If we tensor the short exact sequence \ref{eq_seqdoublebis} with $L^2$ we get
\begin{equation}
0\to L^2\to f_*\O_X\otimes_{\O_Y}L^2\to L\to 0.
\label{eq_seqdoublel2}
\end{equation}
Recall that there exists a natural isomorphism $\phi\colon\P(f_*\O_X\otimes_{\O_Y}L^2)\to\P$ such that $\O_{\P(f_*\O_X\otimes_{\O_Y}L^2)}(1)=\O_\P(1)\otimes_{\P}L^2$ (Equation \ref{eq_diffprojbundle}) and, because $f_*\O_X\otimes_{\O_Y}L^2$ is the central term of a short exact sequence where both the external terms are ample line bundles, Corollary 3.4 of \cite{hartample} shows that $\O_{\P(f_*\O_X\otimes_{\O_Y}L^2)}(1)$ is ample.

As we have shown in Lemma \ref{lemma_doubleprojbund}, the line bundle associated with the divisor $X$ inside $\P$ is $\pi^*L^2\otimes_{\O_\P}\O_\P(2)$.
Moreover $L^2$ admits a section on $Y$ whose zero set $D$ is the branch divisor if the characteristic is different from two and twice the branch divisor if the characteristic is two (cf. Definitions \ref{def_branch} and \ref{def_branch2}).
In particular we can find a divisor $X'$ inside the linear class of the ample line bundle $\pi^*L^4\otimes_{\O_\P}\O_\P(2)$ which consists in the union of $X$ with the pull-back $\pi^*D$ of $D$.
More precisely, denote by $\mathcal{I}$ the ideal sheaf of $X$, by $\mathcal{J}$ the ideal sheaf of $\pi^*D$, and by $Z$ the closed subscheme of $\P$ cut out by the ideal sheaf $\mathcal{I}+\mathcal{J}$ (notice that $Z$, as a closed subscheme of $X$, can be identified with a multiple of the ramification divisor of $f$, hence its reduction is isomorphic to the reduction of $D$ which is the branch divisor of $f$).
By \cite{schwedegl} Corollary 3.9 and the discussion just after it, we know that the push-out $\pi^*D\cup_Z X$\nomenclature{$X\cup_Z Y$}{the push-out of the diagram of schemes $X\leftarrow Z\rightarrow Y$} exists in the category of schemes and its ideal sheaf in $\P$ is $\mathcal{I}\cap\mathcal{J}$ and, by definition, it coincides with $X'$.
Moreover in the following cocartesian diagram
\begin{equation}
\begin{tikzcd}
Z\arrow[hook]{r}{h}\arrow[hook]{d}{k} & X\arrow[hook]{d}{j}\\
\pi^*D\arrow[hook]{r}{l} & X'.
\end{tikzcd}
\label{eq_cartcocart}
\end{equation}
all morphisms are closed inclusion and, by \cite{ferrand} Th\'eor\`eme 5.4 (notice that the hypothesis there required that every finite subset of points is contained in an affine subset is satisfied for every projective scheme over a field thanks to \cite{liu} Proposition 3.3.36(b)), the square \ref{eq_cartcocart} is also Cartesian.

We need to show that the two closed inclusions $j\colon X\to X'$ and $i'\colon X'\to\P$ induce an isomorphism at the level of algebraic fundamental groups.
Applying Theorem \ref{teo_lefhypsec}, we get that $i'_*\colon \pi_1(X',x)\to\pi_1(\P,x)$ is an isomorphism.

Observe that $\restr{\pi}{\pi^*D}\circ k$ is a universal homeomorphism: indeed, thanks to Remark \ref{rem_homeouniv}, in the following commutative diagram
\[
\begin{tikzcd}
Z_{red}\arrow{r}{}\arrow{d}{}& Z\arrow{d}{}\\
D_{red}\arrow{r}{}& D,
\end{tikzcd}
\]
the horizontal arrows are universal homeomorphisms, while the left arrow is clearly an isomorphism.
It follows that $\restr{\pi}{\pi^*D}\circ k$ induces an isomorphism on algebraic fundamental groups (Proposition \ref{propo_fundgrouphomeo}).
In particular, using the fact that $\restr{\pi}{\pi^*D}\pi^*D\to D$ is a $\P^1$ -bundle and Corollary \ref{coro_projbundlefund}, we get that $k$ induces an isomorphism of algebraic fundamental groups and, by Theorem \ref{teo_fundpushfor}, $k^*\colon\Fet_{\pi^*D}\to\Fet_Z$ is an equivalence of categories.

Hence a finite \'etale cover $\widetilde{X}\to X$ induces a finite \'etale cover $\widetilde{\pi^*D}\to\pi^*D$ such that their restrictions on $Z$ coincide, i.e we have the following Cartesian diagrams
\begin{equation}
\begin{tikzcd}
\widetilde{\pi^*D}\arrow{d}{}&\widetilde{Z}\arrow{d}{}\arrow[hook', swap]{l}{\widetilde{k}}\arrow[hook]{r}{\widetilde{h}} & \widetilde{X}\arrow{d}{}\\
\pi^*D & Z\arrow[hook', swap]{l}{k}\arrow[hook]{r}{h} & X.
\end{tikzcd}
\label{eq_cartcocartet}
\end{equation}
Clearly $\widetilde{h}$ and $\widetilde{k}$ are closed inclusions, therefore we can apply Corollary 3.9 of \cite{schwedegl} which ensures that the push-out $\widetilde{\pi^*D}\cup_{\widetilde{Z}}\widetilde{X}$ exists in the category of schemes and $\widetilde{X}$ and $\widetilde{\pi^*D}$ are closed subschemes of it.
Moreover (by \cite{moret} Equation 3.4.1 or \cite{Anantharaman} Remarque 1.1.5(iii)) we have that the natural morphism $\widetilde{\pi^*D}\cup_{\widetilde{Z}}\widetilde{X}\to X'$ induced by the universal property of the push-out is a finite \'etale cover of the same degree of $\widetilde{X}\to X$.
In addition, we can prove that the following two squares are both Cartesian
\begin{equation}
\begin{tikzcd}
\widetilde{\pi^*D}\arrow{d}{}\arrow[hook]{r}{}&\widetilde{\pi^*D}\cup_{\widetilde{Z}}\widetilde{X}\arrow{d}{} & \widetilde{X}\arrow[hook']{l}{}\arrow{d}{}\\
\pi^*D\arrow[hook]{r}{} & X' & X\arrow[hook']{l}{}:
\end{tikzcd}
\label{eq_cartcocartetbis}
\end{equation}
indeed (we prove this just for the right square because the other case is identical) we have the following diagram
\begin{equation}
\begin{tikzcd}
\widetilde{X}\arrow[bend left=20]{rrd}{}\arrow[bend right=20]{ddr}{} \arrow{dr}{} & & \\
&\overline{X} \arrow[hook]{r}{}\arrow{d}{}& \widetilde{\pi^*D}\cup_{\widetilde{Z}}\widetilde{X}\arrow{d}{}\\
&X\arrow[hook]{r}{} & X'
\end{tikzcd}
\label{eq_maverto}
\end{equation}
where $\overline{X}$ is the pull-back of the diagram, $\widetilde{X}\to{X}$ and $\overline{X}\to{X}$ are finite \'etale morphisms of the same degree.
By \cite{liu} Lemma 3.3.15(b) and Exercise 3.3.17(b), we get that $\widetilde{X}\to\overline{X}$ is a finite morphism of degree $1$ and by \cite{egaiv4} Proposition 17.3.4 it is also \'etale, hence an isomorphism.
In particular we have proved that condition 2 of Theorem \ref{teo_fundpushfor} holds.

Now let $\widetilde{X'}\to X'$ be an \'etale finite cover and consider the following diagram
\begin{equation}
\begin{tikzcd}
\widetilde{Z}\arrow{dr}{}\arrow[hook]{rrr}{}\arrow[hook]{ddd}{} & & & \widetilde{X}\arrow{dl}\arrow[hook]{ddd}{}\\
& Z\arrow[hook]{r}{}\arrow[hook]{d}{}& X\arrow[hook]{d}{} & \\
& \pi^*D\arrow[hook]{r}{} & X'&\\
\widetilde{\pi^*D}\arrow[hook]{rrr}{}\arrow{ur}{}& & & \widetilde{X'}\arrow{ul}{}
\end{tikzcd}
\label{eq_orca}
\end{equation}
where $\widetilde{X}$, $\widetilde{\pi^*D}$ and $\widetilde{Z}$ are obtained by base change so that all the squares in it are Cartesian.
Consider the natural morphism $\widetilde{\pi^*D}\cup_{\widetilde{Z}}\widetilde{X}\to\widetilde{X'}$ induced by the universal property of the push-out and consider the following diagram
\begin{equation}
\begin{tikzcd}
\widetilde{\pi^*D}\cup_{\widetilde{Z}}\widetilde{X}\arrow{rr}{}\arrow{dr}{}& &\widetilde{X'}\arrow{dl}{}\\
& X': &
\end{tikzcd}
\label{eq_orcax}
\end{equation}
the two oblique arrows are finite \'etale of the same degree and arguing as above we conclude that the horizontal one is an isomorphism. 
This shows that every finite \'etale cover $\widetilde{X'}\to X'$ is uniquely determined by its pull-back on $X$ and in particular if $\widetilde{X'}$ is connected then so is, using the notation of the diagram \ref{eq_orca}, $\widetilde{X}$ (indeed if it were not the case we would get that $\widetilde{\pi^*D}\cup_{\widetilde{Z}}\widetilde{X}\simeq\widetilde{X'}$ is disconnected, a contradiction).
This means that also condition 1 of Theorem \ref{teo_fundpushfor} holds, i.e. the morphism $j_*\colon\pi_1(X,x)\to\pi_1(X',x)$ is an isomorphism and the proof is completed.
\end{proof}

\begin{corollary}
\label{coro_doubcovamplefundgr}
In the same situation of Theorem \ref{teo_doubcovamplefundgr}, let 
\begin{equation}
\begin{tikzcd}
\widetilde{X}\arrow{d}{\widetilde{f}}\arrow{r}{\phi} & X\arrow{d}{f}\\
\widetilde{Y}\arrow{r}{\psi} & Y
\end{tikzcd}
\label{eq_quad}
\end{equation}
be the canonical resolution of $f\colon X\to Y$ (cf. Section \ref{sec_canres}).
Suppose that all the $m_i$'s appearing are equal to $1$ or, equivalently, that $X$ has only rational double points (e.g. $X$ is the canonical model of $\widetilde{X}$ if $\widetilde{X}$ is a surface of general type).
Then $\phi$ and, consequently, $f\circ\phi$ induce an isomorphism on the algebraic fundamental groups.
\end{corollary}

\begin{proof}
Let $\widetilde{\pi}\colon \widetilde{X}'\to \widetilde{X}$ be a finite \'etale morphism and consider the Stein factorization $\widetilde{X}'\xrightarrow{\phi'}X'\xrightarrow{\pi} X$ of $\phi\circ\widetilde{\pi}$ so that we get the following commutative diagram
\begin{equation}
\begin{tikzcd}
\widetilde{X}'\arrow{r}{\widetilde{\pi}}\arrow{d}{\phi'} & \widetilde{X}\arrow{d}{\phi}\\
X'\arrow{r}{\pi} & X.
\end{tikzcd}
\label{eq_etalefundgr}
\end{equation}

Because $X$ has only rational double points as singularities,  the exceptional divisor $Z$ of $\phi$ (i.e. the divisor contracted by $\phi$) is a disjoint union of  $k$ trees of rational smooth curves (\cite{cossec} Proposition 0.2.4), hence its algebraic fundamental group is trivial (indeed, by Riemann-Hurwitz cf. \cite{hart} Corollary IV.2.4, the only finite \'etale cover of the projective line is the identity).
In particular we have that the \'etale cover $\widetilde{\pi}$ restricted to the pre-image ${Z'}$ of $Z$ consists of a finite union of $n$ disjoint copies of $Z$ where $n$ is the degree of $\widetilde{\pi}$.
Notice in particular that ${Z'}$ is contracted by $\phi'$ into $nk$ distinct points, that $\phi'$ is an isomorphism outside ${Z'}$ and that the resulting singularities of $X'$ are rational double points as well by Artin's contractibility criterion (cf. \cite{badescu} Theorem 3.15). 
We know that $\pi$ is an \'etale morphism when restricted to the smooth locus.
By \cite{redbook} Theorem 3 of Section III.5 we know that, in order to prove that $\pi$ is \'etale everywhere, it is enough to show that the natural morphism $\widehat{\O_{X,x}}\to\widehat{\O_{X',x'}}$ is an isomorphism for every closed singular point $x'$ of $X'$ such that $\pi(x')=x$, where by $\widehat{\O_{X,x}}$ we mean the completion of the local ring $\O_{X,x}$ with respect to its unique maximal ideal $m_x$.
Denote by $\widetilde{X}'_n$ (respectively $\widetilde{X}_n$) the fibre of $\phi'$ (respectively $\phi$) above $\spec(\O_{X',x'}/m_{x'}^n)$ (respectively $\spec(\O_{X,x}/m_x)$): by the Theorem on formal functions (cf. \cite{hart} Theorem III.11.1) we know that $\widehat{\O_{X',x'}}\simeq\varprojlim H^0(\widetilde{X}'_n,\O_{\widetilde{X}'_n})$ (respectively $\widehat{\O_{X,x}}\simeq\varprojlim H^0(\widetilde{X}_n,\O_{\widetilde{X}_n})$).
Because $\widetilde{\pi}$ is an \'etale morphism and $\widetilde{X}_n$ has trivial algebraic fundamental group (recall that the algebraic fundamental group ignores the non-reduced structure, cf. Remark \ref{rem_homeouniv} and Proposition \ref{propo_fundgrouphomeo}, and that the reduction $\widetilde{X}_1$ of $\widetilde{X}_n$ is contained in $Z$ by definition), it induces an isomorphism $\widetilde{X}'_n\simeq\widetilde{X}_n$, hence an isomorphism $ H^0(\widetilde{X}'_n,\O_{\widetilde{X}'_n})\simeq H^0(\widetilde{X}_n,\O_{\widetilde{X}_n})$ for every natural number $n$, therefore on the projective limit.
In particular $\pi$ induces an isomorphism $\widehat{\O_{X,x}}\simeq\widehat{\O_{X',x'}}$, i.e. $\pi$ is \'etale.

Now consider the fibre product $X'\times_X\widetilde{X}$ and the natural morphism $\alpha\colon\widetilde{X}'\to X'\times_X\widetilde{X}$, defined by the commutative diagram \ref{eq_etalefundgr}, which fit in the following commutative diagram:
\begin{equation}
\begin{tikzcd}
\widetilde{X}'\arrow[bend left=20]{rrd}{\widetilde{\pi}}\arrow[bend right=20]{ddr}{\phi'} \arrow{dr}{\alpha} & & \\
&X'\times_X\widetilde{X} \arrow{r}{\beta}\arrow{d}{}& \widetilde{X}\arrow{d}{\phi}\\
&X'\arrow{r}{\pi} & X.
\end{tikzcd}
\label{eq_etalefundgrbis}
\end{equation}
Notice that $\widetilde{\pi}$ and $\beta$ are finite \'etale morphisms of degree $n$ hence $\alpha$ is a finite morphism (cf. \cite{liu} Lemma 3.3.15(b) and Exercise 3.3.17(b)) of degree $1$ and \'etale (cf. \cite{egaiv4} Proposition 17.3.4).
It is then clear that $\alpha$ has to be an isomorphism.

We have proved that condition 2 of Theorem \ref{teo_fundpushfor} holds, while condition 1 is trivially satisfied because $\phi$ is a birational morphism of irreducible varieties.
In particular $\phi$ induces an isomorphism on algebraic fundamental groups.
\end{proof}

\begin{remark}
\label{rem_ultimmss}
Observe that, in the same situation of Corollary \ref{coro_doubcovamplefundgr}, we can assume that $f$ is a purely inseparable double cover over a field of characteristic $2$ and obtain the same result.
Indeed, in this case, Proposition \ref{propo_fundgrouphomeo} applies and shows that $f$ induces an isomorphism on algebraic fundamental groups.
Then the proof that the same holds for $\phi$ is only based on the fact that all the $m_i$'s appearing in the canonical resolution are equal to $1$.
\end{remark}

\section{Inseparable double covers of product surfaces}
\label{sec_insdoubleprod}
In this section we recall and expand a simple constructive process which allows us to give many examples of inseparable double covers of product surfaces  (clearly we are focusing over a field of characteristic $2$).
Indeed a rational vector field on a curve, is simply a rational section of the tangent line bundle and for a surface $Y=A\times B$, where $A$ and $B$ are curves, and rational vector fields $\delta_A$ and $\delta_B$, we can easily construct a "product" rational vector field $\delta_A+\delta_B$ on $Y$ which, thanks to Section \ref{sec_foliations}, gives us a double cover $g\colon Y\to X^{(-1)}$: actually, using the notation of Section \ref{sec_insdouble}, we are interested in the dual cover $f\colon X\to Y$.

This construction is taken from \cite{liedtkeuni}, where it is used to construct a family $\{X_i\}$ of uniruled surfaces of general type all with isomorphic Albanese variety, for which the quantities $h^1(X_i,\O_{X_i})$, $h^0(X,\Omega_{X_i}^1)$ and $h^0(X,\Omega_{X_i}^1)-h^1(X_i,\O_{X_i})$ tend to infinity as far as 
$i\to\infty$.
We will use this construction in Section \ref{sec_exchar2} to construct surfaces of general type with maximal Albanese dimension which lie on the Severi lines over a field of characteristic $2$.

First of all we recall some notable examples of rational vector fields on curves: all of them are essentially taken from Section 2 of \cite{liedtkeuni}.

\begin{example}[Multiplicative and additive rational vector fields on an elliptic curve with only two zeroes and poles]
\label{es_multaddvf1}
Over a field of characteristic $2$ an elliptic curve $E$ can be seen as the projectivization inside $\P^2$ of 
\begin{equation}
y^2+\alpha xy+y+x^3
\label{eq_deuring}
\end{equation} 
where $\alpha\in k$ is such that $\alpha^3\neq 1$: this is the so called Deuring normal form (cf. \cite{silverman} Proposition A.1.3).
In particular the $j$-invariant of $E$ is $j(E)=\frac{\alpha^{12}}{\alpha^3-1}$ and $E$ is supersingular (i.e. it has no non-trivial $2$-torsion points) if and only if $\alpha=0$ (ibidem).

Notice that every rational vector field on $E$ is uniquely determined by its value on the function $x$: indeed we can determine its value on $y$ (hence on all the field $k(E)$ of rational functions on $E$) via the condition that it has to vanish on $y^2+\alpha xy+y+x^3$.
Now let $\partial_x$ be the rational vector field on $E$ which correspond to the derivation  by $x$: from the condition $\partial_x(y^2+\alpha xy+y+x^3)=0$ we get
\begin{equation}
(1+\alpha x)\partial_x(y)=x^2+\alpha y:
\label{eq_dery1}
\end{equation}
we claim that $\delta=(1+\alpha x)\partial_x$ is a never-vanishing vector field on $E$.

Indeed $\delta(x)=1+\alpha x=0$ implies $x=\frac{1}{\alpha}$ while $\delta(y)=x^2+\alpha y=0$ implies that $y=\frac{1}{\alpha^3}$, but $(\frac{1}{\alpha},\frac{1}{\alpha^3})$ is not a point of $E$: in particular, because the tangent sheaf of an elliptic curve is trivial, we can extend this to a regular never-vanishing vector field on $E$.

If we consider the rational vector field $\delta_E'=x\delta$: by the previous discussion we easily see that it has a single zero at $(0,0)$ and at $(0,1)$ and a double pole at $\infty$.
We notice that $\delta_E'$ is a multiplicative rational vector field\index{rational vector field!multiplicative}, indeed
\begin{equation}
\delta_E'^2(x)=\delta_E'(x+\alpha x^2)=x+\alpha x^2=\delta_E'(x)
\label{eq_mult1x}
\end{equation}
and
\begin{equation}
\delta_E'^2(y)=\delta_E'(x^3+\alpha xy)=\delta_E'(y^2+y)=\delta_E'(y).
\label{eq_mult1y}
\end{equation}

On the other hand, if we assume $\alpha\neq 0$, a similar calculation shows that $\delta_E=(1+\alpha x)\delta$ is an additive rational vector\index{rational vector field!additive} field with a double zero at $(\frac{1}{\alpha},\frac{1}{\sqrt{\alpha^3}})$ and a double pole at $\infty$.
If $\alpha=0$, i.e. the elliptic curve is supersingular, we believe it is not possible to find an additive rational vector field with only two (possibly equal) zeroes.

Let $m(x)=\prod_{i=1}^d(x+\gamma_i)$, $n(x)=\prod_{i=1}^d(x+\epsilon_i)$, $\beta(x)=\frac{m^2(x)}{\alpha n^2(x)}$ and $\overline{\beta}(x)=\frac{m^2(x)+xn^2(x)}{n^2(x)}$: we can always choose $m(x)$ and $n(x)$ such that $\gamma_i\neq\frac{1}{\alpha},0\neq \epsilon_k$ for all $i$ and $k$ and that all the zeroes $\mu_i$ for $i=1,\ldots 2d+1$ of $m^2(x)+xn^2(x)$ are non-multiple and mutually different from $\frac{1}{\alpha}, 0$ and $\epsilon_i$ ($i=1,\ldots, d$).
A similar computation as above, shows that $\overline{\delta}_E=\beta(x)\delta_E$ and $\overline{\delta}_E'=\overline{\beta}(x)\partial_x$ are, respectively, an additive and a multiplicative rational vector field on $E$.
Moreover their corresponding divisors are, respectively,
\begin{equation}
\begin{split}
2P_{\frac{1}{\alpha}}-2P_\infty+2G_1+2&G_1'+\ldots+2G_d+\\
+2G_d'-2E_1-2E_1'-&\ldots-2E_d-2E_d'
\end{split}
\label{eq_divdelta}
\end{equation}
and
\begin{equation}
\begin{split}
P_{\infty}-2P_{\frac{1}{\alpha}}+M_1+M_1'+\ldots +&M_{2d+1}+M_{2d+1}'-\\
-2E_1-2E_1'-\ldots-&2E_d-2E_d'
\end{split}
\label{eq_divdelta'}
\end{equation}
where $G_i, G_i'$, $E_i, E_i'$ and $M_i, M_i'$  are the points where $x+\gamma_i$, $x+\epsilon_i$ and $x+\mu_i$ vanish, $P_{\frac{1}{\alpha}}=(\frac{1}{\alpha},\frac{1}{\sqrt{\alpha^3}})$ and $P_\infty$ is the point at infinity (notice that, if $\alpha=0$, then $P_{\frac{1}{\alpha}}=P_{\infty}$).
\end{example}


\begin{example}[Multiplicative and additive rational vector fields on a hyperelliptic curve]
\label{es_multaddvf>1}
Over a field of characteristic $2$, a smooth hyperelliptic curve $C$ of genus $g$ can be described as the projectivization inside the total space of the line bundle $\O_{\P^1}(g+1)$ of
\begin{equation}
y^2+yf(x)+g(x)
\label{eq_hyper2}
\end{equation}
where $f(x)$ is a polynomial of degree $g+1$ and $g(x)$ is a polynomial of degree $2g+2$, such that every zero of $f(x)$ is a simple zero of $g(x)$ (cf. \cite{bhosle} Section 1).
Moreover the zeroes of $f(x)$ represent the branch divisor of the hyperelliptic morphism $C\to \P^1$.
In order to have $C$ as tame as possible, i.e. $C\to \P^1$ has the maximum number of branch points, we require $f(x)$ to have no multiple zeroes.

Then we can rewrite the equation of $C$ as
\begin{equation}
y^2+y\prod_{i=0}^{g}(x+\alpha_i)+\prod_{i=0}^{g}(x+\alpha_i)\prod_{i=0}^{g}(x+\beta_i)
\label{eq_hyper2n}
\end{equation}
where $\alpha_i\neq \alpha_j$ for $i\neq j$ and $\beta_i\neq \alpha_j$ for all $i,j\in\{0,\ldots,g\}$.
Denote by $P_i$ for $i=0,\ldots,g$ the ramification points of $C$ 
and by $P_\infty$, $P_\infty'$ the two points of $C$ lying over $\infty\in\P^1$ and by $h=\prod_{i=0}^{g}(x+\beta_i)$.

Consider, as in Example \ref{es_multaddvf1}, the rational vector field $\delta=\partial_x$ on $C$: the equation defining $C$ tells us that
\begin{equation}
\partial_x(y)=\frac{f'(x)y+f'(x)h(x)+f(x)h'(x)}{f(x)}.
\label{eq_partialyhyper}
\end{equation}
We then easily derive that the divisor corresponding to $(x^2+\alpha)\delta$ is $2P_\alpha+2P_\alpha'-2P_0-\ldots-2P_g$ where $\alpha\neq\alpha_i$ for all $i=0,\ldots,g$ and $P_\alpha$ and $P_\alpha'$ are the two points of $C$ lying over $\alpha\in\P^1$: indeed $f(x)$ vanishes of order $2$ on every $P_i$, $X^2+\alpha$ vanishes of order two on $P_\alpha$ and $P_\alpha'$  and the degree of the tangent line bundle of $C$ is $-2g+2$.
In particular, because the divisor associated with $X^2+\alpha$ is $2P_\alpha+2P_\alpha'-2P_\infty-2P_\infty'$, we have that the divisor associated with $\delta$ is $2P_\infty+2P_\infty'-2P_0-\ldots-2P_g$. 

Let $m(x)=\prod_{i=1}^d(x+\gamma_i)$, $n(x)=\prod_{i=1}^d(x+\epsilon_i)$, $\beta(x)=\frac{m^2(x)}{n^2(x)}$ and $\overline{\beta}(x)=\frac{m^2(x)+xn^2(x)}{n^2(x)}$: we can always choose $m(x)$ and $n(x)$ such that $\gamma_i\neq\alpha_j\neq \epsilon_k$ for all $i$, $j$ and $k$ and that all the zeroes $\mu_i$ for $i=1,\ldots 2d+1$ of $m^2(x)+xn^2(x)$ are non-multiple and mutually different from $\alpha_i$ ($i=0,\ldots g$) and $\epsilon_i$ ($i=1,\ldots, d$).
A similar computation to the one in Example \ref{es_multaddvf1}, shows that $\delta_C=\beta(x)\delta$ and $\delta_C'=\overline{\beta}(x)\delta$ are, respectively, an additive and a multiplicative rational vector field on $C$.
Moreover their corresponding divisors are, respectively,
\begin{equation}
\begin{split}
2P_\infty+2P_\infty'+2G_1+2&G_1'+\ldots+2G_d+\\
+2G_d'-2E_1-2E_1'-\ldots-2E_d&-2E_d'-2P_0-\ldots-2P_g
\end{split}
\label{eq_divdelta}
\end{equation}
and
\begin{equation}
\begin{split}
P_\infty+P_\infty'+M_1+M_1'+\ldots +&M_{2d+1}+M_{2d+1}'-\\
-2E_1-2E_1'-\ldots-2E_d-2E_d'-&2P_0-\ldots-2P_g
\end{split}
\label{eq_divdelta'}
\end{equation}
where $G_i, G_i'$, $E_i, E_i'$ and $M_i, M_i'$  are the points where $x+\gamma_i$, $x+\epsilon_i$ and $x+\mu_i$ vanish.
\end{example}

\begin{proposition}
\label{propo_liedtkeaxb}
Let $Y=A\times B$ be a product surface where $g(A), g(B)\geq 1$, $\delta_A$ and $\delta_B$ be two rational additive (respectively multiplicative) vector fields on $A$ and $B$ whose zeroes and poles have multiplicity at most $2$ both of which have at least one zero.
Then the rational vector field $\delta_A+\delta_B$ naturally induced on $Y$ is still additive (respectively multiplicative).

If we consider the induced inseparable double cover $g\colon Y\to X^{(-1)}=Y/(\delta_A+\delta_B)$ we have that the corresponding 1-foliation $\F$ is isomorphic, as a line bundle, to $\O_{Y}(-\pi_A^*\delta_{A_\infty}-\pi_B^*\delta_{B_\infty})$ where $\delta_{A_\infty}$ (respectively $\delta_{B_\infty}$) is the divisor of poles of $\delta_A$ (respectively $\delta_B$) and $\pi_A$, $\pi_B$ are the projection to $A$ and $B$.
Moreover the singularities of $X^{(-1)}$ lie in the image of the points $(a,b)$ where $a$ and $b$ are both zeroes or both poles of $\delta_A$ and $\delta_B$.

If we consider the dual inseparable double cover $f\colon X\to Y$ we have that its associated line bundle $L$ is isomorphic to $\frac{-K_Y-\F}{2}$ and is effective, in particular $X$ is the canonical model of a surface of general type with maximal Albanese dimension and $q'(X)=q(X)=g(A)+g(B)$.
\end{proposition}

\begin{remark}
Of course $X$ is singular (by hypothesis the foliation $\delta_A+\delta_B$ has isolated zeroes), hence by $\Alb(X)$ we mean the Albanese variety of a minimal smooth model of $X$ and by $q(X)$ we mean, as usual, the dimension of $\Alb(X)$.
\label{rem_alb2qq'}
\end{remark}

\begin{proof}
The fact that $\delta_A+\delta_B$ is additive (respectively multiplicative) is immediate.

It is clear by definition that $\delta_A+\delta_B$ vanishes on points $(a,b)$ where both $a$ and $b$ are vanishing points of $\delta_A$ and $\delta_B$ respectively and has poles on points $(a,b)$ where $a$ or $b$ is a pole of $\delta_A$ or of $\delta_B$ respectively.
In particular via the identification of Proposition \ref{propo_folratvect}, the formula for $\F$ and for singular points of $X^{(-1)}$ easily follows (cf. also \cite{liedtkeuni} Proposition 4.2)

Thanks to the assumption on zeroes and poles of $\delta_A$ and $\delta_B$, $X$ has only rational double points (cf. Remarks \ref{rem_mi2} and \ref{rem_atmost10}). 
The formula for $L$ is simply an application of Equation \ref{eq_lf} and is effective because $\delta_{A_\infty}-K_A>0$ (respectively $\delta_{B_\infty}-K_B>0$) by hypothesis.  
Hence, via Equation \ref{eq_candouble}, we derive that $X$ is the canonical model of a surface of general type.
It is clear by construction and by universal property of the Albanese variety that $\Alb(X)$ lies between $\Alb(Y)$ and $\Alb(Y^{(1)})=\Alb(Y)^{(1)}$ (i.e. there exists a morphism $\Alb(Y)^{(1)}\to\Alb(X)$ and $\Alb(X)\to\Alb(Y)$ such that the composition is the $k$-linear Frobenius) and, similarly, $\Alb(Y)$ lies between $\Alb(X)$ and $\Alb(X)^{(-1)}$.
In particular $X$ has maximal Albanese dimension and $q(X)=g(A)+g(B)$.
Moreover, $q'(X)=g(A)+g(B)$ follows from K\"unneth formula and the long exact sequence in cohomology associated with the short exact sequence \ref{eq_seqdouble} defining the double cover.
\end{proof}

\begin{corollary}
\label{coro_liedtkeaxb}
In the same situation as in Proposition \ref{propo_liedtkeaxb} assume that either the zeroes of $\delta_A$ and $\delta_B$ are simple or they are double.
Let 
\begin{equation}
\begin{tikzcd}
\widetilde{X}\arrow{r}{\phi}\arrow{d}{f_t} & X\arrow{d}{f}\\
Y_t\arrow{r}{\psi} & Y
\end{tikzcd}
\label{eq_canresaxb}
\end{equation}
be the canonical resolution of $f\colon X\to Y$ (cf. Section \ref{sec_canres}).
Then $m_i=1$ for all $i=1,\ldots,t$ and $f_t$ is an inseparable non-splittable double cover.
\end{corollary}

\begin{proof}
We have that $m_i=1$ for all $i=1,\ldots,t$ because all the singularities of $X$ are rational double points.
Because $f$ is inseparable, if $f_t$ were splittable $2L_t$ would have at least a global section (cf. Remark \ref{rem_effdouble}): we are going to show that this is not the case and, in order to do it, we need to better understand the description of $L$.

Let $\delta_{A_0}-\delta_{A_\infty}$ (respectively $\delta_{B_0}-\delta_{B_\infty}$) be the divisor associated with $\delta_A$ (respectively $\delta_B$) where $\delta_{A_0}$, $\delta_{A_\infty}$, $\delta_{B_0}$ and $\delta_{B_\infty}$ are effective with degree $d_{A_0}$, $d_{A_\infty}$, $d_{B_0}$ and $d_{B_\infty}$ respectively.
In Proposition \ref{propo_liedtkeaxb} we have seen that $\F=-\pi_A^*\delta_{A_\infty}-\pi_B^*\delta_{B_\infty}$ and that $2L=-K_Y-\F$.
By definition we have that $\delta_{A_\infty}-K_A=\delta_{A_0}$ and $\delta_{B_\infty}-K_B=\delta_{B_0}$, in particular
\begin{equation}
2L=\pi_A^*\delta_{A_0}+\pi_B^*\delta_{B_0}
\label{eq_2laxb}
\end{equation}
from which we derive $2L.B=d_{B_0}$ and $2L.A=d_{A_0}$ where $A$ and $B$ denote a fibre of $\pi_B$ and $\pi_A$ respectively.

In the canonical resolution of $f\colon X\to Y$, we blow-up all the points $(a,b)$ where both $a$ and $b$ are zeroes or poles of $\delta_A$ and $\delta_B$ and $L_t$ is defined by 
\begin{equation}
L_t=\psi^*L(-\sum_{i=1}^tE_i)
\label{eq_ltcanres}
\end{equation}
because $m_i=1$ for all $i$, where the $E_i$'s are all the exceptional divisors of $Y_t$.
Notice also that if $a$ and $b$ are simple (respectively double) zeroes or poles of $\delta_A$ and $\delta_B$ then the singularity above is an $A_1$ (respectively $D_4$) singularity (cf. \cite{liedtkeuni} Proposition 3.2): in particular one single blow-up is enough to solve it (respectively, after the first blow-up we need three blow-ups in three different points in the exceptional divisor).

Now assume that all the zeroes of $\delta_A$ and $\delta_B$ are simple: if $2L_t$ were effective we would find an effective divisor $D$ linearly equivalent to $2L$ which passes twice through all the couples $(a,b)$ as above.
We can show that, in this case,
there would exist a zero $a_0$ of $\delta_A$ such that $\pi_A^*a_0$ is not contained in $D$.
Indeed, suppose, by contradiction, that for every zero $a$ of $\delta_A$ the fibre $\pi^*_Aa$ is contained in $D$: then, by the fact that $2L.A=d_{A_0}$, we would get that all the other components of $D$ are fibres of $\pi_B\colon A\times B\to B$.
In particular it would follow from $2L.B=d_{B_0}$, that $D=\pi_A^*\delta_{A_0}+\pi_B^*D_B$ with $\deg(D_B)=d_{B_0}$.
It is then immediate that such a divisor cannot pass twice through every couple $(a,b)$ as above.
Hence there exists an $a_0$ such that $\pi_A^*a_0$ is not contained in $D$.
Recall that $D$ should pass through the different $d_{B_0}$ points $(a_0,b_i)$ twice, where $b_i$ are all the simple zeroes of $\delta_B$: hence $D.\pi_A^*a_0\geq2d_{B_0}$ which a contradiction to the existence of such a $D$.

Now assume that all the zeroes of $\delta_A$ and $\delta_B$ are double. If $2L_t$ were effective we would find an effective divisor $D$ linearly equivalent to $2L$ which passes through each point $(a,b)$, with $a$ and $b$ zeroes of $\delta_A$ and $\delta_B$, six times (twice for each of the three different tangents corresponding to the $D_4$ singularity we are solving).
In particular we would have that $D.\pi_A^*a_0\geq 4\frac{d_{B_0}}{2}=2d_{B_0}$, which is a contradiction (we are using that in this situation the number of zeroes of $\delta_A$ and $\delta_B$ are $\frac{d_{A_0}}{2}$ and $\frac{d_{B_0}}{2}$ respectively and that one of the tangents may be vertical, in which case that component could have zero intersection with the fibre).
\end{proof}

\begin{remark}
Recall that if we have two rational vector fields $\delta_1$ and $\delta_2$ we say that they are equivalent if there exists a rational function $\alpha$ such that $\delta_1=\alpha\delta_2$ (Definition \ref{def_ratvect}). 
Moreover we have shown that two equivalent rational vector fields are associated with the same 1-foliation $\F$ and, consequently, to the same inseparable double cover.
Clearly on a curve, all rational vector fields are equivalent and they are associated with the $k$-linear Frobenius morphism.

On the other hand, given equivalent rational vector fields $\delta_A$ and $\delta_A'$ on $A$ and $\delta_B$ and $\delta_B'$ on $B$, then, in general, $\delta_A+\delta_B$ is not equivalent to $\delta_A'+\delta_B'$ on $A\times B$.
\end{remark}

\chapter{Severi type inequalities and surfaces close to the Severi lines over the complex numbers}
\chaptermark{Severi type inequalities over the complex numbers}

\label{chap_severi0}
In this chapter we will discuss surfaces over the complex numbers: we will return to the case of positive characteristic in Chapter \ref{chap_severi+}.

Let $X$ be a minimal surface of general type of maximal Albanese dimension. 
We are interested in characterizing  surfaces which lie on or close to the Severi lines\index{Severi!lines}, i.e. surfaces for which the quantity
\begin{equation}
\label{int}
K_X^2-4\chi(\O_X)-4(q-2)
\end{equation}
vanishes or is "small" provided that $K_X^2<\frac{9}{2}\chi(\O_X)$. This value is strictly related to the so called Severi inequality\index{Severi!inequality} (cf. \cite{par_sev}), which states that a surface of general type of maximal Albanese dimension satisfies
\begin{equation}
K_X^2\geq 4\chi (\O_X).
\label{sev}
\end{equation}
In \cite{barparsto} there is a characterization of surfaces for which the inequality \ref{sev} is indeed an equality, namely these are surfaces whose canonical model is a double cover of its Albanese variety branched over an ample divisor with at most negligible singularities (in particular $q=2$). There are many generalizations of the Severi inequality; in particular Lu and Zuo have proved in \cite{lu} a similar inequality involving also the irregularity $q$:  a surface of general type and maximal Albanese dimension satisfies
\begin{equation}
K_X^2\geq\min\Bigl\{\frac{9}{2}\chi(\O_X),4\chi(\O_X)+4(q-2)\Bigr\}
\label{lusev}
\end{equation}
or, equivalently, if $K_X^2<\frac{9}{2}\chi(\O_X)$ then $K_X^2\geq 4\chi(\O_X)+4(q-2)$.
They also give conditions for a surface to satisfy the equality
\begin{equation}
\frac{9}{2}\chi(\O_X)>K_X^2=4\chi(\O_X)+4(q-2).
\label{uglu}
\end{equation}
The condition $K_X^2<\frac{9}{2}\chi(\O_X)$ is necessary to prove that there exists an involution $i$ for which the Albanese morphism of $X$ is composed with $i$ (cf. \cite{lu} Theorem 3.1) which is central in their argument.
There is a single step, \cite{lu} Lemma 4.4(2), where the condition $K_X^2<\frac{9}{2}\chi(\O_X)$ is really needed in their proof and it is not enough to require that  the Albanese morphism of $X$ is composed with an involution.

The first main result of this chapter, which is treated in Section \ref{sec_onsev}, is the following Theorem.

\begin{theorem}
\label{mioa}
Let $X$ be a minimal surface of general type with maximal Albanese dimension satisfying $q=q(X)\geq 3$ such that $K_X^2<\frac{9}{2}\chi(\O_X)$.
Then 
\begin{equation*}
K_X^2= 4\chi(\O_X) +4(q-2)
\end{equation*}
if and only if the canonical model of $X$ is isomorphic to a double cover of a product elliptic surface $Y=C\times E$ where $E$ is an elliptic curve and  $C$ is a curve of  genus $q-1$, whose branch divisor $R$ has at most negligible singularities and 
\[R\sim_{lin}C_1+C_2+\sum_{i=1}^{2d}E_i,\]
where $E_i$ (respectively $C_i$) is a fibre of the first projection (respectively the second projection) of $C\times E$ and $d> 7(q-2)$. Moreover, we have that $\Alb(X)\simeq\Alb(Y)$ and, in particular, if $q\geq 3$, the Albanese variety of $X$ is not simple.
\end{theorem} 

%

In Example \ref{es1}, we will give a relation between the invariants of $X$ and the number $d$ appearing in the linear class of the branch divisor $R$. Actually we will see that $K_X^2=8(q-2)+4d$ and $\chi(\O_X)=q-2+d$ and we will give a construction of a surface satisfying the equality $K_X^2=4\chi(\O_X)+4(q-2)$ for every $d>7(q-2)$ (this inequality is required to satisfy the hypothesis $K_X^2<\frac{9}{2}\chi(\O_X)$). In particular, for every $q\geq 3$, this gives an infinite set of couples $(a,b)$, for which there exists such a surface with invariants $K_X^2=a$ and $\chi(\O_X)=b$.

The second result, treated in Section \ref{sec_closesev} is the following Theorem.

\begin{theorem}
\label{miob}
Let $X$ be a minimal surface of general type with maximal Albanese dimension with $K_X^2<\frac{9}{2}\chi(\O_X)$.
\begin{enumerate}
\item If $K^2_X>4\chi(\O_X)+4(q-2)$, then $K^2_X\geq 4\chi(\O_X)+8(q-2).$
\item If $q=2$ and $K^2_X>4\chi(\O_X)$, then $K^2_X\geq 4\chi(\O_X)+2.$

\item If $q\geq 3$, equality holds, i.e.
\begin{equation}
\label{eq_8q-2}
K^2_X=4\chi(\O_X)+8(q-2),
\end{equation}
if and only if the canonical model of $X$ is isomorphic to a double cover of a smooth elliptic surface fibration $Y$ over a curve $C$ of genus $q-1$, branched over a divisor $R$ with at most negligible singularities for which $K_Y.R=8(q-2)$. In particular, we have that $\Alb(X)\simeq\Alb(Y)$ and, if $q(X)\geq 3$, we have that the Albanese variety of $X$ is not simple.
\end{enumerate}
\end{theorem}


In Section \ref{sec_ex} we give an exhaustive class of examples of surfaces which satisfy the equalities discussed above with explicit constructions.
In Section \ref{sec_seveq} we recall the argument of \cite{lu} in order to prove their Severi type inequalities.

\section{Examples}
\label{sec_ex}
In this section we give explicit examples of surfaces which satisfy the equalities $K_X^2=4\chi(\O_X)+4(q-2)$ and $K_X^2=4\chi(\O_X)+8(q-2)$, proving that all the inequalities that we are going to prove in this Chapter are sharp. First we give an example of a surface satisfying  $K_X^2=4\chi(\O_X)+4(q-2)$ for $q\geq 3$ (a characterization of surfaces satisfying this equality for $q=2$ is done in \cite{barparsto}).

\begin{example}[double cover of a product elliptic surface]
\label{es1}
We consider an elliptic surface $Y_0=C\times E$ which is the product of an elliptic curve $E$ and a curve $C$ of genus $g>1$.  With an abuse of notation, we call $E$ the class of a fibre of $\pi_C$ in $\ns(Y_0)$\nomenclature{$\ns(X)$}{the N\'eron-Severi group of $X$, i.e. the group of divisors modulo algebraic equivalence} and $C$ the class of a fibre of $\pi_E$ in $\ns(Y_0)$. 

\[
\begin{tikzcd}
	 & [-2em] C_e\arrow[mapsto]{r}{} \arrow[phantom,"\rotatebox{270}{$\subseteq$}"]{d}{} & e \arrow[phantom,"\rotatebox{270}{$\in$}"]{d}{}\\ [-3ex]
	E_c \arrow[phantom,"\subseteq"]{r}{} \arrow[mapsto]{d}{} & C\times E\arrow{r}{\pi_E} \arrow{d}{\pi_C} & E\\
	c \arrow[phantom,"\in"]{r}{} & C & 
\end{tikzcd}
\]
We know that every divisor of even degree on a curve is two-divisible in the Picard group: by this, it follows that $R_0\sim_{hom}2C+2dE$ ($d> 0$) is two-divisible in $\Pic(Y_0)$, i.e. there exists a line bundle $L_0$ such that $R_0=2L_0$, moreover $R_0.E=2$. There certainly exist elements in this homological class that are reduced and have at most negligible singularities: for example it is enough to take different fibres $C_1, C_2$ and $E_1,\ldots,  E_{2d}$, then $C_1+C_2+E_1+\ldots+E_{2d}$ has only double points (Actually, if $d\gg 0$, a general element of the homological class $R_0$ is smooth by Bertini). It is immediate that $2K_{Y_0}.L_0=K_{Y_0}.R_0=4(q-2)$, where $q=q(Y_0)=g+1$.

Thus we obtain a double cover $\pi_0\colon X_0\to Y_0$ and after the canonical resolution (cf. section \ref{sec_canres}) we get a smooth surface $X$ and the following diagram:
\[
\begin{tikzcd}
X\arrow{r}{\phi}\arrow{d}{\pi} & X_0\arrow{d}{\pi_0}\\
Y\arrow{r}{\psi} & Y_0.
\end{tikzcd}
\] 

We know that, if the singularities are at most negligible, $K_X=(\psi\circ\pi)^*(K_{Y_0}+L_0)\sim_{\textit{Num}}(\psi\circ\pi)^*(C+(2q-4+d)E)$. It is easy to see, by the Nakai-Moishezon criterion, that $C+(2q-4+d)E$ is always ample and, from this, it follows that $X$ is of general type. Furthermore $X$ is minimal because its canonical divisor is nef (it is the pull-back of an ample divisor via a generically finite morphism, and its canonical model is $X_0$ because $Y_0$ contains no rational curves).

By Equations \ref{cansq} and \ref{eulchar} we get
\[K_X^2=2(K_{Y_0}+L_0)^2=8(q-2)+4d\]
and
\[\chi(\O_X)=q-2+d.\]
Moreover, because $L_0$ is ample and $(\pi_0)_*\O_{X_0}=\O_{Y_0}\oplus L_0^{-1}$, we have that $q(X)=q$.

In particular 
\[K_X^2-4\chi(\O_X)=4(q-2).\]
Remark \ref{alb} ensures that the Albanese varieties coincide and that the Albanese morphism of $X$ is composed with an involution.

Notice also that 
\[K_X^2-\frac{9}{2}\chi(\O_X)=4(q-2)-\frac{1}{2}(q-2)-\frac{1}{2}d=\frac{1}{2}(7(q-2)-d),\]
which is smaller than zero if and only if $d>7(q-2)$.

Hence we have proved that these surfaces satisfy $K_X^2=4\chi(\O_X)+4(q-2)<\frac{9}{2}\chi(\O_X)$: the next step is to prove that they are the only ones; this will be done in Section \ref{sec_onsev}.
\end{example}

Now we give three examples of surfaces for which $4\chi(\O_X)+4(q-2)<K_X^2<\frac{9}{2}\chi(\O_X)$: in the first two cases we have $q(X)\geq 3$ and $K_X^2=4\chi(\O_X)+8(q-2)$, while in the last example $q(X)=2$ and $K_X^2=4\chi(\O_X)+2$: as we will see in Section \ref{sec_closesev} they are minimal with this property.

\begin{example}
\label{buel}
The easiest possible case is a simple modification of Example \ref{es1}. We take $Y_0$ as before: in this case we just need to take $R_0\sim_{hom}4C+2dE$ and everything is verified in a completely similar way (as before, we need $d> 7(q-2)$). In this case we have $q(X)=q, K_X^2=16(q-2)+8d$ and $\chi(\O_X)=2(q-2)+2d$. Hence $K_X^2-4\chi(\O_X)=8(q-2)$.
\label{es3}
\end{example}

\medskip

Before the next example, we recall some facts that will be useful. It is known that the Jacobian variety $\jac C'$ of a very general curve $C'$ of genus $g'$ is simple (cf. \cite{lange} Theorem 17.5.1). 
Take an \'etale double cover $C\to C'$ and consider its Prym variety $P(C,C')$ (cf. \cite{mumford-prym}): it is known (ibidem) that $\jac C$ is isogenous to $\jac C'\times P(C,C')$.
Because every Jacobian variety is a limit of a family of Prym varieties (\cite{lange} Remark 12.4.3), we can conclude that, for a very general $C'$, also $P(C,C')$ is simple. 
In particular there are no Abelian subvarieties of codimension $1$ of $\jac C$ if $g'>2$ (if $g'=2$, the dimension of $P(C,C')$ is $1$ , in particular $\jac C'$ is an Abelian subvariety of codimension $1$ of $\jac C$). 
So, for every elliptic curve $E$, the set $\Hom(C,E)$, with $C$ chosen as above, contains only constant morphisms.


\begin{example}[Double cover of a non-trivial smooth elliptic surface fibration]
\label{es4}
Here we present  an example of surface $X$ of general type satisfying equality in Theorem \ref{miob}, whose canonical model is a ramified double cover of an elliptic surface which is not a product. We start with $C'$, $C$ and $E$ as above.

Let $G$ be a subgroup of order $2$ of $E$ acting freely on $C$ such that the quotient $C/G$ is $C'$: this action clearly extends diagonally to the product giving a finite morphism of degree two $f\colon C\times E\to Y_0:=(C\times E)/G$. Proposition \ref{propo_piccxe}, together with the non-existence of surjective morphisms from $C$ to $E$, shows that there is no line bundle $L$ fixed by $G$ for which $L.E=1$. By Lemma \ref{lem2} this is enough to prove that $Y_0$ is not a product. We denote by $\pi_{C/G}$ and $\pi_{E/G}$ the two morphisms from $Y_0$ to $C/G$ and $E/G$ respectively, whose generic fibres are $E$ and $C$ respectively. We have the following commutative diagram:
\[
\begin{tikzcd}
E\arrow{r}{f_E} & E/G\\
C\times E\arrow[swap]{u}{\pi_E}\arrow{d}{\pi_C} \arrow{r}{f} & Y_0\arrow[swap]{u}{\pi_{E/G}}\arrow{d}{\pi_{C/G}}\\
C\arrow{r}{f_C} & C/G.
\end{tikzcd}
\]

Recall that, in our case, a line bundle on $C\times E$ descends to a line bundle on $Y_0$ if and only if its class in the Picard group is fixed by $G$  (cf. \cite{drezet} Theorem 2.3). Indeed, when the group $G$ is cyclic, it is always possible to give to a line bundle $L$ fixed by the action of $G$, the structure of a $G$-bundle.

Let $L_C$ and $R_C$ be two line bundles on $C/G$ of degree respectively $d$ and $2d$  such that $2L_C=R_C$. Similarly, let $L_E$ and $R_E$ be two line bundles on $E/G$ of degree respectively $1$ and $2$ such that $2L_E=R_E$. If we take an element $R\in |\pi_{C/G}^*R_C+\pi_{E/G}^*R_E|$ with at most negligible singularities (as in Example \ref{es1}, we can even assume that $R$ is smooth by Bertini if $d\gg 0$) and we denote by $L=\pi_{C/G}^*L_C\otimes \pi_{E/G}^*L_E$, we have $2L=\O_{Y_0}(R)$; hence we get a double cover $\pi_0\colon X_0\to Y_0$, and after the canonical resolution (cf. section \ref{sec_canres}) we get a smooth surface $X$ and the following diagram:
\[
\begin{tikzcd}
X\arrow{r}{\phi}\arrow{d}{\pi} & X_0\arrow{d}{\pi_0}\\
Y\arrow{r}{\psi} & Y_0,
\end{tikzcd}
\]  
where $X$ is smooth and, because $L$ is ample, $q:=q(X)=q(Y_0)$. By Equations \ref{cansq} and \ref{eulchar}
\[K_X^2=2(K_{Y_0}+L)^2=16(q-2)+8d\]
 and 
\[\chi(\O_X)=\frac{1}{2}(K_{Y_0}+L).L=2(q-2)+2d.\]
In particular we have 
\[K^2_X-4\chi(\O_X)=8(q-2).\]
It is obvious (cf. Remark \ref{alb}) that the Albanese morphism of $X$ factors through $Y$ and, as in the previous examples, we require $d>7(q-2)$ in order to have that $K_X^2<\frac{9}{2}\chi(\O_X)$: this concludes our example. 
Notice that the condition $g(C')> 2$ implies that $q(X)=g(C')+1>3$, and we have used it to assures that $\Hom_{c_0}(C,E)=0$. 
Of course, provided that one finds a curve $C'$ of genus $2$ and an elliptic curve $E$ such that $\Hom_{c_0}(C',E)$ is trivial for a point $c_0\in C'$, the same argument as above applies also in this case.
\end{example}

\begin{example}[Double cover of an Abelian variety ramified over a divisor with a quadruple point]
\label{ex_abe}
This is an example of surface of general type with maximal Albanese dimension whose Albanese morphism is composed with an involution with $q(X)=2$ satisfying equality in Theorem \ref{miob}, i.e.
\[K^2_X=4\chi(\O_X)+2.\]

Let $A$ be an Abelian surface and let $H$ be a very ample divisor. Take general elements $H_1, H_2, H_3, H_4, H_5, H_6, H_7, H_8$  inside the linear system  $|H|$ such that they are smooth, they all pass through a point $p$ and $H_i$ intersects $H_j$ for $i\neq j$ transversely at every intersection point. 
Denote by $H_{i,j}=(H_i\cap H_j)\setminus \{p\}$ for $1\leq i<j \leq 8$: we also require that $H_{i,j}\cap H_{i',j'}=\emptyset$ for every $(i,j)\neq(i',j')$. 
Let $D_1=H_1+H_2+H_3+H_4$ and $D_2=H_5+H_6+H_7+H_8$ and consider the pencil $P=\lambda D_1+\mu D_2$: the base locus of $P$ is $\{p\}\sqcup \bigsqcup_{1\leq i,j\leq 4} H_{i,j+4}$. By a Bertini-type argument, the generic element $R\in P$ is smooth away from $p$ and has a quadruple ordinary point at $p$ (actually, it is not hard to see that in the following discussion we just need that all the singularities of $R$ except $p$ are negligible so that we can even take $R=D_1$ or $R=D_2$, because for the computation of numerical invariants of the canonical resolution of a double cover branched over $R$ we can work as if the negligible singularities of $R$ were smooth points). 

It is obvious that $R$ is two-divisible, i.e. there exists $L\in\Pic(A)$ with $2L=\O_A(R)$. Consider the double cover $\pi\colon X\to A$ branched over $R$, and take the canonical resolution $\widetilde{\pi}\colon\widetilde{X}\to\widetilde{A}$ (cf. section \ref{sec_canres}) which consists of a single blow-up of the quadruple point of $R$. 
Because $R$ has a single ordinary quadruple point and $A$ contains no rational curves, $\widetilde{X}$ is the minimal smooth model of $X$ (cf. \cite{barth} III.7). 
Denote by $E$ the exceptional divisor of $\psi\colon\widetilde{A}\to A$: we have that $\widetilde{L}^2=(\psi^*L-2E)^2=L^2-4>0$  if we assume $L$ to be sufficiently ample.
Moreover, up to take a multiple of $H$, we may suppose that the Seshadri constant of $L$ is sufficiently big such that $\widetilde{L}$ intersects positively every curve on $\widetilde{A}$ (cf. \cite{lazarsfeldi} Section 5.1). 
In particular $\widetilde{L}$ is ample by the Nakai-Moishezon criterion (\cite{badescu} Theorem 1.22), from which follows $q(\widetilde{X})=q(A)=2$ and $X$ is of general type. By Equations \ref{cansq} and \ref{eulchar}, we have that
\[K_{\widetilde{X}}^2= 2L^2-2\]
and
\[\chi(\O_{\widetilde{X}})=\frac{1}{2}L^2-1.\]
In particular we have
\[K_{\widetilde{X}}^2-4\chi(\O_{\widetilde{X}})=2\]
and clearly (cf. Remark \ref{alb}) the Albanese morphism of $X$ is the natural morphism to $A$, moreover it is immediate that, for every choice of $H$, $K_X^2<\frac{9}{2}\chi(\O_X)$ holds.
\end{example}

\begin{remark}
Notice that Example \ref{ex_abe} is not exhaustive among examples of surfaces satisfying the equality in Theorem \ref{miob} with $q=2$.
Here we sketch briefly how such a classification would go.

Looking carefully at the construction in Example \ref{ex_abe}, we see that the property of the branch divisor $R$ that we have used in the computation of the invariants of $\widetilde{X}$ is the following: the branch divisor $R$ and all the reduced total transforms of $R$ in the canonical resolution process contains exactly one quadruple or quintuple point and all the other singularities are negligible (notice that the definition of this property is similar to the definition of negligible singularities of curves as in \cite{barth} Section II.8).
We consider quadruple or quintuple points because $\intinf{4}=\intinf{5}=2$. Hence, if there is a single point with this property, there is a single $i$ such that $m_i=2$ in the canonical resolution, in particular $K_{\widetilde{X}}^2=2L^2-2$ and $\chi(\O_{\widetilde{X}})=\frac{1}{2}L^2-1$ (cf.  Equations \ref{cansq} and \ref{eulchar}).
For example, a branch divisors containing a single quintuple ordinary point or a triple point with a single tangent is allowed.

Another property we need is that the surface $\widetilde{X}$ coming from the canonical resolution is minimal, otherwise its minimal model $\overline{X}$ would satisfy $K_{\overline{X}}^2-4\chi(O_{\overline{X}})>2$ because $K_{\overline{X}}^2>K_{\widetilde{X}}^2$.
For example, if we take $R$ having a single quintuple point then $\widetilde{X}$ is minimal (hence we obtain another example similar to the one in Example \ref{ex_abe}) while, if it has a triple point with a single tangent, $\widetilde{X}$ contains a $-1$-curve (hence we exclude this case).

One has the to consider case by case (again, in a similar way as in the study of negligible singularities in \cite{barth} Section II.8)  the resolution of all quadruple and quintuple points and see which of them have the properties just reported.
\label{rem_abe}
\end{remark}

\section{Severi Type inequalities}
\label{sec_seveq}
In this\index{Severi!type inequalities} section we briefly recall the main ideas of Lu and Zuo in their paper \cite{lu} that will be used in our proofs.  As stated there (Theorem 3.1) the condition
\[K_X^2<\frac{9}{2}\chi(\O_X)\] 
is necessary to prove that there exists an involution $i\colon X\to X$ with respect to which the Albanese morphism $\alb_X$ is stable (or, equivalently, $\alb_X$ is composed with $i$), i.e.
\[
\begin{tikzcd}
X \arrow{rr}{i}\arrow[swap]{dr}{\alb_X}& & X\arrow{dl}{\alb_X} \\
 & \Alb(X) &
\end{tikzcd}
\]
is commutative.

We will see a generalization of this result in every characteristic in Section \ref{sec_albfact}.

\begin{remark}
Notice that the condition "$\alb_X$ is composed with an involution" is necessary. Otherwise it is easy to construct a counter example. Take a product of curves $A$ and $B$ with $g(A)=2$ and $g(B)\geq 2$. Then the surface $X=A\times B$, which is of general type, gives the desired counterexample. Indeed, the Albanese morphism is clearly injective and, by K\"unneth formula and the formula of the canonical bundle of a product, we have that 
\[K_X^2=8(g(B)-1) \quad \chi(\O_X)=g(B)-1.\]
Hence 
\[K_X^2-4\chi(\O_X)=4(g(B)-1)=4(q(X)-3)<4(q(X)-2).\]

Notice also that, at least for Theorem \ref{miob},  the condition $K_X^2<\frac{9}{2}\chi(\O_X)$ is also necessary even if we are assuming that the Albanese morphism is composed with an involution.
Let $Y=A\times B$, where $A$ and $B$ are as above with $g(B)=g(A)=2$ and denote by $A$ (respectively $B$) a fibre of the second projection $\pi_2\colon Y\to B$ (respectively of the first projection $\pi_1\colon Y\to A$).
The invariants of $Y$ are:
\begin{itemize}
	\item $q(Y)=4$;
	\item $p_g(Y)=4$,
	\item $\chi(\O_Y)=1$;
	\item $K_Y\sim_{num}2A+2B$;
	\item $K_Y^2=8$.
\end{itemize}
Let $P$ be a point of $A$ such that $h^0(A,2P)=2$ (respectively $Q$ a point of $B$ such that $h^0(B,2Q)=2$) and let $L=A_Q+B_P$ where $A_Q$ is the fibre of the second projection $\pi_2\colon Y\to B$ over $Q$ (respectively $B_P$ is the fibre of the first projection $\pi_1\colon Y\to A$ over $P$) and take $R$ a general element of $|2L|$ (in particular $R\sim_{lin} K_Y$) which, by Bertini, may be assumed to have at most negligible singularities.
Then the desingularization $X$ of the double cover defined by $L^2=\O_Y(R)$ satisfies (cf. Remark \ref{alb} and Equations \ref{cansq} and \ref{eulchar})
\begin{itemize}
	\item $q(X)=q(Y)$ and hence they have the same Albanese variety;
	\item $\chi(\O_X)=5$;
	\item $K_X^2=36$,
\end{itemize}
from which we obtain
\[K_X^2-4\chi(\O_X)=36-20=16=8(q-2)\geq \frac{5}{2}=\frac{1}{2}\chi(\O_X).\]
Moreover the Albanese morphism of $X$ is a double cover onto its image which is not an elliptic surface, which would be a contradiction for Theorem \ref{miob} if we did not have $K_X^2\geq \frac{9}{2}\chi(\O_X)$.
\end{remark}

\medskip

The quotient surface $Y=X/i$ can be  singular, but its singular points are not so bad: they are $A_1$ singularities and they are in one-to-one correspondence with the isolated fixed points of $i$. Let $Y'$ be the resolution obtained by blowing up the singularities and let $X'$ be the blow-up of $X$ at the isolated fixed points of $i$. Denote by $Y_0$ the minimal model of $Y'$ and by $X_0$ the middle term of the Stein factorization of the morphism from $X'$ to $Y_0$. What we get is the following commutative diagram.
\[
\begin{tikzcd}
X \arrow{d}{\pi} & X' \arrow{d}{\pi'}\arrow{l}[swap]{f_X}\arrow{r}{g_X} & X_0 \arrow{d}{\pi_0}\\
Y=X/i 					 & Y' \arrow{r}{g_Y} \arrow{l}[swap]{f_Y}								& Y_0.
\end{tikzcd}
\]
 We know that the double covers $\pi'$ and $\pi_0$ are  given by equations $2L'=\O_{Y'}(R')$ and $2L_0=\O_{Y_0}(R_0)$ respectively where $R'$ and $R_0$ are the branch divisors. Notice that $R_0$ has to be reduced (because $X_0$ is normal), while $R'$ has to be smooth (because $X'$ is smooth). It follows directly from the universal property of the Albanese morphism and the fact that $\alb_X$ factors through $\pi$ that $Y_0$ is a surface of maximal Albanese dimension with $q(Y_0)=q$.

By the classification of minimal surfaces (cf. Sections \ref{sec_surfaces} and \ref{sec_ell}), we know that $Y_0$ has non-negative Kodaira dimension and maximal Albanese dimension and in particular we have the following possibilities :
\begin{itemize}
	\item if $k(Y_0)=0$, then $Y_0$ is an Abelian surface and $q=2$;
	\item if $k(Y_0)=1$, then $Y_0$ is an  elliptic surface only with smooth (but possibly  multiple) fibres over a curve $C$ with genus $g(C)=q-1$ and $q\geq 3$;
	\item if $k(Y_0)=2$, then $Y_0$ is a minimal surface of general  type of maximal Albanese dimension with $q\geq 2$.
\end{itemize}

First, we restrict to the case $k(Y_0)<2$. The surface $X_0$ may not be smooth, so we perform the canonical resolution (cf. Section \ref{sec_canres}). We get the following diagram
\begin{equation}
\label{eq_char0xtoy}
\begin{tikzcd}
X \arrow{d}{\pi} & X_t \arrow{d}{\pi_t}\arrow{l}[swap]{\phi}\arrow{r}{\phi_0} & X_0 \arrow{d}{\pi_0}\\
Y=X/i & Y_t \arrow{r}{\psi_0} \arrow{l}[swap]{\psi} & Y_0.
\end{tikzcd}
\end{equation}
 We notice that $X$ is nothing but the minimal model of $X_t$: thus, there exists an integer $n$ such that $\phi$ is the composition of $n$ blow-ups. In particular $K_X^2=K_{X_t}^2+n$ and $\chi(\O_X)=\chi(\O_{X_t})$.

If $Y_0$ is an elliptic surface over a curve $C$ with $g(C)=q-1$, then, denoting by $F$ a general elliptic fibre of the fibration, $F.R_0=2F.L_0\geq 2$. Indeed if $F.R_0=0$, then we would have that $X$ has an elliptic fibration, which is not the case because $X$ is of general type. Recall that the numerical class of the canonical bundle of such an elliptic surface is (Theorem \ref{teo_canellsur})
\[K_{Y_0}\sim_{\textit{Num}} (2g(C)-2)F+\sum_{j=1}^m(n_j-1)F_j=2(q-2)F+\sum_{j=1}^m(n_j-1)F_j,\]
where $n_jF_j$ are the the multiple fibres with $F_j$ reduced (observe that, with the notation as in Section \ref{sec_ell}, over the complex numbers $T=0$: in particular by Remarks \ref{rem_ellsurf1} and \ref{rem_spectralsfib} we have that $R^1f_*\O_{Y_0}=\L=0$ and $\chi(\O_{Y_0})=0$).
When $Y_0$ is Abelian, we know that the canonical bundle is trivial.

Summarizing we get (thanks to Equations \ref{cansq} and \ref{eulchar})
\begin{equation}
\label{eqell1}
\begin{split}
K_X^2-4\chi(\O_X)=K_{X_t}^2&-4\chi(\O_{X_t})+n=\\
=2(K_{Y_0}^2-4\chi(\O_{Y_0}))+2K_{Y_0}&.L_0+2\sum_{i=1}^{t}(m_i-1)+n=\\
=4(q-2)F.L_0+2\sum_{j=1}^m(n_j-1&)F_j.L_0+2\sum_{i=1}^{t}(m_i-1)+n\geq\\
\geq 4(q-2)+2\sum_{j=1}^m(n_j-1)F_j.L_0&+2\sum_{i=1}^{t}(m_i-1)+n\geq 4(q-2).
\end{split}
\end{equation}

So equality holds if and only if $n=0$ (i.e. $X=X_t$), $m_i=1$ for all $i$, $n_j=1$ for all $j$ (this is required only when $Y_0$ is elliptic) and $K_{Y_0}.L_0=2(q-2)$. The last condition is trivially true in the case $Y_0$ is an Abelian surface, while in the case of an elliptic surface it tells us that there are no multiple fibres  and from this it follows that, after a suitable base change, $Y_0$ is a product of an elliptic curve with a curve of higher genus (cf. Corollary \ref{coro_isoell}). The condition $m_i=1$ implies that the singularities of the branch divisor $R_0$ are at most negligible. Hence the inverse image of  an exceptional curve is a union of $-2$-curves (cf. \cite{barth} III.7 Table 1). This, together with the fact that $Y_0$ has no rational curves, implies that $X_0$ is the canonical model of $X$.

\begin{remark}
\label{uguaglianza}
We stress here what are the necessary numerical conditions on $Y_0$ (in the case it is an elliptic surface) in order to satisfy the equality of Theorem \ref{mioa}. Looking at Equation \ref{eqell1} it is immediate that this happens if and only if 
\begin{itemize}
	\item $F.L_0=1$;
	\item $n_j=1\quad \forall\ j$;
	\item $m_i=1\quad \forall\ i$;
	\item $n=0$.
\end{itemize}
Notice that the first condition together with Corollary \ref{coro_isoell} implies that $Y_0$ is a product elliptic surface.

The same conditions have to be verified in the case $Y_0$ is Abelian without the condition on the multiple fibres.
\end{remark}

Consider the case when $Y_0$ is of general type. By Theorem 1.3 of \cite{lu}, there are two possibilities. First we assume that $K_{Y_0}^2\geq 4\chi(\O_{Y_0})+4(q-2)$. As before we obtain
\begin{equation}
\label{eqgen1}
\begin{split}
K_X^2-4\chi(\O_X)&\geq 2(K_{Y_0}^2-4\chi(\O_{Y_0}))+2K_{Y_0}.L_0+2\sum_{i=1}^{t}(m_i-1)\geq\\
&\geq 8(q-2)\geq 4(q-2).
\end{split}
\end{equation}
Equality would be possible if and only if  $q=2$, $K_{Y_0}^2=4\chi(\O_{Y_0})$, $X=X_t$ is minimal, $m_i=1$ for all $i$ and $K_{Y_0}.R_0=0$.
In particular the last inequality implies that $R_0$ is a, possibly empty, union of $(-2)$-curves (these are the only curves with non-positive intersection with the canonical bundle).
Suppose that $R_0$ is not trivial: then we would have
\begin{equation}
(K_{Y_t}+\frac{1}{2}R_t).\psi_0^*R_0=\Bigl(\psi_0^*(K_{Y_0}+\frac{1}{2}R_0)+\sum_{i=1}^t(1-m_i)E_i\Bigr).\psi_0^*R_0=\frac{1}{2}R_0^2< 0
\label{eq_nonnefcan}
\end{equation}
which would contradict the nefness of the canonical bundle of the minimal surfaces $X=X_t$ (by Equation \ref{eq_candouble}, $K_{X_t}=\pi_t^*(K_{Y_t}+\frac{1}{2}R_t)$ and $\pi_t^*\psi_0^*R_0$ is clearly effective).
By this we have proved $X_0=X_t=X$, $Y_0=Y_t=Y$ and $\pi\colon X\to Y$ is an \'etale double cover such that $\alb_X=\alb_Y\circ\pi$.
In the proof of Theorem 1.3 of \cite{lu} it is proved that this is not possible over the complex numbers using a result similar to Corollary \ref{coro_doubcovamplefundgr} (cf. ibidem Lemma 4.5): we will return back to this argument when we will prove the Severi inequality in positive characteristic in Section \ref{sec_severi+}.

In a completely similar way, we prove that 
\begin{equation}
K_X^2-4\chi(\O_X)=8(q-2)
\label{eq_nuova8q-2}
\end{equation}
if and only if $X=X_t=X_0$, $Y=Y_t=Y_0$ and $\pi\colon X\to Y$ is an \'etale cover and the invariants of $Y$ satisfy $K_Y^2-4\chi(\O_Y)=4(q-2)$.

The other possible case is when $K_{Y_0}^2\geq \frac{9}{2}\chi(\O_{Y_0})$. In this case we have that 
\begin{equation}
\label{eqgen2}
K_X^2-4\chi(\O_X)>64 \max\{q-3,1\}
\end{equation}
(\cite{lu} Lemma 4.4), i.e. we have a much stronger inequality (it is in the proof of this inequality that the condition $K_X^2<\frac{9}{2}\chi(\O_X)$ is really needed).

\section{Surfaces on the Severi lines}
\label{sec_onsev}

In this section\index{Severi!lines, surfaces on}\index{surface!on the Severi lines|see{Severi}} we are going to prove Theorem \ref{mioa}, i.e. we give a characterization of surfaces satisfying $K_X^2=4\chi(\O_X)+4(q-2)<\frac{9}{2}\chi(\O_X)$.

As we have already stressed in Remark \ref{uguaglianza}, we already know that $X$ is a double cover of a surface whose minimal smooth model is a product elliptic surface.
Hence we just need to show that all the possible examples of a $2$-divisible divisor $R$ in an elliptic surface $Y=C\times E$ which intersects the elliptic fibre twice are linearly equivalent to those  in Example  \ref{es1}. The main tools of this part are the See-saw Principle (Theorem \ref{saw}) and the explicit  formula for the Picard group of $C\times E$ (Proposition \ref{propo_piccxe}). 

Let, as usual, $C$ be a curve of genus $g(C)>1$ and $E$ be an elliptic curve. We have the following Lemmas.

\begin{lemma}
Let $Y=C\times E$, then a double cover of $Y$ branched over a divisor $R$ linearly equivalent to
\[\Gamma_f+\Gamma_{f+e}+\sum_{i=1}^{2d}E_i,\]
where $f\in\Hom(C,E)$, is isomorphic to a double cover of $C\times E$ branched over a divisor $R'$ linearly equivalent to
\[C_1+C_2+\sum_{i=1}^{2d}E_i.\]
\label{es2}
\end{lemma}
\begin{proof}
In order to prove this, it is enough to see that there exists an automorphism of $C\times E$ which sends $R$ to $R'$. 
We will prove even more: actually the elliptic fibres of $\pi_C$ will be fixed by this automorphism. Indeed let $\pi_f\colon C\times E\to E$ be the morphism given by $(c,e)\mapsto e-f(c)$, we notice that $\pi_f^{-1}(e)=\Gamma_{f+e}$. 
In particular all the graphs $\Gamma_{f+e}$ are equivalent in $\ns(C\times E)$. 
Observe that if we consider the automorphism $\alpha_f$ of $C\times E$ defined by $(c,e)\mapsto (c,e-f(c))$, this clearly fixes the fibres of $\pi_C$ with respect to which $\pi_f$ is the second projection, i.e. we have the following commutative diagram

\[
\begin{tikzcd}
 & E \\
C\times E \arrow{d}{\pi_C}\arrow[ur, "\pi_f" ] \arrow{r}{\alpha_f}  & C\times E\arrow{d}{\pi_C}\arrow[swap]{u}{\pi_E}\\
C\arrow{r}{Id} &  C:
\end{tikzcd}
\]
this concludes the proof. 
\end{proof}

\begin{remark}
Lemma \ref{es2} says that the result on the branch divisor $R$ in Theorem \ref{mioa} can be given in a seemingly weaker way i.e. we can replace the fibres $C_1$ and $C_2$ by the translated graphs $\Gamma_f$ and $\Gamma_{f+e}$.
\end{remark}

\begin{lemma}
\label{lem1}
Let $Y=C\times E$, $R$ be an effective reduced divisor such that $R.E=2$ and there exists a line bundle $L$ satisfying $2L=R$. Then 
\[R\sim_{lin}\Gamma_{g+e_1}+\Gamma_{g+e_2}+\sum_{i}E_i,\]
where $g\in \Hom_{c_0}(C,E)$ and $e_i$ are elements in $E$.
\end{lemma}

\begin{proof}
Let $R$ be as in the hypothesis and let $e_1, e_2\in E$ be the two points such that $(c_0,e_1), (c_0,e_2)\in R\cap E_{c_0}$ (it could be that $e_1=e_2$). Suppose that $\beta(R)=f$ (cf. Proposition \ref{propo_piccxe}): because $R=2L$ and the map $\beta$ is a group morphism, it follows that there exists an element $g\in \Hom_{c_0}(C,E)$ such that $2g=f$. 

Consider now the divisor 
\[D=R-\Gamma_{g+e_1}-\Gamma_{g+e_2}.\]
It is then obvious that $\beta(D)=0$ and, moreover, $D$ restricted to each fibre $E_c$ is trivial. Indeed we have that
\[0=\beta(D)(c)=i_c^*(D)-i_{c_0}^*(D)=i^*_c(D).\]
Using the See-saw Principle (cf. Theorem \ref{saw}) on $D$ we see that
\[D\sim_{lin}\pi_C^*(B)\]
i.e. 
\[R\sim_{lin}\Gamma_{g+e_1}+\Gamma_{g+e_2}+\pi^*_C(B).\]

It is possible to  show that $B$ has positive degree. Indeed, applying the isomorphism $\alpha_g$ (cf. Lemma \ref{es2}), we may assume that 
\[R\sim_{lin}\pi_E^*(A)+\pi_C^*(B)\]
where $A$ is a degree two divisor of $E$. Then, by K\"unneth formula, we have
\[h^0(R)=h^0(A)h^0(B)\]
which is positive  only if the degree of $B$ is positive. Summing up, we have 
\[R\sim_{lin}\Gamma_{g+e_1}+\Gamma_{g+e_2}+\sum_iE_i.\qedhere\]
\end{proof}

\begin{proof}[Proof of Theorem \ref{mioa}]
This is a direct consequence of Lemma \ref{lem1} and Corollary \ref{coro_isoell}.
By \cite{lu} Theorem 1.3, we know that the canonical model of $X$ is isomorphic to a double cover of a smooth  elliptic surface fibration $Y$, with fibre isomorphic to $E$. Moreover this cover is branched over a divisor $R$, with at most negligible singularities, for which $R.E=2$.
Corollary \ref{coro_isotrivell} proves that $Y$ is a product (the divisor $R$ is two divisible) and Lemma \ref{lem1} gives the linear class of the branch divisor.

Notice that the condition $d>7(q-2)$ is required in order to have $K_X^2-4\chi(\O_X)<\frac{9}{2}\chi(\O_X)$ (cf. Example \ref{es1}).
\end{proof}

\section{Surfaces close to the Severi lines}
\label{sec_closesev}
In this section we prove Theorem \ref{miob}, i.e. we give a characterization of surfaces which lie close to the Severi lines\index{Severi!lines, surfaces close to}\index{surface!close to the Severi lines|see{Severi}} and we will see that the closest satisfy $K_X^2=4\chi(\O_X)+8(q-2)$ (we are always assuming $K_X^2<\frac{9}{2}\chi(\O_X)$).

\begin{proof}[Proof of Theorem \ref{miob}]
Recall, by section \ref{sec_canres}, that we have the following diagram:
\[
\begin{tikzcd}
X \arrow{d}{\pi} & X_t \arrow{d}{\pi_t}\arrow{l}[swap]{\phi_X}\arrow{r}{\phi_t} & X_{t-1}\arrow{d}{\pi_{t-1}}\arrow{r}{\phi_{t-1}} & \ldots \arrow{r}{\phi_2} & X_1\arrow{d}{\pi_1} \arrow{r}{\phi_1} & X_0 \arrow{d}{\pi_0}\\
Y=X/i & Y_t \arrow{r}{\psi_t} \arrow{l}[swap]{\psi_Y} & Y_{t-1}\arrow{r}{\psi_{t-1}} & \ldots \arrow{r}{\psi_2} & Y_1 \arrow{r}{\psi_1} & Y_0
\end{tikzcd}
\]
and recall also that $X$ and $Y_0$ have the same Albanese variety (indeed $\pi$ is a factor of the Albanese morphism $\Alb_X$ of $X$ by construction).

If $Y_0$ is of general type,  the first part of the theorem is proven by equations \ref{eqgen1} and \ref{eqgen2}. In the case $Y_0$ is an Abelian surface the first part is trivial. 

So assume that $Y_0$ is an elliptic surface over a curve $C$ with maximal Albanese dimension. By the classification of surfaces we have that $K_{Y_0}^2=0=\chi(\O_{Y_0})$. The map $\phi_X$ is nothing but a sequence of $n$ blow-ups, in particular $K_X^2=K_{X_t}^2+n$. Moreover the numerical class of the canonical bundle of $Y_0$ is 
\[K_{Y_0}\sim_{\textit{Num}} (2g(C)-2)F+\sum_{j=1}^m(n_j-1)F_j=2(q-2)F+\sum_{j=1}^m(n_j-1)F_j,\]
where $F$ is a general fibre and $n_jF_j$ are the multiple fibres with  $F_j$ reduced.

Rephrasing Equation \ref{eqell1}, we obtain
\begin{equation}
\begin{split}
K_X^2-4\chi(\O_X)=K_{X_t}^2&-4\chi(\O_{X_t})+n=\\
=2(K_{Y_0}^2-4\chi(Y_0))+K_{Y_0}&.R_0+2\sum_{i=1}^{t}(m_i-1)+n=\\
=2(q-2)F.R_0+\sum_{j=1}^m(n_j-1&)F_j.R_0+2\sum_{i=1}^{t}(m_i-1)+n.
\end{split}
\label{eqell2}
\end{equation} 

We already know that $F.R_0$ is divisible by $2$ (recall that there exists $L_0$ such that $2L_0=R_0$) and strictly positive (otherwise $X$ would  be elliptic). As we have already noticed, the conditions in Theorem \ref{mioa}  are equivalent to the following:
\begin{itemize}
	\item $F.R_0=2$;
	\item $n_j=1\quad \forall\ j$;
	\item $m_i=1\quad \forall\ i$;
	\item $n=0$.
\end{itemize}

If we want to increase slightly $K_X^2-4\chi(\O_X)$, we thus have $4$ possibilities.

First we discuss $n$. We know that if all the $m_i=1$, then all the irreducible components of the exceptional curve in the cover surface are $(-2)$-curves (cf. \cite{barth} table 1 page 109). Moreover these are the only possible rational curves on $X_t$ (ibidem). This means that in this case $n=0$. In particular, if $n>0$, then there exists an $i$ such that $m_i>1$.

Now suppose that there exists an $i$ such that $m_i>1$. By the classification of simple singularities of curves (cf. \cite{barth} II.8) we know that we have two possibilities for $R_0$. If $R_0$ has a singular point $x$ of order greater than or equal to $4$, then $(F.R_0)_x\geq 3$ (it may happen that one of the irreducible components of $R_0$ passing through $x$ is a fibre). Hence $F.R_0\geq 4$ because $R_0=2L_0$. The other possibility is that $R_0$ has a triple point $x$ which is not simple. A necessary condition for a  triple point not to be simple is to have a single tangent. If the tangent of $R_0$ in $x$ is transversal to $F$, then $(F.R_0)_x\geq 3$, conversely if it is tangent to $F$, we have $(F.R_0)_x\geq 4$ (it may happen, as before, that one of the irreducible component is $F$ itself). In both cases we have $F.R_0\geq 4$.

Suppose now that $Y_0$ has a multiple fibre $F_j$ with multiplicity $n_j$. In this case we have $F.R_0=n_jF_j.2L_0\geq 2n_j\geq 4$. 

To summarize, whatever quantity we increase, we get $F.R_0\geq 4$: that is to say that whenever 
\[K^2_X>4\chi(\O_X)+4(q-2),\]
we get
\[K^2_X\geq 4\chi(\O_X)+8(q-2)\]
and part 1 is proven.

Now we study the case when $q=2$. First we assume that $Y_0$ is an Abelian surface. In this case we have:
\begin{equation}
\begin{split}
K_X^2-4\chi(\O_X)=K_{X_t}^2&-4\chi(\O_{X_t})+n=\\
=2(K_{Y_0}^2-4\chi(Y_0))+K_{Y_0}&.R_0+2\sum_{i=1}^{t}(m_i-1)+n=\\
=2\sum_{i=1}^{t}(m_i&-1)+n.
\end{split}
\end{equation}
With the same argument as in part 1 of the proof, if $n>0$, then there exists an $i$ such that $m_i>1$. 
Then $K_X^2-4\chi(\O_X)> 0$ implies that there exists an $i$ such that $m_i>1$: in particular $K_X^2-4\chi(\O_X)\geq 2$.

Now suppose that $Y_0$ is of general type. If $K_{Y_0}\geq \frac{9}{2}\chi(\O_{Y_0})$, the proof is immediate thanks to Equation \ref{eqgen2}. In the case $K_{Y_0}^2-4\chi(\O_{Y_0})\geq 4(q-2)$ we have
\begin{equation}
\begin{split}
K_X^2-4\chi(\O_X)=2(K_{Y_0}^2-4\chi(\O_{Y_0}))+K_{Y_0}.R_0+2\sum_{i=1}^{t}(m_i-1)+n=\\
=2\Bigl(K_{Y_0}^2-4\chi(\O_{Y_0})+K_{Y_0}.L_0+\sum_{i=1}^{t}(m_i-1)\Bigr)+n:
\end{split}
\end{equation}
as above, if $n>0$, there exists an $i$ such that $m_i>1$. 
In particular $K_X^2>4\chi(\O_X)$ implies $K_X^2\geq 4\chi(\O_X)+2$ and part two of the Theorem is proven.

Now suppose that equality holds, i.e. $K_X^2=4\chi(\O_X)+8(q-2)<\frac{9}{2}\chi(\O_X)$: it is enough to prove that $Y_0$ is an elliptic surface. Indeed, in this case, it is immediate from Equation \ref{eqell2} that the conditions of the Theorem are necessary and sufficient. 

Suppose by contradiction that $Y_0$ is a surface of general type. For the numerical invariants of $Y_0$ we have two possibilities. First, if $K_{Y_0}^2\geq \frac{9}{2}\chi(\O_{Y_0})$, we know that $K_X^2-4\chi(\O_{X})\geq 64\max\{q-3,1\}$ (this is Equation \ref{eqgen2}): then $K^2_X>4\chi(\O_X)+8(q-2)$, a contradiction. The other possible case is if
\begin{equation}
\label{eq_ultimaaggiunta}
\frac{9}{2}\chi(\O_{Y_0})>K_{Y_0}^2\geq 4\chi(\O_{Y_0})+4(q-2)
\end{equation}
and, from Equations \ref{eqgen1} and \ref{eq_nuova8q-2}, we derive that there $Y_0=Y$, its invariants satisfy
\begin{equation}
\frac{9}{2}\chi(\O_{Y})>K_{Y}^2= 4\chi(\O_{Y})+4(q-2)
\label{eq_ultimaaggiuntabis}
\end{equation} and $\pi$ is an \'etale double cover which is a factor of the Albanese morphism $\Alb_X$ of $X$.
The fact that the invariants of $Y$ satisfy Equation \ref{eq_ultimaaggiuntabis}, implies that (thanks to Theorem \ref{mioa}) the canonical model of $Y$ is a double cover of a product elliptic surface $C\times E$ branched over an ample divisor $R$ and we call $\alpha\colon Y\to C\times E$ the induced generically finite morphism.
By Corollary \ref{coro_doubcovamplefundgr}, we know that $\phi$ induces an isomorphism on the algebraic fundamental groups; in particular there exists an \'etale double cover $\beta\colon \overline{Y}\to C\times E$ and a degree-two morphism $\gamma\colon X\to \overline{Y}$ branched over an ample divisor such that
\begin{equation}
\begin{tikzcd}
X\arrow{r}{\pi}\arrow{d}{\gamma}& Y\arrow{d}{\alpha}\\
\overline{Y}\arrow{r}{\beta} & C\times E
\end{tikzcd}
\label{eq_nuovaagg}
\end{equation}
is a Cartesian diagram (cf. Theorem \ref{teo_fundpushfor}).
If we prove that the Albanese morphism $\Alb_{\overline{Y}}$ of $\overline{Y}$ has degree $1$ onto its image we get a contradiction: indeed in this case, because $\pi$ is a factor of the Albanese morphism $\Alb_X$ of $X$ and the morphism $X\to\Alb_{\overline{Y}}$ has degree two onto its image, we obtain that the $\Alb_X$ has degree two onto its image and, in particular, that $\overline{Y}=Y$ which is not possible because, for example, $Y$ is of general type while $\overline{Y}$ is an elliptic surface.
We know that (cf. Section \ref{sec_double}) $\overline{Y}$ is uniquely determined by an element $L\in\Pic(C\times E)$ of order $2$. 
By Proposition \ref{propo_piccxe} it follows at once that $L\in\Pic^0(C\times E)$: indeed there are  no torsion elements  in $\Hom_{c_0}(C,E)$ and the torsion line bundles of a curve have degree zero.
From this we can derive that there exists an \'etale double cover $A\to\jac(C)\times E$ where $A$ is clearly an Abelian variety and a morphism $\overline{Y}\to A$ such that
\begin{equation}
\begin{tikzcd}
\overline{Y}\arrow{r}{\beta}\arrow{d}{}& C\times E\arrow{d}{}\\
A\arrow{r}{} & \jac(C)\times E
\end{tikzcd}
\label{eq_nonsoh}
\end{equation}
is a Cartesian square.
The morphism $\overline{Y}\to A$ is a closed inclusion because $C\times E\to\jac(C)\times E$ is, hence the Albanese morphism $\Alb_{\overline{Y}}$ of $\overline{Y}$ has degree $1$ and the proof is complete.
\end{proof}

\chapter{Factors of the Albanese morphism and inequalities for surfaces of general type: from characteristic zero to positive characteristic}
\chaptermark{Factors of the Albanese morphism}
\label{chap_from0to+}

 The main result of this Chapter, which holds for every characteristic, is Theorem \ref{teo_albfact2} which states that every surface $X$ of general type of maximal Albanese dimension satisfying $K_X^2<\frac{9}{2}\chi(\O_X)$ admits a morphism of degree two $f\colon X\to Y$ to a normal surface $Y$ which is a factor of the Albanese morphism $\alb_X$.
As we have already seen in the previous chapter, this is the starting point for the Severi type inequalities over the complex numbers and we would like to use it in the next chapter for surfaces over fields of positive characteristic.
This was first proved over the complex numbers in \cite{lu} using  a sequence of surface fibrations $\phi_d\colon \widetilde{X}_d\to \P^1$ and of morphisms of degree two $\psi_d\colon\widetilde{X}_d\to Y_d$ relative to $\phi_d$ related to the surface $X$.
Here $\widetilde{X}_d$ is a blow-up of the \'etale cover $\nu_d\colon X_d\to X$ which comes from the multiplication by $d$ on $\Alb(X)$ and one would like that $\psi_d$ descends to the desired morphism $f\colon X\to Y$.
In characteristic zero, this is done using the linear bound for the automorphism group of $X$ (cf. \cite{xiaobound}) which is known only over the complex numbers.
Such a descent has been generalized in general characteristic in \cite{gusunzhou}, provided that $\psi_d$ is a separable morphism (e.g. $\cha(k)\neq 2$) for infinitely many $d$.

It remains the case where $\alb_X$ is a separable morphism and the morphism $\psi_d$ is inseparable for all but a finite number of  $d$ (the case $\alb_X$ inseparable is trivial, as we will see, thanks to theory of $1$-foliations).
This is done thanks to Theorem \ref{teo_slopegusunzhou}, which is a generalization of Theorem 3.1 of \cite{gusunzhou}, where the generalization is that we allow $c$ to be contained in $[0,\frac{1}{2}]$ instead of $[0,\frac{1}{3}]$.
The proof of this Theorem is very technical and relies on the constructions of, among others, the morphisms $\phi_d$ and $\psi_d$ that we have mentioned above.
Once one has constructed these and other objects  that we will need for the proof of Theorem \ref{teo_slopegusunzhou}, most of the effort in order to get some inequalities such as the slope inequality \cite{xiao}, the Severi inequality \cite{par_sev} and  inequalities for non-hyperelliptic surfaces fibrations \cite{zuo} and \cite{gusunzhou} has been done: that's why, together with the aim to let the reader become more acquainted with these objects, we report the proof of these inequalities here.
We also pay special attention to how all these proofs are still valid in positive characteristic where the main problem is that the general fibre of a surface fibration can be singular.

In Section \ref{sec_slope} we recall the slope inequality with its original proof in \cite{xiao}.
In Section \ref{sec_nonhyper} we see the proof for non-hyperelliptic surface fibrations of \cite{zuo} ($\cha(k)=0$) and \cite{gusunzhou} and we generalize it so that $c$ can take value up to $\frac{1}{2}$.
In Section \ref{sec_sev} we give the proof of the Severi inequality \cite{par_sev}, which relies on the construction of the surface fibrations $\phi_d$.
In Section \ref{sec_albfact} we construct the degree-two morphisms $\psi_d$ and prove the Theorem of the factorization of the Albanese morphism.
Finally, in Section \ref{sec_gusunc12} we generalize the generalized Severi inequality of \cite{gusunzhou} using the fact that $c\in[0,\frac{1}{2}]$.

\section{Surface fibrations and the slope inequality}
\label{sec_slope}

In this section we recall the proof of the slope inequality\index{slope!inequality} for surface fibrations taken from \cite{xiao}: we also see that almost the same proof, paying attention to which conditions are really necessary in some steps, also works in positive characteristic.
The ideas for this generalization of the Xiao's proof of the slope inequality in positive characteristic are taken from \cite{sunsunzhou}  Section 2 and \cite{gu2019slope}  Section 2.1.
Notice that the result here obtained is not completely identical to their:  we will examine closely the differences in Remark \ref{rem_diffslopin}.
%
The main purpose of this section, apart from merely recalling the the result itself, is to give the constructions used in the proof which we will need in the following.
Once this has been explained, the proof does not require much time, hence we report it.

Let $f\colon X\to C$ be a morphism from a smooth surface $X$ to a smooth curve $C$ such that $f_*(\O_X)=\O_C$ and denote by $F$ a general fibre (notice that in positive characteristic $F$ needn't  be smooth but has to be integral: cf. \cite{badescu} Corollary 7.3 and the discussion following Definition 7.5) and by $g$ its arithmetic genus.
We also require $f\colon X\to C$ to be relatively minimal, i.e. there do not exist $(-1)$-curves contracted by $f$.
We are mainly interested in surfaces of general type, hence we assume that $g$ is greater than $1$.

We begin proving, and generalizing, the fundamental Lemma from which the slope inequality is derived.

\begin{lemma}[cf. \cite{xiao} Lemma 2]
\label{lemma_xiao}
Let $D$ be a divisor and suppose there exists a sequence of effective divisors
\begin{equation}
Z_1\geq Z_2\geq\ldots\geq Z_n\geq Z_{n+1}=0
\label{eq_defzi}
\end{equation}
and a sequence of numbers 
\begin{equation}
\mu_1\geq\mu_2\geq \ldots\geq \mu_n, \mu_{n+1}=0
\label{eq_defmui}
\end{equation}
such that for for every $i>1$ we have that $N_i=D-Z_i-\mu_iF$ and $N_{i-1}=D-Z_{i-1}-\mu_{i-1}F$ intersects non-negatively the effective divisor $Z_{i-1}-Z_i$ and $N_1^2\geq 0$.
Then, denoting by $d_i=N_i.F$, we have that
\begin{equation}
D^2\geq\sum_{i=1}^n(d_i+d_{i+1})(\mu_i-\mu_{i+1}).
\label{eq_xiao}
\end{equation}
\end{lemma}

\begin{proof}
We have
\begin{equation*}
\begin{split}
N_i^2=&N_i.(N_{i-1}+(Z_{i-1}-Z_i)+(\mu_{i-1}-\mu_i)F)\geq\\
\geq& N_i(N_{i-1}+(\mu_{i-1}-\mu_i)F)= \\
=& (N_{i-1}+(Z_{i-1}-Z_i)+(\mu_{i-1}-\mu_i)F).(N_{i-1}+(\mu_{i-1}-\mu_i)F)\geq\\
\geq& N_{i-1}^2+(\mu_{i-1}-\mu_i)(2N_{i-1}.F+(Z_{i-1}-Z_i).F)=\\
=&N_{i-1}^2+(\mu_{i-1}-\mu_i)N_{i-1}.F+(\mu_{i-1}-\mu_{i})(N_{i-1}+(Z_{i-1}-Z_i)).F=\\
=&N_{i-1}^2+(\mu_{i-1}-\mu_i)(d_{i-1}+d_i):
\end{split}
\label{eq_xiaoconto}
\end{equation*}
then, by induction, we have
\begin{equation*}
\begin{split}
D^2=N_{n+1}^2\geq N_n^2+(\mu_n-\mu_{n+1})(d_n+d_{n+1})\geq\ldots\geq\\
\geq N_1^2+\sum_{i=1}^n(d_i+d_{i+1})(\mu_i-\mu_{i+1})\geq\sum_{i=1}^n(d_i+d_{i+1})(\mu_i-\mu_{i+1})
\end{split}
\end{equation*}
and  the thesis follows.
\end{proof}

\begin{remark}
\label{rem_funxiao}
Lemma \ref{lemma_xiao} is fundamental in the proof of the slope inequality and we will use it for $D=K_{X/C}$ and $\mu_i$ will be the slopes of the subsequent quotients of the Harder-Narasimhan filtration of $f_*{K_{X/C}}$ so that $\mu_1>\mu_2>\ldots>\mu_n$ while, on the other hand, nothing ensures that $\mu_i\geq\mu_{n+1}=0$. 
Notice that, if $N_i$ is nef for $i\leq n$, $D$ is relatively nef (i.e. $D$ intersect non-negatively every curve contracted by the fibration) and $Z_n$ is  contracted by $f$, then the hypotheses of the Lemma are satisfied.
We will use the Lemma in this situation.
\end{remark}

Now consider the pushforward $f_*\O_X(D)$: by the well known fact that torsion-free modules on a principal ideal domain are free and the Theorems of base change (\cite{hart} Theorem III.12.11) we have that it is a vector bundle of rank $h^0(\restr{\O_X(D)}{F})$ for a general fibre $F$. 
Let $E$ be a non-zero subbundle of $f_*(\O_X(D))$ and consider the natural morphism of sheaves $\Phi\colon f^*E\to \O_X(D)$ which is defined as the composition of the pull-back of the inclusion $i\colon E\to f_*\O_X(D)$ and the counit $\iota\colon f^*f_*\O_X(D)\to\O_X(D)$ of the pushforward-pull-back adjunction: i.e. $\Phi=\iota\circ f^*i$.

\begin{lemma}
\label{lemma_gensur}
If $E\neq0$, the morphism of sheaves $\Phi$ defined above is generically surjective, i.e. there exists an open dense subset $U\subseteq X$ such that $\restr{\Phi}{U}$ is a surjective morphism of sheaves.
\end{lemma} 

\begin{proof}
Thanks to the lemma of Nakayama, it is enough to prove that for a general closed point $x\in X$ we have that the map $f^*E(x)\to \O_X(D)(x)$ is surjective (where, for every sheaf $F$, we denote by $F(x)$ the fibre of $F$ on a closed point $x$, i.e. $F(x)=F\otimes_{\O_X}k$ where $k$ is the algebraically closed field over which our varieties are defined). 

Because $\O_X(D)$ is a line bundle, its fibres are isomorphic to $k$: hence we just need to prove that the morphism of vector spaces $f^*E(x)\to \O_X(D)(x)$ is non-trivial.
Restricting to the open set $V\subseteq C$ where the base change Theorems hold (cf. \cite{hart} Theorem III.12.11) we know that $f^*f_*\O_X(D)(x)=f_*\O_X(D)(f(x))=H^0(F_{f(x)},\restr{\O_X(D)}{F_{f(x)}})$ where $F_{f(x)}$ denotes the fibre of $f$ above $f(x)$.
It follows that $f^*E(x)$ is a linear subspace of $H^0(F_{f(x)},\restr{\O_X(D)}{F_{f(x)}})$ of dimension the rank of $E$ and the map $f^*E(x)\to \O_X(D)(x)$ is given by the evaluation on $x$ of the corresponding section in $H^0(F_{f(x)},\restr{\O_X(D)}{F_{f(x)}})$.
This proves the generic surjectivity of $\Phi$.
\end{proof}

The image of $\Phi$ is isomorphic to a subsheaf of $\O_X(D)$, i.e. $\im(\Phi)=I_Z\otimes_{O_X}\O_X(-Z(E))\otimes_{\O_X}\O_X(D)$ where $I_Z$ is the ideal sheaf of a  subscheme of codimension $2$ and $Z(E)$ is an effective divisor.

\begin{definition}
\label{def_xiaozmn}
We define the fixed part of the linear system $|D|$ with respect to $E$ to be the effective divisor $Z(E)$ defined above.
We also define the moving part of $|D|$ with respect to $E$ to be the divisor $M(E)=D-Z(E)$ and denote by $N(E)=M(E)-\mu_{min}(E)F$ where $\mu_{min}$ is as defined in Definition \ref{def_HN}.
\end{definition}

Notice that $N(E)$ is a rational divisor because $\mu_{min}(E)$ may not be an integer.

\begin{remark}
Here a bit of explanation of the terminology used.
In its original article \cite{xiao}, Xiao defined $Z(E)$ in a slightly different way which makes clearer why we call it the fixed part of $|D|$ with respect to $E$.
In particular if $A$ is a sufficiently ample divisor on $C$ such that $E\otimes_{\O_C} A$ is globally generated, let $P$ be the linear subsystem of $|D+f^*A|$ corresponding to the global sections of $E\otimes_{\O_C} A$. 
Then $Z(E)$ is defined as the fixed part  of this linear system $P$ and clearly $D+f^*A-Z(E)$ is an element of its moving part.

Notice that twisting by $A$ both $E$ and $\O_X(D)$ does not change the local properties of the morphism $\Phi\colon f^*E\to \O_X(D)$: in particular we have that the image of $f^*(E\otimes_{\O_C}A)\to\O_X(D)\otimes_{\O_X}f^*A$ is $\im\Phi\otimes_{\O_X}A=I_Z\otimes_{O_X}\O_X(-Z(E))\otimes_{\O_X}\O_X(D)\otimes_{\O_X}f^*A$.
Because $f^*(E\otimes_{\O_C}A)$ is globally generated, it is clear that $Z(E)$ corresponds to the fixed part of the linear subsystem $P$, consequently the two definitions coincide.

Notice also that, restricting to a general fibre $F$, $\restr{f^*E}{F}$ defines a linear subsystem of $|\restr{D}{F}|$ of dimension equal to the rank of $E$ minus $1$, such that $\restr{M(E)}{F}$ and $\restr{Z(E)}{F}$ are respectively the moving and the fixed part.
Moreover the degree of the linear system $\restr{M(E)}{F}$ is the intersection number $M(E).F=N(E).F$.
\label{rem_xiaoexplanation}
\end{remark}

\begin{proposition}
\label{propo_xiaonef}
If all the quotients of the Harder-Narasimhan filtration of $E$ are strongly semi-stable (e.g. $C=\P^1$ or $\cha(k)=0$), then the line bundle $N(E)$ defined above is nef.
\end{proposition}

\begin{proof}
We know that, essentially by definition of $Z(E)$, $\Phi\colon f^*E\to \O_X(D-Z(E))$ is surjective in codimension $1$: by Proposition \ref{propo_morphtoproj}, this corresponds to a rational map $\phi\colon X\dashrightarrow \P(E)$ such that $\pi\circ \phi=f$, where $\pi\colon\P(E)\to C$ is the natural projection.
Moreover we have that $\O_X(D-Z(E))=\phi^*\O_{\P(E)}(1)$.
We notice that $N(E)=\phi^*(\O_{P(E)}(1)-\mu_{min}(E)F')$ where $F'$ is a general fibre of the projection $\pi\colon \P(E)\to C$.
By Theorem \ref{teo_nefantican} we easily derive the nefness of $N(E)$.
\end{proof}

From now on, let $D$ be the relative canonical divisor\index{canonical divisor!relative} of $f$, i.e. $D=K_{X/C}=K_X-f^*K_C$.
Let $E$ be the pushforward of the relative canonical bundle, i.e. $E=f_*(\O_{X}(K_{X/C}))$, and consider its Harder-Narasimhan filtration
\begin{equation}
0=E_0\subseteq E_1\subseteq\ldots\subseteq E_{n-1}\subseteq E_n=E
\label{eq_HNcan}
\end{equation}
with slopes 
\begin{equation}
\mu_{max}(E)=\mu_1>\mu_2>\ldots>\mu_{n-1}>\mu_n=\mu_{min}(E)
\label{eq_HNmu}
\end{equation}
and ranks
\begin{equation}
r_1< r_2<\ldots< r_{n-1}<r_n.
\label{eq_HNrk}
\end{equation}
We also require that all the quotients $E_i/E_{i-1}$ are strongly semi-stable so that Proposition \ref{propo_xiaonef} applies.
Denote by $Z_i=Z(E_i)$, $M_i=M(E_i)$, $N_i=N(E_i)$, for $i=1,\ldots,n$, as defined in \ref{def_xiaozmn} and let $d_i=N_i.F$, $N_{n+1}=K_{X/C}$, $d_{n+1}=K_{X/C}.F$ and $\mu_{n+1}=0$.
Thinking a little about how $Z_i$ is defined, we can easily derive that $Z_i- Z_{i+1}$ is an effective divisor for $i=1,\ldots,n-1$.

\begin{lemma}
The divisor $Z_n$ is vertical, i.e. it is contracted by $f$.
\label{lemma_vertical}
\end{lemma}

\begin{proof}
Let $c\in C$ such that the fibre $F_c$ is integral: we show that $\restr{Z_n}{F_c}=0$.
Indeed we know that for such an integral curve, its canonical divisor is globally generated (Lemma \ref{lemma_gg}) and the restriction of the relative canonical divisor $K_{X/C}$ to the fibre is its canonical divisor thanks to the adjuction formula (Equation \ref{eq_dual}).
This is enough to prove that $\restr{Z_n}{F_c}=0$ (cf. Remark \ref{rem_xiaoexplanation}).
\end{proof}

\begin{remark}
By Lemma \ref{lemma_vertical} we see that $d_n=N_n.F=(K_{X/C}-Z_n-\mu_nF).F=K_{X/C}.F=2g-2=d_{n+1}$, where, again, the last but one equality follows from the fact that the restriction of $K_{X/C}$ to a fibre is its canonical divisor.
\label{rem_dn}
\end{remark}

\begin{lemma}
We have that $d_i\geq 2r_i-2$.
\label{lemma_diri}
\end{lemma}

\begin{proof}
By Remark \ref{rem_xiaoexplanation}, we know that $d_i$ and $r_i-1$ are the degree and the dimension respectively of a linear subsystem of $|\restr{M_i}{F}|$, which is an effective subcanonical system: hence the hypotheses of  the Clifford's Theorem \ref{teo_cliff} are satisfied and we have the inequality.
\end{proof}

\begin{lemma}
\label{lemma_degcan}
We have that $\deg(f_*\omega_{X/C})=\sum_{i=1}^n r_i(\mu_i-\mu_{i+1})$.
\end{lemma}
\begin{proof}
This is just Lemma \ref{lemma_deghn}.
\end{proof}

\begin{remark}
We would like that $\mu_n\geq 0$, but this is not true in general:  Proposition \ref{propo_nef>} proves it is true provided that $f_*\omega_{X/C}$ is nef.
This is certainly true if the characteristic of the ground field is $0$ (\cite{fujita} Theorem 2.7), but nothing ensures it to be true in positive characteristic.
Another nice property would be that the relative canonical bundle of $f\colon X\to C$ is nef: this is known to be true if $X$ is relatively minimal and the general fibre is smooth (Appendix of \cite{debarre} ) or if $C=\P^1$ (in this case we have that $K_{X/C}$ is numerically equivalent to $K_X+2F$ where $F$ is a fibre of $f$, which is clearly nef for a relatively minimal surface fibration).
These properties are used in the original proof of Xiao, but we will see that they are not really necessary.
\end{remark}

\begin{remark}
Notice that we can apply Lemma \ref{lemma_xiao} for $D=K_{X/C}$ and $Z_i$, $\mu_i$ and $N_i$ for $i=1,\ldots,n+1$ defined by the Harder-Narasimhan filtration of $f_*K_{X/C}$ provided that the subsequent quotients $E_i/E_{i-1}$ are strongly semi-stable.
By Remark \ref{rem_funxiao} it is enough to prove that the $N_i$'s are nef (which is done in Proposition \ref{propo_xiaonef}), that $Z_n$ is a vertical divisor (Lemma \ref{lemma_vertical}) and that $K_{X/C}$ is a relatively nef divisor (which holds because the fibration is relatively minimal).

In the next Proposition we see how we can apply Lemma \ref{lemma_xiao} even to a subset of the indexes $\{1,\ldots,n+1\}$.
\label{rem_nuovoxiao}
\end{remark}

\begin{proposition}
\label{propo_xiaorelcan}
In the same situation as in Remark \ref{rem_nuovoxiao}, for any set of indexes $1\leq i_1\leq \ldots\leq i_{k-1}\leq i_k=n$ we have
\begin{equation}
K_{X/C}^2\geq \sum_{j=1}^k (d_{i_j}+d_{i_{j+1}})(\mu_{i_j}-\mu_{i_{j+1}}),
\label{eq_xiaopos}
\end{equation}
where $i_{k+1}=n+1$ and  $\mu_{n+1}$ is defined to be zero.
If $K_{X/C}$ is nef (e.g. the characteristic is zero or $C=\P^1$), then the same inequality is true for any set of indexes $1\leq i_1\leq \ldots\leq i_{k-1}\leq i_k\leq n$.
\end{proposition}
\begin{proof}
This is an application of Lemma \ref{lemma_xiao}: as noticed in Remark \ref{rem_nuovoxiao} we have that $N_{i_j}$ is nef for every $j=1,\ldots,k$, $Z_n=Z_{i_k}$ is vertical and $K_{X/C}$ is relatively nef, hence we can apply Lemma \ref{lemma_xiao}. 
We stress that for the general case we really need that $i_k=n$ because we need $Z_{i_k}=Z_n$ so that it is vertical and $K_{X/C}.Z_{i_k}\geq 0$ which is needed in the proof of Lemma \ref{lemma_xiao}.
If we further assume that $K_{X/C}$ is nef we have that $K_{X/C}.Z_{i_k}\geq 0$ even if $Z_{i_k}$ is not vertical, hence Lemma \ref{lemma_xiao} applies for every subset of indexes of $\{1,\ldots,n\}$.
\end{proof}

\begin{definition}
We define the slope\index{slope! of a surface fibration} of a surface fibration $f\colon X\to C$ to be $\lambda(f)=\frac{K_{X/C}^2}{\deg(f_*\omega_{X/C})}$ when $\deg(f_*\omega_{X/C})\neq 0$.
\end{definition}

\begin{theorem}[Slope inequality: \cite{xiao} Theorem 2]
\label{teo_slope}
Let $f\colon X\to C$ a relatively minimal surface fibration\index{slope!inequality} with general fibre of arithmetic genus $g\geq 2$.
Suppose that all the subsequent quotients of the Harder-Narasimhan filtration of $f_*K_{X/C}$ are strongly semi-stable.
Then we have that $K_{X/C}^2\geq\frac{4g-4}{g}\deg(f_*\omega_{S/C})$.
\end{theorem}

\begin{remark}
\label{rem_nolambda}
One could also write the thesis of the Theorem as $\lambda(f)\geq\frac{4g-4}{g}$, but, in order to do so, one has to be sure that $\deg(f_*\omega_{X/C})\neq 0$: this is known to be true if the fibration is smooth and isotrivial (i.e. all the fibres are isomorphic one another).
The original statement over the complex numbers by Xiao uses the slope and require that the fibration is not smooth isotrivial: we will not make this further assumption because it is not necessary in our formulation. 
\end{remark}

\begin{proof}
As before let 
\begin{equation}
0\subseteq E_1\subseteq E_2\subseteq \ldots\subseteq E_n=E=f_*\omega_{X/C}
\label{eq_HNslope}
\end{equation}
be the Harder-Narasimhan filtration with relative slopes $\mu_i$ and let $M_i$, $Z_i$, $N_i$, $r_i$ and $d_i$ as defined above.
Observe that $r_{i+1}\geq r_i+1$ for $i=1,\dots,n-1$ by the properties of the Harder-Narasimhan filtration but $r_{n+1}=r_n$ by definition.

Suppose that 
\begin{equation}
\mu_1+\mu_n>\frac{2}{g}\deg(E):
\label{eq_suppose1}
\end{equation}
applying Proposition \ref{propo_xiaorelcan} to the set of indexes $1,n$ we get:
\begin{equation}
K_{X/C}^2\geq(d_1+d_n)(\mu_1-\mu_n)+(d_n+d_{n+1})(\mu_n-\mu_{n+1}).
\label{eq_suppose2}
\end{equation}
By the properties of the Harder-Narasimhan filtration we know that $\mu_1-\mu_n>0$, by definition $\mu_{n+1}=0$, by Remark \ref{rem_dn} $d_n=d_{n+1}=2g-2$: putting all of these together, along with the assumption \ref{eq_suppose1},
we obtain
\begin{equation}
\begin{split}
K_{X/C}^2\geq d_1(\mu_1-\mu_n&)+d_n\mu_1+d_{n+1}\mu_n\geq\\
\geq(2g-2)(\mu_1+&\mu_n)>\frac{4g-4}{g}\deg(E).
\end{split}
\label{eq_suppose3}
\end{equation}

Now suppose that 
\begin{equation}
\mu_1+\mu_n\leq\frac{2}{g}\deg(E):
\label{eq_suppose1'}
\end{equation}
applying Proposition \ref{propo_xiaorelcan} to the set of indexes $1,2,\ldots,n-1,n$ and   Lemmas \ref{lemma_diri} and \ref{lemma_degcan} we get:
\begin{equation}
\begin{split}
K_{X/C}^2\geq\sum_{i=1}^{n}(d_i+&d_{i+1})(\mu_i-\mu_{i+1})\geq\\
\geq \sum_{i=1}^{n-1}(4r_i-2)(&\mu_i-\mu_{i+1})+4(g-1)\mu_n=\\
=4\sum_{i=1}^nr_i(\mu_i-\mu_{i+1})-2(\mu_1+&\mu_n)=4\deg(E)-2(\mu_1+\mu_n)\geq \\
\geq 4\deg(E)-\frac{4}{g}&\deg(E)=\frac{4g-4}{g}\deg(E) 
\end{split}
\label{eq_suppose2'}
\end{equation}
and the Theorem is proven.
\end{proof}

\begin{remark}
\label{rem_diffslopin}
As we have already noticed at the beginning of this Section, even though the ideas used here in order to prove the slope inequality are taken from \cite{sunsunzhou} and \cite{gu2019slope}, we have presented a slightly different result.

In particular we have shown that if we assume that all the subsequent quotients of the Harder-Narasimhan filtration of $f_*(K_{X/C})$ for a relatively minimal surface fibration $f\colon X\to C$ are strongly semi-stable, then Xiao's proof of slope inequality still holds over any algebraically closed field.
On the other hand in \cite{sunsunzhou} and \cite{gu2019slope} a result of Langer (which ensures that there exists an integer $k$ such that the subsequent quotients of the Harder-Narasimhan filtration of $(F^{k})^*E$ are strongly semi-stable, where $F^{k}$ is the $k$-th iterated absolute Frobenius morphism and $E$ is a vector bundle on a smooth curve, cf. \cite{langer} Theorem 2.7) is used in order to avoid the hypothesis of strongly semi-stability provided that the generic fibre of $f$ is smooth or hyperelliptic or the genus of $C$ is smaller than or equal to $1$ (cf. \cite{sunsunzhou} Corollary 1 and \cite{gu2019slope} Theorem 2.7).

The problem that arises when all the fibres are singular (so that we really need to require the strongly semi-stability hypothesis) is that the generic fibre of the pull-back fibration $f'\colon X'\to C$ via the $k$-th iterated absolute Frobenius morphism may not be normal.
More precisely, consider the following Cartesian diagram
\begin{equation}
\begin{tikzcd}
X'\arrow{d}{f'}\arrow{r}{}&X\arrow{d}{f}\\
C\arrow{r}{F^{k}} & C
\end{tikzcd}
\label{eq_cartdignonnorm}
\end{equation}
where $F$ is the absolute Frobenius morphism of $C$ and $k$ is as above:
we know that the generic fibre $X_\eta$ of $f$ is regular (cf. \cite{liedtke+} Section 5.1) even if it is not smooth; moreover over a reduced scheme of dimension $1$, regularity and normality coincide (cf. for example \cite{egaiv2} Remarques 5.7.8, Proposition 5.8.5 and Th\'eor\`eme 5.8.6), so that $X_\eta$ is normal.
On the other hand, if we assume that $X_\eta$ is not smooth i.e. all the fibres of $f$ are singular, it can happen that the generic fibre $X'_\eta$ is not regular, hence not normal.
Observe that this fact is an obstruction to extend the proof of slope inequality given in \cite{sunsunzhou} and \cite{gu2019slope} unless the generic fibre $X_\eta$ is smooth or hyperelliptic (cf. \cite{gu2019slope} Theorem 2.7 points (a) and (b)) or $k=1$ (partially ibidem Theorem 2.7(c) and more generally the case here treated) and that in general the slope inequality can be false (\cite{gu2019slope} Section 4.1).

Notice that in general the condition that all the subsequent quotients of the Harder-Narasimhan filtration are strongly semi-stable  is quite technical and generally difficult to check while the conditions of Theorem 2.7 of \cite{gu2019slope} (that the generic fibre $X_\eta$ is smooth or hyperelliptic or that $g(C)\leq1$) are more natural and of geometric type.
We decided to use this condition because in the case we are really interested in (i.e. when $C=\P^1$, cf. the proof of Theorem \ref{teo_sev}) the fact that all the subsequent quotients of the Harder-Narasimhan filtration of $f_*{K_{X/C}}$ are strongly semi-stable is easily verified (cf. Remark \ref{rem_finitestable}) and the proof is easier. 
\end{remark}

\section{Non-hyperelliptic surface fibrations}
\label{sec_nonhyper}
The idea of using the second multiplication map $f_*\omega_{X/C}\otimes_{\O_C}f_*\omega_{X/C}\to f_*(\omega_{X/C}^2)$ in the study of the slope of a surface fibration is due to Lu and Zuo \cite{zuo}.
In \cite{gusunzhou} there is an extension of these results to the case of positive characteristic.
Theorem \ref{teo_slopegusunzhou} is a generalization of Theorem 3.1 of \cite{gusunzhou}, which allows us to generalize Theorem 3.1 of \cite{lu} in positive characteristic (Theorem \ref{teo_albfact2}) which gives us conditions for the Albanese morphism of a surface of general type to be composed with a morphism of degree two.

Throughout this section we consider a surface fibration $f\colon X\to C$ such that the relative canonical bundle is nef:  this is always true if the fibration is relatively minimal (i.e. there are no $(-1)$-curves contracted by $f$) and the general fibre is smooth (cf. Theorem d'Arakelov in the appendix of \cite{debarre}) or if $C=\P^1$.
We also require that the general fibre of $f$ is not hyperelliptic.
Another hypothesis that we need, is that the successive quotients of the Harder-Narasimhan filtration of $f_*(\omega_{X/C})$ are strongly semi-stable, which is induced, for example, by $g(C)\leq 1$.
As already said, the results of this section will be used in the proof of Theorem \ref{teo_albfact2}. 
The fibrations appearing in the proof of this Theorem are all over the projective line, hence they satisfy the hypotheses we are requiring in this section.

Consider a non-hyperelliptic surface fibration $f\colon X\to C$ with nef relative canonical bundle and let $g$ be the arithmetic genus of a generic fibre.
As in Section \ref{sec_slope}, let $E=f_*\omega_{X/C}$,
\begin{equation}
0=E_0\subseteq E_1\subseteq\ldots\subseteq E_{n-1}\subseteq E_n=E
\label{eq_HNcan2}
\end{equation}
be the Harder-Narasimhan filtration of $E$, suppose that all the quotients $E_i/E_{i-1}$ are strongly semi-stable and define $\mu_i$, $r_i$, $d_i$, $Z_i$, $M_i$ and $N_i$ as in section \ref{sec_slope}.

We have seen that (cf. proof of Proposition \ref{propo_xiaonef} and Proposition \ref{propo_morphtoproj})  $f^*E_i\to \omega_{X/C}$ is generically surjective and canonically defines a rational map $\phi_i\colon X\dashrightarrow\P(E_i)$.
It is clear, by the fact that restriction $\restr{\phi_i}{F}\colon F\to \P^{r_i-1}$ is non degenerate, that, if $r_i>1$ the image of $\phi_i$ is a surface: denote by $X_i$ a minimal smooth model of this surface.
Moreover let $\delta_i$ be the degree of $\phi_i\colon X\dashrightarrow X_i$ and let $g_i$ be the arithmetic genus of the generic fibre of the surface fibration $X_i\to C$ induced by the projection $\pi_i\colon \P(E_i)\to C$.
It is immediate by construction that $\phi_i$ factors through $\phi_{i+1}$, in particular $\delta_{i+1}|\delta_i$.
If $r_1=1$ (notice that $r_i>r_{i-1}\geq 1$, so that only $r_1$ can be equal to one), then it is clear that $d_1=0$ (indeed in this case we have $\P(E_1)=C$ and $N_1$ is the pull-back of a line bundle on $C$, cf. proof of Proposition \ref{propo_xiaonef}) and we define $\delta_1=\infty$ and $g_1=\infty$.

Now consider the second multiplication map
\begin{equation}
\sigma\colon f_*{\omega_{X/C}}\otimes_{\O_C} f_*{\omega_{X/C}}\to f_*{\omega_{X/C}^2}
\label{eq_secondmult}
\end{equation}
and denote by $\F$ its image: by Max Noether's Theorem \ref{teo_max} it is generically surjective.
In particular we have $\deg(\F)\leq \deg(f_*\omega_{X/C}^2)$.

\begin{lemma}
\label{lemma_degs2f}
Let $f\colon X\to C$ a non-hyperelliptic relatively minimal surface fibration.
We have 
\begin{equation}
K_{X/C}^2\geq \deg(\F)-\deg(f_*\omega_{X/C}).
\end{equation}
\end{lemma}

\begin{proof}
First of all we notice that $R^1f_*\omega_{X/C}^2=0$ where, for a coherent sheaf $E$, $R^1f_*E$ denotes the first higher pushforward of $E$: by the base change Theorem (cf. \cite{hart} Corollary III.12.9) we know that it is enough  to prove that $h^1(F,\restr{\omega_{X/C}^2}{F})=0$ for every fibre $F$.
By Riemann-Roch Theorem \ref{teo_RR} and the adjunction Formula \ref{eq_dual}, this is equivalent to the vanishing of $h^0(F,-K_F)$.
Suppose that $h\in H^0(F,-K_F)$ is a section and let $A$ be the maximal closed subscheme of $F$ such that $\restr{h}{A}\equiv 0$, i.e. $h\in H^0(B,\O_B(-K_F-A))$ with $F=A+B$ (cf. Equation \ref{eq_decseq}).
Because the fibration $f$ is relatively minimal, $B$ can contain only $(-2)$-curves: indeed on all the other components the degree of the restriction of $-K_F$ is strictly negative and the restriction of $h$ has to be identically zero.
Moreover $-K_F.A<0$ because otherwise $f$ would be a genus-one fibration, which is excluded by hypothesis, and $A$ and $B$ can not be linearly equivalent.

In particular $h$ defines the following exact sequence
\begin{equation}
0\to\O_B\to\O_B(-K_F-A)\to Q\to 0
\label{eq_exactrelativelyminimal}
\end{equation}
where $Q$ is a finitely supported sheaf, which means $\chi(Q)=h^0(B,Q)\geq0$.
By Riemann-Roch \ref{teo_RR} we have
\begin{equation}
-A.B=\deg(\O_B(-K_F-A))=\chi(\O_B(-K_F-A))-\chi(\O_B)=h^0(B,Q)\geq 0
\label{eq_rrabf}
\end{equation}
which, unless $B=0$, is in contradiction with Zariski's Lemma for surface fibrations (\cite{badescu} Corollary 2.6).
In particular $h$ has to be identically zero.

By the Leray spectral sequence (cf. \cite{badescu} page 103), relative duality (ibidem) and Riemann-Roch Theorem \ref{teo_RR}, we have that 
\begin{equation}
\begin{split}
\chi(\O_X)=\chi(f_*\O_X)&-\chi(R^1f_*\O_X)=\\
=-\deg(R^1f_*\O_X)-(\rk(R^1&f_*\O_X)-1)\chi(\O_C)=\\
=-\deg(R^1f_*\O_X/T(R^1f_*\O_X))&-l-(g-1)\chi(\O_C)=\\
=\deg(f_*\omega_{X/C})-(g-1)&\chi(\O_C)-l
\end{split}
\label{eq_pushforstructure}
\end{equation}
(where $l$ denotes the dimension of the torsion subsheaf $T(R^1f_*\O_X)$ of $R^1f_*\O_X$) and

\begin{equation}
\chi(\omega_{X/C}^2)=\chi(f_*\omega_{X/C}^2)=\deg(f_*\omega_{X/C}^2)+3(g-1)\chi(\O_C).
\label{eq_pushforsym}
\end{equation}

Using Riemann-Roch Theorem \ref{teo_rrss} on $X$ we have
\begin{equation}
\begin{split}
\chi(\omega_{X/C}^2)=\chi(\O_X)&+K_{X/C}.(2K_{X/C}-K_X)=\\
=\chi(\O_X)+K_{X/C}^2-K_{X/C}.f^*K_C&=\chi(\O_X)+K_{X/C}^2+4(g-1)\chi(\O_C).\hspace{-10pt}
\end{split}
\label{eq_pushforrr}
\end{equation}

Combining Equations \ref{eq_pushforstructure}, \ref{eq_pushforsym} and \ref{eq_pushforrr} we get
\begin{equation}
\begin{split}
K_{X/C}^2=\chi(\omega_{X/C}^2)-\chi(\O_X)-4(g-1)\chi(\O_C)=\\
=\deg(f_*\omega_{X/C}^2)-\deg(f_*\omega_{X/C})+l\geq \deg(\F)-\deg(f_*\omega_{X/C})
\end{split}
\label{eq_pushforfin}
\end{equation}
and the Lemma is proven.
\end{proof}

Let's consider the following filtration for $\F$:
\begin{equation}
0=\F_0\subseteq \F_1\subseteq \F_2\subseteq\ldots\subseteq\F_{n-1}\subseteq\F_n=\F
\label{eq_filtrf}
\end{equation}
where $\F_i$ is the image of $\sigma$ restricted to $E_i\otimes E_i$.

Observe that, by Corollary \ref{coro_luzuo2} and Remarks \ref{rem_luzuo2} and \ref{rem_lineargg}, we have
\begin{equation}
\rk(\F_i)\geq\min\{3r_i-3,2r_i+g_i-1\}\geq 2r_i-1
\label{eq_rkfi2}
\end{equation} 
and, if $\phi_i$ is birational, $\rk(\F_i)\geq 3r_i-3$ (recall that $\restr{\phi_i}{F}$ is induced by a subcanonical system on a generic fibre $F$). 
Notice that if $r_1=1$, it is not true that $\min\{3r_1-3,2r_1+g_1-1\}\geq 2r_1-1$, but it remains true that $\rk(\F_1)\geq 2r_1-1=1$.
Denote by $I_1=\{i\ |\ \delta_i=1\}$, $I_2=\{i\ |\ \delta_i=2,\ r_i\geq g_i+2\}$, $I_3=\{ i |\ \delta_i=2,\ r_i< g_i+2\}$  and $I_4=\{i\ |\ \delta_i\geq 3\}$. 
If $r_1=1$ we have defined $\delta_1=\infty$: in this case $1$ is an element of $I_4$.
Notice also that $n\in I_1$: indeed, by construction, we have that $\restr{\phi_n}{F}$ is the canonical morphism of $F$ which is an embedding (cf. Lemma \ref{lemma_canbirnonhyp}), in particular $\phi_n$  is birational.

We rephrase Equation \ref{eq_rkfi2} as
\begin{equation}
\rk(\F_i)\geq
\begin{cases}
2r_i-1 & i\in I_4;\\
3r_i-3 &  i\in I_3;\\
2r_i+g_i-1 &  i\in I_2;\\
3r_i-3 & i\in I_1.
\end{cases}
\label{eq_rkfi}
\end{equation}
where for $i\in I_4$ we are using the weaker inequality $\rk(\F_i)\geq 2r_i-1$, for $i \in I_3\cup I_2$ we are using the inequalities which come from the definition of $I_3$ and $I_2$ and for $i\in I_1$ we are using the stronger inequality $\rk(\F_i)\geq 3r_i-3$ which holds because for these indexes $\phi_i$ is birational and subcanonical on $F$.
Notice that we have $\rk(\F_n)=\rk(\F)=3r_n-3=3g-3$: this is clear because of the generic surjectivity of the morphism in \ref{eq_secondmult} and is needed for the application of Corollary \ref{coro_belloneq} in order to estimate the degree of $\F$.
Denote  by
\begin{equation}
\widetilde{r}_i=
\begin{cases}
2r_i-1 & i\in I_4\\
3r_i-3 & i\in I_3;\\
2r_i+g_i-1 & i\in I_2;\\
3r_i-3 & i\in I_1.
\end{cases}
\label{eq_ritilde}
\end{equation}

Using Corollary \ref{coro_luzuo} and Remarks \ref{rem_luzuo} and \ref{rem_lineargg}, we see that 
\begin{equation}
\label{eq_casesdegrk}
d_i\geq 
\begin{cases}
\delta_i(r_i-1)\geq3(r_i-1) & i\in I_4;\\
2(2r_i-2) & i\in I_3;\\
2(r_i+g_i-1) & i\in I_2;\\
2(r_i-1) & i\in I_1.
\end{cases}
\end{equation}
where for $i\in I_4$ we are using the weaker inequality $d_i\geq \delta_i(r_i-1)$, for $i\in I_3\cup I_2$ we are using the inequalities which come from the definition of $I_3$ and $I_2$ and for $i\in I_1$ we are using the stronger inequality $d_i\geq 2r_i-2$ which holds because for these indexes $\phi_i$ is birational and $\restr{M_i}{F}$ is subcanonical.
Observe that this inequality is still valid when $r_1=1$ and hence $d_1=0$.

Combining Corollary \ref{coro_tensorss} with Remark \ref{rem_minmaxss}, we know that 
\begin{equation}
\label{eq_degf2}
\mu_{min}({\F_i})\geq\mu_{min}(E_i\otimes E_i)= 2\mu_i.
\end{equation}

Now, thanks to Corollary \ref{coro_belloneq}, where, in the notation of that Corollary, $\widetilde{\mu}_i=2\mu_i$ and $\widetilde{r}_i$ as defined in Equation \ref{eq_ritilde},  we get
\begin{equation}
\begin{split}
\deg(\F)\geq \sum_{i\in I_4}(4r_i-2)(\mu_i-\mu_{i+1})+\sum_{i\in I_3}(6r_i-6)(\mu_i-\mu_{i+1})+\\
+\sum_{i\in I_2}(4r_i+2g_i-2)(\mu_i-\mu_{i+1})+\sum_{i\in I_1}(6r_i-6)(\mu_i-\mu_{i+1}).
\end{split}
\label{eq_nuovanew10}
\end{equation}

\begin{theorem}[cf. Theorem 1.2 \cite{zuo} for 1. over the complex numbers and Theorem 3.1 \cite{gusunzhou} for 2. with $0\leq c\leq\frac{1}{3}$]
\label{teo_slopegusunzhou}
Let $f\colon X\to C$ be a non-hyperelliptic  surface fibration such that $K_{X/C}$ is nef and all the quotients appearing in the Harder-Narasimhan filtration of $f_*\omega_{X/C}$ are strongly semi-stable.
\begin{enumerate}
	\item Assume that there exists an integer $\delta>1$ such that $\delta_i=1$ or $\delta_i>\delta$ for all $i$: then
	\begin{equation}
	K_{X/C}^2\geq \Bigl(5-\frac{1}{\delta}\Bigr)\frac{g-1}{g+2}\deg(f_*\omega_{X/C}).
	\label{eq_gusunzhou1}
	\end{equation}
	\item Assume that there exists a constant $0\leq c\leq\frac{1}{2}$ such that $g_i\geq cg$ whenever $\delta_i=2$: then
	\begin{equation}
	K_{X/C}^2\geq (4+c)\frac{g-1}{g+2}\deg(f_*\omega_{X/C}).
	\label{eq_gusunzhou2}
	\end{equation}
\end{enumerate}
\end{theorem}

\begin{proof}
First of all we notice that if $\deg(f_*\omega_{X/C})$ is non-positive, the thesis is trivial in both cases.
Hence we may assume that $\deg(f_*\omega_{X/C})>0$ from which, essentially by definition, we derive that $\mu_1>0$.

By Lemmas \ref{lemma_degs2f} and  \ref{lemma_degcan} we have
\begin{equation}
K_{X/C}^2\geq \deg(\F)-\deg(f_*\omega_{X/C})=\deg(\F)-\sum_{i=1}^nr_i(\mu_i-\mu_{i+1}).
\label{eq_nuovanew1}
\end{equation}


Combining Equations \ref{eq_nuovanew1} and \ref{eq_nuovanew10}, we get
\begin{equation}
\begin{split}
K_{X/C}^2\geq \sum_{i\in I_4}(3r_i-2)(\mu_i-\mu_{i+1})+\sum_{i\in I_3}(5r_i-6)(\mu_i-\mu_{i+1})+\\
+\sum_{i\in I_2}(3r_i+2g_i-2)(\mu_i-\mu_{i+1})+\sum_{i\in I_1}(5r_i-6)(\mu_i-\mu_{i+1}).
\end{split}
\label{eq_slopezhou1}
\end{equation}

In case 1, we have that $I_3$ and $I_2$ are empty.
By Equation \ref{eq_casesdegrk}, for $i\in I_4$, we have
\begin{equation}
d_i+d_{i+1}\geq 2d_i\geq 2\delta_i(r_i-1)\geq (2\delta+2)(r_i-1),
\end{equation}
 for $i\in I_1$ we have
\begin{equation}
d_i+d_{i+1}\geq 2d_i\geq 4r_i-4.
\end{equation}

Combining these Equations with Proposition \ref{propo_xiaorelcan}, with set of indexes $1,\ldots, n$, we have
\begin{equation}
K_{X/C}^2\geq\sum_{i\in I_4}(2\delta+2)(r_i-1)(\mu_i-\mu_{i+1})+\sum_{i\in I_1}(4r_i-4)(\mu_i-\mu_{i+1}).
\label{eq_slopezhou8}
\end{equation}

In particular $\frac{\delta-1}{\delta}\ref{eq_slopezhou1}+\frac{1}{\delta}\ref{eq_slopezhou8}$ says
\begin{equation}
\begin{split}
K_{X/C}^2\geq\Bigl(5-\frac{1}{\delta}\Bigr)&\sum_{i=1}^{n}r_i(\mu_i-\mu_{i+1})+\\
-4\sum_{i\in I_4}(\mu_i-\mu_{i+1})+&\frac{2-6\delta}{\delta}\sum_{i\in I_1}(\mu_i-\mu_{i+1})\geq\\
\geq \Bigl(5-\frac{1}{\delta}\Bigr)\deg(&f_*\omega_{X/C})-6\mu_1
\end{split}
\label{eq_slopezhou3}
\end{equation}
where the last inequality follows because $-4\geq\frac{2-6\delta}{\delta}\geq-6$, $\mu_i-\mu_{i+1}\geq 0$ for all $i\in I_4$ and $\mu_1>0$ by the assumption we made at the beginning of this proof.

Equation \ref{eq_slopezhou3} becomes
\begin{equation}
K_{X/C}^2\geq\Bigl(5-\frac{1}{\delta}\Bigr)\deg(f_*\omega_{X/C})-6\mu_1.
\label{eq_slopezhou11}
\end{equation}
Combining this equation with $K_{X/C}^2\geq(2g-2)\mu_1$ (for this, we are using Proposition \ref{propo_xiaorelcan} with set of indexes $\{1\}$: here it is needed that $K_{X/C}$ is nef) we obtain
\begin{equation}
K_{X/C}^2\geq\Bigl(5-\frac{1}{\delta}\Bigr)\frac{g-1}{g+2}\deg(f_*\omega_{X/C}):
\label{eq_slopezhou10}
\end{equation}
this finishes the first part of the proof.

Now consider the second part of the Theorem.
By  Equation \ref{eq_casesdegrk} we get
\begin{equation}
\label{eq_boh1}
d_i+d_{i+1}\geq 2d_i\geq 6(r_i-1)
\end{equation}
for $i\in I_4$ (this also works for $r_1=1$).
Similarly we obtain 
\begin{equation}
d_i+d_{i+1}\geq 2(2r_i-2)+2(r_{i+1}-1)\geq6r_i-4\geq 4r_i-4
\end{equation}
for $i\in I_3$,
\begin{equation}
d_i+d_{i+1}\geq 2(r_i+g_i-1)+2(r_{i+1}-1)\geq 4r_i+2g_i-2
\end{equation}
for $i\in I_2$ and
\begin{equation}
d_i+d_{i+1}\geq 2(r_i-1)+2(r_{i+1}-1)\geq 4r_i-2
\end{equation}
for $i\in I_1\setminus\{n\}$.
Moreover $d_n+d_{n+1}=4g-4=4r_n-4$.
Notice that for $d_{i+1}$ we are always using the weaker inequality $d_{i+1}\geq 2r_{i+1}-2\geq 2r_i$ of \ref{eq_casesdegrk} which holds for every index $i$ except for $n$ which we are excluding from $I_1$ because it is the unique index $i$ for which $r_{i+1}>r_i$ does not hold.

Combining these equations with Proposition \ref{propo_xiaorelcan} with set of indexes $1,\ldots,n$, we get
\begin{equation}
\begin{split}
K_{X/C}^2\geq \sum_{i\in I_4}(6r_i-6)(\mu_i-\mu_{i+1})+\sum_{i\in I_3\cup\{n\}}(4r_i-4)(\mu_i-\mu_{i+1})\\
+\sum_{i\in I_2}(4r_i+2g_i-2)(\mu_i-\mu_{i+1})+\sum_{i\in I_1\setminus\{n\}}(4r_i-2)(\mu_i-\mu_{i+1}).
\end{split}
\label{eq_slopezhou2}
\end{equation} 

Calculating $c(\ref{eq_slopezhou1})+(1-c)(\ref{eq_slopezhou2})$ we have (recall that $g=r_n\geq r_i$ and $g_i\geq cg$)
\begin{equation}
\label{eq_slopezhou4}
\begin{split}
K_{X/C}^2\geq (4+c)\sum_{i=1}^{n}r_i(\mu_i-\mu_{i+1})+(2-4c)\sum_{i\in I_4}r_i(\mu_i-\mu_{i+1})\\
-(6-4c)\sum_{i\in I_4}(\mu_i-\mu_{i+1})-(4+2c)\sum_{i\in I_3}(\mu_i-\mu_{i+1})\\
-2\sum_{i\in I_2}(\mu_i-\mu_{i+1})-(2+4c)\sum_{i\in I_1\setminus\{n\}}(\mu_i-\mu_{i+1})-(4+2c)\mu_n.
\end{split}
\end{equation}

We notice that, for $1/3\leq c\leq 1/2$, we have $4+2c\geq 6-4c$, $4+2c\geq 2$, $4+2c\geq 2+4c$ and $2-4c\geq 0$: in particular
\begin{equation}
\begin{split}
K_{X/C}^2\geq (4+c)\deg(f_*\omega_{X/C})&-(4+2c)(\mu_1-\mu_n)-(4+2c)\mu_n\geq\\
\geq(4+c)\deg(f_*&\omega_{X/C})-(4+2c)\mu_1
\end{split}
\end{equation}
(we are again using that $\mu_i-\mu_{i+1}\geq 0$ for $i\neq n$).

Combining this equation with $K_{X/C}^2\geq(2g-2)\mu_1$ (here we are using Proposition \ref{propo_xiaorelcan} with set of indexes $\{1\}$: here it is needed that $K_{X/C}$ is nef) we obtain
\begin{equation}
K_{X/C}^2\geq(4+c)\frac{g-1}{g+1+c}\deg(f_*\omega_{X/C})\geq(4+c)\frac{g-1}{g+2}\deg(f_*\omega_{X/C})
\end{equation}
where we are using the assumption $\deg(f_*\omega_{X/C})>0$ made at the beginning of this proof.
Hence we have proven 2 for a constant $1/3\leq c\leq 1/2$.

One can substitute Equation \ref{eq_boh1} by
\begin{equation}
d_i+d_{i+1}\geq 5r_i-4
\end{equation}
for $i\in I_4$: indeed, for $i\in I_4$, if $r_i>1$ we clearly have $6r_i-6\geq 5r_i-4$, while if $r_1=1$ we have that $d_1+d_2=d_2\geq 1=5r_1-4$.
In this case Proposition   \ref{propo_xiaorelcan}, with set of indexes $1,\ldots,n$, gives
\begin{equation}
\begin{split}
K_{X/C}^2\geq \sum_{i\in I_4}(5r_i-4)(\mu_i-\mu_{i+1})+\sum_{i\in I_3\cup\{n\}}(4r_i-4)(\mu_i-\mu_{i+1})\\
+\sum_{i\in I_2}(4r_i+2g_i-2)(\mu_i-\mu_{i+1})+\sum_{i\in I_1\setminus\{n\}}(4r_i-2)(\mu_i-\mu_{i+1})
\end{split}
\label{eq_slopezhou5}
\end{equation} 
and  $c(\ref{eq_slopezhou1})+(1-c)(\ref{eq_slopezhou5})$ becomes
\begin{equation}
\label{eq_slopezhou6}
\begin{split}
K_{X/C}^2\geq (4+c)\sum_{i=1}^{n}r_i(\mu_i-\mu_{i+1})+(1-3c)\sum_{i\in I_4}r_i(\mu_i-\mu_{i+1})\\
-(4-2c)\sum_{i\in I_4}(\mu_i-\mu_{i+1})-(4+2c)\sum_{i\in I_3}(\mu_i-\mu_{i+1})\\
-2\sum_{i\in I_2}(\mu_i-\mu_{i+1})-(2+4c)\sum_{i\in I_1\setminus\{n\}}(\mu_i-\mu_{i+1})-(4+2c)\mu_n.
\end{split}
\end{equation}

We notice that, for $0\leq c\leq 1/3$, we have $4+2c\geq 4-2c$, $4+2c\geq 2$, $4+2c\geq 2+4c$ and $1-3c\geq 0$: in particular
\begin{equation}
\begin{split}
K_{X/C}^2\geq (4+c)\deg(f_*\omega_{X/C})&-(4+2c)(\mu_1-\mu_n)-(4+2c)\mu_n\geq\\
\geq(4+c)\deg(f_*&\omega_{X/C})-(4+2c)\mu_1.
\end{split}
\end{equation}
Combining this equation with $K_{X/C}^2\geq(2g-2)\mu_1$ (here we are using Proposition \ref{propo_xiaorelcan} with set of indexes $\{1\}$: here it is needed that $K_{X/C}$ is nef) we obtain
\begin{equation}
K_{X/C}^2\geq(4+c)\frac{g-1}{g+1+c}\deg(f_*\omega_{X/C})\geq(4+c)\frac{g-1}{g+2}\deg(f_*\omega_{X/C})
\end{equation}
where we are using the assumption $\deg(f_*\omega_{X/C})>0$ made at the beginning of this proof.
This concludes the proof.
\end{proof}

%

\begin{corollary}[cf. \cite{lu} Theorem 2.1 for the same result over the complex numbers]
Let $f\colon X\to C$ a surface fibration  such that $K_{X/C}$ is nef and all the quotients appearing in the Harder-Narasimhan filtration of $f_*\omega_{X/C}$ are strongly semi-stable and suppose that $K_{X/C}^2<\frac{9(g-1)}{2(g+2)}\deg(f_*\omega_{X/C})$. 
Then $f$ factors through a morphism $\pi\colon X\to Y$ of degree $2$, i.e. the following diagram commutes
\begin{equation}
\begin{tikzcd}
X\arrow{rr}{\pi}\arrow{dr}{f} & & Y\arrow{ld}{}\\
& C. & 
\end{tikzcd}
\end{equation}
\label{coro_slopegusunzhou2}
\end{corollary}

\begin{proof}
Suppose by contradiction that there is not such a morphism $f\colon X\to Y$: then we can apply part 1 of Theorem \ref{teo_slopegusunzhou} with $\delta=2$ and we get a contradiction.
\end{proof}

\section{Severi inequality}
\label{sec_sev}
In this section we give the proof of the Severi inequality\index{Severi!inequality} for surfaces of general type with maximal Albanese dimension. 
The original proof in \cite{par_sev} is just for complex surfaces and it is strongly based on the slope inequality of certain surface fibrations over the projective line constructed via a pencil that comes from an ample line bundle on the Albanese variety: as noticed in \cite{gusunzhou}, because, as we have seen, it is possible to extend the proof of the slope inequality in positive characteristic, the same proof works in general.
In \cite{yuan} it is possible to find the first proof of Severi inequality in positive characteristic which uses different methods.

Throughout this section we construct a surface fibration $\phi_d\colon\widetilde{X}_d\to\P^1$ which depends on the choice of a pencil inside $|L|$ for a fixed very ample line bundle $L$ in $\Alb(X)$: we decided not to mention the dependence of $\phi_d$ on the choice of such a pencil in order  to avoid a too heavy notation.
The reader should always keep this in mind.

Let $X$ be a minimal surface of general type with maximal Albanese dimension, let $\alb_X\colon X\to\Alb(X)$ be its Albanese morphism, denote by $q=\dim(\Alb(X))$ and let $L$ be a very ample line bundle on $\Alb(X)$ with $h^0(\Alb_X,L)\gg 0$.
For any integer $d$ coprime with the characteristic, consider the multiplication by $d$ morphism $\mu_d\colon \Alb(X)\to\Alb(X)$: it is known that (\cite{mumford} Proposition at the end of Section 6 and Corollary 1 Section 7)  it is an \'etale morphism of degree $d^{2q}$ and that (ibidem Section 8.(iv))
\begin{equation}
\mu_d^*(L)\sim_{num}d^2L.
\label{eq_'etalb}
\end{equation}

Now consider the following Cartesian diagram:
\begin{equation}
\label{eq_cartalb}
\begin{tikzcd}
X_d\arrow{r}{\nu_d} \arrow[dr, phantom, "\square"] \arrow{d}{a_d} & X\arrow{d}{\alb_X}\\
\Alb(X)\arrow{r}{\mu_d} & \Alb(X).
\end{tikzcd}
\end{equation}
It is clear that $\nu_d$ is \'etale and that $X_d$ is a connected minimal surface of general type of maximal Albanese dimension with invariants:
\begin{equation}
K_{X_d}^2=d^{2q}K_X^2,\qquad \chi(\O_{X_d})=d^{2q}\chi(\O_X).
\label{eq_invcartalb}
\end{equation}

Denote by $L_X=\alb_X^*L$ and $L_{X_d}=a_d^*L$: we have $L_{X_d}\sim_{num}d^{-2}\nu_d^*L_X$, in particular
\begin{equation}
L_{X_d}^2=d^{2q-4}L_X^2, \qquad K_{X_d}.L_{X_d}=d^{2q-2}K_X.L_X.
\label{eq_multcartalb}
\end{equation}

\begin{lemma}[cf. \cite{gusunzhou} Lemmas 4.1 and 5.3]
\label{lemma_bertinialb}
A general element $D$ of $a_d^*|L|$ is integral and non-hyperelliptic.
\end{lemma}

\begin{proof}[Sketch of proof]
Let $Z=\{(x,H)\in X_d\times \P\ |\ a_d(x)\in H\}$, where by $\P=(\P^{h^0(L)-1})^*$ we mean the dual projective space, i.e. the projective space whose points are the hyperplanes of $\P^{h^0(L)-1}$, and denote by $\pi_{X_d}\colon Z\to X_d$ and $\pi_{\P}\colon Z\to\P$ the two natural projections.
It is clear that $\pi_{X_d}$ is a $\P^{h^0(L)-2}$-bundle and the fibre of $\pi_{\P}$ over $H$ coincides with $a_d^*H$.

 By Theorem I.6.10(3) of \cite{jouanolou}, it follows that, for a general hyperplane $H$, $a_d^*H$ is irreducible.
In Theorem 2.1 of \cite{gusunzhou}, it is proved that, if for every hyperplane $H$, the divisor $a_d^*H$ were non-reduced then the morphism $a_d^*\Omega^1_{\Alb(X)}\to\Omega^1_{X_d}$ would be zero and, similarly, it would follow that $\alb_X^*\Omega^1_{\Alb(X)}\to\Omega^1_{X}$ is zero as well: this is not possible thanks to Remark \ref{rem_albfrob}.

Now suppose, by contradiction, that the general element of $a_d^*|L|$ is hyperelliptic and let $D_1,D_2\in a_d^*|L|$ be two general elements.
Consider the induced surface fibration $\widetilde{\phi}_d\colon\widetilde{X}_d\to\P^1$ where $\widetilde{X}_d$ is obtained by successive blow-ups of the intersection points of the strict transforms of $D_1$ and $D_2$.
Clearly the generic fibre cannot be inseparably hyperelliptic, otherwise it would be rational contradicting the hypothesis of maximality of the Albanese dimension of $X$.

Now consider $(\widetilde{\phi}_d)_*K_{\widetilde{X}_d/\P^1}$: this is a vector bundle of degree the arithmetic genus $g$ of a general fibre of $\widetilde{\phi}_d$.
Let $\P=\P( (\widetilde{\phi}_d)_*K_{\widetilde{X}_d/\P^1})$ be the projective bundle associated with $(\widetilde{\phi}_d)_*K_{\widetilde{X}_d/\P^1}$ and denote by $g\colon \widetilde{X}_d\dashrightarrow\P$ the rational map induced by the generically surjective morphism of sheaves $\widetilde{\phi}_d^*(\widetilde{\phi}_d)_*K_{\widetilde{X}_d/\P^1}\to K_{\widetilde{X}_d/\P^1}$ (cf. Proposition \ref{propo_morphtoproj} and the proof of Proposition \ref{propo_xiaonef}).
Because the restriction of $K_{\widetilde{X}_d/\P^1}$ to a general fibre is its canonical divisor, it is immediate that $g$ restricted to a general fibre is its canonical morphism.
In particular, it follows at once that $g$ is generically finite of degree $2$ and is induced by a birational involution $\widetilde{\sigma}_d$ on $\widetilde{X}_d$ such that $\widetilde{\phi}_d\circ\widetilde{\sigma}_d=\widetilde{\phi}_d$.

Because the birational automorphism group of a surface coincides with the automorphism group of its minimal model (Theorem \ref{teo_kxgeq0}), we have that there exists an automorphism $\sigma_d$ of $X_d$ such that $\phi_d\circ\sigma_d=\phi_d$ where, by $\phi_d$, we mean the rational fibration induced on $X_d$.
From the general choice of $D_1$ and $D_2$ and the ampleness of $L$, we have that $a_d\circ\sigma_d=a_d$ which is a contradiction, because in this case $a_d$ would contracts the fibres of $\phi_d$. 
\end{proof}

\begin{corollary}[\cite{gusunzhou} Corollary 4.2]
\label{coro_hypergu}
For a general element $(D_1,D_2)\in a_d^*|L|\times a_d^*|L|$ we have that $D_1$ and $D_2$ are integral and the intersection points of $D_1$ with $D_2$ lie in their smooth locus.
\end{corollary}

Now consider the pencil generated by $D_1=a_d^*H_1$ and $D_2=a_d^*H_2$ satisfying the conditions of Lemma \ref{lemma_bertinialb} and of Corollary \ref{coro_hypergu} and the corresponding rational map $X_d\dashrightarrow \P^1$.
Let $\phi_d\colon \widetilde{X}_d\to \P^1$ be its resolution of indeterminacy  obtained by successive blow-ups of the intersection points of the strict transforms of $D_1$ and $D_2$ and denote by $\pi_d\colon \widetilde{X}_d\to X_d$ the natural morphism.

\begin{lemma}[\cite{gusunzhou} Proposition 4.4]
Let $\pi_d\colon \widetilde{X}_d\to X_d$ be the morphism defined above and denote by $p_1,\ldots,p_r$ the intersection points of $D_1$ with $D_2$. 
Then we have
\begin{enumerate}
	\item $K_{\widetilde{X}_d}=\pi_d^*K_{X_d}+\sum_{i=1}^{r}(E_{i_1}+2E_{i_2}+\ldots+m_iE_{i_{m_i}})$ where $E_{i_j}$ denotes the strict transform of the exceptional divisor of the $j$-th blow-up over $p_i$;
	\item $E_{i_j}$ are $(-2)$-curves for $j<m_i$ and $E_{i_{m_i}}$ is a $(-1)$-curve;
	\item the strict transforms $\widetilde{D}_1$ and $\widetilde{D}_2$ of $D_1$ and $D_2$ are fibres of $\phi_d$, in particular the fibration has connected fibres;
	\item $E_{i_{m_i}}$ is a section of $\phi_d$ and $E_{i_j}$ is vertical for every $j<m_i$, in particular the fibration has no multiple fibres;
	\item the relative canonical bundle $K_{\widetilde{X}_d}-\phi_d^*K_{\P^1}$ is nef.
\end{enumerate}
\end{lemma}

%
%
%

\begin{remark}
\label{rem_notorsionhigher}
The fact that the fibration $\phi_d$ has a section implies that it has no multiple fibres.
Let $F$ be  a fibre of $\phi_d$.
For every decomposition $F=A+B$ the intersection $A.B$ is greater than $0$: indeed $0=F^2=A^2+2A.B+B^2<2A.B$ where the last inequality follows by the Zariski's Lemma for fibrations (\cite{badescu} Corollary 2.6).
We have thus shown that $F$ is $1$-connected (cf. \cite{barth} Section II.12), in particular $h^0(F,\O_F)=1$ for every fibre $F$ of $\phi_d$ (cf. ibidem Lemma II.12.2).
Combining this with the base change Theorem (\cite{hart} Corollary III.12.9) we can derive that $R^1f_{d*}\O_{\widetilde{X}_d}$ is locally free: This is always true for a surface fibration over the complex numbers (cf. \cite{barth} Corollary III.11.2), but it is not always the case for a surface fibration in positive characteristic (cf. \cite{badescu} Chapter 7). 
\end{remark}

Because the blown-up points on $X_d$ lie above the smooth locus of $D_1$, its arithmetic genus does not vary: in particular the arithmetic genus of a fibre of $\phi_d$ is
\begin{equation}
g_d=1+\frac{K_{X_d}.L_{X_d}+L_{X_d}^2}{2}=1+\frac{d^{2q-2}K_X.L_X+d^{2q-4}L_X^2}{2}.
\label{eq_genfibr}
\end{equation}
Notice that, clearly $K_X.L_X>0$ because $X$ is minimal and $L_X$ is the pull-back of an ample divisor via the Albanese morphism, hence it can not be a union of $-2$-curves: for $d\to\infty$ we have $g_d\to\infty$.

Moreover
\begin{equation}
\begin{split}
K_{\widetilde{X}_d/\P^1}^2=K_{\widetilde{X}_d}^2+4F.K_{\widetilde{X}_d}&=K_{X_d}^2+E^2+8(g_d-1)=\\
=K_{X_d}^2-L_{X_d}^2+4K&_{X_d}.L_{X_d}+4L_{X_d}^2=\\
=d^{2q}K_X^2+4d^{2q-2}K_X.&L_X+3d^{2q-4}L_X^2
\end{split}
\label{eq_cansqfib}
\end{equation}
where $E=\sum_{i=1}^r(E_{i_1}+2E_{i_2}+\ldots+E_{i_{m_i}})$.

Combining Equation \ref{eq_pushforstructure} with Remark \ref{rem_notorsionhigher} we get
\begin{equation}
\begin{split}
\deg((\phi_d)_*\omega_{\widetilde{X}_d/\P^1})=\chi(\O_{\widetilde{X}_d})+g_d-1=\\
=d^{2q}\chi(\O_X)+\frac{d^{2q-2}K_X.L_X+d^{2q-4}L_X^2}{2}.
\end{split}
\label{eq_deghigherpushf}
\end{equation}

Notice that, thanks to Theorem \ref{teo_shepbarr}, we have that $\chi(\O_X)>0$ which implies that $\deg((\phi_d)_*\omega_{\widetilde{X}_d/\P^1})>0$ for $d\gg 0$.
In particular we can define the slope of $\phi_d$ and we have that
\begin{equation}
\lambda(\phi_d)=\frac{2d^{2q}K_X^2+8d^{2q-2}K_X.L_X+6d^{2q-4}L_X^2}{2d^{2q}\chi(\O_X)+d^{2q-2}K_X.L_X+d^{2q-4}L_X^2}\xrightarrow{d\to\infty}\frac{K_X^2}{\chi(\O_X)}.
\label{eq_limitslope}
\end{equation}
 
So far we have proven (cf. \cite{par_sev} Proposition 2.3):
\begin{proposition}
\label{propo_luzuopard}
Let $X$ be a minimal surface of general type with maximal Albanese dimension.
Then there exists a sequence of surface fibrations $\phi_d\colon \widetilde{X}_d\to\P^1$ for $d>0$ and coprime to the characteristic of the ground field $k$ such that:
\begin{enumerate}
	\item $\lim_{d\to\infty}g_d=\infty$;
	\item $\lim_{d\to\infty}\lambda(\phi_d)=\frac{K_X^2}{\chi(\O_X)}$.
\end{enumerate}
\end{proposition}

\begin{figure}%
\centering
\begin{tikzpicture}[shorten >=1pt,node distance=2cm,on grid,auto, xscale=0.8, yscale=0.8, every node/.style={scale=1}]
    \begin{scope}
        \begin{axis}[ticks=none, xmin=0, xmax=7, ymin=0, ymax=30, axis lines=middle, xlabel=$\chi(\O_S)$, ylabel=$K^2_S$, enlargelimits]

	\draw (0,0) -- (3.4,30.6) node [sloped, midway, above] { \quad BMY $K_X^2\leq 9\chi(\O_X)$};

	\draw (0,0) -- (7.5,30) node [sloped,midway,below] {Severi ($q\geq2$) $K_X^2\geq 4\chi(\O_X)$};
	
	\filldraw  (1,4) circle (0.3pt);
	\filldraw  (1,5) circle (0.3pt);
	\filldraw  (1,6) circle (0.3pt);
	\filldraw  (1,7) circle (0.3pt);
	\filldraw  (1,8) circle (0.3pt);
	\filldraw  (1,9) circle (0.3pt);
	\filldraw  (2,8) circle (0.3pt);
	\filldraw  (2,9) circle (0.3pt);
	\filldraw  (2,10) circle (0.3pt);
	\filldraw  (2,11) circle (0.3pt);
	\filldraw  (2,12) circle (0.3pt);
	\filldraw  (2,13) circle (0.3pt);
	\filldraw  (2,14) circle (0.3pt);
	\filldraw  (2,15) circle (0.3pt);
	\filldraw  (2,16) circle (0.3pt);
	\filldraw  (2,17) circle (0.3pt);
	\filldraw  (2,18) circle (0.3pt);
	\filldraw  (3,12) circle (0.3pt);
	\filldraw  (3,13) circle (0.3pt);
	\filldraw  (3,14) circle (0.3pt);
	\filldraw  (3,15) circle (0.3pt);
	\filldraw  (3,16) circle (0.3pt);
	\filldraw  (3,17) circle (0.3pt);
	\filldraw  (3,18) circle (0.3pt);
	\filldraw  (3,19) circle (0.3pt);
	\filldraw  (3,20) circle (0.3pt);
	\filldraw  (3,21) circle (0.3pt);
	\filldraw  (3,22) circle (0.3pt);
	\filldraw  (3,23) circle (0.3pt);
	\filldraw  (3,24) circle (0.3pt);
	\filldraw  (3,25) circle (0.3pt);
	\filldraw  (3,26) circle (0.3pt);
	\filldraw  (3,27) circle (0.3pt);
	\filldraw  (4,16) circle (0.3pt);
	\filldraw  (4,17) circle (0.3pt);
	\filldraw  (4,18) circle (0.3pt);
	\filldraw  (4,19) circle (0.3pt);
	\filldraw  (4,20) circle (0.3pt);
	\filldraw  (4,21) circle (0.3pt);
	\filldraw  (4,22) circle (0.3pt);
	\filldraw  (4,23) circle (0.3pt);
	\filldraw  (4,24) circle (0.3pt);
	\filldraw  (4,25) circle (0.3pt);
	\filldraw  (4,26) circle (0.3pt);
	\filldraw  (4,27) circle (0.3pt);
	\filldraw  (4,28) circle (0.3pt);
	\filldraw  (4,29) circle (0.3pt);
	\filldraw  (4,30) circle (0.3pt);
	\filldraw  (5,20) circle (0.3pt);
	\filldraw  (5,21) circle (0.3pt);
	\filldraw  (5,22) circle (0.3pt);
	\filldraw  (5,23) circle (0.3pt);
	\filldraw  (5,24) circle (0.3pt);
	\filldraw  (5,25) circle (0.3pt);
	\filldraw  (5,26) circle (0.3pt);
	\filldraw  (5,27) circle (0.3pt);
	\filldraw  (5,28) circle (0.3pt);
	\filldraw  (5,29) circle (0.3pt);
	\filldraw  (5,30) circle (0.3pt);
	\filldraw  (6,24) circle (0.3pt);
	\filldraw  (6,25) circle (0.3pt);
	\filldraw  (6,26) circle (0.3pt);
	\filldraw  (6,27) circle (0.3pt);
	\filldraw  (6,28) circle (0.3pt);
	\filldraw  (6,29) circle (0.3pt);
	\filldraw  (6,30) circle (0.3pt);
	\filldraw  (7,28) circle (0.3pt);
	\filldraw  (7,29) circle (0.3pt);
	\filldraw  (7,30) circle (0.3pt);

	\end{axis}
    \end{scope}

    \begin{scope}[xshift=7.5cm]
        \begin{axis}[ticks=none, xmin=-1, xmax=7, ymin=0, ymax=30, axis lines=middle, xlabel=$\chi(\O_S)$, ylabel=$K^2_S$, enlargelimits]

	\draw (0,0) -- (2.5,30) node [sloped, midway, above] {\quad \ $K_X^2<12\chi(\O_X)$};
	
	\draw (0,0) -- (7.5,30) node [sloped,midway,below] {Severi ($q\geq2$) $K_X^2\geq 4\chi(\O_X)$};

  \filldraw  (1,4) circle (0.3pt);
	\filldraw  (1,5) circle (0.3pt);
	\filldraw  (1,6) circle (0.3pt);
	\filldraw  (1,7) circle (0.3pt);
	\filldraw  (1,8) circle (0.3pt);
	\filldraw  (1,9) circle (0.3pt);
	\filldraw  (1,10) circle (0.3pt);
	\filldraw  (1,11) circle (0.3pt);
	\filldraw  (2,8) circle (0.3pt);
	\filldraw  (2,9) circle (0.3pt);
	\filldraw  (2,10) circle (0.3pt);
	\filldraw  (2,11) circle (0.3pt);
	\filldraw  (2,12) circle (0.3pt);
	\filldraw  (2,13) circle (0.3pt);
	\filldraw  (2,14) circle (0.3pt);
	\filldraw  (2,15) circle (0.3pt);
	\filldraw  (2,16) circle (0.3pt);
	\filldraw  (2,17) circle (0.3pt);
	\filldraw  (2,18) circle (0.3pt);
	\filldraw  (2,19) circle (0.3pt);
	\filldraw  (2,20) circle (0.3pt);
	\filldraw  (2,21) circle (0.3pt);
	\filldraw  (2,22) circle (0.3pt);
	\filldraw  (2,23) circle (0.3pt);
	\filldraw  (3,12) circle (0.3pt);
	\filldraw  (3,13) circle (0.3pt);
	\filldraw  (3,14) circle (0.3pt);
	\filldraw  (3,15) circle (0.3pt);
	\filldraw  (3,16) circle (0.3pt);
	\filldraw  (3,17) circle (0.3pt);
	\filldraw  (3,18) circle (0.3pt);
	\filldraw  (3,19) circle (0.3pt);
	\filldraw  (3,20) circle (0.3pt);
	\filldraw  (3,21) circle (0.3pt);
	\filldraw  (3,22) circle (0.3pt);
	\filldraw  (3,23) circle (0.3pt);
	\filldraw  (3,24) circle (0.3pt);
	\filldraw  (3,25) circle (0.3pt);
	\filldraw  (3,26) circle (0.3pt);
	\filldraw  (3,27) circle (0.3pt);
	\filldraw  (3,28) circle (0.3pt);
	\filldraw  (3,29) circle (0.3pt);
	\filldraw  (3,30) circle (0.3pt);
	\filldraw  (4,16) circle (0.3pt);
	\filldraw  (4,17) circle (0.3pt);
	\filldraw  (4,18) circle (0.3pt);
	\filldraw  (4,19) circle (0.3pt);
	\filldraw  (4,20) circle (0.3pt);
	\filldraw  (4,21) circle (0.3pt);
	\filldraw  (4,22) circle (0.3pt);
	\filldraw  (4,23) circle (0.3pt);
	\filldraw  (4,24) circle (0.3pt);
	\filldraw  (4,25) circle (0.3pt);
	\filldraw  (4,26) circle (0.3pt);
	\filldraw  (4,27) circle (0.3pt);
	\filldraw  (4,28) circle (0.3pt);
	\filldraw  (4,29) circle (0.3pt);
	\filldraw  (4,30) circle (0.3pt);
	\filldraw  (5,20) circle (0.3pt);
	\filldraw  (5,21) circle (0.3pt);
	\filldraw  (5,22) circle (0.3pt);
	\filldraw  (5,23) circle (0.3pt);
	\filldraw  (5,24) circle (0.3pt);
	\filldraw  (5,25) circle (0.3pt);
	\filldraw  (5,26) circle (0.3pt);
	\filldraw  (5,27) circle (0.3pt);
	\filldraw  (5,28) circle (0.3pt);
	\filldraw  (5,29) circle (0.3pt);
	\filldraw  (5,30) circle (0.3pt);
	\filldraw  (6,24) circle (0.3pt);
	\filldraw  (6,25) circle (0.3pt);
	\filldraw  (6,26) circle (0.3pt);
	\filldraw  (6,27) circle (0.3pt);
	\filldraw  (6,28) circle (0.3pt);
	\filldraw  (6,29) circle (0.3pt);
	\filldraw  (6,30) circle (0.3pt);
	\filldraw  (7,28) circle (0.3pt);
	\filldraw  (7,29) circle (0.3pt);
	\filldraw  (7,30) circle (0.3pt);
	
	\end{axis}
    \end{scope}
\end{tikzpicture}\caption{A comparison for the geography of Surfaces of general type with maximal Albanese dimension for $\cha(k)=0$ (on the left) and $\cha(k)>0$ (on the right). Every dot symbolizes a couple $(a,b)$ in the plane for which it may be possible to find a surface of general type with maximal Albanese dimension $X$ with invariants $K_X^2=a$ and $\chi(\O_X)=b$.}
\label{geogen}%
\end{figure}
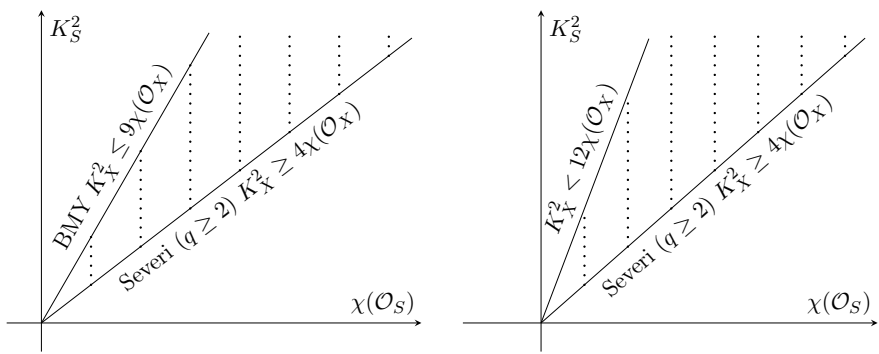

\begin{theorem}[Severi inequality, cf. \cite{par_sev} Theorem 2.1]
\label{teo_sev}
Let $X$ be a minimal surface of general type of maximal Albanese dimension.
Then $K_X^2\geq 4\chi(\O_X)$.
\end{theorem}
\begin{proof}
Consider the surface fibrations $\phi_d\colon \widetilde{X}_d\to\P^1$ constructed above.
By the Slope inequality \ref{teo_slope} we have that $\lambda(\phi_d)\geq\frac{4g_d-4}{g_d}$, i.e. (combining this with Equations \ref{eq_genfibr} and \ref{eq_limitslope})
\begin{equation}
\begin{split}
\lambda(\phi_d)=\frac{2d^{2q}K_X^2+8d^{2q-2}K_X.L_X+6d^{2q-4}L_X^2}{2d^{2q}\chi(\O_X)+d^{2q-2}K_X.L_X+d^{2q-4}L_X^2}\geq \\
\geq\frac{4d^{2q-2}K_X.L_X+4d^{2q-4}L_X^2}{d^{2q-2}K_X.L_X+d^{2q-4}L_X^2+2}
\end{split}
\label{eq_sevfinale}
\end{equation}
and, for $d\to\infty$ we get the thesis.
\end{proof}

\section{Factors of the Albanese morphism}
\label{sec_albfact}
In this section we are going to prove that if a surface of general type $X$ with maximal Albanese\index{Albanese!morphism, factors of} dimension satisfies $K_X^2<\frac{9}{2}\chi(\O_X)$, then there exists a morphism of degree two $f\colon X\to Y$ such that the Albanese morphism of $X$ factors through $f$ (Theorem \ref{teo_albfact2}).
If the characteristic is $2$ and the Albanese morphism is inseparable, this directly follows from Lemma \ref{lemma_albinsfact}.
If the characteristic is zero, this result has been proven in \cite{lu} Theorem 3.1 using the linear bound for the automorphism group of a surface of general type of Xiao (\cite{xiaobound} Theorem 1) which does not hold in positive characteristic. 
In positive characteristic this is some how hidden in \cite{gusunzhou}, provided that one refines their Theorem 3.1(2) as we have done in the second part of Theorem \ref{teo_slopegusunzhou}. 

Throughout this section we still use the fibrations $\phi_d$ constructed in the previous section: as before, the reader should keep in mind their dependence on the choice of a pencil, which is not stressed by the notation.
The same is still valid for the rational degree two maps $f_d$ that we are going to construct here. 

\begin{theorem}
Let $X$ be a surface of general type and let $\phi_d$ be the surface fibration constructed in Section \ref{sec_sev} for $d\gg 0$ and coprime with the characteristic of the ground field.
If $K_X^2<\frac{9}{2}\chi(\O_X)$ then there exists a rational map $f_d\colon X_d\to Y_d$ to a surface $Y_d$ of degree $2$ such that
\begin{equation}
\begin{tikzcd}
X_d\arrow[dashed]{dr}{\phi_d}\arrow[dashed]{rr}{f_d} & & Y_d\arrow{dl}{\psi_d}\\
& \P^1 &
\end{tikzcd}
\end{equation}
commutes.
\label{teo_lualb2}
\end{theorem}

\begin{proof}
It is clear by hypothesis that $K_X^2\leq\frac{9}{2}\chi(\O_X)-\frac{1}{2}$.
Recall that by Proposition \ref{propo_luzuopard} we have that $\lambda(\phi_d)\to\frac{K_X^2}{\chi(\O_X)}$ and $g_d\to\infty$ as $d\to\infty$ where $g_d$ is the arithmetic genus of a generic fibre of $\phi_d$.
From this we derive that, for a sufficiently large $d$,
\begin{equation}
\lambda(\phi_d)<\frac{K_X^2}{\chi(\O_X)}+\frac{1}{4\chi(\O_X)}\leq\frac{9}{2}-\frac{1}{4\chi(\O_X)}<\frac{9}{2}\frac{g_d-1}{g_d+2}.
\end{equation}

Applying Corollary \ref{coro_slopegusunzhou2} we get that there necessarily exists a rational map of degree two $f_d\colon {X}_d\dashrightarrow Y_d$ satisfying the thesis of the Theorem.
Notice that $Y_d\to\P^1$ is an actual morphism by the construction of $Y_d$: recall that it is the image of a rational map from $X_d$ to a projective bundle over $\P^1$ (cf. introduction to Section \ref{sec_nonhyper}).
\end{proof}

\begin{remark}
It is possible to prove also that $f_d$ is an actual morphism if it is separable.
Indeed, in this case $f_d$ is obtained as a quotient of a rational involution $\sigma_d$ of $X_d$, which extends to an actual involution because $X_d$ is minimal (Theorem \ref{teo_kxgeq0}).

On the other hand, up to passing to the normalization $\nu\colon\widetilde{Y}_d\to Y_d$, if $f_d$ is inseparable we obtain that $f_d\circ\pi_d\colon\widetilde{X}_d\to \widetilde{Y}_d$ is the normalized relative Frobenius morphism (cf. \cite{gusunzhou} Example 2.5), i.e. 
\begin{equation}
\begin{tikzcd}
\widetilde{X}_d\arrow{r}{f_d\circ\pi_d}\arrow[bend left=30]{rr}{F_{X_d/\P^1}}\arrow{dr}{\phi_d} & \widetilde{Y}_d\arrow{d}{\psi_d}\arrow{r}{\nu}& \P^1\times_{\P^1}\widetilde{X}_d\arrow{d}{}\arrow{r}{}& \widetilde{X}_d\arrow{d}{\phi_d}\\
& \P^1\arrow[equal]{r}{} & \P^1\arrow{r}{F_k}&\P^1.
\end{tikzcd}
\end{equation}
\end{remark}

\begin{lemma}[cf. \cite{lu} Lemma 3.2]
Suppose that $f_d$ is separable for a general pencil in $a_d^*|L|$ and assume, as usual, that $d\gg 0$ and coprime with the characteristic of the ground field.
Then, for a general pencil in $a_d^*|L|$, the involution $\sigma_d$ of $X_d$ corresponding to $f_d$ satisfies $a_d\circ\sigma_d=a_d$.
\end{lemma}

\begin{proof}
The proof given for Lemma 3.2 \cite{lu} is characteristic-free if we assume that $f_d$ is separable.
\end{proof}

\begin{proposition}[cf. \cite{gusunzhou} Proposition 5.10]
\label{propo_gusunzhoufactalb}
Suppose that for infinitely many prime numbers $l$ we have a suitable subpencil of $a_l^*|L|$ such that the corresponding degree-two morphism $f_l$ is separable and let $\sigma_l$ be its corresponding involution.
Then $f_l$ descends to a separable degree-two morphism $f\colon X\to Y$ induced by an involution $\sigma$ of $X$, such that $\alb_X\circ\sigma=\alb_X$.
\label{propo_infiniteseparable2}
\end{proposition}

For the proof of this Proposition we need some preliminary Lemmas.

\begin{lemma}[cf. \cite{gusunzhou} Lemma 5.11]
Consider the following cartesian diagram
\begin{equation}
\begin{tikzcd}
V_l\arrow[dr, phantom, "\square"]\arrow{r}{\nu_l}\arrow{d}{}&V\arrow{d}{}\\
A\arrow{r}{\mu_l}&A
\end{tikzcd}
\label{eq_gusunzhouvarious}
\end{equation}
where $V$ is an integral scheme, $A$ is an Abelian variety, $\mu_l$ is the multiplication by $l$ morphism where $l$ is a prime different from the characteristic of the ground field.
Suppose there exists a generically finite separated morphism $\pi\colon V'\to V$ where $V'$ is an integral scheme and that for infinitely many primes $l$ the morphism $\nu_l$ factors through $\pi$: then $\pi$ is an isomorphism.
\label{lemma_gusunzhouvarious}
\end{lemma}

\begin{proof}
By assumption we know that there exists an $l>|k(V'):k(V)|$\nomenclature{$!k_1:k_2|$}{the degree of an extension of fields $k_2\subseteq k_2$} where $k(V)$, respectively $k(V')$, represents the field of rational functions of $V$, respectively $V'$, $|k(V'):k(V)|$ the finite degree of the extension, such that $\nu_l$ factors through $\pi$.
In particular there exists a connected component $V_l'$ of $V_l$ such that
\begin{equation}
\begin{tikzcd}
V_l'\arrow{r}{\tau}\arrow{dr}{\nu_l'}& V'\arrow{d}{\pi}\\
& V
\end{tikzcd}
\label{eq_non serve}
\end{equation}
commutes and $\nu_l'$ is given by a quotient of $V_l'$ via the free action of a subgroup $\Gamma$ of $(\Z/l\Z)^{2q}$ where $q$ is the dimension of $A$.

In particular we have that 
\begin{equation}
|k(V'):k(V)||k(V_l'):k(V')|=|\Gamma|=l^a
\label{eq_maj}
\end{equation}
with $a\leq 2q$.
From this and the assumption on $l$, we derive that $K(V)=K(V')$, i.e. $\pi$ is a birational map, hence $\tau$ is dominant.

By \cite{liu} Lemma 3.3.15(b) and Exercise 3.3.17(b) we have that $\tau$ is finite hence it is also surjective (indeed any finite morphism is proper, hence closed, ibidem Exercise 3.3.17(d)).
From this it follows that $\pi$ is quasi-finite (i.e. it has finite fibres).
By \cite{liu} Proposition 3.3.16(f) we have that $\pi$ is proper and, by \cite{egaiv4} Corollaire 18.12.4, it is also finite.

We have that $\O_V\subseteq \pi_*\O_{V'}$ because $\pi$ is a dominant morphism and $\pi_*\O_{V'}\subseteq (\nu_l')_*\O_{V_l'}^\Gamma=\O_V$ where $\O_{V_l'}^\Gamma$ denotes the subsheaf of $\O_{V_l'}$ of sections fixed by the action of $\Gamma$.
The last inclusion follows from the fact that $V'$ is birational to $V$, in particular the pull-back of every function on $V'$ is fixed by the action of $\Gamma$ on an open dense subset, hence it is fixed everywhere.
Clearly a finite morphism which induces an isomorphism on the structure sheaves is an isomorphism and we are done.
\end{proof}

Now for every scheme $X$ over a base $S$ and another scheme $T$ over the same base, we denote by $A_{X/S}(T)=\aut_T(X_T)$\nomenclature{$A_{X/S}$}{the functor of the automorphism of an $S$-scheme $X$}where $X_T=X\times_S T$, $\aut_T(X_T)$ denotes the set of automorphism of $X\times_S T$ seen as a $T$-scheme. 
For every morphism of schemes $\tau\colon T'\to T$ we can naturally define a map $f^*\colon A_{X/S}(T)\to A_{X/S}(T')$. 
In particular $A_{X/S}$ is a contravariant functor from the category of schemes over $S$ to the category of groups.
It is known that (cf. \cite{fantechi} Theorem 5.23 and the Exercise immediately following) that this functor is represented by a group scheme (which we still call $A_{X/S}$) separated over $S$ provided that $X\to S$ is a projective flat morphism.
In particular it is representable if $S=\spec(k)$: in this case we can even prove (cf. \cite{matsura} Theorem 3.7) that this functor is represented by a group scheme locally of finite type over $k$ which we call $A_X$.
Moreover, if $X$ is a surface of general type, $A_X$ is a finite group scheme because $\aut(X)$ is a finite group (cf. Theorem \ref{teo_autogeneral}).

\begin{lemma}
\label{lemma_functorofpoints}
Let $f\colon X\to V$ be a surjective morphism from a smooth surface of general type to a projective integral surface. 
Let $f'\colon X'\to V'$ the restriction of $f$ to the locus where it is flat (recall that $X'$ and $V'$ are open dense subsets of $X$ and $V$ respectively by \cite{egaiv2} 6.9.1) and denote by $A_{X'/V'}[2]$\nomenclature{$A_{X/S}[2]$}{the subfunctor of involutions of $A_{X/S}$} the subfunctor of $A_{X'/V'}$ which, to a $V'$-scheme $T$, associates the set of involution of $X'\times_{V'} T$ as a scheme over $T$.
Then $A_{X'/V'}[2]$ is represented by a closed subscheme of $A_{X'/V'}$, which we still denote by $A_{X'/V'}[2]$, that is generically finite and separated over $V'$.
\end{lemma}

\begin{proof}
Observe that $A_{X'/V'}$ is representable because $X'\to V'$ is projective and flat by definition.
Let $S\to A_{X'/V'}$ be a morphism of schemes: by Yoneda's Lemma (\cite{fantechi} 2.1.2) this is equivalent to take an element $\sigma_S\in A_{X'/V'}(S)$.
Denote by $S_{\sigma_S}$ the closed subscheme of $S$ over which $\sigma_S^2=1$.
It is then clear that the Cartesian product $A_{X'/V'}[2]\times_{A_{X'/V'}} S$ is isomorphic to $S_{\sigma_S}$: indeed, for every scheme $T$, one has 
\begin{equation*}
\begin{split}
(A_{X'/V'}[2]&\times_{A_{X'/V'}} S)(T)=\\
=\{(\sigma_T,\tau)\in A_{X'/V'}(T)&\times\Hom(T,S)\ |\  \sigma_T^2=Id, \tau^*\sigma_S=\sigma_T\}=\\
=&\Hom(T,S_{\sigma_S}).
\end{split}
\label{eq_functorclosed}
\end{equation*}
From this, for $S=A_{X'/V'}$, it follows at once that $A_{X'/V'}[2]$ is represented by a closed subscheme of $A_{X'/V'}$, in particular it is separated over $V$ (cf. \cite{hart} Corollary II.4.6(a)(b)).

Observe that $A_{X'/V'}\to V'$ is generically finite (hence also $A_{X'/V'}[2]\to V'$ is) because $X\to V$ is: indeed for a general point $p\in V'$ we have that the fibre of $f'$ is finite, in particular the automorphism group of the fibre is finite.
\end{proof}

\begin{proof}[Proof of Proposition \ref{propo_gusunzhoufactalb}]
Denote by $V$ the image of the Albanese morphism $\alb_X\colon X\to \Alb(X)$ and by $X'$ and $V'$ the open dense subsets where $\alb_X$ is flat. 
By Lemma \ref{lemma_functorofpoints} we know that $M=A_{X'/V'}[2]$ is a generically finite separated scheme over $V'$.
By hypothesis, there is an irreducible (and reduced) component $M_i$ such that for infinitely many primes $l$ we have a morphism $V_l'\to M_i$ (where $V_l'$ is the pull-back of $V'$ via the multiplication by $l$ on $\Alb(X)$): indeed, by  Yoneda's Lemma to construct a morphism to $M$, it is enough to give an involution $\sigma_l$ on $X_l'$ (where $X_l'$ is the pull-back via $V_l'$ of $X$) with $a_l\circ\sigma_l=a_l$ (and we have that such an involution exists on $X_l$ by hypothesis for $l$ sufficiently large) and clearly, because $M$ has finitely many irreducible component, there is at least one $M_i$ such that the image of $V_l'$ is contained in $M_i$.

Using Lemma \ref{lemma_gusunzhouvarious} we have that $M_i=V'$.
In particular we have a section $V'\to M\to V'$ which, again by Yoneda's Lemma, proves that there exists a rational involution $\sigma$ of $X$ satisfying $\alb_X\circ \sigma=\alb_X$ and the Proposition follows from the fact that on a minimal surface of general type every birational map is an isomorphism (Theorem \ref{teo_kxgeq0}).
\end{proof}

There is still one case left before the proof of the main Theorem of this Section.
If the characteristic of the ground field is $2$  and the Albanese morphism is separable, it may happen that, for all but finite $d$  coprime with $2$ and for a general choice of a pencil in $a_d^*|L|$, the degree two rational map $f_d$ is inseparable: in this situation it is not possible to apply Proposition \ref{propo_gusunzhoufactalb}.
We are going to prove that in this case we have $K_X^2\geq\frac{9}{2}\chi(\O_X)$.

\begin{lemma}[\cite{gusunzhou} Lemma 5.8] 
\label{lemma_sepalbchar2}
Suppose that $\alb_X\colon X\to\Alb(X)$ is a separable morphism, that the characteristic of the ground field is $2$ and that for a general choice of a pencil in $a_d^*|L|$ for a fixed $d$ coprime with $2$, we have that the corresponding rational map $f_d\colon X_d\to Y_d$ is inseparable.
Then for a suitable pencil we have
\begin{equation}
\rk((\phi_d)_*(T(\Omega_{\widetilde{X}_d/\P^1}))\leq c_1(\Omega_{X_d/\Alb(X)}).L_{X_d}
\label{eq_sepalbchar2}
\end{equation}
where $\phi_d$ is the same as above and, for a sheaf $E$, by $T(E)$ we mean the torsion subsheaf of $E$.
\end{lemma}

With the notation as in Theorem \ref{teo_lualb2}, let $\widetilde{Y}_d$ be the minimal model of the resolution of singularities of $Y_d$ and let $g_d'$ be the arithmetic genus of a general fibre of $\psi_d\colon \widetilde{Y}_d\to\P^1$.
We are going to give an estimate of $g_d'/g_d$ in order to apply  Theorem \ref{teo_slopegusunzhou}.

\begin{proposition}[cf. \cite{gusunzhou} Propositions 5.2 and 5.13] 
Suppose that the hypotheses of Lemma \ref{lemma_sepalbchar2} are satisfied for all but finitely many $d$.
Then we have
\begin{equation}
\lim_{d\to\infty}\frac{g_d'}{g_d}\geq\frac{1}{2}.
\label{nonmiserve}
\end{equation}
\label{propo_nonmiserve}
\end{proposition}

\begin{proof}
By \cite{gu_2016} Proposition 2.2 (cf. also \cite{gusunzhou} proof of Proposition 5.13) we have that $g_d'=g_d-\frac{1}{4}\rk\Bigl((\phi_d)_*\bigl(T(\Omega_{\widetilde{X}_d/\P^1}^1)\bigr)\Bigr)$.
So, using Lemma \ref{lemma_sepalbchar2} and Equation \ref{eq_genfibr} we get
\begin{equation}
\label{eq_naltra}
\begin{split}
g_d'\geq g_d-\frac{L_{X_d}.c_1\bigl(\Omega_{X_d/\Alb(X)}^1\bigr)}{4}=&\\
=\frac{d^{2q-2}K_X.L_X+d^{2q-4}L_X^2}{2}-\frac{d^{2q-2}c_1(\Omega_{X/\Alb(X)}^1).L_X}{4}&+1=\\
=\frac{d^{2q-2}(2K_X-c_1(\Omega_{X/\Alb(X)}^1)).L_X+2d^{2q-4}L_X^2}{4}&+1.
\end{split}
\end{equation}
 
By \cite{gusunzhou} Proposition 5.2, we have  
\begin{equation}
\label{eq_gusunzhouc1}
(2K_X-c_1(\Omega_{X/\Alb(X)}^1)).L_X\geq K_X.L_X.
\end{equation}
Hence 
\begin{equation}
\lim_{d\to\infty}\frac{g_d'}{g_d}\geq\lim_{d\to\infty}\frac{2(d^{2q-2}K_X.L_X+2d^{2q-4}L_X^2)+8}{4(d^{2q-2}K_X.L_X+2d^{2q-4}L_X^2)+8}=\frac{1}{2}
\label{eq}
\end{equation}
and the Proposition is proven.
\end{proof}

\begin{theorem}
Let $X$ be a surface of general type with maximal Albanese dimension and suppose that $K_X^2<\frac{9}{2}\chi(\O_X)$. 
Then there exists a morphism of degree two $f\colon X\to Y$ to a normal surface $Y$ such that the Albanese morphism of $X$ factors through $f$, i.e. the following diagram commutes:
\begin{equation}
\begin{tikzcd}
X\arrow{rr}{f}\arrow{dr}{\alb_X} & & Y\arrow{ld}{}\\
& \Alb(X). & 
\end{tikzcd}
\end{equation}
\label{teo_albfact2}
\end{theorem}

\begin{proof}
Suppose that the characteristic of the ground field is different from $2$.
Then the Theorem directly follows from Proposition \ref{propo_gusunzhoufactalb}.

Suppose that the characteristic of the ground field is two.
If the Albanese morphism is inseparable, we have already noticed at the beginning of this section that the assert of the Theorem is exactly Lemma \ref{lemma_albinsfact}.

Now suppose that the Albanese morphism is separable.
Thanks to Proposition \ref{propo_nonmiserve} and part two of Theorem \ref{teo_slopegusunzhou}, we can derive that if all but finitely many $f_d$ are inseparable, then $K_X^2\geq\frac{9}{2}\chi(\O_X)$ which is excluded by the hypothesis.
Hence, also in this case the thesis is a consequence of Proposition \ref{propo_gusunzhoufactalb}.
\end{proof}

\section{Refined Severi inequality}
\label{sec_gusunc12}
In this\index{Severi!inequality, refined} Section we  outline how to generalize Theorem 5.14 of \cite{gusunzhou} thanks to the generalization we have given in Theorem \ref{teo_slopegusunzhou}.
For the details we refer to their Section 5, even though most of their effort have been recalled throughout this Chapter.

Let $X$ be a minimal surface of general type with maximal Albanese dimension and denote by $\alb_X\colon X\to \Alb(X)$ its Albanese morphism.
Notice that there are finitely many rational maps of degree two $\pi_i\colon X\dashrightarrow Y_i$ with $Y_i$ smooth and minimal such that there exists a morphism $b_i\colon Y_i\to \Alb(X)$ for which $b_i\circ\pi_i=\alb_X$.
Indeed the separable morphisms $\pi_i$ are in natural correspondence with a subset of the group of automorphism of $X$ which is finite (cf. Theorem \ref{teo_autogeneral}), while there exists a unique inseparable morphism $\pi_i$ of degree two when $\alb_X$ is inseparable and the characteristic of the ground field is two and it is uniquely determined by the 1-foliation $\T_{X/\Alb(X)}$ (cf. Lemma \ref{lemma_albinsfact}).

\begin{definition}
Let $\pi_i\colon X\dashrightarrow Y_i$ as above and let $L$ be a very ample line bundle on $\Alb(X)$.
Denote by $L_X=\alb_X^*L$ and $L_{Y_i}=b_i^*L$.
We define 
\begin{equation}
c_i(X,L)=\frac{K_{Y_i}.L_{Y_i}}{K_X.L_X}.
\label{eq_cixl}
\end{equation}
If the characteristic of the ground field is two and $\alb_X$ is separable we define
\begin{equation}
c_0(X,L)=\frac{(2K_X-c_1(\Omega_{X/\Alb(X)}^1)).L_X}{2K_X.L_X}.
\label{eq_c0xl}
\end{equation}

Now take 
\begin{equation}
c(X,L)=\min_{i>0}\{c_i(X,L)\}
\label{eq_cxl}
\end{equation}
 and
\begin{equation}
c(X)=\sup\{c(X,L)\}
\label{eq_cx}
\end{equation}
for all possible very ample line bundles $L$ on $\Alb(X)$.
\label{def_cxl}
\end{definition}

It is clear by definition that $c_i(X,L)=0$ if and only if $Y_i$ is an Abelian surface.

\begin{remark}
Observe that the definition of all $c_i(X,L)$ and, consequently the definition of $c(X,L)$, does not change if we substitute $L$ with a multiple of $L$.
\label{rem_cxl}
\end{remark}

\begin{theorem}[cf. \cite{gusunzhou} Theorem 5.14]
Let $X$ be a minimal surface of general type of maximal Albanese dimension.
Then we have
\begin{equation}
K_X^2\geq\Bigl(4+\min \{c(X),\frac{1}{2}\}\Bigr)\chi(\O_X)
\label{eq_cxl}
\end{equation}
where $c(X)$ is as in Definition \ref{def_cxl}.
\label{teo_cxl}
\end{theorem}

\begin{proof}
Recall that in Section \ref{sec_albfact} we constructed the following commutative diagram
\begin{equation}
\begin{tikzcd}
& X_d\arrow[dashed]{dd}{f_d}\arrow[dashed]{dl}{\phi_d}\arrow{dr}{a_d}& \\
\P^1 & & \Alb(X)\\
& \widetilde{Y}_d \arrow{ul}{\psi_d}\arrow{ur}{b_d}& 
\end{tikzcd}
\label{eq_nonsocomexcl}
\end{equation}
where $a_d\colon X_d\to \Alb(X)$ is the pull-back via the multiplication by an integer $d$ coprime with the characteristic of the ground field of the Albanese morphism $\alb_X\colon X\to \Alb(X)$ (cf. Equation \ref{eq_invcartalb}), $f_d$ is a degree-two rational map and $\widetilde{Y}_d$ can be assumed to be smooth and minimal.
We have defined $g_d$ to be the arithmetic genus of a general fibre of $\phi_d$ and $g'_d$ to be the arithemtic genus of a general fibre of $\psi_d$.
Recall also that this construction depends on the choice of  a pencil in $a_d^*|L|$ where $L$ is a very ample line bundle on $\Alb(X)$.
Arguing as in Lemma \ref{lemma_bertinialb}, we can restrict the choice of the pencil so that the general fibre of $\psi_d$ is integral too.

If the characteristic of the ground field is $2$ and $\Alb(X)$ is inseparable or if the the hypotheses of Proposition \ref{propo_gusunzhoufactalb} are satisfied, then the right hand side of diagram \ref{eq_nonsocomexcl} descends to
\begin{equation}
\begin{tikzcd}
X\arrow[dashed]{rr}{\pi_i}\arrow{dr}{\alb_X}& & Y_i\arrow{dl}{b_i}\\
& \Alb(X) &
\label{eq_nonsoxcl}
\end{tikzcd}
\end{equation}
where the horizontal arrow has to be one of the $\pi_i$ defined above.
In this case, a formula similar to the one for $g_d$ (Equation \ref{eq_genfibr}) holds also for $g_d'$, that is
\begin{equation}
g_d'=1+\frac{d^{2q-2}K_{Y_i}.L_{Y_i}+d^{2q-4}L_{Y_i}^2}{2}
\label{eq_genfibr'}
\end{equation}
Combining Equations \ref{eq_genfibr} and \ref{eq_genfibr'}, we easily derive that 
\begin{equation}
\lim_{d\to\infty}\frac{g_d'}{g_d}=c_i(X,L);
\label{eq_ultimaboh}
\end{equation}
in particular, for $d$ sufficiently large, Theorem \ref{teo_slopegusunzhou}(2) applies and, for $d\to\infty$, we get (cf. Proposition \ref{propo_luzuopard})
\begin{equation}
K_X^2\geq \Bigl(4+\min \{c_i(X,L),\frac{1}{2}\}\Bigr)\chi(\O_X).
\label{eq_mahult}
\end{equation}

If the characteristic of the ground field is two and $\alb_X$ is a separable morphism and $f_d$ does not descend to a degree two morphism $X\to Y_i$, combining Equations \ref{eq_genfibr} and \ref{eq_naltra} we get
\begin{equation}
\lim_{d\to\infty}\frac{g_d'}{g_d}\geq c_0(X,L).
\label{eq_anctra}
\end{equation}
Notice that Equation \ref{eq_gusunzhouc1} states that $c_0(X,L)\geq\frac{1}{2}$.
Again, applying Theorem \ref{teo_slopegusunzhou}(2), we get
\begin{equation}
K_X^2\geq \Bigl(4+\min \{c_0(X,L),\frac{1}{2}\}\Bigr)\chi(\O_X)=\frac{9}{2}\chi(\O_X)
\label{eq_mahult'}
\end{equation}
which concludes the proof.
\end{proof}

\begin{remark}
In \cite{gusunzhou} Corollary 5.15, Theorem \ref{teo_cxl} is used to characterize minimal surfaces of general type with maximal Albanese dimension as double covers of Abelian varieties: notice that the double cover can be purely inseparable.
\label{rem_=sevineq}
\end{remark}

\chapter{Severi type inequalities and surfaces close to the Severi lines in positive characteristic}
\chaptermark{Severi type inequalities in positive characteristic}
\label{chap_severi+}

In this chapter we are going to extend the results of Chapter \ref{chap_severi0} over fields of positive characteristic.
Because most of the arguments shown there rely on the theory of double covers, everything proceeds smoothly, up to some minor changes, when the characteristic of the ground field is different from two.
In particular in this case we have to pay a particular attention to the fact that in positive characteristic it may happen that for a surface $X$ we have $q(X)\neq q'(X)$ and, if $X$ is an elliptic surface, this is related to the presence of exceptional fibres.

When the characteristic of the ground field is $2$, the theory of double covers gets more wild and we are just able to give some partial results in that context.

In the first two Section we restrict to the case of characteristic different from two.
In particular in Section \ref{sec_severi+} we prove Severi type inequalities  and we give a characterization of surfaces for which equality holds while in Section \ref{sec_closesev+} we characterize surfaces lying close to the Severi lines: we see that almost everything works as over the complex numbers.
In the last two Sections we consider the case of characteristic $2$.
In Section \ref{sec_exchar2} we give some constructive Examples of Surfaces lying over the Severi lines and in Section \ref{sec_char2} we collect all the partial results we can obtain regarding Severi type inequalities.


\section{Severi type inequalities in characteristic different from $2$}
\label{sec_severi+}

In this Section we\index{Severi!type inequalities} generalize the results obtained over the complex numbers in Chapter \ref{chap_severi0}.
Because most of the work done in that Chapter depends on the theory of double covers of surfaces, which is essentially the same for fields of characteristic different from $2$, we will see that, up to some minor changes, the same proof works also in characteristic different from $2$.

\begin{theorem}
\label{teo_ultimmio+}
Let $X$ be a minimal surface of general type of maximal Albanese dimension over a field of characteristic different from $2$ and assume that $K_X^2<\frac{9}{2}\chi(\O_X)$.
Then we have that 
\begin{equation}
K_X^2\geq 4\chi(\O_X)+4(q-2)
\label{eq:luzuoneq3}
\end{equation}
where $q$ is the dimension of the Albanese variety of $X$.

In particular equality\index{Severi!lines, surfaces on}  holds, i.e.
\begin{equation}
K_X^2=4\chi(\O_X)+4(q-2)
\label{eq_teomioa+}
\end{equation}
if and only if the canonical model of $X$ is isomorphic to a double cover of a product elliptic surface ($q\geq3$) $Y=C\times E$ where $E$ is an elliptic curve and $C$ is a curve of genus $q-1$, whose branch divisor $R$ has a most negligible singularities and
\begin{equation}
R\sim_{lin} C_1+C_2+\sum_{i=1}^{2d}E_i
\label{eq_teomioabranch}
\end{equation}
where $E_i$ (respectively $C_i$) is a fibre of the first projection (respectively the second projection) of $C\times E$ and $d>7(q-2)$ or the canonical model of $X$ is a double cover of an Abelian surfaces branched over an ample divisor ($q=2$). 
Moreover we have that $\Alb(X)=\Alb(Y)$ and $q(X)=q'(X)$, i.e. the Picard scheme of $X$ is reduced.
In particular, if $q\geq 3$, the Albanese variety of $X$ is not simple.
\end{theorem}

The proof of the inequality of this  Theorem is very similar to the case over the complex numbers which can be found in \cite{lu} and which we have outlined in Section \ref{sec_seveq}.
The proof of the equality is almost the same as the one over the complex numbers we have given in Section \ref{sec_onsev}.
We analyse here what problems arise in characteristic different from $2$ and we see how to solve them.

Let $X$ be a minimal surface of general type of maximal Albanese dimension satisfying $K_X<\frac{9}{2}\chi(\O_X)$.
By Theorem \ref{teo_albfact2}, we know that there exists a  morphism $\pi\colon X\to Y$ of degree two where $Y$ is a normal surface induced by an involution $i\colon X\to X$, i.e. $Y=X/i$.
Recall that $\pi$ is a factor of the Albanese morphism i.e. the following diagram commutes:
\[
\begin{tikzcd}
X \arrow{rr}{\pi}\arrow[swap]{dr}{\alb_X}& & Y\arrow{dl}{} \\
 & \Alb(X). &
\end{tikzcd}
\]

As in Section \ref{sec_seveq} let $Y'$ be the resolution obtained by blowing up the singularities and let $X'$ be the base change of $X$ over $Y$. Denote by $Y_0$ the minimal model of $Y'$ and by $X_0$ the middle term of the Stein factorization of the morphism from $X'$ to $Y_0$. 
It follows directly from the universal property of the Albanese morphism and the fact that $\alb_X$ factors through $\pi$ that $Y_0$ is a surface of maximal Albanese dimension with $q(Y_0)=q$ and $\Alb(X)=\Alb(Y_0)$.

The surface $X_0$ may not be smooth, so we perform the canonical resolution (cf. Section \ref{sec_canres}). We get the following diagram
\begin{equation}
\begin{tikzcd}
X \arrow{d}{\pi} & X_t \arrow{d}{\pi_t}\arrow{l}[swap]{\phi_X}\arrow{r}{\phi_t} & X_{t-1}\arrow{d}{\pi_{t-1}}\arrow{r}{\phi_{t-1}} & \ldots \arrow{r}{\phi_2} & X_1\arrow{d}{\pi_1} \arrow{r}{\phi_1} & X_0 \arrow{d}{\pi_0}\\
Y=X/i & Y_t \arrow{r}{\psi_t} \arrow{l}[swap]{\psi_Y} & Y_{t-1}\arrow{r}{\psi_{t-1}} & \ldots \arrow{r}{\psi_2} & Y_1 \arrow{r}{\psi_1} & Y_0.
\end{tikzcd}
\label{eq_diagpos+}
\end{equation}
where every $\psi_i$ is a single blow-up.
Recall that we denote by $L_i$ the line bundle associated with the double cover $\pi_i\colon X_i\to Y_i$ and that $L_{i}=\psi_i^*L_{i-1}(-m_iE_i)$ where $E_i$ is the exceptional divisor of $\psi_i$ and $m_i=\intinf{d_i/2}$ where $d_i$ is the multiplicity of the branch divisor $R_{i-1}$ of $\pi_{i-1}$ at the point blown-up by $\psi_i$.

By the classification of minimal surfaces (cf. Sections \ref{sec_surfaces} and \ref{sec_ell}, we know that $Y_0$ has non-negative Kodaira dimension and maximal Albanese dimension and in particular we have the following possibilities:
\begin{itemize}
	\item if $k(Y_0)=0$, then $Y_0$ is an Abelian surface and $q=2$;
	\item if $k(Y_0)=1$, then $Y_0$ is an  elliptic surface with only smooth (but possibly  multiple) fibres over a curve $C$ with genus $g(C)=q-1$ and $q\geq 3$ (cf. Lemma \ref{lemma_ellalbjac});
	\item if $k(Y_0)=2$, then $Y_0$ is a minimal surface of general  type of maximal Albanese dimension with $q\geq 2$.
\end{itemize}

If $Y_0$ is an Abelian surface, we see that every argument of Chapter \ref{chap_severi0} works and needs no changes.

If $Y_0$ is an elliptic surface we have (cf. Theorem \ref{teo_canellsur})
\begin{equation}
\begin{split}
K_{Y_0}\sim_{num}(2g(C)-2-&\deg(\L_0))F+\sum_{j=1}^m(a_j-1)F_j=\\
= 2(q-2)F-\deg(\L_0)F+&\sum_{j=1}^m(a_j-1)F_j
\end{split}
\label{eq_candegell234}
\end{equation}
(recall that $\L_0$ is the locally-free part of $R^1(\alpha_0)_*\O_{Y_0}$ where $\alpha_0\colon Y_0\to C_0$ is the elliptic fibration and that $-\deg(\L_0)\geq0$, cf. Equation \ref{eq_bohdeg}).
Moreover, we have $F.L_0>0$ because $X$ is of general type. 
Indeed, assume by contradiction that $F.L_0\leq 0$.
In particular it follows that $\psi^*F.L_t=\psi^*F.(\psi^*L_0(-\sum_{i=1}^tm_iE_i))=F.L_0\leq 0$ (where $\psi=\psi_1\circ\ldots\circ\psi_t$) which is a contradiction because $K_{X_t}=\pi_t^*(K_{Y_t}+L_t)$ and, numerically, $K_{Y_t}$ is a multiple of a fibre.

Using Equations \ref{cansq} and \ref{eulchar}, we get
\begin{equation}
\label{eq_ell1}
\begin{split}
K_X^2-4\chi(\O_X)=K_{X_t}^2&-4\chi(\O_{X_t})+n=\\
=2(K_{Y_0}^2-4\chi(\O_{Y_0}))+2K_{Y_0}&.L_0+2\sum_{i=1}^{t}(m_i-1)+n=\\
=4(q-2)F.L_0-&2\deg(\L_0)F.L_0+\\
+2\sum_{j=1}^m(a_j-1)F_j.L_0+&2\sum_{i=1}^{t}(m_i-1)+n\geq\\
\geq 4(q-2)+2\sum_{j=1}^m(a_j-1)F_j.L_0&+2\sum_{i=1}^{t}(m_i-1)+n\geq 4(q-2).
\end{split}
\end{equation}

It is then clear that equality holds if and only if $n=0$ (i.e. $X=X_t$), $m_i=1$ for all $i$, $F.L_0=1$, $a_j=1$ for all $j$ and $\deg(\L_0)=0$.
As noticed in Remark \ref{rem_spectralsfib}, $\deg(\L_0)=0$ implies that there are no exceptional fibres and in particular that $a_j=n_j$ is the multiplicity of the multiple fibres.
In particular we can substitute the conditions $a_j=1$ for all $j$ and $\deg(\L_0)=0$ with $n_j=1$ for all $j$, as over the complex numbers, in order to characterize surfaces for which equality holds.

Notice that, in this case, we are in the situation of Corollary \ref{coro_isotrivell} with $d=1$,
in particular we know that $Y_0$ is a product elliptic surface. 

Observe that the proofs of Lemmas \ref{es2} and \ref{lem1} are characteristic-free.

The fact that, in case of equality, we have $q'(X)=q'(Y)=q(Y)=q(X)$ follows easily from K\"unneth's formula and the explicit formula of the branch divisor when $Y_0$ is an elliptic surface, while, when $Y_0$ is Abelian it is even more immediate by the fact that $L_0$ has to be ample.

Now suppose that $Y_0$ is a surface of general type.
If $K_{Y_0}^2\geq \frac{9}{2}\chi(\O_{Y_0})$ we have to generalize the proof of Lemma 4.4(2) of \cite{lu} in positive characteristic (there it is used that $h^0(Y_0,K_{Y_0})\geq 2q'(Y_0)-4$ which holds only over the complex numbers).
We have the following Lemma.

\begin{lemma}
In the above setting, if $K_{Y_0}^2\geq\frac{9}{2}\chi(\O_{Y_0})$ then $K_X^2-4\chi(\O_X)> 28 (q(X)-2)$.
\label{lemma_luzuogen+}
\end{lemma}

\begin{proof}
If $K_{Y_0}^2\geq\frac{9}{2}\chi(\O_{Y_0})$ then, thanks to Equations \ref{cansq} and \ref{eulchar} together with the assumption $K_X^2<\frac{9}{2}\chi(\O_X)$ (which holds throughout all this Section) and  using the notation as in diagram \ref{eq_diagpos+}, we have that 
\begin{equation}
\begin{split}
0>K_X^2-\frac{9}{2}\chi(\O_X)\geq& K_{X_t}^2-\frac{9}{2}\chi(\O_{X_t})=\\
=2\Bigl(K_{Y_0}^2-\frac{9}{2}\chi(\O_{Y_0})\Bigr)+\frac{7}{4}K_{Y_0}L_0-&\frac{1}{4}L_0^2+\frac{1}{4}\sum_{i=1}^t(m_i-1)(m_i+8)\geq\\
\geq \frac{1}{4}\Bigl(7K_{Y_0}.&L_0-L_0^2\Bigr).
\end{split}
\label{eq_>2gentype}
\end{equation}
In particular we have $L_0^2>7K_{Y_0}.L_0\geq 0$ where the last equality follows from the fact that $2L_0=R_0$ is effective and $Y_0$ is minimal.

By the Hodge index Theorem we have that $(K_{Y_0}.L_0)^2\geq K_{Y_0}^2L_0^2$ which, combined with Equation \ref{eq_>2gentype}, guarantees that
\begin{equation}
L_0^2 (K_{Y_0}.L_0)>7(K_{Y_0}.L_0)^2\geq 7K_{Y_0}^2L_0^2
\label{eq_>2genhodge}
\end{equation}
i.e.
\begin{equation}
K_{Y_0}.L_0>7K_{Y_0}^2.
\label{eq_>2genhodge2}
\end{equation}

Therefore the first two lines of Equation \ref{eq_ell1} becomes
\begin{equation}
\begin{split}
K_X^2-4\chi(\O_X)\geq 2(K_{Y_0}^2&-4\chi(\O_{Y_0}))+2K_{Y_0}.L_0>\\
>\chi(O_{Y_0})+14K_{Y_0}^2> 28(p_g(Y_0&)-2)\geq 28(q'(Y_0)-2)\geq\\
\geq 28(q(Y_0)-2)&=28(q(X)-2)
\end{split}
\label{eq_y0gentype}
\end{equation}
where we are using Noether's inequality \ref{eq_noetine}, Theorem \ref{teo_shepbarr} and Remark \ref{rem_qq'}.
\end{proof}

\begin{remark}
Notice that the proof of Lemma \ref{lemma_luzuogen+} follows step by step the proof of Lemma 4.4(2) of \cite{lu} (whose result is Equation \ref{eqgen2}) and the only difference is that we are using Noether's inequality instead of $h^0(Y_0,\O_{Y_0})\geq2q'(Y_0)-4$ in order to obtain an estimate of $K_X^2-4\chi(\O_X)$.
We stress here once again that it is here that we really need to assume that $K_X^2<\frac{9}{2}\chi(\O_X)$ (and now it is clear why we need it) and it is not enough to require that the Albanese morphism of $X$ factors through a morphism of degree two. 

Observe also that, even though the estimate given in Lemma \ref{lemma_luzuogen+} is weaker than the one in Lemma 4.4(2) of \cite{lu}, it is enough in the discussion in Chapter \ref{chap_severi0} for the characterization of surfaces which satisfy $K_X^2-4\chi(\O_X)=4(q-2)$ or $K_X^2-4\chi(\O_X)=8(q-2)$: indeed we just need that in this case $K_X^2-4\chi(\O_X)>28(q-2)\geq8(q-2)$.
\end{remark}

If $K_{Y_0}^2\geq 4\chi(\O_{Y_0})+4(q-2)$ we get, as over the complex  numbers (cf. Equation \ref{eqgen1} and the discussion therein) that equality holds if and only if $Y_0=Y$, $\pi\colon X\to Y$ is an \'etale double cover and the invariants of $Y$ satisfy $K_Y^2=4\chi(\O_Y)$ and $q(Y)=q(X)=2$.

The proof for the inequality in case $Y_0$ is of general type then proceeds as in \cite{lu} Theorem 1.3: if $K_{Y_0}^2\geq \frac{9}{2}\chi(\O_{Y_0})$ we use Lemma \ref{lemma_luzuogen+}, otherwise, thanks to Theorem \ref{teo_lualb2}, we know that there exists an involution on $Y_0$ relative to the Albanese morphism $\Alb_{Y_0}$ of $Y_0$.
Proceeding by induction on the power of $2$ dividing the degree of $\alb_X$, we get a sequence of surfaces $Y_i$ with $i=1,\dots,n$ with rational maps of degree two $\phi_i\colon Y_i\dashrightarrow Y_{i+1}$ such that the Albanese morphism $\alb_X$ of $X$ factors through the $\phi_i$, $Y_i$ is a surface of general type whose invariants satisfy $K_{Y_i}^2<\frac{9}{2}\chi(\O_{Y_i})$ for $i=1,\ldots,n-1$ and $Y_n$ is either a surface of general type with $K_{Y_n}^2\geq\frac{9}{2}\chi(\O_{Y_n})$ or a surface not of general type.
In both cases, through the above discussion, we have that $K_{Y_{n-1}}^2\geq 4\chi(\O_{Y_{n-1}})+4(q-2)$ and proceeding by induction and using Equation \ref{eqgen1} we conclude that $K_X^2\geq 4\chi(\O_X)+4(q-2)$.

The fact that if $Y_0$ is of general type we can not have equality, proceeds as done over the complex numbers in \cite{lu}.
Indeed, if this happens, we have that $Y=Y_0$, $\pi\colon X\to Y$ is an \'etale double cover, $K_Y^2=4\chi(\O_Y)$ and $q=2$ as we have pointed out above.
Proceeding by induction on the power of $2$ dividing the degree of $\alb_X$, we obtain, similarly as above, a sequence
\begin{equation}
X=Z_0\xrightarrow{\phi_1=\pi}Y=Z_1\xrightarrow{\phi_2}Z_2\xrightarrow{\phi_3}\ldots\xrightarrow{\phi_{n-1}}Z_{n-1}\xrightarrow{\phi_n}Z_n=\Alb(X)
\label{eq_mavero}
\end{equation}
where $n\geq2$, $Z_i$ is a surface of general type for $i=0\ldots,n-1$, $\phi_i$ is an \'etale double cover for $i=1,\ldots,n-1$, $Z_n$ is the Albanese variety of $X$ and $\phi_n$ is the resolution of singularities of a flat double cover branched over an ample divisor.
By Corollary \ref{coro_doubcovamplefundgr}, we know that $\phi_{n}$ induces an isomorphism of algebraic fundamental groups. 
In particular there exists an \'etale double cover $A\to \Alb(X)$ where $A$ is an Abelian surface such that 
\begin{equation}
\begin{tikzcd}
Z_{n-2}\arrow{r}{\phi_{n-1}}\arrow{d}{}&Z_{n-1}\arrow{d}{\phi_n}\\
A\arrow{r}{}& \Alb(X)
\end{tikzcd}
\label{eq_mahvero}
\end{equation}
is a Cartesian diagram (cf. Theorem \ref{teo_fundpushfor}).
But this is a contradiction to the fact that $\phi_n\circ\ldots\circ\phi_1=\alb_X\colon X\to \Alb(X)$.
The proof of Theorem \ref{teo_ultimmio+} is then complete.

\section{Surfaces close to the Severi lines in characteristic different from $2$}
\label{sec_closesev+}
In this section\index{Severi!lines, surfaces close to} we are going to generalize Theorem \ref{miob} over a field of characteristic different from $2$.
This is almost immediate after the discussion of Section \ref{sec_severi+}.

\begin{theorem}
\label{teo_closesev+}
Let $X$ be a minimal surface of general type with maximal Albanese dimension with $K_X^2<\frac{9}{2}\chi(\O_X)$ and $q$ be the dimension of the Albanese variety $\Alb(X)$ of $X$.
\begin{enumerate}
\item If $K^2_X>4\chi(\O_X)+4(q-2)$, then $K^2_X\geq 4\chi(\O_X)+8(q-2).$
\item If $q=2$ and $K^2_X>4\chi(\O_X)$, then $K^2_X\geq 4\chi(\O_X)+2.$
\item If $q\geq 3$, equality holds, i.e.
\begin{equation}
\label{eq_8q-2}
K^2_X=4\chi(\O_X)+8(q-2),
\end{equation}
if and only if the canonical model of $X$ is isomorphic to a double cover of a smooth isotrivial elliptic surface $Y$ over a curve $C$ of genus $q-1$, branched over a divisor $R$ with at most negligible singularities for which $K_Y.R=8(q-2)$. In particular, we have that $\Alb(X)\simeq\Alb(Y)$ and $q(X)=q'(X)$, i.e. the Picard variety of $X$ is reduced, and the Albanese variety of $X$ is not simple.
\end{enumerate}
\end{theorem}

\begin{proof}
The proof is completely similar to the one of Theorem \ref{miob}, we just stress here the differences.

If $Y_0$ is Abelian the proof is identical as over the complex numbers.

If $Y_0$ is an elliptic surface, by Equation \ref{eq_ell1} we have that 
\begin{equation}
\begin{split}
K_X^2-4\chi(\O_X)=4(q-2)F&.L_0-2\deg(\L_0)F.L_0+\\
+2\sum_{j=1}^m(a_j-1)F_j.L_0+&2\sum_{i=1}^{t}(m_i-1)+n.
\end{split}
\label{eq_ell1b}
\end{equation}
As already observed in Section \ref{sec_severi+}, $K_X^2-4\chi(\O_X)=4(q-2)$ holds if and only if $F.L_0=1$, $\deg(\L_0)=0$, $a_j=1$ for all $j$, $m_i=1$ for all $i$ and $n=0$.
Notice that if $\deg(\L_0)<0$ or there exists an index $j$ such that $a_j>1$ then there exists at least a multiple fibre in the elliptic surface $Y_0$.
Hence we can argue as over the complex numbers: whatever of the above quantities we increase we obtain that $F.L_0\geq 2$ and equality holds if and only if $F.L_0=2$,  $\deg(\L_0)=0$, $a_j=1$ for all $j$, $m_i=1$ for all $i$ and $n=0$.

If $Y_0$ is a surface of general type we argue as over the complex numbers using Lemma \ref{lemma_luzuogen+} instead of Equation \ref{eqgen2} and we are done.

It remains to show that $q(X)=q'(X)$, which we prove in the following Lemma.
\end{proof}

\begin{lemma}
Let $\alpha\colon Y\to C$ be a smooth elliptic surface fibration with maximal Albanese dimension and let $X$ be the minimal smooth model of the double cover given by the equation $2L=\O_Y(R)$ where $R$ has at most  negligible singularities, denote by $\pi\colon X\to Y$ the induced morphism.
Suppose moreover that $X$ is a surface of general type, $q(X)=q(Y)$, $L.F=2$ where $F$ is a fibre of $Y\to C$ and $K_X^2<\frac{9}{2}\chi(\O_X)$.
Then we have that $q'(X)=q(X)$
\label{lemma_q'=q}
\end{lemma}

\begin{proof}
By 
Remark \ref{rem_ellsurf1}, we know that $q(Y)=q'(Y)=g(C)+1$ and by Corollary \ref{coro_isoell} we know that $Y\to C$ becomes trivial after an \'etale cover $\delta\colon\widetilde{C}\to C$ of degree $d$ coprime with the characteristic of the ground field (here we are using $L.F=2$ and $\cha(k)\neq 2$).
In particular we have the following Cartesian diagram
\begin{equation}
\begin{tikzcd}
\widetilde{X}\arrow{r}{\beta}\arrow{d}{\widetilde{\pi}}& X\arrow{d}{\pi}\\
\widetilde{Y}=\widetilde{C}\times F\arrow{d}{\widetilde{\alpha}}\arrow{r}{\gamma} & Y\arrow{d}{\alpha}\\
\widetilde{C}\arrow{r}{\delta} & C.
\end{tikzcd}
\label{eq_hofinitoleidee}
\end{equation} 

If we assume that $L$ is ample we conclude by
\begin{equation}
\begin{split}
q'(X)=h^1(X,\O_X)&=h^1(Y,\pi_*\O_X)=\\
=h^1(Y,\O_{Y})+h^1(Y,L^{-1})=&q'(Y)=q(Y)=q(X)
\end{split}
\label{eq_q'=q}
\end{equation}
where $h^1(Y,L^{-1})=0$ follows by Lemma \ref{lemma_kodvan}.

Notice that $L$ is ample if and only if $\widetilde{L}=\gamma^*L$ is, which is the line bundle associated with the double cover $\widetilde{\pi}$.
We know that $\widetilde{L}$ is effective up to a multiple (twice $\widetilde{L}$ is the ramification divisor of $\widetilde{\pi}$) and, by the formula of the Picard group for $\widetilde{C}\times F$ (cf. Proposition \ref{propo_piccxe}), we know it is numerically equivalent to $\Gamma_f+\widetilde{C}+lF$ where $l\geq0$, $F$ and $\widetilde{C}$ are fibres of $\widetilde{\alpha}$ and $\widetilde{C}\times F\to F$ respectively and $f\colon\widetilde{C}\to F$ is a morphism.
If $f$ is not a constant function, it is clear by the Nakai-Moishezon criterion (cf. \cite{badescu} Theorem 1.22) that $\widetilde{L}$ is an ample line bundle (recall that $\Gamma_f^2=0$, cf. Section \ref{sec_piccxe}).
So suppose that $f$ is constant, i.e. $\widetilde{L}$ is numerically equivalent to $2\widetilde{C}+lF$. 
If $l=0$ we would have that $\widetilde{L}^2=0$ and $\widetilde{L}.K_{\widetilde{Y}}>0$.
In particular Equations \ref{cansq} and \ref{eulchar} and the fact that for \'etale covers the Euler characteristic and the self intersection of the canonical divisor are multiplicative tell us that
\begin{equation}
\begin{split}
0>K_X^2-\frac{9}{2}\chi(\O_X)=\frac{1}{d}(K_{\widetilde{X}}^2-\frac{9}{2}\chi(\O_{\widetilde{X}}))=\\
=\frac{2}{d}(K_{\widetilde{Y}}^2-\frac{9}{2}\chi(\O_{\widetilde{Y}}))+\frac{7}{4d}\widetilde{L}.K_{\widetilde{Y}}-\frac{1}{4d}\widetilde{L}^2>0
\end{split}
\label{eq_bohultcha+}
\end{equation}
which is a contradiction.
So $l>0$ and, again by the Nakai-Moishezon criterion, $\widetilde{L}$ is ample also in this case.
\end{proof}

\begin{remark}
Observe that all the Examples over the complex numbers given in Section \ref{sec_ex} are still valid with the same constructions over a field of characteristic different from $2$.
In particular Examples \ref{es1} and \ref{buel} are clearly still valid in this context;
for Example \ref{es4} we have to notice that $\jac C'$ of a very general curve $C'$ of genus $g'$ is simple even in positive characteristic (cf. Appendix B of \cite{moonen}) and that every Jacobian variety is a limit of a family of Prym varieties even in positive characteristic different from $2$ (cf. \cite{beauville-prym} or \cite{farkas} Section 3.1) so that even in this case we can conclude that the very general Prym variety $P(C,C')$ is simple.

In Example \ref{ex_abe} we used a version of Bertini Theorem which still holds in positive characteristic because we considered very ample divisors (cf. \cite{cumino} Corollary 1 for the most general version of Bertini Theorem in positive characteristic).
\end{remark}

\section{Examples in characteristic $2$}
\label{sec_exchar2}
In this Section we give explicit examples of surfaces of general type with maximal Albanese dimension over a field of characteristic $2$ lying over the Severi lines\index{Severi!lines, surfaces on}, i.e. which satisfy $K_X^2-4\chi(\O_X)=4(q-2)<\frac{9}{2}\chi(\O_X)$.
For all these examples we will have that the Picard scheme is reduced, even if this is not ensured by Lemma \ref{lemma_specialtype+} of the next Section.
Examples \ref{es_addelles} and \ref{es_addhyperelles} have purely inseparable Albanese morphism of degree $2$ , while Examples \ref{es_splitsepdoubleabe} and \ref{es_splitsepdoubleprodell} have separable Albanese morphism of degree $2$.

\begin{example}[Double cover of a product Abelian surface via a rational vector field]

\label{es_addelles}

This is a generalization of Example 7.1 of \cite{gusunzhou} and Example 6.2 of \cite{takeda} (there it is considered only the case of a double cover of a product elliptic supersingular elliptic surface via a multiplicative rational vector field).
Let $E_1$ and $E_2$ be two elliptic curves and let $\delta_{E_1}'$ and $\delta_{E_2}'$ be the multiplicative rational vector fields defined in Example \ref{es_multaddvf1}.
Consider the quotient $g_0\colon Y_0=E_1\times E_2\to X_0^{(-1)}$ via the rational vector field $\delta_{E_1}'+\delta_{E_2}'$ and its dual morphism $\pi_0\colon X_0\to Y_0$.
After the canonical resolution we get the following diagram 
\begin{equation}
\begin{tikzcd}
X\arrow{d}{\pi}\arrow{r}{\phi} & X_0\arrow{d}{\pi_0}\\
Y\arrow{r}{\psi} & Y_0=E_1\times E_2:
\end{tikzcd}
\label{eq_canres2es}
\end{equation}
by Proposition \ref{propo_liedtkeaxb} and Corollary \ref{coro_liedtkeaxb} we know that $X$ is a surface of general type of maximal Albanese dimension with $q'(X)=q(X)=g(E_1)+g(E_2)=2$ whose canonical model is $X_0$. 
Moreover (ibidem) the $1$-foliation $\F_0$ associated with $\delta_{E_1}'+\delta_{E_2}'$ is linearly equivalent to $-2E_1'-2E_2'$ and the line bundle $L_0$ associated with $\pi_0$ is linearly equivalent to $E_1+E_2$ where $E_2$, $E_2'$ and $E_1$, $E_1'$ are fibres of the first and the second projection of $E_1\times E_2$ and $\pi_0$ is an inseparable double cover.

Using Equations \ref{cansq} and \ref{eulchar} we get
\begin{equation}
K_X^2=2L_0^2=4
\label{eq_mah11}
\end{equation}
and
\begin{equation}
\chi(\O_X)=\frac{1}{2}L_0^2=1
\label{eq_max12}
\end{equation}
from which
\begin{equation}
K_X^2-4\chi(\O_X)=0
\label{eq_char2elell}
\end{equation}
follows.
It is clear that the Albanese morphism of $X$ is $\psi\circ\pi$ because otherwise we would have that $X$ is birational to an Abelian variety.

If we assume that $E_1$ and $E_2$ are not supersingular and we take the two additive rational vector fields $\delta_{E_1}$ and $\delta_{E_2}$ as in Example \ref{es_multaddvf1} we obtain another construction of a surface of general type lying over the Severi line and everything works as above.

With similar calculation as above, we obtain a surface $X$ satisfying $K_X^2=4\chi(\O_X)$ if in the above discussion we substitute the rational vector fields $\delta_{E_1}$, $\delta_{E_2}$, $\delta_{E_1}'$ and $\delta_{E_2}'$ with $\overline{\delta}_{E_1}$, $\overline{\delta}_{E_2}$, $\overline{\delta}_{E_1}'$ and $\overline{\delta}_{E_2}'$.
\end{example}

\begin{example}[Splittable separable double cover of an Abelian surface]
\label{es_splitsepdoubleabe}
Let $Y_0$ be an Abelian surface and $L_0$ be an ample line bundle on $Y_0$.
Let $R_1$ and $R_2$ be two effective divisor such that $\O_{Y_0}(R_1)=L_0$ and $\O_{Y_0}(R_2)=L_0^{\otimes2}$ and every singular point of $R_1$ is a smooth point of $R_2$ and no irreducible component of $R_1$ is contained twice in $R_2$.
Notice that it is not difficult to find such divisors in an Abelian surface: if we further assume $L_0$ to be very ample, we can even take $R_1$ to be smooth by Bertini's Theorem (cf. \cite{hart} Theorem II.8.18) and $R_2$ whatever effective divisor in the linear series of $2R_1$ different from $2R_1$ itself.

Consider the splittable separable double cover $\pi_0\colon X_0\to Y_0$ associated with $R_1$ and $R_2$ (cf. Remark \ref{rem_splitdoublecover2}): we know that it has only rational double points as singularities (outside the singular points of $R$ we derive this from Lemma \ref{lemma_branchan} while on singular points of $R$ it is a consequence of Remark \ref{rem_normsep2}) and, after the canonical resolution, we get
\begin{equation}
\begin{tikzcd}
X\arrow{d}{\pi}\arrow{r}{\phi} & X_0\arrow{d}{\pi_0}\\
Y\arrow{r}{\psi} & Y_0
\end{tikzcd}
\label{eq_quad}
\end{equation}
where $X_0$ is the canonical model of $X$.
By the ampleness of $L_0$, it is clear that $X$ is a surface of general type (indeed $K_{X_0}=\pi_0^*L_0$) and that 
\begin{equation}
q'(X)=h^1(Y_0,\O_{Y_0})+h^1(Y_0,L_0^{-1})=2
\label{eq_stellaq}
\end{equation}
(cf. \cite{mumford} The Vanishing Theorem Section III.16 for the vanishing of $h^1(Y_0,L_0^{-1})$). 
Moreover the Albanese morphism of $X$ is $\psi\circ\pi$ because otherwise $X$ would be birational to an Abelian surface.

Using Equations \ref{cansq} and \ref{eulchar} we get
\begin{equation}
K_X^2=2L_0^2
\label{eq_mah11bis}
\end{equation}
and
\begin{equation}
\chi(\O_X)=\frac{1}{2}L_0^2
\label{eq_max12bis}
\end{equation}
from which
\begin{equation}
K_X^2-4\chi(\O_X)=0
\label{eq_char2elellbis}
\end{equation}
follows.
\end{example}

\begin{example}[Double cover of a product  elliptic surface via a rational vector field]

\label{es_addhyperelles}

Let $C$ be a hyperelliptic curve of genus $g$ and $E$ be an elliptic curve and let $\delta_{E}'$ and $\delta_{C}'$ be the multiplicative rational vector fields defined in Examples \ref{es_multaddvf1} and \ref{es_multaddvf>1}.
Consider the quotient $g_0\colon Y_0=C\times E\to X_0^{(-1)}$ via the rational vector field $\delta_{C}'+\delta_{E}'$ and its dual morphism $\pi_0\colon X_0\to Y_0$.
After the canonical resolution we get the following diagram 
\begin{equation}
\begin{tikzcd}
X\arrow{d}{\pi}\arrow{r}{\phi} & X_0\arrow{d}{\pi_0}\\
Y\arrow{r}{\psi} & Y_0=C\times E:
\end{tikzcd}
\label{eq_canres2es}
\end{equation}
by Proposition \ref{propo_liedtkeaxb} and Corollary \ref{coro_liedtkeaxb} we know that $X$ is a surface of general type of maximal Albanese dimension with $q'(X)=q(X)=g(C)+g(E)=g+1$ (which we denote simply by $q$) whose canonical model is $X_0$.
Moreover (ibidem) the $1$-foliation $\F_0$ associated with $\delta_{C}'+\delta_{E}'$ is linearly equivalent to $-2C'-\sum_{i=1}^{2d+g+1}2E_i'$ and the line bundle $L_0$ associated with $\pi_0$ is linearly equivalent to $C+\sum_{i=1}^{2d+2}E_i$ where $E_i$, $E_i'$ and $C$, $C'$ are fibres of the first and the second projection of $C\times E$, $d$ is as in the definition of the rational vector field $\delta_C'$ and $\pi_0$ is an inseparable double cover.

Using Equations \ref{cansq} and \ref{eulchar} we get
\begin{equation}
K_X^2=2K_{Y_0}^2+4K_{Y_0}.L_0+2L_0^2=8(g-1)+8(d+1)=8(q-2)+8(d+1)
\label{eq_mah21}
\end{equation}
and
\begin{equation}
\chi(\O_X)=2\chi(\O_{Y_0})+\frac{1}{2}K_{Y_0}.L_0+\frac{1}{2}L_0^2=q-2+2(d+1)
\label{eq_mah22}
\end{equation}
from which we derive
\begin{equation}
K_X^2-4\chi(\O_X)=4(q-2)
\label{eq_char2elellbis}
\end{equation}
and
\begin{equation}
K_X^2-\frac{9}{2}\chi(\O_X)=\frac{7}{2}(q-2)-d-1.
\label{eq_mahboh}
\end{equation}
It easily follows that $K_X^2-\frac{9}{2}\chi(\O_X)<0$ if and only if $d>\frac{7}{2}(q-2)-1$: in this case we can apply Theorem \ref{teo_albfact2} which implies that $\psi\circ\pi\circ\alb_{C\times E}$ is the Albanese morphism of $X$. 

If we assume that $E$ is not not supersingular and we take the two additive rational vector fields $\delta_{C}$ and $\delta_{E}$ as in Examples \ref{es_multaddvf1} and \ref{es_multaddvf>1} we obtain another construction of a surface of general type lying over the Severi lines as above.
\end{example}

\begin{example}[Splittable separable double cover of a product elliptic surface]
\label{es_splitsepdoubleprodell}
Let $C$ be a curve of genus $g\geq 2$, $E$ be an elliptic curve and denote by $Y_0=C\times E$.
Let $R_1$ and $R_2$ be two effective divisors such that $\O_{Y_0}(R_1)=L_0$ and $\O_{Y_0}(R_2)=L_0^{\otimes2}$ such that every singular point of $R_1$ is a smooth point of $R_2$ and no irreducible components of $R_1$ are contained twice in $R_2$ where $L_0\sim_{hom}\O_{Y_0}(C+dE)$ with $d>7(q-2)$ and $q=g+1$.
For example we can take $R_1=C_1+\sum_{i=1}^dE_i$ and $R_2=C_1+C_2+\sum_{i=d+1}^{3d}E_i$ where $C_i$ (respectively $E_i$) are fibres of the natural projection $\pi_E\colon Y_0\to E$ (respectively $\pi_C\colon Y_0\to C$) and $C_1\neq C_2$ (respectively $E_i\neq E_j$ for $i\neq j$).
Notice that, by definition, $C_1$ is the only curve contained in both $R_1$ and $R_2$.

Consider the splittable separable double cover $\pi_0\colon X_0\to Y_0$ associated with $R_1$ and $R_2$ (cf. Remark \ref{rem_splitdoublecover2}): we know that it has only rational double points as singularities (outside the singular points of $R$ we derive this from Lemma \ref{lemma_branchan} while on singular points of $R$ it is a consequence of Remark \ref{rem_normsep2}) and after the canonical resolution we get
\begin{equation}
\begin{tikzcd}
X\arrow{d}{\pi}\arrow{r}{\phi} & X_0\arrow{d}{\pi_0}\\
Y\arrow{r}{\psi} & Y_0=C\times E
\end{tikzcd}
\label{eq_quad}
\end{equation}
where $X_0$ is the canonical model of $X$.
By the ampleness of $L_0$, it is clear that $X$ is a surface of general type (indeed $K_{X_0}=\pi_0^*(L_0+K_{Y_0})$) and that 
\begin{equation}
q'(X)=h^1(Y_0,\O_{Y_0})+h^1(Y_0,L_0^{-1})=q.
\label{eq_stellaq}
\end{equation}

By Equations \ref{cansq} and \ref{eulchar} we get
\begin{equation}
K_X^2=2K_{Y_0}^2+4K_{Y_0}.L_0+2L_0^2=8(q-2)+4d
\label{eq_mambo1}
\end{equation}
and 
\begin{equation}
\chi(\O_X)=\chi(\O_{Y_0})+\frac{1}{2}K_{Y_0}.L_0+\frac{1}{2}L_0^2=q-2+d
\label{eq_mambo2}
\end{equation}
from which we easily derive
\begin{equation}
K_X^2-4\chi(\O_X)=4(q-2)
\label{eq_mambo3}
\end{equation}
and
\begin{equation}
K_X^2-\frac{9}{2}\chi(\O_X)=\frac{7}{2}(q-2)-\frac{1}{2}d<0
\label{eq_mambo4}
\end{equation}
where the last inequality follows from the assumption $d>7(q-2)$.

Then we can apply Theorem \ref{teo_albfact2} which implies that $\psi\circ\pi\circ\alb_{Y_0}$ is the Albanese morphism of $X$. 
\end{example}

\section{Partial results in characteristic $2$}
\label{sec_char2}
In this section we are going to see what can be said about our results when the characteristic of the ground field is $2$.
The results of this thesis are based on the theory of double covers of surfaces and that's why things get wilder when the characteristic of the ground field is exactly $2$.
In chapter \ref{chap_double} we have shown the theory of double covers, with a particular attention to the case of characteristic $2$: even though this case is more complicated we have some similarities.

First of all we noticed that for every double cover $f\colon X\to Y$ ther is naturally associated a line bundle $L$ via the following short exact sequence:
\begin{equation}
0\to \O_Y\to f_*\O_X\to L^{-1}\to 0.
\label{eq:}
\end{equation}

If the characteristic is different from two, we have seen that this exact sequence splits and the divisor associated with $L^2$ is the ramification divisor of $f$: in particular we have that $L$ is an effective divisor up to a multiple.
We have used this effectiveness result for $L$  to show that $K_Y.L\geq 0$ when $Y$ is a minimal non-ruled surface, and this is the real result we would like to extend in characteristic two.

As we have already noticed in Remark \ref{rem_effdouble}, the line bundle $L$ is effective up to a multiple in most cases even when the characteristic of the ground field is $2$.
There is just one single case where nothing assures $L$ to be effective up to a multiple, i.e. for non-splittable inseparable double covers.

Now suppose that $f\colon X\to Y$ is an inseparable double cover such that $f$ is a factor of the Albanese morphism $\alb_X\colon X\to \Alb(X)$. 
For simplicity we also assume that both $X$ and $Y$ are smooth.
As we have seen in Section \ref{sec_insdouble}, it is canonically defined another double cover $g^{(1)}\colon Y^{(1)}\to X$. 

As we have seen in Remark \ref{rem_albinsfact}, the line bundle $M^{-1}$, which is the inverse of the line bundle associated with $g^{(1)}\colon Y^{(1)}\to X$, has many sections.
At first glance, we thought it were possible to prove that $L$ has to be effective up to a multiple using this and the relations of Lemma \ref{lemma_flgm} between $L$, $M$ and the associated foliations on $Y^{(1)}$ and $X$, but in the end we were not able to do so.

Without this, we can only collect here partial results concerning Severi type inequalities when the characteristic of the ground field is $2$.

Let $X$ be a surface of general type with maximal Albanese dimension over a field of characteristic $2$ satisfying $K_X^2<\frac{9}{2}\chi(\O_X)$: by Theorem \ref{teo_albfact2} we know that there exists a morphism of degree two $\pi\colon X\to Y$ with $Y$ normal such that
\[
\begin{tikzcd}
X \arrow{rr}{\pi}\arrow[swap]{dr}{\alb_X}& & Y\arrow{dl}{} \\
 & \Alb(X) &
\end{tikzcd}
\]
commutes.

As we have done over the complex numbers (cf. Equation \ref{eq_char0xtoy}) or for field of characteristic different from two (cf. Equation \ref{eq_diagpos+}) we obtain the following diagram
\begin{equation}
\begin{tikzcd}
X\arrow{d}{\pi}& X_t\arrow{l}{\phi}\arrow{d}{\pi_t}\arrow{r}{\phi_0} & X_0\arrow{d}{\pi_0}\\
Y & Y_t\arrow{l}{\psi}\arrow{r}{\psi_0} & Y_0
\end{tikzcd}
\label{eq_diagpos2}
\end{equation}
where $Y_0$ is the minimal smooth model of $Y$ and $\pi_t\colon X_t\to Y_t$ is the canonical resolution of $\pi_0\colon X_0\to Y_0$ (cf. Section \ref{sec_canres}).

\begin{lemma}
\label{lemma_specialtype+}
Suppose that $Y_0$ has Kodaira dimension strictly smaller than $2$.
Then we have that 
\begin{equation}
K_X^2\geq 4\chi(\O_X)+4(q-2)
\label{eq_+sevegen}
\end{equation}
and equality holds if and only if the canonical model of $X$ is isomorphic to a double cover of a product elliptic surface $C\times E$ ($q\geq 3$) or of an Abelian surface ($q=2$).

Suppose moreover that $\pi_0$ is separable and $q\geq3$: in this case we can further prove that the line bundle $L_0$ associated with $\pi_0$ is effective and is linearly equivalent to 
\begin{equation}
L\sim_{lin}C+\sum_{i=1}^dE_i
\label{eq_lineqchar2}
\end{equation}
where, as usual, $C$ (respectively $E_i$) is a fibre of the first projection (respectively are fibres of  the second projection) of $C\times E$  and $d>7(q-2)$ and $q'(X)=q(X)$, i.e. the Picard scheme of $X$ is reduced.
\end{lemma}

\begin{proof}
If $Y_0$ is an Abelian surface it is clear that everything works as over a field of characteristic different from two.

If $Y_0$ is an elliptic surface we have seen in Section \ref{sec_severi+}, we have that $F.L_0>0$ and the proof there given only relies on the fact that $X$ is a surface of general type and the canonical bundle of $Y_0$ is numerically equivalent to a rational multiple of $F$, so it works even if we do not know whether $L_0$ is effective up to a multiple or not.
Once this has been noticed, the Lemma is proven by Equation \ref{eq_ell1} and Corollary \ref{coro_isotrivell}. 
%
%
%
If $\pi_0$ is separable, the linear class of $L$ is proved in a similar and simpler, way as in Lemma \ref{es2}.
The fact that $d>7(q-2)$ is required in order to have $K_X^2<\frac{9}{2}\chi(\O_X)$ (cf. Example \ref{es1}) which is an assumption for all this Section.
Then K\"unneth formula and the long exact sequence in cohomology associated with the short exact sequence associated with the double covers $\pi_0\colon X_0\to Y_0$ shows that $q'(X)=g(E)+g(C)=q(X)$.
\end{proof}

\begin{definition}
Let $X$ be a surface of general type of maximal Albanese dimension.
Suppose that the Albanese morphism $\alb_X$ of $X$ factorizes as
\begin{equation}
X\xrightarrow{f_1} X_1\xrightarrow{f_2}\ldots \xrightarrow{f_n} X_n \xrightarrow{g} \Alb(X)
\label{eq_star}
\end{equation}
where the minimal models $\widetilde{X}_1, \ldots, \widetilde{X}_n$ of $X_1, \ldots, X_{n-1}$ are surfaces of general type, $f_i$ are finite morphism of degree $2$ and either $\widetilde{X}_n$ is not of general type or  $\widetilde{X}_n$ is of general type and $g\colon X_n\to\Alb(X)$ does not factor through any morphism of degree two.
In this situation, we say that $X$ has the property (*) if it satisfies the following condition: the line bundle $\widetilde{L}_i$ associated with the double cover $\widetilde{f}_i\colon\widetilde{X}_{i-1}\to\widetilde{X}_i$ intersects non negatively the canonical bundle of $\widetilde{X}_i$ for every $i$ such that the minimal model of $X_i$ is of general type.
\label{def_star}
\end{definition}

\begin{remark}
Observe that the condition (*) is very technical, hence we give some stronger conditions which assure that it holds. 
If all the $f_i$ defined above are separable or splittable and inseparable (e.g. $\alb_X$ is a separable morphism), we clearly have that $X$ satisfies (*): indeed in this case we have seen that $\widetilde{L}_i$ is effective up to a multiple (cf. Remark \ref{rem_effdouble}).
\label{rem_star}
\end{remark}

\begin{theorem}
Let $X$ be a surface of general type with maximal Albanese dimension satisfying $K_X^2<\frac{9}{2}\chi(\O_X)$ over a field of characteristic $2$ for which condition (*) holds. 
Then we have $K_X^2-4\chi(\O_X)\geq 4(q-2)$ and equality holds if and only if it is a flat double cover of the blow-up of a product elliptic surface $C\times E$ ($q\geq 3$) or of an Abelian surface ($q=2$).
\label{teo_star2}
\end{theorem}

\begin{proof}
The condition (*) is exactly the one that assures that the same proof over a field of characteristic different from two works also in this setting.
Observe that, using the notation of Equation \ref{eq_mavero}, the morphism $\phi_n\colon Z_{n-1}\to Z_n$ can be the resolution of singularities of a purely inseparable flat double cover: also in this case, thanks to Remark \ref{rem_ultimmss}, we have that $\phi_n$ induces an isomorphism on algebraic fundamental groups and the same argument used there works.
\end{proof}

\backmatter 

\phantomsection
\addcontentsline{toc}{chapter}{Bibliography} 
\nocite{*} 
\bibliographystyle{alpha} 
\bibliography{bibliografia} 

\phantomsection 
\addcontentsline{toc}{chapter}{\indexname} 
\printindex 

\end{otherlanguage} 

\end{document}